\documentclass{article}
\usepackage{amsthm,amsmath,amsfonts,amssymb}
\usepackage{natbib}
\usepackage[colorlinks,citecolor=blue,urlcolor=blue]{hyperref}
\usepackage{graphicx}

\usepackage{enumitem,epsfig}
\usepackage[dvipsnames]{xcolor}
\usepackage{multirow}
\usepackage{multicol}
\usepackage{cleveref}
\usepackage{subcaption}
\usepackage[hmargin={1.0in, 1.0in}, vmargin={1.2in, 1.2in}]{geometry}
\usepackage{booktabs}
\usepackage{array}
\newcolumntype{Z}{>{\setbox0=\hbox\bgroup}c<{\egroup}@{\hspace*{-\tabcolsep}}}
\usepackage[title]{appendix}

\usepackage{setspace}
\usepackage{diagbox}
\usepackage{authblk}
\usepackage{comment}

\allowdisplaybreaks[4]

\newtheorem{example}{Example} 
\newtheorem{theorem}{Theorem}
\newtheorem{lemma}[theorem]{Lemma} 
\newtheorem{proposition}[theorem]{Proposition} 
\newtheorem{remark}[theorem]{Remark}

\newtheorem{definition}[theorem]{Definition}

\newtheorem{assumpt}{Assumption}
\usepackage[title]{appendix}

\newenvironment{assumptbis}[1]
{\renewcommand{\theassumpt}{\ref{#1}$'$}%
 \addtocounter{assumpt}{-1}%
 \begin{assumpt}}
{\end{assumpt}}

\newcommand{\ca}{\mathcal{A}}
\newcommand{\cb}{\mathcal{B}}
\newcommand{\cc}{\mathcal{C}}
\newcommand{\cd}{\mathcal{D}}
\newcommand{\ce}{\mathcal{E}}
\newcommand{\cf}{\mathcal{F}}

\newcommand{\ch}{\mathcal{H}}

\newcommand{\ck}{\mathcal{K}}

\newcommand{\cn}{\mathcal{N}}
\newcommand{\cp}{\mathcal{P}}

\newcommand{\cs}{\mathcal{S}}
\newcommand{\ct}{\mathcal{T}}
\newcommand{\cu}{\mathcal{U}}
\newcommand{\cv}{\mathcal{V}}
\newcommand{\cw}{\mathcal{W}}

\newcommand{\be}{\mathbb{E}}
\newcommand{\bn}{\mathbb{N}}
\newcommand{\bp}{\mathbb{P}}

\newcommand{\br}{\mathbb{R}}
\newcommand{\bz}{\mathbb{Z}}
\newcommand{\bone}{\mathbb{I}}

\newcommand{\fX}{\mathbf{X}}
\newcommand{\fY}{\mathbf{Y}}
\newcommand{\fZ}{\mathbf{Z}}

\newcommand{\fone}{\mathbf{I}}

\newcommand{\fzeta}{\boldsymbol{\zeta}}

\newcommand{\Var}{\mbox{Var}}

\newcommand{\Cov}{\mbox{Cov}}

\newcommand{\cum}{\mbox{cum}}
\newcommand{\MMD}{\mbox{MMD}}

\newcommand{\stsim}[1]{\stackrel{#1}{\sim}}
\newcommand{\beft}[2]{\mathbb{E}_{#1}[{#2}]}
\newcommand{\bbe}[1]{\mathbb{E}\big[{#1}\big]}
\newcommand{\bbeft}[2]{\mathbb{E}_{#1}\big[{#2}\big]}
\newcommand{\blre}[1]{\mathbb{E}\left[{#1}\right]}
\newcommand{\blreft}[2]{\mathbb{E}_{#1}\left[{#2}\right]}

\newcommand{\blrp}[1]{\mathbb{P}\left({#1}\right)}

\newcommand{\bigp}[1]{\big({#1}\big)}
\newcommand{\lrp}[1]{\left({#1}\right)}
\newcommand{\bigcp}[1]{\big\{{#1}\big\}}
\newcommand{\lrcp}[1]{\left\{{#1}\right\}}

\newcommand{\lrbk}[1]{\left[{#1}\right]}
\newcommand{\lrag}[1]{\langle{#1}\rangle}
\newcommand{\lrabs}[1]{\left|{#1}\right|}
\newcommand{\lrfloor}[1]{\lfloor{#1}\rfloor}
\newcommand{\lrceil}[1]{\lceil{#1}\rceil}

\newcommand{\sumin}{\sum\limits_{i=1}^{n}}

\newcommand{\sumjm}{\sum\limits_{j=1}^{m}}

\newcommand{\BE}{\begin{equation}}
\newcommand{\EE}{\end{equation}}
\newcommand{\BEqn}{\begin{eqnarray*}}
\newcommand{\EEqn}{\end{eqnarray*}}

\newcommand{\bin}[2]{
  \left(
        \begin{array}{@{}c@{}}
          #1 \\ #2
        \end{array}
  \right)}

\makeatletter
\newcommand*{\addFileDependency}[1]{
  \typeout{(#1)}
  \@addtofilelist{#1}
  \IfFileExists{#1}{}{\typeout{No file #1.}}
}
\newcommand*\showfontsize{\f@size{} point}
\makeatother

\newcommand{\blind}{1}

\begin{document}
 
\if1\blind
{
    \title{Two Sample Testing in High Dimension via \\ Maximum Mean Discrepancy}
		
    \author[1]{Hanjia Gao}
    \author[1]{Xiaofeng Shao}
    \affil[1]{Department of Statistics, University of Illinois at Urbana-Champaign}
	\date{}	
	\maketitle
} \fi


\begin{abstract}
Maximum Mean Discrepancy (MMD) has been widely used in the areas of machine learning and statistics to quantify the distance between two distributions in the $p$-dimensional Euclidean space. The asymptotic property of the sample MMD has been well studied when the dimension $p$ is fixed using the theory of U-statistic. As motivated by the frequent use of MMD test for data of moderate/high dimension, we propose to investigate the behavior of the sample MMD in a high-dimensional environment and develop a new studentized test statistic. Specifically, we obtain the central limit theorems for the studentized sample MMD as both the dimension $p$ and sample sizes $n,m$ diverge to infinity. Our results hold for a wide range of kernels, including popular Gaussian and Laplacian kernels, and also cover energy distance as a special case. We also derive the explicit rate of convergence under mild assumptions and our results suggest that the accuracy of normal approximation can improve with dimensionality. Additionally, we provide a general theory on the power analysis under the alternative hypothesis and show that our proposed test can detect difference between two distributions in the moderately high dimensional regime. Numerical simulations demonstrate the effectiveness of our proposed test statistic and normal approximation.
\end{abstract}



\section{Introduction}
\label{Sec:Intro}

Testing whether two samples are drawn from the same distribution is  a classical problem in statistics. Mathematically speaking, given independent and identically distributed (iid) $p$-dimensional  samples $X_1,\dots,X_n$ from the distribution $F_X$ and $Y_1,\dots,Y_m$ from the distribution $F_Y$, we aim to test the hypothesis $H_0:F_X=F_Y$ versus $H_A:F_X\not=F_Y$. There is a rich literature for the two-sample testing and  well-known tests include Kolmogorov-Smirnov test [\cite{kolmogorov1933sulla}, \cite{smirnov1939estimation}], Cramer von-Mises test [\cite{cramer1928composition}] and Anderson-Darling test [\cite{anderson1952asymptotic}]. Other notable ones include Wald-Wolfowitz runs test [\cite{wald1940test}], Mann-Whitney test [\cite{mann1947test}] for univariate distributions and their multivariate generalizations  [\cite{friedman1979multivariate}], among others.

In this article, we focus on the test based on maximum mean discrepancy
(MMD, hereafter) [\cite{gretton2012kernel}], which is defined as the largest difference in expectations over functions in the unit ball of a reproducing kernel Hilbert space (RKHS). Since its introduction in the machine learning literature, it has gained growing popularity in both statistics and machine learning and found numerous real-world applications,  ranging from biological data integration [\cite{borgwardt2006integrating}], to neural networks training [\cite{dziugaite2015training}], to the evaluation of a generative model in generative adversarial networks (GAN) [\cite{arbel2018gradient}, \cite{binkowski2018demystifying}].

As a distance metric that measures the closeness of two distributions, MMD belongs to the category of interpoint distance based metric. In this category, a notable member is energy distance (ED, hereafter) [\cite{szekely2004testing}, \cite{szekely2013energy}], which can be viewed as a special case of MMD [\cite{sejdinovic2013equivalence}]. 
ED has been applied to many statistical problems, including two sample testing [\cite{szekely2004testing}, \cite{zhu2021interpoint}], change-point detection [\cite{matteson2014nonparametric}], hierarchical clustering [\cite{szekely2005hierarchical}], assessment of the quality of probabilistic forecasts via new scoring rules [\cite{gneiting2007strictly}], and covariate balancing in causal inference [\cite{huling2020energy}].

Motivated by the increasing use of MMD test for data of moderate and high dimension [\cite{borgwardt2006integrating},  \cite{zhu2017maximum}, \cite{zhao2019identification}], we propose to study the behavior of sample MMD in the high-dimensional setting, which seems relatively less explored. To the best of our knowledge, we are only aware of recent contributions from \cite{zhu2021interpoint} and \cite{chakraborty2021new}.  
In \cite{zhu2021interpoint}, they showed that under the setting  $p\gg\max(n,m)$, the MMD permutation tests are inconsistent when the two high dimensional distributions correspond to the same marginal distributions but differ in other aspects of the distributions in that the ED and MMD tests mainly target the differences between marginal means and sum of componentwise variances; see \cite{chakraborty2021new} for similar findings. Note that the computational complexity of MMD permutation test is $O((n+m)^2pB)$ with $B$ being the number of permutations employed and the computational cost is expensive for large scale data, whereas that of our proposed method is $O((n+m)^2p)$.

As close relatives of ED, distance covariance (dcov, hereafter) and its standardized version distance correlation (dcor, hereafter) were proposed by \cite{szekely2007measuring} to measure the dependence between two random vectors $X\in \br^p$ and $Y\in \br^q$ of arbitrary dimensions. The high-dimensional behavior of sample dcov has been studied in \cite{zhu2020distance} and \cite{gao2019asymptotic}. In \cite{zhu2020distance}, they showed that under the setting $\min(p,q)\gg n$, the dcov is unable to capture full nonlinear dependence between $X$ and $Y$ and it is only capable of capturing componentwise cross-covariance, a phenomenon reminiscent of the one in \cite{zhu2021interpoint} and \cite{chakraborty2021new}. Additionally, their results have been shown to hold for sample HSIC (Hilbert-Schmidt Independence Criterion), which can be viewed as a {\color{black} kernelized} version of sample dcov; see \cite{sejdinovic2013equivalence}. On the other hand, \cite{gao2019asymptotic} showed that a rescaled sample dcor is capable of detecting full nonlinear dependence as long as $p=q=o(\sqrt{n})$ and other regularity conditions hold. Thus the results in \cite{zhu2020distance} and \cite{gao2019asymptotic}  complement each other and suggest that there are several interesting regimes for the asymptotic behavior of sample dcov and sample dcor. \cite{han2021generalized} derived the first non-null central limit theorem {\color{black} (CLT, hereafter)} for the sample distance covariance, as well as the more general sample HSIC in high dimensions, and their results were obtained primarily in the Gaussian case. 


Despite the aforementioned recent advances, the asymptotic theory for sample MMD under the null hypothesis $H_0$ in general case of $n,m$ and $p$  diverging in an arbitrary fashion remains unexplored. Our first main contribution is to obtain central limit theorems for a studentized sample MMD. We also obtain the explicit rates of convergence to the limiting standard normal distribution. As another important contribution, we provide a general theory for the power analysis for our studentized sample MMD and provide several non-overlapping cases to discuss when the power of our MMD test is asymptotically one. One of the main findings is that in the moderately high-dimensional regime, the proposed studentized test statistic is able to detect the difference between two distributions with high power. The difference can lie in the means, marginal variances, componentwise covariances, and higher-order features associated with two high-dimensional distributions. The theoretical results are new to the literature and can be considered as substantial extensions over those obtained in 
{\color{black} \cite{zhu2020distance},   \cite{zhu2021interpoint}, and
\cite{gao2019asymptotic}. As compared to \cite{zhu2021interpoint}, who focused on the behavior of MMD-based permutation test in both High-dimensional Low Sample Size (HDLSS) and High-dimensional Medium Sample Size (HDMSS) settings, we aim to derive a simple studentized test statistic with standard normal limiting null distribution, under less stringent restrictions on the growth rate of $p$ as a function of $n$. Some detailed comparisons with their power results are deferred to Section~\ref{subsec:example}. 

As two sample testing and independence testing are very much related, our work is also inspired by the dcov-based testing in high dimensional setting in \cite{zhu2020distance} and \cite{gao2019asymptotic}. In particular, since our work and \cite{gao2019asymptotic} share some technical arguments (say, Berry-Esseen bound for martingale), it pays to highlight the main difference  between these two papers. First, the main U-statistic ({\color{black} that is,} sample dcov) in  \cite{gao2019asymptotic} is based on a one-sample kernel of order four, whereas we need to deal with a  two-sample kernel of order $(2,2)$.  Consequently, some new theoretical tools need to be developed, such as the moment inequality for the two-sample U-statistic. Second, to form the studentized test statistic, we estimate the variance of sample MMD under the null using the pooled sample. The asymptotic behavior of this variance estimate is studied under both the null and the alternative. In particular, we have shown that it is a ratio-consistent estimator of  HSIC of a mixture distribution with itself under some mild conditions.
Lastly, our asymptotic theory is developed for a large class of kernels,  including the $L_2$ norm as well as the Gaussian kernel, the Laplacian kernel, and many other kernels used in the machine learning literature. This generality is achieved by substantial new technical developments and very involved asymptotic analysis. 
}

Recently, \cite{yan2023kernel} have obtained some related results for MMD-based test in high dimension. Specifically, they propose a studentized MMD-based test statistic under a specific model structure and establish the null CLT as well as the non-null CLTs under fixed and local alternatives for an (infeasible) standardized statistic. Though both papers consider the two-sample {\color{black} MMD-based} testing problem when both $(n,m)$ and $p$ diverge and propose a studentized  statistic, there are significant differences in terms of settings, technical tools and theoretical results. Firstly, the problem set-ups are different. \cite{yan2023kernel} consider a special factor-like model which has been adopted in high-dimensional two sample mean testing [\cite{chen2010two}]. All of our theory, including the CLT, the general Berry-Esseen bound and the power results, are established with no specific model constraints, and are thus applicable to a broader set of data generating processes. 
Secondly, the technical tools and primary results established in the two papers are very different. The most striking contribution in \cite{yan2023kernel} is the non-null CLTs for the standardized statistic, which are very interesting and seem only achievable under the specific model assumption,  whereas we only present the null CLT for our studentized test statistic but additionally derive a Berry-Esseen bound under the null.
Thirdly, the power results and the regimes under which the power approaches one are very different, and more discussion can be found in Section~\ref{subsec:power}. Overall, we view the results in  \cite{yan2023kernel} and our paper complementary to each other. Together they provide a more complete {\color{black} portrayal} of the high-dimensional behavior of MMD-based statistics.

The rest of this paper is organized as follows. Section \ref{Sec:MMD} introduces the maximum mean discrepancy, its sample version as a two sample U-statistic and its Hoeffeding decomposition. The distributional properties when the dimension $p$ is fixed is  also  described. We propose a studentized test statistic and present the main theorems in Section \ref{Sec:AsyDist}. To be specific, we present the CLT for the studentized MMD and obtain the rates of convergence under the null. We also provide a general theory for the power under the alternative in this section. Finite sample performance is examined via simulations in  Section \ref{Sec:NumStudies}. In Section \ref{Sec:Discussion}, we summarize our results and discuss some potential extensions. {\color{black} Some illustrative examples, all the technical details, and some additional simulation results are presented in the online appendices; see \url{https://arxiv.org/abs/2109.14913}.}

Let $c,d$ be any positive integers and $\phi(x_1,\dots,x_c,y_1,\dots,y_d)$ denote a two-sample kernel function. For any $0 \le c' \le c$, $0 \le d' \le d$, and subsets $\{i_1,\dots,i_{c'}\}$, $\{j_1,\dots, j_{d'}\}$, define 
\begin{equation*}
    \blreft{X_{i_1},\dots,X_{i_{c'}},Y_{j_1},\dots,Y_{j_{d'}}}{\phi(X_1,\dots,X_c,Y_1,\dots,Y_d)} = \displaystyle{\int\dots\int \phi(X_1,\dots,X_c,Y_1,\dots,Y_d) \prod\limits_{s=1}^{c'}dF_{X_{i_s}} \prod\limits_{r=1}^{d'}dF_{Y_{j_r}}}.
\end{equation*}
For simplicity, we write
\begin{equation*}
    \blreft{X_1,\dots,X_c,Y_1,\dots,Y_d}{\phi(X_1,\dots,X_c,Y_1,\dots,Y_d)} 
  = \blre{\phi(X_1,\dots,X_c,Y_1,\dots,Y_d)}.
\end{equation*}
For two random vectors $V_1,V_2$, the notation $V_1 =^d V_2$ means that they are identically distributed. We use $\rightarrow^d$ and $\rightarrow^p$ to denote convergence in distribution and in probability respectively. For two real-valued sequences $a_n,b_n$, we say $a_n = O(b_n)$ or $a_n \lesssim b_n$ if there exist $M,C>0$, such that $a_n \le C b_n$ for $n>M$. If there exist $M,C_1,C_2>0$, such that $C_1 b_n \le a_n \le C_2 b_n$ for $n>M$, then we say $a_n = O_s(b_n)$. In addition, we say $a_n = o(b_n)$ { or $a_n \prec b_n$} if $a_n/b_n \rightarrow 0$ as $n\rightarrow\infty$. For any $p$-dimensional vectors $a,b$, we use $|a-b|$ to denote the Euclidean distance between $a$ and $b$. For a function $f$, we use $f_i$ to denote its $i$-th order derivative, and $\tilde{f}$ to denote its centered version, that is, $\tilde{f}(V_1,\dots,V_k) = f(V_1,\dots,V_k)-\bbe{f(V_1,\dots,V_k)}$. We use $C(u_1,\dots,u_k)$ to denote a positive and finite constant that depends only on the parameters $u_1,\dots,u_k$ and the values of $C(u_1,\dots,u_k)$ may vary from line to line. Additionally, we use $\cum(x_1,\cdots,x_k)$ to denote the joint  cumulant of the random variables $x_1,\cdots,x_k$.


\section{Maximum Mean Discrepancy and Its Properties}
\label{Sec:MMD}

\subsection{The Definition of Maximum Mean Discrepancy}

{\color{black} We follow Definition 2 in \cite{gretton2012kernel} to provide a formal definition of MMD.
\begin{definition}\label{Def:MMD-original}
    Let $X\sim P_1$ and $Y\sim P_2$ be independent random vectors in $\br^p$ and let $\cf_0$ be a class of functions $f: \br^p \rightarrow \br$. We define the maximum mean discrepancy (MMD) as
    \begin{equation}\label{Equ:MMD-original}
        \MMD(P_1,P_2) = \sup\limits_{f\in\cf_0} \lrcp{\bbe{f(X)} - \bbe{f(Y)}}.
    \end{equation}
\end{definition}
With properly selected function class $\cf_0$, $\MMD(P_1,P_2)$ has some special properties. To facilitate the subsequent discussion, we follow the introduction in Section 2.2 of \cite{gretton2012kernel} to provide some basic properties of the reproducing kernel Hilbert space (RKHS). 

Specifically, let $\cf$ be an RKHS on the separable metric space $(\br^{p},\cp)$, where $\cp$ denotes the set of Borel probability measures on $\br^p$. By the property of the RKHS and the Riesz representation theorem, there is a feature mapping $\phi(x): \br^p \rightarrow \br$ such that $f(x) = \lrag{f,\phi(x)}_\cf$ for any $f\in\cf$. Furthermore, there exists a symmetric and positive definite kernel $\bar{k}$ associated with $\cf$ such that $\phi$ takes the canonical form $\phi(x) = \bar{k}(x,\cdot)$. It follows that $\bar{k}(x,x')=\lrag{\phi(x),\phi(x')}_\cf$ for any $x,x'\in\br^p$.

For any distribution $P\in\cp$, we define the mean embedding of $\mu_P\in\cf$ as the function satisfying that $\be[f(x)] = \lrag{f,\mu_P}_F$ for any $f\in\cf$. It is shown in Lemma 3 and Lemma 4 of \cite{gretton2012kernel} that, when the aforementioned kernel $\bar{k}$ is measurable and satisfies $\be[\sqrt{\bar{k}(X,X)}] < \infty$, $\be[\sqrt{\bar{k}(Y,Y)}] < \infty$, then MMD can be expressed as the distance in $\cf$ between mean embeddings, that is, $\MMD^2(P_1,P_2):=\|\mu_{P_1}-\mu_{P_2}\|_{\cf}^2$.} Equivalently, as stated in Lemma 6 of \cite{gretton2012kernel},  $\MMD$ can be expressed through the kernel as
\begin{equation}\label{Equ:MMD}
    \MMD(P_1,P_2) = \MMD(X,Y|\bar{k}):=\lrp{-2\blre{\bar{k}(X,Y)} + \blre{\bar{k}(X,X')} + \blre{\bar{k}(Y,Y')}}^{1/2},
\end{equation}
where $X',Y'$ are independent and identical copies of $X \sim P_1$ and $Y \sim P_2$, respectively. 

{\color{black} When $\cf_0$ in Equation (\ref{Equ:MMD-original}) is the unit ball in the RKHS $(\cf,\bar{k})$, \cite{gretton2012kernel} has shown that $\MMD(P_1,P_2)$ is a nonnegative metric and $\MMD(P_1,P_2)=0$ if and only if $P_1 =^d P_2$. Similar results have been generalized by using the equivalent definition of MMD. In particular,} if $\bar{k}$ {\color{black} in Equation (\ref{Equ:MMD})} is characteristic on $\br^p$ (i.e., the corresponding mean map $\mu_{P}$ is injective), then the associated MMD is a metric on $\cp$, which satisfies $\MMD(P_1,P_2)=0$ if and only if $P_1=P_2$ [\cite{fukumizu2007kernel}, \cite{sejdinovic2013equivalence}]. Many commonly used kernels are shown to be characteristic kernels on $\br^p$, including the Gaussian kernel and Laplacian kernel [\cite{fukumizu2007kernel}]. 

{\color{black} We note that when $k(x,y)=|x-y|$, Equation (\ref{Equ:MMD}) coincides with the formulation of ED (\cite{szekely2004testing}).
\begin{definition}
    Let $X,X',Y,Y'$ be independence random vectors in $\br^p$ that satisfies $X,X'\sim P_1$ and $Y,Y'\sim P_2$, we define the energy distance (ED) as
    \begin{equation}\label{Equ:ED}
        \mbox{ED}(P_1,P_2) = \mbox{ED}(X, Y) = \lrp{2\bbe{|X-Y|} - \bbe{|X-X'|} - \bbe{|Y-Y'|}}^{1/2}.
    \end{equation}
\end{definition}
ED is a nonnegative metric and it holds that $ED(P_1,P_2)=0$ if any only if $P_1=P_2$.}

In this paper, we aim to provide a unified treatment of ED and MMD, so we follow the approach in \cite{zhu2021interpoint} and mimic the definition of energy distance in \cite{szekely2013energy} and \cite{huang2017efficient} to give the definition of {\color{black} MMD} with a general kernel $k$. Though kernel is commonly used to measure similarity in the machine learning literature,  we use kernel throughout this article to refer to a broader range of metrics of dissimilarities, which include both a semimetric of strong negative type on $\br^{p}$ (Definition 1 and Definition 28 of \cite{sejdinovic2013equivalence}) and a characteristic kernel multiplied by $-1$, as formally stated in Definition \ref{Def:kernel} below. {\color{black} For notational simplicity, we shall use MMD (ED) instead of MMD$^2$ (ED$^2$) as in their original definitions (\ref{Equ:MMD}) and (\ref{Equ:ED}), and the same is done for dcov later.}

\begin{definition}\label{Def:kernel}
  Define $k:\br^p\times \br^p\rightarrow [0,\infty)$ to be a kernel that satisfies either of the following conditions:
  \begin{enumerate}[label=(\roman*)]
      \item \label{Def:kernel-1}
      for any $x,y\in\br^{p}$, it holds that $k(x,y)=k(y,x)$ and $k(x,y)=0$ if and only if $x=y$, and additionally, for any Borel probability measures $P,Q$ on $\br^{p}$ satisfying {\color{black} $\int k(z,z) dP(z)<\infty$ and $\int k(z,z) dQ(z) < \infty$}, $P\neq Q$ implies that $\int k d([P-Q]\times[P-Q])<0$.
      \item  \label{Def:kernel-2}
      $(\cf,-k)$ is an RKHS on $(\br^{p},\cp)$ and the kernel $-k$ is characteristic.
  \end{enumerate}
\end{definition}
{\color{black} To ease the reading, we mimic Table 1 of \cite{zhu2021interpoint} to summarize a few kernels covered by Definition \ref{Def:kernel} in the following table.
\begin{table}[!h]
    \centering
    \scalebox{0.9}{
    \begin{tabular}{c|c|c}
    \hline\hline
        Kernel $k$ & Expression of $k$ & Condition satisfied \\
        \hline
        Euclidean distance & $k(x,y)=|x-y|$ & $k$ satisfies \ref{Def:kernel-1} \\
        Gaussian kernel (multiplied by -1) & $k(x,y)=-\exp\lrp{-|x-y|^2/(2\gamma^2)}$ & $k$ satisfies \ref{Def:kernel-2} \\
        Laplacian kernel (multiplied by -1) & $k(x,y)=-\exp\lrp{-|x-y|/\gamma}$ & $k$ satisfies \ref{Def:kernel-2} \\
    \hline\hline
    \end{tabular}}
    \caption{Examples of kernel $k$ covered by Definition \ref{Def:kernel}.}
    \label{Tab:Example}
\end{table}}

{\color{black} Then we are ready to propose the unified definition of ED and MMD.}
\begin{definition}\label{Def:MMD}
  {\color{black} Let $k$ denote a kernel defined as Definition~\ref{Def:kernel}. Suppose that $X,Y\in\br^p$ are two independent random vectors satisfying that $\blre{\lrabs{k(X,X')}} + \blre{\lrabs{k(X,Y)}} + \blre{\lrabs{k(Y,Y')}}<\infty$, then we define}
  \begin{equation}
      \ce^k(X,Y) 
    = 2\blre{k(X,Y)} - \blre{k(X,X')} - \blre{k(Y,Y')},
  \end{equation}
  where $X',Y'$ are independent and identical copies of $X$ and $Y$, respectively.
\end{definition}

As shown in \cite{sejdinovic2013equivalence}, $\ce^k(X,Y)$ is always non-negative and is zero if and only if $X =^d Y$. Similar to \cite{huang2017efficient}, who expressed sample {\color{black} ED} as a U-statistic, we can find an unbiased estimator of $\ce^k(X,Y)$  via a U-statistic with a two-sample kernel.
\begin{proposition}\label{Prop:Enm-k}
  Define the two-sample kernel
  \begin{equation}\label{Equ:h^k}
      h^k(X_1,X_2,Y_1,Y_2) 
    = \frac{1}{2} \sum_{i=1}^2 \sum_{j=1}^2 k(X_i, Y_j) 
    - k(X_1, X_2) 
    - k(Y_1, Y_2),
  \end{equation}
  which satisfies $\blre{h^k(X_1,X_2,Y_1,Y_2)}=\ce^k(X,Y)$. Then an unbiased estimator of $\ce^k(X,Y)$ can be defined as
  \BEqn
    \ce_{n,m}^k(X,Y)
  &=& \bin{n}{2}^{-1} \bin{m}{2}^{-1} \sum\limits_{1 \le i_1 < i_2 \le n} \sum\limits_{1 \le j_1 < j_2 \le m} h^k(X_{i_1},X_{i_2},Y_{j_1},Y_{j_2}) \\
  &=& \frac{2}{mn} \sum\limits_{i=1}^{n} \sum\limits_{j=1}^{m} k(X_i, Y_j)
     - \bin{n}{2}^{-1} \sum\limits_{1 \leq i < j \leq n} k(X_i, X_j) 
     - \bin{m}{2}^{-1} \sum\limits_{1 \leq i < j \leq m} k(Y_i, Y_j).
  \normalsize
  \EEqn
\end{proposition}

It follows from the Hoeffding decomposition that two-sample U-statistic $\ce_{n,m}^{k}(X,Y)$ can be decomposed into the sum of a leading term and a remainder term.
{\color{black} In particular, let $G_x$ denote the distribution function of a single point mass at $x$, and for $0\le c,d \le 2$, define 
\begin{align*}
     h^{(c,d)}(X_1,\dots,X_c;Y_1,\dots,Y_d)  
  =& \int\dots\int h^{k}(u_1,u_2,v_1,v_2) \prod\limits_{i=1}^{c}\lrp{dG_{X_i}(u_i) - dF_X(u_i)} \prod\limits_{i=c+1}^{2}dF_X(u_i) \\
   &\quad \times \prod\limits_{j=1}^{d}\lrp{dG_{Y_j}(v_j) - dF_Y(v_j)} \prod\limits_{j=d+1}^{2}dF_Y(v_j).
\end{align*}
Then it holds that $\ce_{n,m}^{k}(X,Y) = L_{n,m}^{k}(X,Y) + R_{n,m}^{k}(X,Y)$, where
\begin{align*}
    L_{n,m}^k(X,Y)
  =& \ce^k(X,Y) 
     + \frac{2}{n}\sum\limits_{i=1}^n h^{(1,0)}(X_i) 
     + \frac{2}{m}\sum\limits_{j=1}^{m} h^{(0,1)}(Y_j) \\
   & + \frac{2}{n(n-1)} \sum\limits_{1 \le i_1 < i_2 \le n} h^{(2,0)}(X_{i_1},X_{i_2})
     + \frac{4}{nm} \sumin \sumjm h^{(1,1)}(X_i,Y_j) \nonumber \\
   & + \frac{2}{m(m-1)} \sum\limits_{1 \le j_1 < j_2 \le m} h^{(0,2)}(Y_{j_1},Y_{j_2}), \\
    R_{n,m}^k(X,Y)
  =& \frac{4}{n(n-1)m} \sum\limits_{1 \le i_1 < i_2 \le n} \sumjm h^{(2,1)}(X_{i_1},X_{i_2},Y_j) \\
  & + \frac{4}{nm(m-1)} \sumin \sum\limits_{1 \le j_1 < j_2 \le m} h^{(1,2)}(X_i,Y_{j_1},Y_{j_2}) \\
  & + \frac{4}{n(n-1)m(m-1)} \sum\limits_{1 \le i_1 < i_2 \le n} \sum\limits_{1 \le j_1 < j_2 \le m} h^{(2,2)}(X_{i_1},X_{i_2},Y_{j_1},Y_{j_2}),
\end{align*}}
By generalizing some results established in \cite{huang2017efficient}, the expressions of $L_{n,m}^{k}(X,Y)$ and $R_{n,m}^{k}(X,Y)$ can be greatly simplified, which are stated in Proposition \ref{Prop:est_MMD_simplified}.

\begin{proposition}\label{Prop:est_MMD_simplified}
Assume that $\blre{\lrabs{k(X,X')}} + \blre{\lrabs{k(X,Y)}} + \blre{\lrabs{k(Y,Y')}}<\infty$, then it holds that $R_{n,m}^{k}(X,Y) = 0$, and
\begin{equation*}
\begin{array}{lcl}
L_{n,m}^{k}(X,Y) & = &   3\ce^k(X,Y) 
                       - \frac{4}{n}\sum\limits_{i=1}^n h_{10}(X_i) 
                       - \frac{4}{m}\sum\limits_{j=1}^{m} h_{01}(Y_j) 
                       + \frac{4}{nm} \sumin \sumjm h_{11}(X_i,Y_j) \vspace{1mm}\\
                 &   & + \bin{n}{2}^{-1} \sum\limits_{1 \le i_1 < i_2 \le n} h_{20}(X_{i_1},X_{i_2})
                       + \bin{m}{2}^{-1} \sum\limits_{1 \le j_1 < j_2 \le m} h_{02}(Y_{j_1},Y_{j_2}),
\end{array}
\end{equation*}
\normalsize
{\color{black}
where for $1 \le i_1 < i_2 \le n$ and $1 \le j_1 < j_2 \le m$, we define
\BEqn
&  h_{10}(X_{i_1}) = \beft{X_{i_2},Y_{j_1},Y_{j_2}}{h^k(X_{i_1},X_{i_2},Y_{j_1},Y_{j_2})}, 
&  h_{01}(Y_{j_1}) = \beft{X_{i_1},X_{i_2},Y_{j_2}}{h^k(X_{i_1},X_{i_2},Y_{j_1},Y_{j_2})}, \\
&  h_{20}(X_{i_1},X_{i_2}) = \beft{Y_{j_1},Y_{j_2}}{h^k(X_{i_1},X_{i_2},Y_{j_1},Y_{j_2})}, 
&  h_{02}(Y_{j_1},Y_{j_2}) = \beft{X_{i_1},X_{i_2}}{h^k(X_{i_1},X_{i_2},Y_{j_1},Y_{j_2})},
\EEqn
and $h_{11}(X_{i_1},Y_{j_1}) = \beft{X_{i_2},Y_{j_2}}{h^k(X_{i_1},X_{i_2},Y_{j_1},Y_{j_2})}$.}

If additionally $X =^d Y$, then $\ce_{n,m}^{k}(X,Y) = L_{n,m}^{k}(X,Y)$ with {\color{black} the simplified expression}
\begin{equation*}
\hspace*{-1em}
\begin{array}{l}
    L_{n,m}^k(X,Y)
=   \bin{n}{2}^{-1} \sum\limits_{1 \le i_1 < i_2 \le n} h_{20}(X_{i_1},X_{i_2})
  + \frac{4}{nm} \sumin \sumjm h_{11}(X_i,Y_j)
  + \bin{m}{2}^{-1} \sum\limits_{1 \le j_1 < j_2 \le m} h_{02}(Y_{j_1},Y_{j_2}).
\end{array}
\end{equation*}
\end{proposition}

Proposition \ref{Prop:est_MMD_simplified} simplifies $\ce_{n,m}^{k}(X,Y)$ and facilitates the subsequent analysis.

\subsection{Distributional Properties of $\ce_{n,m}^k(X,Y)$ with Fixed $p$}

The asymptotic behavior of sample {\color{black} MMD (and sample ED)} has been well studied when the dimension $p$ is fixed; see \cite{szekely2004testing}, \cite{gretton2009fast} and \cite{gretton2012kernel}. In particular, the asymptotic distribution of $\ce_{n,m}^k(X,Y)$ (in the case of sample  {\color{black} MMD}) under the null is established in Theorem 12 of \cite{gretton2012kernel}. 
Define
\begin{equation}\label{Equ:d^k}
    d^k(X_1,X_2) 
  = k(X_1,X_2) - \blreft{X_1}{k(X_1,X_2)} - \blreft{X_2}{k(X_1,X_2)} + \blre{k(X_1,X_2)}.
\end{equation}
Assume that $\blre{k^2(X,X')} < \infty$, and that $\lim\limits_{n,m\rightarrow\infty}\frac{n}{n+m}=\rho$ for some fixed $0<\rho<1$, then under the null, $\ce_{n,m}^k(X,Y)$ converges in distribution according to
\begin{equation*}
    (n+m)\ce_{n,m}^k(X,Y)
    \rightarrow^{d}
    \sum\limits_{\ell=1}^{\infty} \lambda_\ell \lrp{\lrp{\rho(1-\rho)}^{-1} - \lrp{\rho^{-1/2}a_\ell - (1-\rho)^{-1/2}b_\ell}^2},
\end{equation*}
where
{\color{black} $\{a_\ell\}_{\ell\geq1}, \{b_{\ell}\}_{\ell\geq1} \sim \cn(0,1)$ are two independent sequences of iid Gaussian random variables, and $\{\lambda_i\}_{i=1}^{\infty}$ and $\{\Psi_i(x)\}_{i=1}^{\infty}$ are respectively the eigenvalues and the eigenfunctions of the equation $\blreft{X}{d^k(X,X')\Psi_i(X)} = \lambda_i \Psi_i(X')$.}
Note that the above limiting null distribution is not pivotal, so critical values are not directly available. Several approximation methods have been developed in the special case of $m=n$ in \cite{gretton2009fast}.


\section{Studentized Statistic and Asymptotic Theory}
\label{Sec:AsyDist}

\subsection{Studentized Statistic}

To develop our studentized statistic, we need to find the variance of $L_{n,m}^k$ under the null, which can be shown to have a strong connection with the Hilbert-Schmidt independence criterion (HSIC) [\cite{gretton2007kernel}]. We thus start from the definition and some basic results of HSIC, and then move on to deriving the variance of $L_{n,m}^k$ under the null before proposing a studentized test statistic.

\begin{definition}\label{Def:DC-k}
  Suppose $X\in\br^p$ and $Y\in\br^q$ with integers $p,q\geq1$. Let $k$ be a kernel defined as Definition \ref{Def:kernel}, then a generalized Hilbert-Schmidt independence criterion (HSIC) $\cv_k^2(X,Y)$ between the distributions $X$ and $Y$ is defined as
  \begin{equation}\label{Equ:DC-k}
      \cv_k^2(X,Y) 
    = \bbe{k(X_1,X_2) k(Y_1,Y_2)} - 2\bbe{k(X_1,X_2) k(Y_1,Y_3)} + \bbe{k(X_1,X_2)} \bbe{k(Y_1,Y_2)}, 
  \end{equation}
  where $(X_1,Y_1)$, $(X_2,Y_2)$ and $(X_3,Y_3)$ are independent copies of $(X,Y)$.
\end{definition}

\begin{remark}
This definition directly follows from Lemma 1 of \cite{gretton2005measuring} when $-k(\cdot,\cdot)$ is a { characteristic kernel w.r.t. $(\br^{p},\cp)$}. When $k$ is a semimetric of strong negative type according to Definition \ref{Def:kernel}, the resulting metric also {\color{black} fully quantifies nonlinear dependence in the sense that $\cv_k^2(X,Y) \geq 0$ and equals to zero if and only if $X, Y$ are independent; see \cite{sejdinovic2013equivalence}}. 
{\color{black} In fact, $\cv_k^2(X,Y)$ can be viewed as a generalization of distance covariance (dcov), as the expression of $\cv_k^2(X,Y)$ with $k(x,y)=|x-y|$ coincides with the well-known dcov [\cite{szekely2007measuring}].}
For simplicity, we call $\cv_k^2(X,Y)$ the HSIC in both cases.
\end{remark}

Next, we introduce a mixture distribution of $X$ and $Y$, defined by 
\begin{equation}\label{Equ:mixZ}
    Z 
  = \left\{
      \begin{array}{ll}
         X, & \text{with probability } \rho, \\
         Y, & \text{with probability } 1-\rho,
      \end{array}
    \right.
\end{equation}
{\color{black} where $\rho$ is the limit of the size proportion $n/(n+m)$. By introducing the mixture distribution $Z$, we are able to aggregate $X$ and $Y$ with respect to their occurrence frequencies in the pooled sample. Consequently,  we can directly establish some unified results in terms of $Z$, which are more succinct than establishing the counterparts in terms of $X$ and $Y$ respectively.}

Let $Z_1,Z_2,Z_3$ be three independent copies of $Z$, then HSIC of $Z$ with itself is given by
\begin{equation}\label{Equ:Vk}
    \cv_k^2(Z)
  = \blre{k^2(Z_1,Z_2)} - 2\blre{k(Z_1,Z_2) k(Z_1,Z_3)} + \lrp{\blre{k(Z_1,Z_2)}}^2.
\end{equation}
{\color{black} With $\cv_k^2(X,Y)$ being a generalization of distance covariance, we can also view $\cv_k^2(Z)$ as a generalization of distance variance of $Z$.} It is trivial that $X =^d Y =^d Z$ under the null. Furthermore, the variance of $L_{n,m}^k$ under the null can be written in terms of $\cv_k^2(Z)$.

\begin{proposition}\label{Prop:Var(Lnm)}
  If $X$ and $Y$ are identically distributed, and $\blre{k^2(X_1,X_2)}< \infty$,
  then it holds that $\Var(L_{n,m}^k(X,Y)) = c_{n,m} \cv_k^2(Z)$, where $c_{n,m} = \frac{2}{n(n-1)} + \frac{4}{nm} + \frac{2}{m(m-1)}$.
\end{proposition}

If given $n$ independent and identically distributed observations $\fZ = (Z_1,\dots,Z_n)$ from the mixture distribution $Z$, an unbiased estimator of $\cv_k^2(Z)$ with $k(x,y)=|x-y|$ can be obtained using the  $U$-centering approach in \cite{szekely2013distance} and \cite{szekely2014partial}. However, the mixture distribution $Z$ is unobserved, and we only have two independent random samples $\fX =  (X_1,\dots,X_n)$ from the distribution of $X$ and $\fY = (Y_1,\dots,Y_m)$ from the distribution of $Y$. Let $N = n+m$ denote the total sample size. Throughout, we assume that there exists some constant $0 < \rho < 1$, such that   $n/N \rightarrow \rho$ as $\min\{n,m\}\rightarrow\infty$.
We propose to use the pooled sample to estimate  $\cv_k^2(Z)$  as follows.

\begin{proposition}\label{Prop:estDC_k}
  { For any fixed $p$ and kernel $k:\br^p \times \br^p \mapsto \br$}, assume that $\blre{k^2(X_1,X_2)}$, $\blre{k^2(X_1,Y_1)}$ and $\blre{k^2(Y_1,Y_2)}$ are all finite, and $k(X,X) = a_0^k$ is a finite constant independent of $X$. For $1\le s,t\le N$, define
  \begin{equation}
      a_{s,t}^{k}
    = \left\{
        \begin{array}{ll}
          k(X_s,X_t),        & 1 \le s,t \le n \vspace{3pt}\\
          k(X_s,Y_{t-n}),    & 1 \le s \le n < t \le N \vspace{3pt}\\
          k(X_t,Y_{s-n}),    & 1 \le t \le n < s \le N \vspace{3pt}\\
          k(Y_{s-n},Y_{t-n}), & n+1 \le s,t \le N
        \end{array}
      \right.
  \end{equation}
  Define the $\cu$-centered distances with kernel $k$ as
  $A_{s,t}^{k\ast} = a_{s,t}^{k} - \tilde{a}_{\cdot t}^{k} - \tilde{a}_{s \cdot}^{k} + \tilde{a}_{\cdot \cdot}^{k}$,  where
  \begin{equation*}
      \tilde{a}_{\cdot t}^{k}
    = \frac{1}{N-2} \sum\limits_{i=1}^{N} a_{i,t}^{k},\quad
      \tilde{a}_{s \cdot}^{k} 
    = \frac{1}{N-2} \sum\limits_{j=1}^{N} a_{s,j}^{k},\quad
      \tilde{a}_{\cdot \cdot}^{k}
    = \frac{1}{(N-1)(N-2)} \sum\limits_{i,j=1}^{N} a_{i,j}^{k}.
  \end{equation*}
  Then under the null, it holds for any fixed $p$ and kernel $k$ that,
  \begin{equation}
    \cv_{n,m}^{k\ast}(X,Y) 
  = \frac{1}{N(N-3)} \sum\limits_{s \neq t} \lrp{A_{s,t}^{k\ast}}^2
  - \frac{(a_0^k)^2}{(N-1)(N-3)}
  \end{equation}
  is an unbiased estimator of $\cv_k^2(Z)$. Furthermore, under the alternative, $\cv_{n,m}^{k\ast}(X,Y)$ is asymptotically unbiased of $\cv_k^2(Z)$ for any fixed $p$ and kernel $k$, that is, $\blre{\cv_{n,m}^{k\ast}(X,Y)} \rightarrow \cv_k^2(Z)$ as $n,m\rightarrow\infty$.
\end{proposition}

\begin{remark}
  If the kernel $k$ is chosen to be the $L_2$ norm, we have $a_0^k = 0$ and the estimate $\cv_{n,m}^{k\ast}(X,Y)$ reduces to the {\color{black} traditional} $\cu$-centering based sample distance variance based on the pooled sample. However, for a general kernel $k$, $a_0^k$ may be nonzero, and the correction term $-\frac{(a_0^k)^2}{(N-1)(N-3)}$ is necessary to obtain the unbiasedness. This bias correction is important for Gaussian and Laplacian kernels as the use of biased variance estimate leads to noticeable size distortion in the small sample in our (unreported) simulations. 
\end{remark}

{\color{black} To our best knowledge, the proposed estimate of $\cv_k^2(Z)$ based on the pooled sample is a new addition to the literature, and it is different from the studentizers proposed in \cite{chakraborty2021new} and \cite{yan2023kernel}; see Remark \ref{Rmk:CompareTestStat-1} and Remark \ref{Rmk:CompareTestStat-2}.}


In terms of computational complexity, the computation of all $a_{s,t}^{k}$'s is of order $O((n+m)^2p)$. Since the $\cu$-centering only requires $O((n+m)^2)$ computation,  the computational complexity of $\cv_{n,m}^{k\ast}(X,Y)$ and the studentized statistic $T_{n,m,p}^k$ defined below is of order $O((n+m)^2p)$. By contrast, the computational complexity for the permutation based test {\color{black} in \cite{zhu2021interpoint}} is of order $O((n+m)^2pB)$, where $B$ is the number of {\color{black} permutations}.

To test $H_0: X =^d Y$ against $H_a: X \neq^d Y$, it is natural to use the following studentized test statistic:
\begin{equation}
    T_{n,m,p}^k 
  = \frac{\ce_{n,m}^k(X,Y)}{\sqrt{c_{n,m} \cv_{n,m}^{k\ast}(X,Y)}},
\end{equation}
where $c_{n,m}$ is defined in Proposition \ref{Prop:Var(Lnm)}.

{\color{black} 
Similar test statistics for the two-sample problem have been previously discussed in other existing papers; see \cite{chakraborty2021new} and \cite{yan2023kernel}. We conjecture that all three studentizers are asymptotically equivalent. Some additional discussions can be found in the following remarks.}

{\color{black} 
\begin{remark}\label{Rmk:CompareTestStat-1}
In \cite{chakraborty2021new}, a different studentized test statistic is proposed in the form of $\displaystyle{\tilde{T}_{n,m,p}^k=\frac{\ce_{n,m}^k(X, Y)}{\sqrt{c_{n,m}S_{n,m}/2}}}$. The major difference between $T_{n,m,p}^k$ and $\tilde{T}_{n,m,p}^k$ is the variance estimator of $\ce_{n,m}^{k}(X,Y)$ in the denominator. Specifically, 
\begin{equation*}
    S_{n,m}
  = \frac{4(n-1)(m-1)cdCov^2_{n,m}(X,Y)+4v_n \cv_{n}^{k\ast}(X) + 4v_m \cv_{m}^{k\ast}(Y) }{(n-1)(m-1) + n(n-3)/2 + m(m-3)/2},
\end{equation*}
where $\cv_{n}^{k\ast}(X),\cv_{m}^{k\ast}(Y)$ are respectively the $\cu$-centering based unbiased estimators of $\cv_k^2(X),\cv_k^2(Y)$, and $cdCov^2_{n,m}(X,Y)$ is the cross distance covariance between $X$ and $Y$, given by
\begin{equation*}
    cdCov^2_{n,m}(X,Y)=\frac{1}{(n-1)(m-1)}\sum_{k=1}^n\sum_{l=1}^m\hat{k}(X_k,Y_l)^2
\end{equation*}
with $\hat{k}(X_k,Y_l)=k(X_k,Y_l)-\frac{1}{n}\sum_{i=1}^{n}k(X_i,Y_l)-\frac{1}{m}\sum_{j=1}^{m}k(X_k,Y_j)+\frac{1}{nm}\sum_{i=1}^{n}\sum_{j=1}^{m}k(X_i,Y_j)$. In Theorem 4.2 of \cite{chakraborty2021new}, they derived the limiting distribution of $\tilde{T}_{n,m,p}^k$ under both the null and alternatives when $p\rightarrow\infty$ whereas $(n,m)$ are fixed. 
\end{remark}
}


{\color{black} 
\begin{remark}\label{Rmk:CompareTestStat-2}
In a very recent paper by \cite{yan2023kernel}, they also proposed a studentized {\color{black} MMD test}, and their studentizer is based on a linearization argument and differs from ours and the one in \cite{chakraborty2021new}. However, the CLT results in \cite{yan2023kernel} are established for the standardized statistic instead of the studentized statistic under both the null and alternative when both the dimension and sample size diverge, and the standardizer is actually infeasible.
\end{remark}
}


As we present below, we will be investigating the asymptotic behavior of our studentized statistic $ T_{n,m,p}^k $ under the setting $\min(n,m,p)\rightarrow\infty$ using a different set of technical arguments and our results are complementary to those in theirs.

\subsection{Asymptotic Distributions}

For each $p\in\bn$, let $k^{(p)}: \br^p\times\br^p \mapsto \br$ be a kernel defined as in Definition \ref{Def:kernel} and $(k^{(p)}:p\in\bn)$ thus forms a sequence of kernels. Throughout the paper, we let $\cc = \lrcp{(k^{(p)}: p\in\bn)}$ be the set of all the kernel sequences of interest. Again, we drop the symbol $(p)$ for simplicity when we are focusing on a specific kernel given a fixed $p$.

{\color{black} Let $\tilde{k}(Z_1,Z_2)=k(Z_1,Z_2)-\bbe{k(Z_1,Z_2)}$ denote} the centered version of $k$. We have already shown that $\cv_{n,m}^{k\ast}(X,Y)$ is {\color{black} an unbiased estimator of $\cv_k^2(Z)$ under the null and is asymptotically unbiased under the alternative, then} we are ready to state that $\cv_{n,m}^{k\ast}(X,Y)$ is ratio-consistent for $\cv_k^2(Z)$ {\color{black} under both the null and the alternative with some conditions}.

\begin{proposition}\label{Prop:ratio-const-k-H0}
  Assume that $\blre{k^4(Z_1,Z_2)}<\infty$ { for each $k=k^{(p)}\in\cc$} and { $n/N \rightarrow \rho$ as $n,m\rightarrow\infty$ for some $0 < \rho < 1$.} Under the null {\color{black} when $X =^d Y =^d Z$}, if for some constant $0 < \tau \le 1$, it is satisfied when $N,p\rightarrow\infty$ that
  \begin{equation}\label{Equ:consistency_H0}
      \frac{\blre{\lrabs{\tilde{k}(Z_1,Z_2)}^{2+2\tau}}}{N^{\tau}\lrp{\cv_k^2(Z)}^{1+\tau}} \rightarrow 0,
  \end{equation}
  then we have $\displaystyle{\frac{\cv_{n,m}^{k\ast}(X,Y)}{\cv_k^2(Z)} \longrightarrow^{p} 1}$.
  As a direct consequence, under the null we have that $\displaystyle{\frac{c_{n,m} \cv_{n,m}^{k\ast}(X,Y)}{\Var(L_{n,m}^k(X,Y))} \rightarrow^p 1}$.
\end{proposition}

Condition (\ref{Equ:consistency_H0}) is directly from the use of Markov's inequality. In fact, to show that $\cv_{n,m}^{k\ast}(X,Y)$ is ratio-consistent for $\cv_k^2(Z)$, it suffices to find an upper bound of $\frac{\blre{\lrabs{\cv_{n,m}^{k\ast}(X,Y) - \cv_k^2(Z)}^{1+\tau}}}{\lrp{\cv_k^2(Z)}^{1+\tau}}$, where $\tau\in (0,1]$. As shown in Proposition \ref{Prop:estDC_k}, $\cv_{n,m}^{k\ast}(X,Y)$ is unbiased of $\cv_k^2(Z)$ under the null. Furthermore, $\cv_{n,m}^{k\ast}(X,Y)-\cv_k^2(Z)$ can be decomposed as a summation of multiple U-statistics with mean zero. By applying a moment inequality for the U-statistics, we can show that the deviation $\bbe{|\cv_{n,m}^{k\ast}(X,Y) - \cv_k^2(Z)|^{1+\tau}}$ is bounded by $\bbe{|\tilde{k}(Z_1,Z_2)|^{2+2\tau}}/N^{\tau}$ from above. 
This leads to condition (\ref{Equ:consistency_H0}).

As a counterpart of Proposition \ref{Prop:ratio-const-k-H0}, the ratio-consistency of the sample estimate $\cv_{n,m}^{k\ast}(X,Y)$ under the alternative is established in Proposition \ref{Prop:ratio-const-k-H1}.

\begin{proposition}\label{Prop:ratio-const-k-H1}
  Assume that $\blre{k^4(Z_1,Z_2)}<\infty$ for each $k=k^{(p)}\in\cc$ and $n/N = \rho + O(1/N^s)$ as $n,m\rightarrow\infty$ for some $0 < \rho < 1$ and $s>0$.   
  Under the alternative, if for some constant $0 < \tau \le 1$, it is satisfied that
  \begin{equation}\label{Equ:consistency_1}
      \frac{\lrabs{\ce^k(X,Y)}^{2+2\tau}}{N^{\tau} \lrp{\cv_k^2(Z)}^{1+\tau}} \rightarrow 0, \quad\quad\quad
      { \frac{\bbe{k^2(Z_1,Z_2)}}{N^s\cv_k^2(Z)} \rightarrow 0,}
  \end{equation}
  and
  \begin{equation}\label{Equ:consistency_2}
      \frac{\blre{\lrabs{\tilde{k}(X_1,X_2)}^{2+2\tau} + \lrabs{\tilde{k}(X_1,Y_1)}^{2+2\tau} + \lrabs{\tilde{k}(Y_1,Y_2)}^{2+2\tau}}}{N^{\tau} \lrp{\cv_k^2(Z)}^{1+\tau}} \rightarrow 0,
  \end{equation}
  then it holds that $\displaystyle{\frac{\cv_{n,m}^{k\ast}(X,Y)}{\cv_k^2(Z)} \longrightarrow^{p} 1}$ as { $N,p \rightarrow \infty$}. 
\end{proposition}

The argument to show Proposition \ref{Prop:ratio-const-k-H1} is quite similar to that of Proposition \ref{Prop:ratio-const-k-H0}, with the main difference being attributed to the fact that $\cv_{n,m}^{k\ast}(X,Y)$ is not an unbiased estimator of $\cv_k^2(Z)$ under the alternative. In this case, to bound $\bbe{|\cv_{n,m}^{k\ast}(X,Y) - \cv_k^2(Z)|^{1+\tau}}$, we break it into two parts, namely, $\bbe{|\cv_{n,m}^{k\ast}(X,Y)-\be[\cv_{n,m}^{k\ast}(X,Y)]|^{1+\tau}}$ and $|\be[\cv_{n,m}^{k\ast}(X,Y)] - \cv_k^2(Z)|^{1+\tau}$. Note that $\cv_{n,m}^{k\ast}(X,Y)-\bbe{\cv_{n,m}^{k\ast}(X,Y)}$ can be decomposed as a combination of multiple $U$-statistics and its upper bound is obtained by a 
moment inequality. This is manifested in the first condition in (\ref{Equ:consistency_1}) and condition (\ref{Equ:consistency_2}). A major difference from the condition under the null is that, to make the pooled-sample estimate $\cv_{n,m}^{k\ast}(X,Y)$ ratio-consistent for $\cv_k^2(Z)$, the discrepancy between the distributions of $X$ and $Y$, as quantified by $\ce^k(X,Y)$, cannot be too large, as regulated by the first condition in (\ref{Equ:consistency_1}).

The upper bound of $|\bbe{\cv_{n,m}^{k\ast}(X,Y)} - \cv_k^2(Z)|^{1+\tau}$ corresponds to the second condition in (\ref{Equ:consistency_1}). In fact, it follows from some simple calculations that, under the assumption $n/N = \rho + O\lrp{1/N^{s}}$, the bias $\bbe{\cv_{n,m}^{k\ast}(X,Y)} - \cv_k^2(Z)$ can be bounded by $\bbe{k^2(Z_1, Z_2)}/N^s$ up to a multiplicative constant, where the convergence rate of $n/N$ is involved.

To establish the central limit theorem for the proposed test, we define the functionals
\begin{equation}\label{Def:gk}
  g^k(X_1,X_2,X_3,X_4)
= d^k(X_1,X_2) d^k(X_1,X_3) d^k(X_2,X_4) d^k(X_3,X_4),
\end{equation}
where $d^k$ is defined as (\ref{Equ:d^k}). We can obtain the following central limit theorem for the proposed test statistic under the null.

\begin{theorem}\label{Thm:clt}
Assume that $\blre{k^4(Z_1,Z_2)}<\infty$ for each $k=k^{(p)}\in\cc$ and $n/N \rightarrow \rho$ as $n,m\rightarrow\infty$ for some $0 < \rho < 1$. Under the null {\color{black} when $X =^d Y =^d Z$}, if for some constant $0 < \tau \le 1$, it is satisfied for some $\{(k^{(p)})\}\in\cc$ that
\begin{equation}\label{Equ:clt_1}
    \frac{\blre{\lrabs{\tilde{k}(Z_1,Z_2)}^{2+2\tau}}}{N^{\tau} \lrp{\cv_k^2(Z)}^{1+\tau}} \rightarrow 0,
\end{equation}
and
\begin{equation}\label{Equ:clt_2}
   \frac{\lrp{\blre{g^k(Z_1,Z_2,Z_3,Z_4)}}^{(1+\tau)/2}}{\lrp{\cv_k^2(Z)}^{1+\tau}} 
\rightarrow 0,
\end{equation}
when $N,p\rightarrow\infty$, then it holds for this sequence $\{k^{(p)}\}$ that $\displaystyle{T_{n,m,p}^k \longrightarrow^d \cn(0,1)}$.
\end{theorem}

This theorem can be viewed as a counterpart of Theorem 1 in \cite{gao2019asymptotic} but is stated for a general kernel. Under the null, it follows from Proposition \ref{Prop:ratio-const-k-H0} that $\cv_{n,m}^{k\ast}(X,Y)$ is ratio-consistent for $\cv_k^2(Z)$ under $H_0$ when (\ref{Equ:clt_1}) is satisfied, and it is derived in Proposition \ref{Prop:est_MMD_simplified} that $\ce_{n,m}^{k}(X,Y)=L_{n,m}^{k}(X,Y)$. To derive the central limit theorem of $T_{n,m,p}^{k}$, it suffices to investigate the asymptotic behavior of $\frac{L_{n,m}^k(X,Y)}{\sqrt{c_{n,m}\cv_k^2(Z)}}$ using the martingale central limit theorem, since $L_{n,m}^{k}(X,Y)$ forms a martingale. Condition (\ref{Equ:clt_2}) is basically Lyapunov-type condition in the use of martingale central limit theorem.

Note that condition (\ref{Equ:clt_2}) only depends on $p$ and is free of the sample size $N$, whereas condition (\ref{Equ:clt_1}) depend on both $N$ and $p$ and thus might impose some implicit constraints between the divergence rate of $N$ and $p$. For any fixed $\tau$, if the order of $\bbe{|\tilde{k}(Z_1,Z_2)|}^{2+2\tau}$ does not exceed that of $(\cv_k^2(Z))^{1+\tau}$, then the quantity in (\ref{Equ:clt_1}) naturally goes to zero as long as $N$ diverges without additional restrictions between $N$ and $p$. As it turns out, it can be shown that the orders of $\bbe{|\tilde{k}(Z_1,Z_2)|}^{2+2\tau}$ and $(\cv_k^2(Z))^{1+\tau}$ are the same for the Gaussian kernel, the Laplacian kernel as well as the $L_2$ norm, hence the first term is independent of $p$ for these kernels. Furthermore, later in this paper, we show that this is true as long as the kernel satisfies some technical conditions.
 
In the literature, \cite{zhu2021interpoint} obtained the asymptotic distribution for the MMD permutation test statistics under the HDLSS (high-dimensional low sample size, where $p$ grows to infinity and $(n,m)$ is fixed) and HDMSS (high-dimensional medium sample size, where $(p,n,m)$ all grow to infinity but $p$ grows faster than $N$). The asymptotic results for the studentized test proposed in \cite{chakraborty2021new} are also limited to the HDLSS setting. \cite{yan2023kernel} obtained the CLT of a standardized {\color{black} MMD statistic} for a factor-like model allowing $(p,n,m)$ to diverge without constraints.

\subsection{Rate of Convergence}

We can further obtain the rate of convergence of the test statistic under the null using the Berry-Esseen bound for martingales, which has been used in \cite{gao2019asymptotic}. Here we first follow their steps to find an upper bound of $\sup\limits_{x\in\br} \lrabs{\blrp{T_{n,m,p}^k \le x} - \Phi(x)}$.

Note that under the null, for any $0<\gamma<1$, we have
\BEqn
& & \sup\limits_{x\in\br} \lrabs{\blrp{T_{n,m,p}^k \le x} - \Phi(x)}
 =  \sup\limits_{x\in\br} \lrabs{\blrp{\frac{\ce_{n,m}^k}{\sqrt{c_{n,m} \cv_{n,m}^{k\ast}(X,Y)}} \le x} - \Phi(x)} \\
&\le& 2 \sup\limits_{x\in\br} \lrabs{\blrp{\frac{L_{n,m}^k}{\sqrt{\Var(L_{n,m}^k)}} \le x} - \Phi(x)}
+ \sup\limits_{x\in\br} \lrabs{\Phi(x) - \Phi(x\sqrt{1+\gamma})} \\
& &+\sup\limits_{x\in\br} \lrabs{\Phi(x) - \Phi(x\sqrt{1-\gamma})}
+ 2 \blrp{\lrabs{\frac{c_{n,m} \cv_{n,m}^{k\ast}(X,Y)}{\Var(L_{n,m}^k)} - 1} > \gamma} =:
2P_1 + P_2 + P_3 + 2P_4.
\EEqn

By upper bounding each $P_i$, we obtain the following theorem.
\begin{theorem}\label{Thm:BEbd-k}
Let $Z$ denote the mixture distribution of $X$ and $Y$ defined as (\ref{Equ:mixZ}). Assume that $\blre{k^4(Z_1,Z_2)}<\infty$ for each $k=k^{(p)}\in\cc$ and $n/N \rightarrow \rho$ as $n,m\rightarrow\infty$ for some $0 < \rho < 1$. Under the null, it holds for any $n,m,p$ and $0<\tau\le1$ that
\begin{equation*}
    \sup\limits_{x\in\br} \lrabs{\blrp{T_{n,m,p}^k \le x} - \Phi(x)}
\le C(\rho,\tau)
    \lrcp{\frac{\blre{\lrabs{\tilde{k}(Z_1,Z_2)}^{2+2\tau}}}{N^\tau(\cv_k^2(Z))^{1+\tau}}
  + \frac{\lrp{\blre{g^k(Z_1,Z_2,Z_3,Z_4)}}^{\frac{1+\tau}{2}}}{(\cv_k^2(Z))^{1+\tau}}}^{\frac{1}{3+2\tau}}.
\end{equation*}
\end{theorem}

Theorem \ref{Thm:BEbd-k} states a non-asymptotic Berry-Esseen bound of the proposed test statistic. The two terms in the bound, $\frac{\bbe{|\tilde{k}(Z_1,Z_2)|}^{2+2\tau}}{N^\tau(\cv_k^2(Z))^{1+\tau}}$ and $\frac{\bigp{\bbe{g^k(Z_1,Z_2,Z_3,Z_4)}}^{(1+\tau)/2}}{(\cv_k^2(Z))^{1+\tau}}$ jointly determine the accuracy of normal approximation. As we have mentioned after Theorem \ref{Thm:clt}, the second term is solely determined by $p$, whereas the first term might depend on both $N$ and $p$. Although the bound established in Theorem \ref{Thm:BEbd-k} is valid for any $n,m,p$, the accuracy of normal approximation is guaranteed only when both quantities are close to zero, which might put some restrictions on the way the dimension $p$ diverges with respect to $N$. Such restriction is implicit for a general kernel $k$, but under some assumptions we can explicitly calculate the order of each term on the right-hand side, which enables us to derive the specific regime where the bound goes to zero. To this end, we first present a computational formula for $\bbe{g^k(Z_1,Z_2,Z_3,Z_4)}$ in the following proposition.

\begin{proposition}\label{Prop:gk}
  Assume that $\blre{k^4(Z_1,Z_2)}<\infty$, it holds that
  \begin{equation*}
      \blre{g^k(Z_1,Z_2,Z_3,Z_4)}
    = G_1 + G_2 + G_3 + G_4,
  \end{equation*}
  where 
  $ G_1 = \blre{k(Z_1,Z_2) k(Z_1,Z_3) k(Z_2,Z_4) k(Z_3,Z_4)}
         -4\blre{k(Z_1,Z_2) k(Z_1,Z_3) k(Z_2,Z_4) k(Z_4,Z_5)}  
         +2\blre{k(Z_1,Z_2) k(Z_1,Z_3)}^2$,
  $ G_2 = 4\blre{k(Z_1,Z_2)}\blre{k(Z_1,Z_2) k(Z_1,Z_3) k(Z_2,Z_4)}$,
  $ G_3 = -4\blre{k(Z_1,Z_2)}^2\\ \times\blre{k(Z_1,Z_2) k(Z_1,Z_3)}$, and
  $ G_4 = \blre{k(Z_1,Z_2)}^4$.
\end{proposition}

If we restrict our attention to the kernels of the form $k(x,y) = f(|x-y|)$ for some smooth function $f$, we can derive the explicit rate of convergence. To this end, in the following  we state the technical assumptions on  $f$ and the distributions of $X$ and $Y$.

\begin{assumpt}\label{Assumpt:uniform-kernel}
  Assume that for each $k^{(p)}\in\cc$, there exists some function $f^{(p)}$, such that $k^{(p)}(x,y) = f^{(p)}(|x-y|)$ for any $x,y\in\br^p$. Let $D$ be the domain of $f^{(p)}$ and $D_0 \subseteq D$ be a set that contains $A_0=\blre{|Z_1-Z_2|^2}^{1/2}$ and $A_0^{XY}=\blre{|X_1-Y_1|^2}^{1/2}$.
  
  Additionally, assume that 
  \begin{enumerate}[label=(\roman*)]
      \item \label{Assumpt:uniform-kernel-1} 
      for each $f^{(p)}$ and any $s\in D$ and $s_0 \in D_0$, it holds that
      \begin{equation*}
          f^{(p)}(s)
        = \sum\limits_{i=0}^{6} \frac{1}{i!} f_i^{(p)}(s_0) (s-s_0)^i 
        + f_7^{(p)}(\xi(s,s_0)) (s-s_0)^7,
      \end{equation*}
      where $f_i^{(p)}$ denotes the $i$-th order derivative of $f^{(p)}$, and $\xi(s,s_0)$ denotes some point between $s$ and $s_0$.
      
      \item \label{Assumpt:uniform-kernel-2}
      these exists a positive constant $\tilde{M}<\infty$, such that for any $f^{(p)}$ and any $s_0 \in D_0$, it holds that
      \begin{equation*}
          \max\limits_{1 \le i \le 7} \sup\limits_{s\in D} |f_i^{(p)}(s)|\cdot|s_0^i| \le \tilde{M} |f_0^{(p)}(s_0)|.
      \end{equation*}
      
      \item \label{Assumpt:uniform-kernel-3}
      there exists a positive constant $\hat{M}<\infty$, such that for any $f^{(p)}$ and any $s_0 \in D_0$, it holds that
      \begin{equation*}
          |f_0^{(p)}(s)| \le \hat{M} \min\{|\frac{1}{2}f_1^{(p)}(s)s|, |-\frac{1}{8}f_1^{(p)}(s)s + \frac{1}{8}f_2^{(p)}(s)s^2|\}.
      \end{equation*}
      
      \item \label{Assumpt:uniform-kernel-4}
      there exists a positive constant $\hat{M}<\infty$, such that for any $f^{(p)}$ and any $s_0 \in D_0$, it holds that
      \BEqn
          |f_0^{(p)}(s)| 
      &\le& \hat{M} \min\{|\frac{1}{16}f_1^{(p)}(s)s - \frac{1}{16}f_2^{(p)}(s)s^2 + \frac{1}{48}f_3^{(p)}(s)s^3|, \\
      & & \hspace{3.5em}
          |-\frac{5}{128}f_1^{(p)}(s)s + \frac{5}{128}f_2^{(p)}(s)s^2 - \frac{1}{64}f_3^{(p)}(s)s^3 + \frac{1}{384}f_4^{(p)}(s)s^4|\}.
      \EEqn
  \end{enumerate}
\end{assumpt}

Assumption \ref{Assumpt:uniform-kernel}\ref{Assumpt:uniform-kernel-1} is mild and it only requires that the function $f$ is smooth enough and has continuous derivatives up to the 7-th order. Assumption \ref{Assumpt:uniform-kernel}\ref{Assumpt:uniform-kernel-2}-\ref{Assumpt:uniform-kernel-4} further regulates the smoothness of the derivatives of $f$ and will be used to determine the exact orders of $\ce^k(X,Y)$ and $\cv_k^2(Z)$. 
{\color{black} Later in Section \ref{subsec:example}, we will use the Gaussian kernel as a special example to demonstrate the verification of Assumption \ref{Assumpt:uniform-kernel}. Additional examples for $L_2$ norm and the Laplacian kernel can be found in the online supplement.}

Before stating the next assumption, we introduce some useful notations. Define $A=\blre{|Z_1-Z_2|^2}$, $A^{X}=\blre{|X_1-X_2|^2}$, $A^{XY}=\blre{|X_1-Y_1|^2}$ and $A^{Y}=\blre{|Y_1-Y_2|^2}$. Let $A_0 = A^{1/2}$ and define $A_0^{X},A_0^{XY}$ and $A_0^{Y}$ in the same way. For each $X$ and $Y$, let $\mu_X = \blre{X}$ and $\mu_Y = \blre{Y}$ be the mean vectors, and use $\tilde{X}=X-\mu_X$ and $\tilde{Y}=Y-\mu_Y$ to denote the centered version of $X,Y$ respectively. Additionally, we use $\Delta = \mu_X - \mu_Y$  to denote the mean difference between $X$ and $Y$ and we denote the covariance matrices of $X,Y$ by $\Sigma_X = (\sigma^2_{X,j_1 j_2})_{j_1,j_2}$ and $\Sigma_Y = (\sigma^2_{Y,j_1 j_2})_{j_1,j_2}$. 

\begin{assumpt}\label{Assumpt:component-dept}
  For any fixed $p$ and $X = (x_1,\dots,x_p)^{\top}$ and $Y = (y_1,\dots,y_p)^{\top}$, assume that
  \begin{enumerate}[label=(\roman*)]
      \item \label{Assumpt:component-dept-1}
      there exists an integer $1 \le \alpha(p) \prec p$, such that $X$ and $Y$ have $\alpha(p)$-dependent components, respectively. Specifically, for any $1 \le j \le p-\alpha(p)-1$ and $\ell > \alpha(p)$, $\{x_i\}_{i \le j}$ is independent of $\{x_i\}_{i \geq j+\ell}$, and $\{y_i\}_{i \le j}$ is independent of $\{y_i\}_{i \geq j+\ell}$.   
      \item \label{Assumpt:component-dept-2}
      there exists a constant $0 < U^{\ast} < \infty$, such that
      \begin{equation*}
          \max\limits_{1 \le j \le p} \max\limits_{1 \le r \le 128} \{\be[|x_j|^r], \be[|y_j|^r]\} < U^{\ast}.
      \end{equation*}
      \item \label{Assumpt:component-dept-3}
      there exists some positive constants $0<L_0,U_0<\infty$, such that 
      \begin{equation*}
          L_0 p \le \min\{A^{X},A^{XY},A^{Y}\} \le \max\{A^{X},A^{XY},A^{Y}\} \le U_0 p.
      \end{equation*}
      \item \label{Assumpt:component-dept-4}
      there exists some positive constants $0<L_0^{\ast},U_0^{\ast}<\infty$, such that 
      \begin{equation*}
          L_0^{\ast} \alpha(p) p 
      \le \min\{\|\Sigma_X\|_F^2, \|\Sigma_Y\|_F^2\} 
      \le \max\{\|\Sigma_X\|_F^2, \|\Sigma_Y\|_F^2\} 
      \le U_0^{\ast} \alpha(p) p.
      \end{equation*}
  \end{enumerate}
\end{assumpt}
Assumption \ref{Assumpt:component-dept}\ref{Assumpt:component-dept-1} imposes some condition on the weak componentwise dependence within $X$ and $Y$ and it only needs to hold for some permutation of components of $X$ and $Y$, as our test statistic is permutation-invariant when the kernel $k(x,y)=f(|x-y|)$.  It is worth noting that $\alpha(p)$ may vary w.r.t. $p$ and thus the range of dependence is allowed to grow when $p$ increases. Assumption \ref{Assumpt:component-dept}\ref{Assumpt:component-dept-2} requires a uniform bound of the componentwise moments of both distributions, which  can be relaxed at the expense of lengthy proofs. Assumption \ref{Assumpt:component-dept}\ref{Assumpt:component-dept-3} requires both $\bbe{|\tilde{X}|^2}$ and $\bbe{|\tilde{Y}|^2}$ are strictly of order $p$, which is a mild condition. Finally, Assumption \ref{Assumpt:component-dept}\ref{Assumpt:component-dept-4} specifies the order of $\|\Sigma_X\|_F^2$ and $\|\Sigma_Y\|_F^2$, which seems reasonable in views of the $\alpha(p)$-dependent assumption. With Assumption \ref{Assumpt:component-dept}, we are able to calculate the orders of the quantities involved in our main theorems, which lead to a specific convergence rate of normal approximation and some explicit power results to be stated in the next section. Note that it is not our intention to showcase the convergence rate of normal approximation under the weakest possible assumption, as that is at the expense of very complicated arguments.  Assumption~\ref{Assumpt:component-dept} is quite reasonable to illustrate the convergence rate in a case of broad interest.

\begin{proposition}\label{Prop:BEbd-f}
  Assume that $\blre{k^4(Z_1,Z_2)}<\infty$ and $n/N \rightarrow \rho$ as $n,m\rightarrow\infty$ for some $0 < \rho < 1$.
  Suppose that Assumptions \ref{Assumpt:uniform-kernel}\ref{Assumpt:uniform-kernel-1}-\ref{Assumpt:uniform-kernel-2} and Assumptions \ref{Assumpt:component-dept}\ref{Assumpt:component-dept-1}-\ref{Assumpt:component-dept-2} hold, then there exists some $p_0=p_0(\tilde{M},\hat{M},U^{\ast},L_0,U_0,L_0^{\ast},U_0^{\ast})$, such that for any $p \geq p_0$, it holds under the null that for any $k = k^{(p)}\in\cc$, 
  \begin{equation}
      \sup\limits_{x\in\br} \lrabs{\blrp{T_{n,m,p}^{k} \le x} - \Phi(x)}
      \le
      C(\tilde{M},\hat{M},U^{\ast},L_0,U_0,L_0^{\ast},U_0^{\ast}) \lrp{\frac{1}{N} + \frac{\alpha(p)}{p}}^{1/5},
  \end{equation}
  where $\tilde{M},\hat{M}$ are defined in Assumption \ref{Assumpt:uniform-kernel} and $U^{\ast},L_0,U_0,L_0^{\ast},U_0^{\ast}$ are defined in Assumption \ref{Assumpt:component-dept}.
\end{proposition}

Proposition \ref{Prop:BEbd-f} provides a uniform explicit rate of convergence for a class of kernels and for $X$ and $Y$ with weakly dependent components. In fact, the rate of convergence is determined only by $N,p$, and the parameters $\tilde{M},\hat{M}$ from Assumption \ref{Assumpt:uniform-kernel}, as well as $\alpha(p),U^{\ast},L_0,U_0,L_0^{\ast},U_0^{\ast}$ from Assumption \ref{Assumpt:component-dept}. One implication of Proposition \ref{Prop:BEbd-f} is that, the empirical distribution of the proposed test statistic 
can be accurately approximated by the standard Gaussian distribution only when both $N$ and $p$ diverge to infinity, though no constraint is required regarding the divergence rate between $N$ and $p$. Another implication is that, the dependence within $X$ and $Y$ is allowed to grow as $p$ increases, but at the sacrifice of the accuracy of normal approximation. 
When the dependence within $X$ and $Y$ gets stronger, accurate normal approximation can only be obtained with larger $p$. This theoretical phenomenon is consistent with our empirical finding in Section~\ref{Sec:NumStudies}.

The main theoretical tool we use to obtain the rate of convergence is the Berry-Esseen bound for martingale [\cite{haeusler1988rate}], as also used in \cite{gao2019asymptotic}. One important difference between \cite{gao2019asymptotic} and our work is that the denominator of our test statistic is estimated over the pooled sample, and its leading term is a combination of several two-sample U-statistics, and the tools provided in \cite{gao2019asymptotic} are not sufficient for our theory. To this end, we generalize the moment inequality for the one-sample U-statistic to the two-sample U-statistic. Furthermore, Berry-Esseen bound obtained here is valid for a general kernel, and the rate of convergence can be explicitly derived under some mild conditions as shown in Proposition \ref{Prop:BEbd-f}.

\subsection{Power Analysis}
\label{subsec:power}

Next we look into the power behavior of the studentized test statistic. In the following theorem, we can show that the power of the proposed test is asymptotically one under some conditions.

\begin{theorem}\label{Thm:power-k}
  Assume that $\blre{k^4(Z_1,Z_2)}<\infty$ and { $n/N = \rho + O(1/N^s)$ as $n,m\rightarrow\infty$ for some $0 < \rho < 1$ and $s>0$}. If for some constant $0 < \tau \le 1$, it holds that
  \begin{eqnarray}
 & \displaystyle{\frac{\lrabs{\ce^k(X,Y)}^{2+2\tau}}{N^{\tau} \lrp{\cv_k^2(Z)}^{1+\tau}} \rightarrow 0},
    \hspace{7em}
   { \displaystyle{\frac{\bbe{k^2(Z_1,Z_2)}}{N^s \cv_k^2(Z)} \rightarrow 0}},
 & \label{Condt:power-1} \\[3mm] 
 & \displaystyle{\frac{\blre{\lrabs{\tilde{k}(X_1,X_2)}^{2+2\tau} + \lrabs{\tilde{k}(X_1,Y_1)}^{2+2\tau} + \lrabs{\tilde{k}(Y_1,Y_2)}^{2+2\tau}}}{N^{\tau} \lrp{\cv_k^2(Z)}^2} \rightarrow 0},
 & \label{Condt:power-2} \\[3mm]
 & \displaystyle{\frac{N \lrp{\ce^k(X,Y)}^2}{\blre{\lrp{h^k(X_1,X_2,Y_1,Y_2)}^2}}} \rightarrow \infty,
    \hspace{3em}
    \displaystyle{\frac{N \ce^k(X,Y)}{\sqrt{\cv_k^2(Z)}} \rightarrow \infty},
 & \label{Condt:power-3} 
  \end{eqnarray}
where $\tilde{k}$ denotes the centered version of $k$. Then for any $C>0$, we have $\blrp{T_{n,m,p}^{k} > C}\rightarrow 1$ as $n,m,p \rightarrow \infty$.
\end{theorem}


Theorem \ref{Thm:power-k} gives the conditions under which the power of the test can be asymptotically one for a general kernel. Note that conditions (\ref{Condt:power-1}) and (\ref{Condt:power-2}) are introduced in Proposition \ref{Prop:ratio-const-k-H1} to ensure the ratio-consistency of the pooled-sample estimate $\cv_{n,m}^{k\ast}(X,Y)$. In the proof,  we show that when $\frac{N \lrp{\ce^k(X,Y)}^2}{\blre{\lrp{h^k(X_1,X_2,Y_1,Y_2)}^2}} \rightarrow \infty$, the sample estimate $\ce_{n,m}^k(X,Y)$ closely approximates its population counterpart $\ce^k(X,Y)$ and the asymptotic divergence of $T_{n,m,p}^{k}(X,Y)$ is equivalent to  $\frac{\ce^k(X,Y)}{\sqrt{c_{n,m}\cv_k^2(Z)}}$ diverging to infinity as $N,p$ increases. It is then not difficult to see that the asymptotic power one of the proposed test can be achieved under the condition $\frac{N \ce^k(X,Y)}{\sqrt{\cv_k^2(Z)}} \rightarrow \infty$.

The conditions presented in Theorem \ref{Thm:power-k} are sufficient but may not be necessary due to the technical arguments we employed. Nevertheless, Theorem \ref{Thm:power-k} can provide us some interesting insights of the regimes where our proposed test has nontrivial power. Below we shall discuss multiple scenarios based on the leading terms 
of $\ce^k(X,Y)$ and $\cv_k^2(Z)$. For the sake of readability, we only present the results when $\alpha(p) = O(1)$ (i.e. fixed) and $s=1$ and leave the general results when $\alpha(p)=p^{\delta_0}$ with $0\le\delta_0<1$ and $s>0$ in online appendices.

\begin{assumpt}\label{Assumpt:PowerAnalysis-1}
  For any fixed $p$ and $X=(x_1,\dots,x_p)^{\top}$ and $Y=(y_1,\dots,y_p)^{\top}$, assume that there exists some positive constants $L_1,U_1<\infty$, such that
  \begin{equation*}
      L_1 p 
  \le \max\{|\Delta|^2, \lrabs{\bbe{|\tilde{X}_1|^2} - \bbe{|\tilde{Y}_1|^2}} \}
  \le U_1 p.
  \end{equation*}
\end{assumpt}

Assumption \ref{Assumpt:PowerAnalysis-1} focuses on the scenario when at least one of $|\Delta|^2$ and $|\bbe{|\tilde{X}_1|^2} - \bbe{|\tilde{Y}_1|^2}|$ is strictly of order $p$. It holds under Assumption \ref{Assumpt:component-dept}\ref{Assumpt:component-dept-2} that $|\Delta|^2=\sum_{j=1}^{p} (\bbe{x_j}-\bbe{y_j})^2 \le 2\sum_{j=1}^{p} (\bbe{x_j}^2 + \bbe{y_j}^2) \le 4U^{\ast} p$ and $\lrabs{\bbe{|\tilde{X}_1|^2} - \bbe{|\tilde{Y}_1|^2}} = \lrabs{\sum_{j=1}^{p} (\Var(x_j)-\Var(y_j))} \le \sum_{j=1}^{p} (\Var(x_j)+\Var(y_j)) \le 2U^{\ast}p$, then Assumption \ref{Assumpt:PowerAnalysis-1} implies that the differences in componentwise mean or variance attain the highest possible order. 

\begin{proposition}\label{Prop:PowerAnalysis-1}
  Suppose that Assumption \ref{Assumpt:uniform-kernel}\ref{Assumpt:uniform-kernel-1}-\ref{Assumpt:uniform-kernel-2}, Assumption \ref{Assumpt:component-dept}\ref{Assumpt:component-dept-1}-\ref{Assumpt:component-dept-3} and Assumption \ref{Assumpt:PowerAnalysis-1} hold, and additionally, for any $k=k^{(p)}\in\cc$, assume that $\blre{k^4(Z_1,Z_2)}<\infty$ and $n/N = \rho + O(1/N)$ as $n,m\rightarrow\infty$ for some $0 < \rho < 1$. When there exists some positive constant $L^{\ast}<\infty$, such that
  \begin{equation}\label{Condition:PowerAnalysis-1}
     \lrabs{2f(A_0^{XY}) - f(A_0^{X}) - f(A_0^{Y})}
  \geq L^{\ast} |f(A_0^{XY})|,
  \end{equation}
  then it holds that $\blrp{T_{n,m,p}^{k}>C} \rightarrow 1$ as $n,m,p\rightarrow\infty$.
\end{proposition}

The condition (\ref{Condition:PowerAnalysis-1}) requires that the leading term of $\lrabs{2f(A_0^{XY}) - f(A_0^{X}) - f(A_0^{Y})}$ can be lower bounded by $|f(A_0^{XY})|$ up to a multiplicative constant, which is a mild condition and can be satisfied by many kernel functions; {\color{black} see Section \ref{subsec:example} for its verification of the Gaussian kernel and the online supplement for the verifications of the $L_2$ norm and the Laplacian kernel.}
Under the assumptions in Proposition \ref{Prop:PowerAnalysis-1}, $2f(A_0^{XY}) - f(A_0^{X}) - f(A_0^{Y})$ is the leading term of $\ce^k(X,Y)$ and $(2f(A_0^{XY}) - f(A_0^{X}) - f(A_0^{Y}))^2$ is that of $\cv_k^2(Z)$. It follows that $(\ce^k(X,Y))^2$ and $\cv_k^2(Z)$ are of the same order, and 
both of them dominate $\bbe{k^2(Z_1,Z_2)}$ and $\bbe{(h^k(X_1,X_2,Y_1,Y_2))^2}$. Additionally, $(\cv_k^2(Z))^2$ dominated the numerator of condition (\ref{Condt:power-2}). Consequently, all the conditions in Theorem \ref{Thm:power-k} are naturally satisfied and the nontrivial power is obtained with no constraints on the order of $p$ relative to $N$, which seems reasonable in view of significant differences in either the means and/or the sum of marginal variances.
In comparison, \cite{zhu2021interpoint} obtained the asymptotic power one result for the MMD permutation test under the HDLSS  and HDMSS settings only. The asymptotic power function for the studentized test proposed in \cite{chakraborty2021new} is also derived only under the HDLSS setting. 
{\color{black} Additional comparison with \cite{zhu2021interpoint} under the special case when $X$ and $Y$ have either identical means or identical covariance matrices are discussed for the Gaussian kernel in Section \ref{subsec:example}; see Remark \ref{Rmk:ComparePower}.}

Next, we further investigate the scenarios where the differences in marginal mean or variance are weaker. 

\begin{assumpt}\label{Assumpt:PowerAnalysis-2}
  For any fixed $p$ and $X=(x_1,\dots,x_p)^{\top}$ and $Y=(y_1,\dots,y_p)^{\top}$, assume that there exists some positive constants $L_2,U_2<\infty$ and $0\le\delta_1<1$, $-1\le\delta_2<1$, such that $L_2\alpha(p)p \le \|\rho\Sigma_X+(1-\rho)\Sigma_Y\|_F^2 \le U_2\alpha(p)p$, and
  \begin{equation*}
      L_2 p^{\delta_1} \le |\Delta|^2 \le U_2 p^{\delta_1} 
  \quad\mbox{and}\qquad
      L_2 p^{(1+\delta_2)/2} \le \lrabs{\bbe{|\tilde{X}_1|^2} - \bbe{|\tilde{Y}_1|^2}} \le U_2 p^{(1+\delta_2)/2}.
  \end{equation*}
\end{assumpt}

Assumption \ref{Assumpt:PowerAnalysis-2} considers the case where the orders of both $|\Delta|^2$ and $\lrabs{\bbe{|\tilde{X}_1|^2} - \bbe{|\tilde{Y}_1|^2}}$ are strictly smaller than $p$ but are no smaller than a constant. In this case, the differences in the marginal mean and variance still dominate those in higher moments as long as $\max\{\delta_1,\delta_2\}>0$, resulting in high power under certain rate constraints on $p$. The condition $L_2\alpha(p)p \le \|\rho\Sigma_X+(1-\rho)\Sigma_Y\|_F^2 \le U_2\alpha(p)p$ is mild under Assumption \ref{Assumpt:component-dept}\ref{Assumpt:component-dept-4}.

\begin{proposition}\label{Prop:PowerAnalysis-2}
  Suppose that Assumption \ref{Assumpt:uniform-kernel}\ref{Assumpt:uniform-kernel-1}-\ref{Assumpt:uniform-kernel-3} and Assumption \ref{Assumpt:component-dept}\ref{Assumpt:component-dept-1}-\ref{Assumpt:component-dept-4} hold, Assumption \ref{Assumpt:PowerAnalysis-2} holds with $\delta_1\neq\delta_2$ and $\max\{\delta_1,\delta_2\}>0$, and additionally, for any $k=k^{(p)}\in\cc$, assume that $\blre{k^4(Z_1,Z_2)}<\infty$ and $n/N = \rho + O(1/N)$ as $n,m\rightarrow\infty$ for some $0 < \rho < 1$, then it holds that $\blrp{T_{n,m,p}^{k}>C} \rightarrow 1$ when $n,m,p\rightarrow\infty$ as long as $p = o\lrp{N^{1/(2-2\max\{\delta_1,\delta_2,1/2\})}}$.
\end{proposition}

When the leading difference between the distributions of $X$ and $Y$ lies in marginal mean or variance, the $\ce^{k}(X,Y)$ is of order $|f_0(A_0^{XY})|p^{\max\{\delta_1,\delta_2\}-1}$, while $\cv_k^2(Z)$ is of order $f_0^2(Z_0^{XY}) p^{2\max\{\delta_1,\delta_2,1/2\}-2}$. When $\max\{\delta_1,\delta_2\} > 1/2$, $(\ce^k(X,Y))^2$ has the same order as of $\cv_k^2(Z)$, thus condition (\ref{Condt:power-2})-(\ref{Condt:power-3}) are naturally satisfied. To have $\frac{\bbe{k^2(Z_1,Z_2)}}{N \cv_k^2(Z)} \rightarrow 0$, we need additional constraint between $N$ and $p$ since $\cv_k^2(Z)$ is dominated by $\bbe{k^2(Z_1,Z_2)}$. For the case that $\max\{\delta_1,\delta_2\}\le 1/2$, the order of $\cv_k^2(Z)$ becomes $f_0^2(A_0^{XY})p^{-1}$ but that of $\ce^k(Z_1,Z_2)$ remains unchanged. Hence $(\ce_k^2(X,Y))^2$ is no longer capable of dominating $\cv_k^2(Z)$ and $\bbe{(h^k(X_1,X_2,Y_1,Y_2))^2}$, which leads to the constraint
 $p=o(N)$ to make condition (\ref{Condt:power-3}) hold. 

Intuitively, as the disparities in marginal mean and variance between $X$ and $Y$ weakens to the point $\max(\delta_1,\delta_2)\le 1/2$, our proposed test has nontrivial power only when the growth rate of $p$ is strictly smaller than that  of $N$. When $\delta_1=\delta_2$, similar power results can be attained given a specific kernel function following some lengthy analysis, but we exclude this case for simplicity.


Next we investigate the scenario when the differences in the marginal mean and variance between $X$ and $Y$ further diminish.

\begin{assumpt}\label{Assumpt:PowerAnalysis-3}
  For any fixed $p$ and $X=(x_1,\dots,x_p)^{\top}$ and $Y=(y_1,\dots,y_p)^{\top}$, assume that $|\Delta|=0$ and there exists some positive constants $L_3,U_3<\infty$ and $0\le \delta_3,\delta_4<1$, such that
  \begin{equation*}
      L_3 p^{\delta_3/2} \le \lrabs{\bbe{|\tilde{X}_1|^2} - \bbe{|\tilde{Y}_1|^2}} \le U_3 p^{\delta_3/2}
  \quad\mbox{and}\quad
      L_3 p^{\delta_4/2} \le \|\Sigma_X-\Sigma_Y\|_F \le U_3 p^{\delta_4/2}.
  \end{equation*}
\end{assumpt}

Assumption \ref{Assumpt:PowerAnalysis-3} targets at the case when $X$ and $Y$ have the identical mean, and their leading disparities fall within the covariances. Note that with $\alpha(p)=O(1)$ and Assumption \ref{Assumpt:component-dept}\ref{Assumpt:component-dept-4}, the order of $\|\Sigma_X-\Sigma_Y\|_F^2$ won't exceed $p$, and under Assumption \ref{Assumpt:PowerAnalysis-3} we set it to be $p^{\delta_4}$, where $\delta_4\in (0,1)$.

\begin{proposition}\label{Prop:PowerAnalysis-3}
  Suppose that Assumption \ref{Assumpt:uniform-kernel}\ref{Assumpt:uniform-kernel-1}-\ref{Assumpt:uniform-kernel-3} and Assumption \ref{Assumpt:component-dept}\ref{Assumpt:component-dept-1}-\ref{Assumpt:component-dept-4} hold, Assumption \ref{Assumpt:PowerAnalysis-3} holds with $\delta_3\neq\delta_4$ and  $\max\{\delta_3,\delta_4\}>0$, and additionally, for any $k=k^{(p)}\in\cc$, assume that $\blre{k^4(Z_1,Z_2)}<\infty$ and $n/N = \rho + O(1/N)$ as $n,m\rightarrow\infty$ for some $0 < \rho < 1$, then it holds that $\blrp{T_{n,m,p}^{k}>C} \rightarrow 1$ when $n,m,p\rightarrow\infty$ as long as $p = o\lrp{N^{1/(3-2\max\{\delta_3,\delta_4\})}}$.
\end{proposition}

Under Assumption \ref{Assumpt:PowerAnalysis-3}, the order of $\ce^k(X,Y)$ decreases to $|f_0(A_0^{XY})| p^{\max\{\delta_3,\delta_4\}-2}$, as compared to the second scenario (under Assumption~\ref{Assumpt:PowerAnalysis-2}). Recall that $\ce^k(X,Y)$ characterizes the disparity between the distributions $X$ and $Y$, then it is not surprising that its order decreases as the leading disparities between $X$ and $Y$ move to some higher moment quantities. In this case, $\ce^k(X,Y)$ may not dominate $\bbe{(h^k(X_1,X_2,Y_1,Y_2))^2}$, and furthermore, $\ce^k(X,Y)$ is dominated by $\cv_k^2(Z)$, whose order stays at $f_0^2(A_0^{XY}) p^{-1}$. Therefore, additional constraints on $p$ are required to satisfy condition (\ref{Condt:power-3}).

\cite{zhu2021interpoint} showed that in the HDMSS setting, when 
$|\Delta|^2=o(\sqrt{p}/N)$, $\lrabs{\bbe{|\tilde{X}_1|^2} - \bbe{|\tilde{Y}_1|^2}} = o(\sqrt{p}/N)$ and $\|\Sigma_X-\Sigma_Y\|_F=o(\sqrt{p})$, MMD permutation test has trivial power.
In the special case $\delta_3=0, \delta_4\in (0,1)$, it is easy to see that both the condition in \cite{zhu2021interpoint} and our Assumption~\ref{Assumpt:PowerAnalysis-3} can be 
satisfied for different sets of $(p,N)$. 
 The resulting power phenomenon is strikingly different with  
the MMD permutation test being powerless and our studentized test being power one asymptotically. This difference is not a contradiction but is mainly attributed to the different regimes, since HDMSS setting implies $p\gg N$ whereas our Proposition~\ref{Prop:PowerAnalysis-3} requires $p\ll N$. 
This is an example that shows that even for the same alternative, the order of $p$ relative to $N$ can play an important role in determining the power behavior. 

Finally, we look into the scenario when $X$ and $Y$ have identical means and covariance matrices to complete the discussions in this section. With $\cum(\cdot)$ denoting the cumulant, we propose the following assumption.

\begin{assumpt}\label{Assumpt:PowerAnalysis-4}
  For any fixed $p$ and $X=(x_1,\dots,x_p)^{\top}$ and $Y=(y_1,\dots,y_p)^{\top}$, assume that $|\Delta|=\|\Sigma_X-\Sigma_Y\|_F=0$ and there exists some positive constants $L_4,U_4<\infty$ and $0\le \delta_5,\delta_6,\delta_7<1$, such 
  \begin{equation*}
      \begin{array}{l}
         L_4 p^{\delta_5} \le \sum\limits_{j_1,j_2,j_3=1}^{p}\lrp{\cum(\tilde{x}_{1 j_1},\tilde{x}_{1 j_2},\tilde{x}_{1 j_3}) - \cum(\tilde{y}_{1 j_1},\tilde{y}_{1 j_2},\tilde{y}_{1 j_3})}^2 \le U_4 p^{\delta_5},\\
         L_4 p^{\delta_6} \le \sum\limits_{j_2=1}^{p}\lrp{\sum\limits_{j_1=1}^{p}\{Cov(\tilde{x}_{1 j_1}^2-\sigma_{X, j_1},\tilde{x}_{1 j_2}) - \Cov(\tilde{y}_{1 j_1}^2-\sigma_{Y, j_1},\tilde{y}_{1 j_2})
    \}   }^2
         \le U_4 p^{\delta_6},\\
         L_4 p^{2\delta_7} \le \lrp{\sum\limits_{j_1,j_2=1}^{p}
    \{   \Cov(\tilde{x}_{1 j_1}^2-\sigma_{X, j_1},\tilde{x}_{1 j_2}^2-\sigma_{X, j_2}) - \Cov(\tilde{y}_{1 j_1}^2-\sigma_{Y, j_1},\tilde{y}_{1 j_2}^2-\sigma_{Y, j_2})  \}
    }^2 \le U_4 p^{2\delta_7}.\\
      \end{array}
  \end{equation*}
\end{assumpt}

\begin{proposition}\label{Prop:PowerAnalysis-4}
  Suppose that Assumption \ref{Assumpt:uniform-kernel}\ref{Assumpt:uniform-kernel-1}-\ref{Assumpt:uniform-kernel-4} and Assumption \ref{Assumpt:component-dept}\ref{Assumpt:component-dept-1}-\ref{Assumpt:component-dept-4} hold, Assumption \ref{Assumpt:PowerAnalysis-4} holds with $1+\max\{\delta_5,\delta_6\} \neq 2\delta_7$,
  $\max\{\delta_5,\delta_6,\delta_7\}>0$, and additionally, for any $k=k^{(p)}\in\cc$, assume that $\blre{k^4(Z_1,Z_2)}<\infty$ and $n/N = \rho + O(1/N)$ as $n,m\rightarrow\infty$ for some $0 < \rho < 1$, then it holds that $\blrp{T_{n,m,p}^{k}>C} \rightarrow 1$ when $n,m,p\rightarrow\infty$ as long as $p = o\lrp{N^{1/(7-2\max\{1+\max\{\delta_5,\delta_6\},2\delta_7\})}}$.
\end{proposition}

Proposition \ref{Prop:PowerAnalysis-4} implies that, when $\delta_7 < (1+\max\{\delta_5,\delta_6\})/2$, nontrivial power against the alternative is obtained when $p = o\lrp{N^{1/(5-2(\delta_5\vee\delta_6))}}$ and otherwise the corresponding regime is $p = o\lrp{N^{1/(7-4\delta_7)}}$. In fact, the order of $\cv_k^2(Z)$ remains $f_0^2(A_0^{XY}) p^{-1}$ while that of $\ce^k(X,Y)$ drops to $|f_0(A_0^{XY})| p^{\max\{\delta_5,\delta_6,2\delta_7-1\}-3}$, and following some similar arguments as in the previous scenario, we obtain the constraint between $N$ and $p$ for this case.

To summarize, Proposition \ref{Prop:PowerAnalysis-1}-Proposition \ref{Prop:PowerAnalysis-3} jointly investigate the cases when the discrepancy between two distributions is dominated by their differences in the mean and/or  covariance matrices, which correspond to S1 in Section 4 of \cite{yan2023kernel}. Proposition \ref{Prop:PowerAnalysis-4} corresponds to  the scenario where two distributions have identical first and second moments, and the difference lies  in the third and/or fourth moments/cumulants.
This scenario corresponds to S2 with $\ell=3$  in Section 4 of \cite{yan2023kernel}. The latter authors provided a comprehensive description of when their test  has trivial power, nontrivial power and asymptotic power one based on non-null CLT obtained. In general, we feel it is difficult to directly compare the power results in  \cite{yan2023kernel} with ours due to the different settings and regimes we explored. In particular, we mainly focus on the regime $p\ll N^{\omega_1}$ for some $\omega_1>0$ as stated in Proposition~\ref{Prop:PowerAnalysis-2}- Proposition~\ref{Prop:PowerAnalysis-4}, whereas \cite{yan2023kernel} focus on the regime where $N\ll p^{\omega_2}$ for some $\omega_2\ge 1/2$. The two regimes may have overlap (i.e., the intersection is nonzero),  their power one results and ours complement each other, and both contribute to the understanding of the space of alternatives for which  the MMD-based test has high power.

As revealed by the four propositions above,  our test is powerful against a wide range of alternatives, including the differences in means, variances, covariances and high-order features associated with the distributions. 

The intuition behind all these propositions is that, the disparities that fall within lower moments between $X$ and $Y$ are easier to be detected by our proposed test. When the leading differences move to higher moment quantities, stricter constraints between $N$ and $p$ are required to make the test powerful. This phenomenon is consistent with that found by \cite{yan2023kernel}, who provided an asymptotic exact power analysis and revealed a delicate interplay between the detectable moment discrepancy and the dimension-and-sample orders (see Table 1 therein).

{\color{black}
\subsection{An Illustrative Example with the Gaussian Kernel}\label{subsec:example}

As shown in Table \ref{Tab:Example}, a special case covered by our setup is the Gaussian kernel multiplied by -1, that is, $k(x,y)=-\exp\lrp{-|x-y|^2/(2\gamma^2)}$, where $\gamma$ is a pre-specified tuning parameter. We note that many technical assumptions and theoretical results are presented in the previous sections, which may be difficult to digest. In this section, we use the Gaussian kernel as a special example to demonstrate the verification of Assumption \ref{Assumpt:uniform-kernel} and condition (\ref{Condition:PowerAnalysis-1}) in the previous section.

We define $D=[0,\infty)$ and $D_0=\lrbk{\be[|Z_1-Z_2|^2]^{1/2}, \be[|X_1-Y_1|^2]^{1/2}}$. For each fixed $p$ and the tuning parameter $\gamma$ that depends on $p$, we consider the Gaussian kernel $k^{(p)}(x,y) = -\exp\lrp{-\frac{|x-y|^2}{2\gamma^2}}$. Here, different choices of $\gamma^2$ lead to different Gaussian kernels, and we restrict our interest to $k^{(p)}$ with specific $\gamma$, that is,
\begin{equation}\label{Equ:c}
    \cc = \cc(\ell,u)
 := \lrcp{k^{(p)}: \mbox{for each }p, \frac{\blre{|X_1-Y_1|^2}^{1/2}}{u} \le \gamma \le \frac{\blre{|Z_1-Z_2|^2}^{1/2}}{\ell}},
\end{equation}
where $0<\ell,u<\infty$ are some specified constants such that $\cc(\ell,u)$ is well defined. Note that $\cc$ is a set of Gaussian kernel sequences with growing $p$, for each $k^{(p)} \in \cc$, we define $f^{(p)}(s) = -\exp\lrp{-\frac{s^2}{2\gamma^2}}$ to be the unique smooth function associated with $k^{(p)}$ and for simplicity, we drop the superscript hereafter. With the explicit expression of $f$, we obtain the derivatives of $f$ up to the 7th order, that is
\begin{equation*}
\begin{array}{ll}
   \displaystyle
   f_0(s) = - \exp\lrp{-\frac{s^2}{2\gamma^2}},
&  \displaystyle
   f_1(s) = \frac{s}{\gamma^2} \exp\lrp{-\frac{s^2}{2\gamma^2}}, \\[5mm]
   \displaystyle
   f_2(s) = \lrp{\frac{1}{\gamma^2} - \frac{s^2}{\gamma^4}} \exp\lrp{-\frac{s^2}{2\gamma^2}}, 
&  \displaystyle
   f_3(s) = \lrp{-\frac{3s}{\gamma^4} + \frac{s^3}{\gamma^6}} \exp\lrp{-\frac{s^2}{2\gamma^2}}, \\[5mm]
   \displaystyle
   f_4(s) = \lrp{- \frac{3}{\gamma^4} + \frac{6s^2}{\gamma^6} - \frac{s^4}{\gamma^8}} \exp\lrp{-\frac{s^2}{2\gamma^2}}, 
&  \displaystyle
   f_5(s) = \lrp{\frac{15s}{\gamma^6} - \frac{10s^3}{\gamma^8} + \frac{s^5}{\gamma^{10}}} \exp\lrp{-\frac{s^2}{2\gamma^2}}, \\[5mm]   
\end{array}
\end{equation*}
\begin{equation*}
\begin{array}{l}
   \displaystyle
   f_6(s) = \lrp{\frac{15}{\gamma^6} - \frac{45s^2}{\gamma^8} + \frac{15s^4}{\gamma^{10}} - \frac{s^6}{\gamma^{12}}} \exp\lrp{-\frac{s^2}{2\gamma^2}}, \\[5mm]
   \displaystyle
   f_7(s) = \lrp{-\frac{105s}{\gamma^8} + \frac{105s^3}{\gamma^{10}} - \frac{21s^5}{\gamma^{12}} + \frac{s^7}{\gamma^{14}}} \exp\lrp{-\frac{s^2}{2\gamma^2}}.
\end{array}
\end{equation*}
It follows from the Taylor theorem with the Lagrange form of remainder that $\cc$ satisfies Assumption \ref{Assumpt:uniform-kernel}\ref{Assumpt:uniform-kernel-1} 

To verify Assumption \ref{Assumpt:uniform-kernel}\ref{Assumpt:uniform-kernel-2}, we note that 
\begin{equation*}
    \sup\limits_{s\in D}|f_1(s)| 
  = \sup\limits_{s\geq0} \lrabs{\frac{s}{\gamma^2}\exp\lrp{-\frac{s^2}{\gamma^2}}} = \frac{1}{\gamma} \sup\limits_{t\geq0} |t\exp(-t^2/2)| = \frac{0.607}{\gamma}.
\end{equation*}
It follows from similar steps that
\begin{equation*}
\begin{array}{lll}
    \displaystyle
    \sup\limits_{s\in D}|f_2(s)| = \frac{1}{\gamma^2},
  & \displaystyle
    \sup\limits_{s\in D}|f_3(s)| = \frac{1.38}{\gamma^3},
  & \displaystyle
    \sup\limits_{s\in D}|f_4(s)| = \frac{3}{\gamma^4}, \\[5mm]
    \displaystyle
    \sup\limits_{s\in D}|f_5(s)| = \frac{5.783}{\gamma^5},
  & \displaystyle
    \sup\limits_{s\in D}|f_6(s)| = \frac{15}{\gamma^6},
  & \displaystyle
    \sup\limits_{s\in D}|f_7(s)| = \frac{35.539}{\gamma^7}.
\end{array}
\end{equation*}
Recall that $\cc$ has restrictions on the tuning parameter $\gamma$ associated with $k^{(p)}$ such that $\blre{|X_1-Y_1|^2}^{1/2}/u \le \gamma \le \blre{|Z_1-Z_2|^2}^{1/2}/\ell$, then it holds that $\sup\limits_{s_0\in D_0} \frac{s_0}{\gamma} \le u$, thus $\max\limits_{1\le i\le 7}\sup\limits_{s\in D}|f_i(s)|\cdot|s_0^i| \le \tilde{M} |f_0(s_0)|$ holds for any $s_0\in D_0$ when 
\BEqn
    \tilde{M}
&=& \max\{0.607 u\exp(u^2/2),
        u^2\exp(u^2/2),
        1.38 u^3\exp(u^2/2),
        3 u^4\exp(u^2/2), \\
& & \hspace{2.5em}
        5.783 u^5\exp(u^2/2),
        15 u^6\exp(u^2/2),
        35.539 u^7\exp(u^2/2)\},
\EEqn
which completes the verification of Assumption \ref{Assumpt:uniform-kernel}\ref{Assumpt:uniform-kernel-2}.

As for Assumption \ref{Assumpt:uniform-kernel}\ref{Assumpt:uniform-kernel-3}, it follows from direct computation that $|f_0(s)| = \exp\lrp{-\frac{1}{2}\lrp{\frac{s}{\gamma}}^2}$, $\frac{1}{2}|f_1(s)s| = \frac{1}{2} \lrp{\frac{s}{\gamma}}^2 \exp\lrp{-\frac{1}{2}\lrp{\frac{s}{\gamma}}^2}$, and
\begin{equation*}
    |-\frac{1}{8}f_1(s)s + \frac{1}{8}f_2(s)s^2| 
 =  \frac{1}{8} \lrp{\frac{s}{\gamma}}^4 \exp\lrp{-\frac{1}{2}\lrp{\frac{s}{\gamma}}^2}.
\end{equation*}
Again, it follows from the definition of $\cc$ that $\inf\limits_{s_0\in D_0} (s_0/\gamma)\geq\ell$. Therefore, Assumption \ref{Assumpt:uniform-kernel}\ref{Assumpt:uniform-kernel-3} holds with $\hat{M}=\max\{2/\ell^2,8/\ell^4\}$.

Lastly, we look into Assumption \ref{Assumpt:uniform-kernel}\ref{Assumpt:uniform-kernel-4}. Note that $|f_0(s)| = \exp\lrp{-\frac{1}{2}\lrp{\frac{s}{\gamma}}^2}$, and
\BEqn
& & |\frac{1}{16}f_1(s)s - \frac{1}{16}f_2(s)s^2 + \frac{1}{48}f_3(s)s^3| 
 =  \frac{1}{48} \lrp{\frac{s}{\gamma}}^6 \exp\lrp{-\frac{1}{2}\lrp{\frac{s}{\gamma}}^2}, \\
& & |-\frac{5}{128}f_1(s)s + \frac{5}{128}f_2(s)s^2 - \frac{1}{64}f_3(s)s^3 + \frac{1}{384}f_4(s)s^4|
 =  \frac{1}{384} \lrp{\frac{s}{\gamma}}^8 \exp\lrp{-\frac{1}{2}\lrp{\frac{s}{\gamma}}^2},
\EEqn
then it is trivial Assumption \ref{Assumpt:uniform-kernel}\ref{Assumpt:uniform-kernel-4} is satisfied with $\hat{M} = \max\{48/\ell^6,384/\ell^8\}$.

To conclude, we present the results in the following proposition.
\begin{proposition}\label{Prop:Example-Gaussian}
    Let $\cc$ denote the set of Gaussian kernel sequences as defined in Equation (\ref{Equ:c}), it holds that $\cc$ satisfies Assumption \ref{Assumpt:uniform-kernel}.
\end{proposition}

Next, we verify condition \ref{Condition:PowerAnalysis-1} in Proposition \ref{Prop:PowerAnalysis-1}, that is
\begin{equation*}
    \lrabs{2f(A_0^{XY}) - f(A_0^{X}) - f(A_0^{Y})}
\geq L^{\ast} |f(A_0^{XY})|,
\end{equation*}
Again, we restrict the analysis to the set $\cc$. For each $k^{(p)}\in\cc$, we define $f(s)=-\exp\lrp{-\frac{s^2}{2\gamma^2}}$, and it follows from the definition of $\cc$ that $\bbe{|X_1-Y_1|^2}^{1/2}/u \le \gamma \le \bbe{|Z_1-Z_2|^2}^{1/2}/\ell$.
Note that
\begin{equation*}
    \lrp{2f(A_0^{XY}) - f(A_0^{X}) - f(A_0^{Y})}/f(A_0^{XY}) \\
 =  2 - \exp\lrp{-\frac{A^{X}-A^{XY}}{2\gamma^2}} - \exp\lrp{-\frac{A^{Y}-A^{XY}}{2\gamma^2}},
\end{equation*}
where 
\BEqn
& & \exp\lrp{-\frac{A^{X}-A^{XY}}{2\gamma^2}} = \exp\lrp{\frac{|\Delta|^2}{2\gamma^2}}\exp\lrp{-\frac{\bbe{|\tilde{X}_1|^2} - \bbe{|\tilde{Y}_1|^2}}{2\gamma^2}}, \\
& & \exp\lrp{-\frac{A^{Y}-A^{XY}}{2\gamma^2}} = \exp\lrp{\frac{|\Delta|^2}{2\gamma^2}}\exp\lrp{\frac{\bbe{|\tilde{X}_1|^2} - \bbe{|\tilde{Y}_1|^2}}{2\gamma^2}},
\EEqn
then it follows from the fact $\exp(s)>1$ and $\exp(s)+\exp(-s)>2$ for any $s>0$ that
\BEqn
& & \exp\lrp{-\frac{A^{X}-A^{XY}}{2\gamma^2}} + \exp\lrp{-\frac{A^{Y}-A^{XY}}{2\gamma^2}} \\
&=& \exp\lrp{\frac{|\Delta|^2}{2\gamma^2}}
    \lrp{\exp\lrp{-\frac{\bbe{|\tilde{X}_1|^2} - \bbe{|\tilde{Y}_1|^2}}{2\gamma^2}}
         + \exp\lrp{\frac{\bbe{|\tilde{X}_1|^2} - \bbe{|\tilde{Y}_1|^2}}{2\gamma^2}}} \\
&>& 2.
\EEqn
Consequently, the condition is naturally satisfied with
\BEqn
& & L^{\ast} \\ 
&=& \inf\limits_{A_0^{XY}/u \le \gamma \le A_0/\ell} 
    \lrcp{\exp\lrp{\frac{|\Delta|^2}{2\gamma^2}}
    \lrp{\exp\lrp{-\frac{\bbe{|\tilde{X}_1|^2} - \bbe{|\tilde{Y}_1|^2}}{2\gamma^2}}
         + \exp\lrp{\frac{\bbe{|\tilde{X}_1|^2} - \bbe{|\tilde{Y}_1|^2}}{2\gamma^2}}} -2}.
\EEqn

Suppose that Assumption \ref{Assumpt:component-dept} holds, and we have $\be[k^4(Z_1,Z_2)]<\infty$ for each $k=k^{(p)}\in\cc$, we summarize a few regimes where the asymptotic power of the Gaussian kernels in $\cc$ is one. We want to emphasize that Table \ref{Tab:Example-power} only includes a few special cases, whereas our proposed test is guaranteed to obtain full power asymptotically across a wider range of regimes. 

\begin{table}[!h]
    \centering
    \begin{tabular}{c||c|c|c|c}
    \hline\hline
        $|\Delta|^2$ & $O_s(p)$ & & $O_s(p^{1/2})$ & 0 \\[1mm] 
        $\lrabs{\bbe{|\tilde{X}_1|^2} - \bbe{|\tilde{Y}_1|^2}}$ & & $O_s(p)$ & $O_s(p^{1/2})$ & $O_s(p^{1/4})$ \\[2mm] 
        $\|\rho\Sigma_X+(1-\rho)\Sigma_Y\|_F^2$ & & & $O_s(\alpha(p)p)$ \\[2mm] 
        $\|\Sigma_X-\Sigma_Y\|_F$ & & & & $O_s(p^{1/4})$ \\[1mm] 
        \hline\hline 
        Regime & $N,p\rightarrow\infty$ & $N,p\rightarrow\infty$ & $p=o(N)$ & $p=o(N^{1/2})$ \\
    \hline\hline
    \end{tabular}
    \caption{Selected regimes where the power of the Gaussian kernels is asymptotically one.}
    \label{Tab:Example-power}
\end{table}

In the following remark, we compare the sufficient conditions for asymptotically power one derived in \cite{zhu2021interpoint} and in this article under the special case when $X$ and $Y$ have either identical means or identical covariance matrices.

\begin{remark}\label{Rmk:ComparePower}
    Both \cite{zhu2021interpoint} and our work aim to test for the distributional discrepancy, that is, to test for $H_0: X =^d Y$ versus $H_1: X \neq^d Y$. The discussion in \cite{zhu2021interpoint} is limited to MMD with a user-specified kernel $\hat{k}$ of the following expression
    \begin{equation*}
        \hat{k}(X,Y)
      = \varphi\lrp{\frac{1}{p}\sum\limits_{j=1}^{p}\psi(x_j, y_j)},
    \end{equation*}
    where $\psi\geq0$ and $\varphi$ has continuous second order derivative on $(0,\infty)$. It is trivial that the Gaussian kernel is covered by the set of $\hat{k}$. 
    
    Both the high dimensional low sample size setting (HDLSS) when $n,m$ are fixed but $p\rightarrow\infty$ and the high dimensional medium sample size (HDMSS) setting when $p\rightarrow\infty$ and $n:=n(p)\rightarrow\infty$ are investigated in \cite{zhu2021interpoint}, but here we only focus on HDMSS setting. It is shown that the permutation test in \cite{zhu2021interpoint} is consistent within the following consistency space $\ch$:
    \begin{equation*}
        \ch_c
      = \lrcp{(X,Y): ~2\varphi(e_{XY}) \neq \varphi(e_X) + \varphi(e_Y)},
    \end{equation*}
    where
    \begin{equation*}
        e_X = \lim\limits_{p\rightarrow\infty} \be[\bar{\psi}(X,X')], \qquad
        e_Y = \lim\limits_{p\rightarrow\infty} \be[\bar{\psi}(Y,Y')], \qquad
        e_{XY} = \lim\limits_{p\rightarrow\infty} \be[\bar{\psi}(X,Y)],
    \end{equation*}
    are all assumed to exist, with $X', Y'$ being independent copies of $X, Y$ and $\bar{\psi}$ denoting the average distance over components
    \begin{equation*}
        \bar\psi(Z_i,Z_j) = \frac{1}{p} \sum\limits_{s=1}^{p} \psi(z_{is},z_{js}).
    \end{equation*}
    For the Gaussian kernels, $\ch_c$ can be further characterized as 
    \begin{equation*}
        \ch_c 
      = \lrcp{(X,Y): ~\sum\limits_{j=1}^{p}\lrp{\be[x_j]-\be[y_j]}^2=o(p), ~|\sum\limits_{j=1}^{p}\lrp{\Var(x_j)-\Var(y_j)}| = o(p)}^c.
    \end{equation*}
    
    As shown in Theorem 3.5 of \cite{zhu2021interpoint}, it holds under the HDMSS that
    \begin{equation*}
        \lim\limits_{p\rightarrow\infty} \bp_{\ch_c}\lrp{\ce_{n,m}^k(X,Y) > c} = 1
    \end{equation*}
    for any $c\in\{Q_{\hat{R},1-\alpha}, Q_{\tilde{R},1-\alpha}\}$, where $Q_{\hat{R},1-\alpha}$ is the critical value obtained from $(n+m)!$ permutations with $\hat{R}$ being the randomization distribution of $\ce_{n,m}^k(X,Y)$, and $Q_{\hat{R},1-\alpha}$ is the critical value obtained from a fixed number $S$ of permutations with $\tilde{R}$ being the counterpart of $\hat{R}$ with only $S$ permutations. In other words, it is shown in \cite{zhu2021interpoint} that the asymptotic power of their permutation-based test is one within $\ch_c$. 
    
    Now we are able to compare the sufficient conditions for consistent power derived in \cite{zhu2021interpoint} and ours when $X$ and $Y$ have either identical means or identical covariance matrices. When $\mu_X=\mu_Y$, the condition in \cite{zhu2021interpoint} reduces to $|\sum\limits_{j=1}^{p}\lrp{\Var(x_j)-\Var(y_j)}| \gtrsim p$, which is equivalent to $\lrabs{\be[|\tilde{X_1}|^2] - \be[|\tilde{Y_1}|^2]} = O_s(p)$ in our article. When $\Sigma_X=\Sigma_Y$, the condition in \cite{zhu2021interpoint} is simplified as $\sum\limits_{j=1}^{p}\lrp{\be[x_j]-\be[y_j]}^2 \gtrsim p$, which exactly matches the condition $|\mu_X-\mu_Y|^2 = O_s(p)$ in our work. In summary, the sufficient conditions derived in both articles are equivalent under the special case that $X$ and $Y$ have either the same means or the same covariance matrices. However, it is worth mentioning that both works require some additional regularity conditions, which are not enumerated here.
\end{remark}
}


\section{Numerical Experiments}
\label{Sec:NumStudies}

In this section, we carry out several simulation studies to examine the finite-sample performance of the proposed test statistics and compare with permutation-based counterparts.

\subsection{Normal Approximation Accuracy}

We generate two independent random samples $\fX = \{X_1,\dots,X_n\}$ and $\fY = \{Y_1,\dots,Y_m\}$ as follows.
\begin{example}\label{Ex:Hist-R1}
  Generate independent samples: $X_1,\dots,X_n \stsim{iid} \cn(0,\Sigma)$, $Y_1,\dots,Y_m \stsim{iid} \cn(0,\Sigma)$, where $\Sigma = \lrp{\sigma_{ij}} \in \br^{p \times p}$ with $\sigma_{ij} = \rho^{|i-j|}$ and $\rho = 0.5$. We set the sample size ratio $m/n=1$, and consider the setting that $n \in \{25,50,100,200,400\}$ and the data dimensionality $p \in \{25,50,100,200\}$. As for the kernel $k$, we consider the $L_2$-norm $k_{L_2}(x,y) = |x-y|$, the Gaussian kernel multiplied by -1, that is,, $k_G(x,y) = -\exp\lrp{-|x-y|^2/(2\gamma^2)}$ with $\gamma^2 = \text{Median}\{|X_{i_1}-X_{i_2}|^2, |X_i-Y_j|^2, |Y_{j_1}-Y_{j_2}|^2\}$, and the Laplacian kernel multiplied by -1, that is,, $k_L(x,y) = -\exp\lrp{-|x-y|/\gamma}$ with $\gamma = \text{Median}\{|X_{i_1}-X_{i_2}|, |X_i-Y_j|, |Y_{j_1}-Y_{j_2}|\}$. The median heuristic is a  popular way of choosing $\gamma$; see \cite{gretton2012kernel}. 
\end{example}

\begin{figure}[h!]
  \begin{center}
    \includegraphics[width=0.9\textwidth]{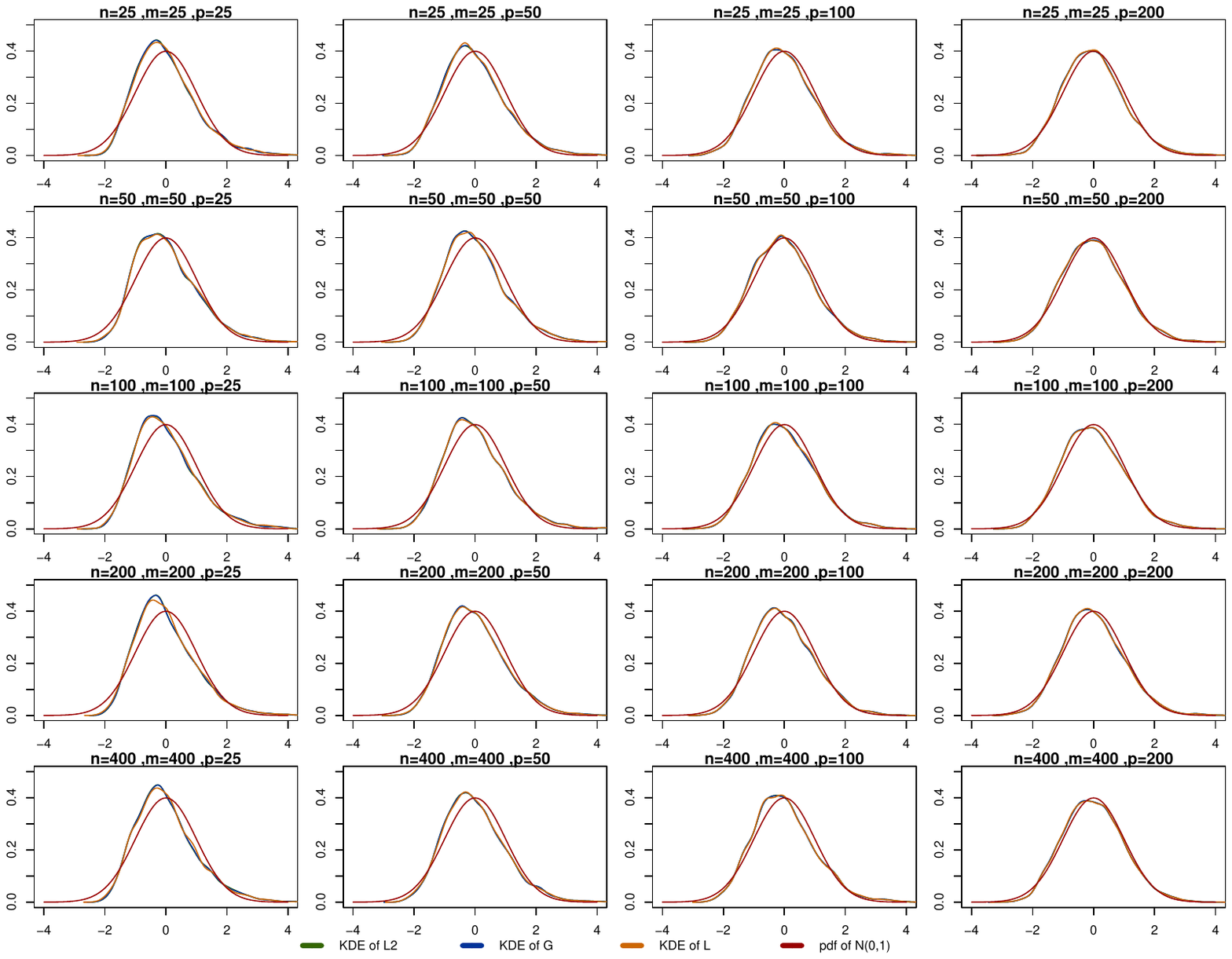}
    \caption{Kernel density estimates of the studentized test statistic $T_{n,m,p}^k$ with different kernels for Example \ref{Ex:Hist-R1} when $m/n=1$. The four columns correspond to different $p$'s and the five rows correspond to different pairs of $(n,m)$.}
    \label{Fig:Hist-R1-r1}
  \end{center}
\end{figure}

Throughout the simulations, our proposed methods are averaged over $5000$ Monte Carlo replications, whereas those of the permutation tests are averaged over $1000$ Monte Carlo replications owing to the high computational cost; see subsequent section. $300$ permutations are conducted for each replication. Given the 5000 replicates of the studentized test statistic $T_{n,m,p}^{k}$, we plot the kernel density estimates (KDE) for the three kernels and the standard normal density function for each combination of sample size and dimension, see Figure \ref{Fig:Hist-R1-r1}. Each row of Figure \ref{Fig:Hist-R1-r1} corresponds to a fixed pair of $(n,m)$ whereas each column represents a fixed choice of $p$. 

As shown in Figure \ref{Fig:Hist-R1-r1}, when $p$ is fixed, the improvement of  normal approximation accuracy is minimal as $N$ increases. However, we do observe significant improvement in the accuracy as $p$ grows for fixed $n=m$. The three kernels correspond to 
very similar empirical distributions suggesting the insensitivity to the kernel choice in terms of size. 
It is worth noticing that  normal approximation is already quite accurate when the sample size and the dimensionality are relatively small, say $N=50$ and $p=100$, and higher accuracy is achieved with larger $N$ and $p$. Such requirements of $N$ and $p$ are usually not demanding in real-world applications, which shows the applicability of the proposed test.

Additional simulation results regarding normal approximation accuracy can be found in online appendices, including the results when the sample sizes $n,m$ are unequal with the difference beyond a constant, and the Kolmogorov-Smirnov distance as well the Wasserstein distance between the standard normal distribution and the empirical distribution of our proposed test statistic under the null. The overall finding from the unbalanced setting is qualitatively similar to what we observe here.

\subsection{Empirical Size}

Under the significance level $\alpha = 0.05$, we reject the null hypothesis if $T_{n,m,p}^k > \Phi(1-\alpha)$. As a comparison, we also consider the permutation test based on the sample MMD studied in \cite{zhu2021interpoint}. For our test statistic, we consider the three kernels introduced in Example \ref{Ex:Hist-R1}, while for the permutation test, we additionally consider the $L_1$-norm $k_{L_1}(x,y) = |x-y|_1$, which is advocated in \cite{zhu2021interpoint}.

In this section, we consider a simulated example that mimics Example 4.1 in \cite{zhu2021interpoint}.

\begin{example}\label{Ex:Size-1}
  Generate independent samples: $X_1,\dots,X_n \stsim{iid} (V^{1/2} \Sigma V^{1/2})^{1/2} Z_X$ and $Y_1,\dots,Y_m \stsim{iid} (V^{1/2} \Sigma V^{1/2})^{1/2} Z_Y$, where $\Sigma = \lrp{\sigma_{ij}} \in \br^{p \times p}$ with $\sigma_{ij} = \rho^{|i-j|}$ and $\rho \in \{0.4, 0.7\}$. 
    We consider the setting that $(n,m) \in \{(25,25), (50,50), (50,100), (100,100), \\ (200,200)\}$ and $p \in \{50, 100\}$. Here, $V$ is { a diagonal matrix with $V_{ii}^{1/2} = 1$ or uniformly drawn from the interval $(1, 5)$}. $Z_X, Z_Y$ are iid copies of $Z$ drawn from the following two distributions:
  \begin{enumerate}[label=(\roman*)]
      \item $Z = \lrp{z_1,\dots,z_p}$ with $z_1,\dots,z_p \stsim{iid} \cn(0,1)$. \label{Ex:Size-Z1}
      \item $Z = \lrp{z_1-1,\dots,z_p-1}$ with $z_1,\dots,z_p \stsim{iid} \text{Exponential}(1)$. \label{Ex:Size-Z2}
  \end{enumerate}
\end{example}

As reported in Table \ref{Tab:Size-1-Cr}, our test exhibit some mild size distortion due to 
the inaccuracy of normal approximation in finite sample. 
 
However, even when $n=m=25$ and $p=50$, the rejection rate is only slightly higher than the nominal level $5\%$, suggesting that our test is practically useful for small sample and moderate dimensional setting. The size distortion tends to increase when the componentwise dependence gets stronger, which matches the theory developed in previous sections; see Proposition~\ref{Prop:BEbd-f}. As we expect, the permutation tests exhibit accurate size, but at the cost of computation. Comparing all three kernels, it seems that no kernel dominates the other in terms of size accuracy.

\begin{table}[h!]
    \noindent
    \scalebox{0.75}{
    \begin{tabular}{c|c|c|c||ccc|cccc||ccc|cccc}
    \hline\hline
        \multirow{3}{*}{$(n,m)$} & \multirow{3}{*}{$p$} & \multirow{3}{*}{$\rho$} & \multirow{3}{*}{$V$} & \multicolumn{7}{c||}{Example \ref{Ex:Size-1}\ref{Ex:Size-Z1}} & \multicolumn{7}{c}{Example \ref{Ex:Size-1}\ref{Ex:Size-Z2}} \\ \cline{5-18}
        & & & & \multicolumn{3}{c|}{Proposed} & \multicolumn{4}{c||}{Permutation} & \multicolumn{3}{c|}{Proposed} & \multicolumn{4}{c}{Permutation} \\ \cline{5-18} 
        & & & & $L_2$ & G & L & $L_2$ & G & L & $L_1$ & $L_2$ & G & L & $L_2$ & G & L & $L_1$ \\ \hline
        \multirow{8}{*}{(25,25)} & \multirow{4}{*}{50} & \multirow{2}{*}{0.4} & Id & 6.98 & 6.92 & 6.80 & 4.60 & 4.60 & 4.20 & 3.80 & 5.66 & 5.80 & 5.80 & 4.50 & 4.60 & 4.50 & 4.90 \\
        & & & Unif & 6.54 & 6.54 & 6.50 & 4.60 & 4.60 & 4.10 & 4.50 & 6.72 & 6.72 & 6.60 & 5.40 & 5.40 & 5.50 & 5.50 \\ \cline{3-18}
        & & \multirow{2}{*}{0.7} & Id & 5.86 & 5.80 & 6.00 & 5.20 & 5.00 & 5.00 & 5.20 & 6.38 & 6.26 & 6.00 & 6.00 & 5.80 & 4.90 & 4.80 \\
        & & & Unif & 6.50 & 6.46 & 6.10 & 5.00 & 4.90 & 4.90 & 4.80 & 6.98 & 6.78 & 6.66 & 4.90 & 5.00 & 5.20 & 5.10 \\ \cline{2-18}
        & \multirow{4}{*}{100} & \multirow{2}{*}{0.4} & Id & 5.52 & 5.54 & 5.52 & 4.90 & 4.90 & 4.60 & 4.20 & 6.12 & 6.04 & 6.08 & 5.00 & 5.00 & 5.10 & 5.70 \\
        & & & Unif & 6.42 & 6.36 & 6.34 & 5.50 & 5.40 & 5.70 & 5.10 & 6.06 & 6.04 & 6.08 & 4.90 & 4.80 & 5.40 & 4.60 \\ \cline{3-18} 
        & & \multirow{2}{*}{0.7} & Id & 6.90 & 6.86 & 7.10 & 5.80 & 5.70 & 5.10 & 4.90 & 6.52 & 6.56 & 6.46 & 5.80 & 5.70 & 5.40 & 5.50 \\
        & & & Unif & 6.60 & 6.60 & 6.74 & 5.30 & 5.40 & 5.50 & 6.30 & 6.62 & 6.62 & 6.70 & 5.50 & 5.40 & 5.50 & 5.40 \\ \hline
        \multirow{8}{*}{(50,50)} & \multirow{4}{*}{50} & \multirow{2}{*}{0.4} & Id & 6.06 & 6.06 & 6.00 & 5.80 & 5.60 & 5.50 & 5.20 & 6.36 & 6.44 & 6.64 & 4.00 & 4.10 & 4.20 & 4.10 \\
        & & & Unif & 5.86 & 5.94 & 5.82 & 4.00 & 3.90 & 3.80 & 4.00 & 6.46 & 6.30 & 6.26 & 4.00 & 3.90 & 3.80 & 3.80 \\ \cline{3-18}
        & & \multirow{2}{*}{0.7} & Id & 7.36 & 7.28 & 7.00 & 5.50 & 5.30 & 5.50 & 5.20 & 7.00 & 6.88 & 6.96 & 5.10 & 5.30 & 5.10 & 5.00 \\
        & & & Unif & 6.56 & 6.66 & 6.68 & 5.90 & 6.30 & 6.30 & 6.50 & 6.46 & 6.36 & 6.22 & 5.10 & 5.30 & 5.40 & 5.20 \\ \cline{2-18}
        & \multirow{4}{*}{100} & \multirow{2}{*}{0.4} & Id & 5.98 & 6.00 & 6.18 & 5.10 & 5.00 & 5.30 & 5.30 & 5.82 & 5.72 & 5.74 & 4.30 & 4.30 & 4.10 & 4.00 \\
        & & & Unif & 6.42 & 6.40 & 6.42 & 4.90 & 4.90 & 4.70 & 4.60 & 5.62 & 5.50 & 5.50 & 5.40 & 5.50 & 5.10 & 5.30 \\ \cline{3-18}
        & & \multirow{2}{*}{0.7} & Id & 6.90 & 6.86 & 6.74 & 5.30 & 5.30 & 5.60 & 4.80 & 6.38 & 6.46 & 6.68 & 5.60 & 5.80 & 6.00 & 6.10 \\
        & & & Unif & 6.68 & 6.68 & 6.50 & 6.30 & 6.00 & 5.50 & 5.70 & 6.10 & 6.12 & 6.10 & 4.20 & 4.20 & 4.60 & 4.60 \\ \hline
        \multirow{8}{*}{(50,100)} & \multirow{4}{*}{50} & \multirow{2}{*}{0.4} & Id & 6.34 & 6.26 & 6.16 & 5.10 & 5.10 & 4.80 & 4.90 & 6.38 & 6.38 & 6.04 & 5.40 & 5.20 & 5.00 & 5.00 \\
        & & & Unif & 6.66 & 6.76 & 6.80 & 6.00 & 6.00 & 5.50 & 6.00 & 6.08 & 6.20 & 5.92 & 4.90 & 5.10 & 5.40 & 4.10 \\ \cline{3-18}
        & & \multirow{2}{*}{0.7} & Id & 6.44 & 6.40 & 6.38 & 5.20 & 5.10 & 4.70 & 5.10 & 6.50 & 6.36 & 6.30 & 4.30 & 4.10 & 4.10 & 4.00 \\
        & & & Unif & 6.74 & 6.62 & 6.30 & 5.20 & 5.20 & 4.90 & 5.10 & 5.98 & 5.94 & 6.02 & 3.30 & 3.30 & 3.50 & 3.80 \\ \cline{2-18}
        & \multirow{4}{*}{100} & \multirow{2}{*}{0.4} & Id & 5.84 & 5.86 & 5.92 & 5.30 & 5.30 & 5.00 & 4.60 & 5.86 & 5.84 & 5.70 & 4.90 & 4.80 & 5.00 & 4.30 \\
        & & & Unif & 6.14 & 6.22 & 6.24 & 5.50 & 5.60 & 6.20 & 6.80 & 6.04 & 6.06 & 6.14 & 5.70 & 5.70 & 5.50 & 5.10 \\ \cline{3-18}
        & & \multirow{2}{*}{0.7} & Id & 6.50 & 6.54 & 6.30 & 5.40 & 5.30 & 4.20 & 5.30 & 6.16 & 6.22 & 6.08 & 5.30 & 5.30 & 5.30 & 5.20 \\
        & & & Unif & 6.48 & 6.38 & 6.42 & 4.70 & 4.40 & 4.70 & 4.10 & 6.64 & 6.72 & 6.72 & 5.20 & 5.40 & 5.40 & 5.30 \\ \hline
        \multirow{8}{*}{(100,100)} & \multirow{4}{*}{50} & \multirow{2}{*}{0.4} & Id & 6.24 & 6.14 & 6.16 & 5.30 & 5.30 & 5.80 & 5.10 & 5.54 & 5.56 & 5.68 & 4.80 & 4.50 & 4.90 & 4.40 \\
        & & & Unif & 6.12 & 6.08 & 6.06 & 5.50 & 5.50 & 5.20 & 4.90 & 7.04 & 7.04 & 6.90 & 6.40 & 6.40 & 6.40 & 5.60 \\ \cline{3-18}
        & & \multirow{2}{*}{0.7} & Id & 6.56 & 6.66 & 6.44 & 4.00 & 4.10 & 4.00 & 4.40 & 6.00 & 6.02 & 6.14 & 5.10 & 5.10 & 5.00 & 5.00 \\
        & & & Unif & 6.40 & 6.36 & 6.42 & 5.10 & 5.00 & 5.30 & 4.90 & 6.66 & 6.60 & 6.70 & 4.50 & 4.30 & 4.70 & 4.40 \\ \cline{2-18}
        & \multirow{4}{*}{100} & \multirow{2}{*}{0.4} & Id & 5.98 & 5.94 & 6.08 & 5.50 & 5.30 & 5.60 & 5.00 & 5.92 & 5.86 & 5.76 & 4.40 & 4.60 & 4.20 & 3.90 \\
        & & & Unif & 6.24 & 6.22 & 6.10 & 5.00 & 5.10 & 5.50 & 5.70 & 6.04 & 6.00 & 5.94 & 4.10 & 4.10 & 3.50 & 3.40 \\ \cline{3-18}
        & & \multirow{2}{*}{0.7} & Id & 6.14 & 6.08 & 6.06 & 4.90 & 4.90 & 5.20 & 5.00 & 6.46 & 6.36 & 6.44 & 4.50 & 4.40 & 4.40 & 4.40 \\
        & & & Unif & 6.60 & 6.66 & 6.54 & 4.50 & 4.30 & 4.30 & 4.70 & 6.76 & 6.88 & 6.82 & 5.30 & 5.40 &  5.50 & 4.60 \\ \hline
        \multirow{8}{*}{(200,200)} & \multirow{4}{*}{50} & \multirow{2}{*}{0.4} & Id & 5.60 & 5.60 & 5.42 & 4.30 & 4.20 & 4.30 & 4.30 & 6.12 & 6.22 & 6.28 & 6.20 & 6.60 & 5.70 & 5.60 \\
        & & & Unif & 6.38 & 6.52 & 6.32 & 5.20 & 5.20 & 4.80 & 5.10 & 6.14 & 6.24 & 6.26 & 5.00 & 5.20 & 5.00 & 5.20 \\ \cline{3-18}
        & & \multirow{2}{*}{0.7} & Id & 6.78 & 6.84 & 6.52 & 4.50 & 4.40 & 4.10 & 4.10 & 6.66 & 6.60 & 6.56 & 5.60 & 5.50 & 4.80 & 5.40 \\
        & & & Unif & 6.56 & 6.68 & 6.64 & 4.10 & 4.00 & 4.40 & 4.20 & 6.94 & 6.88 & 6.74 & 6.10 & 5.90 & 6.10 & 5.80 \\ \cline{2-18}
        & \multirow{4}{*}{100} & \multirow{2}{*}{0.4} & Id & 6.18 & 6.18 & 6.06 & 5.30 & 5.30 & 4.80 & 5.00 & 5.60 & 5.60 & 5.66 & 5.40 & 5.40 & 5.50 & 5.80 \\
        & & & Unif & 6.74 & 6.74 & 6.66 & 5.10 & 5.10 & 5.00 & 4.70 & 5.50 & 5.52 & 5.56 & 4.30 & 4.50 & 4.10 & 3.90 \\ \cline{3-18}
        & & \multirow{2}{*}{0.7} & Id & 7.22 & 7.24 & 7.10 & 5.70 & 5.70 & 5.90 & 5.80 & 6.42 & 6.50 & 6.16 & 4.40 & 4.30 & 3.60 & 3.40 \\
        & & & Unif & 6.56 & 6.50 & 6.32 & 5.70 & 5.40 & 5.40 & 4.40 & 5.56 & 5.54 & 5.76 & 3.90 & 4.10 & 4.30 & 4.60 \\ \hline\hline
    \end{tabular}}
    \caption{Size comparison for Example \ref{Ex:Size-1}. All the empirical sizes are reported in percentage.}\label{Tab:Size-1-Cr}
\end{table}

\subsection{Power Behavior}

Next we investigate the power behavior. The simulated example is adopted from the setting of Example 4.2 in \cite{zhu2021interpoint}. We present the simulation results for the alternative of mean difference in this section, and defer to online appendices the results when two distributions differ in the covariance matrices.

\begin{example}\label{Ex:Power-1}
  Generate independent samples: $X_1,\dots,X_n \stsim{iid} (V^{1/2} \Sigma V^{1/2})^{1/2} Z_X$ and $Y_1,\dots,Y_m \stsim{iid} (V^{1/2} \Sigma V^{1/2})^{1/2} Z_Y + \lrp{0.15\times\mathbf{1}_{\beta p},\mathbf{0}_{(1-\beta)p}}$, where $\Sigma,Z_X,Z_Y$ are defined the same as in Example \ref{Ex:Size-1} and $V$ is chosen as the identity matrix. Here, we fix $\rho = 0.5$, consider $(n,m) \in \{(25,25), (50,50), (100,100), (200,200)\}$, $p \in \{50, 100\}$ and $\beta \in \{0, 0.1, \dots, 1\}$.
\end{example}

\begin{figure}[h!]
  \begin{center}
    \includegraphics[width=0.85\textwidth]{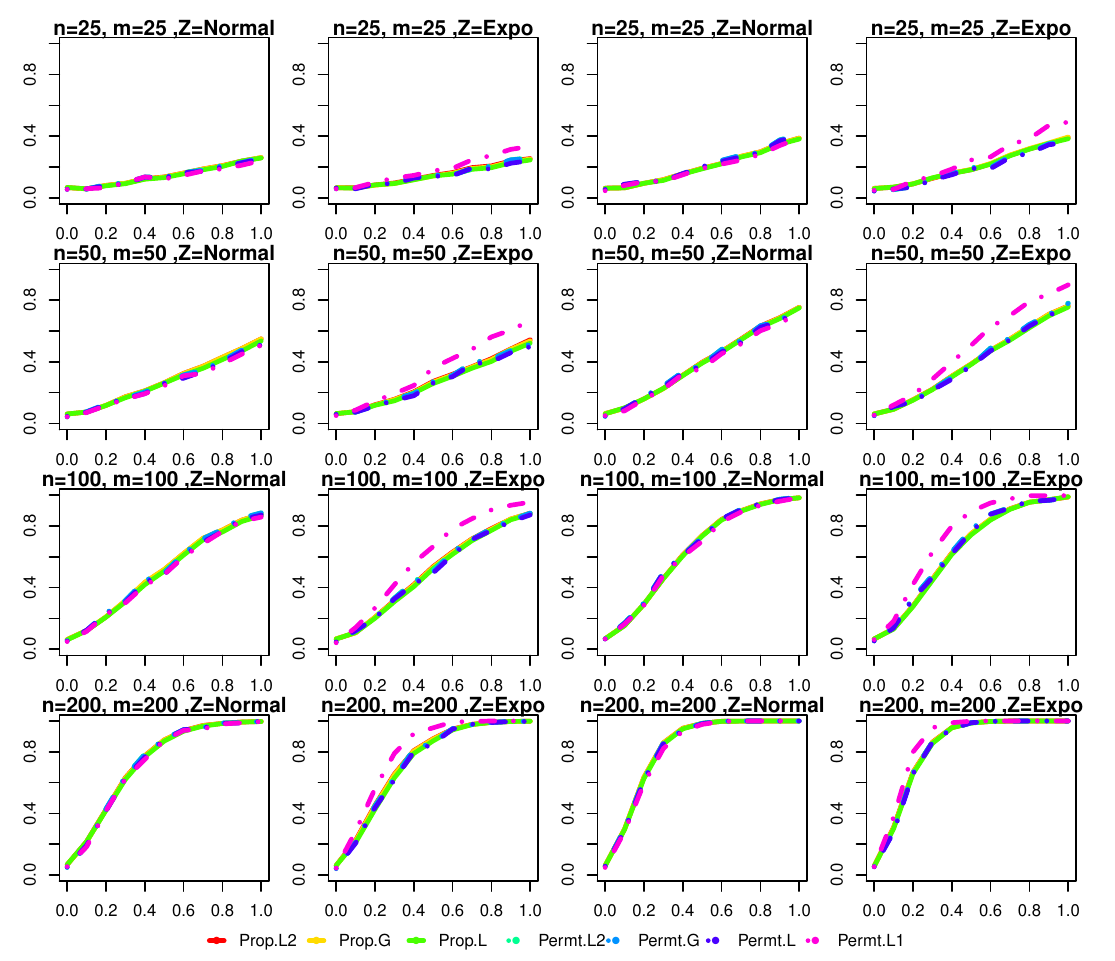}
    \caption{Size-adjusted Power Curves for Example  \ref{Ex:Power-1}. The first two columns correspond to $p=50$ while the last two columns correspond to $p=100$.}\label{Fig:Power-1}
  \end{center}
\end{figure}

We plot the size-adjusted power curves against $\beta$ in Figure \ref{Fig:Power-1}. Note that we only made critical value adjustment in calculating the size-adjusted power for our method, as there is little distortion for permutation-based test. 
As can be seen from Figure \ref{Fig:Power-1}, when there is a mean shift, our test statistic and permutation-based counterpart have almost identical power for all kernels. The use of $L_1$ norm brings some power gain in some cases. As $N$ increases, we do observe a significant improvement in power, regardless of the choice of $p$, which is consistent with our intuition. When $p$ increases from $50$ to $100$, the power increases noticeably for fixed $(n,m)$, as the alternative gets farther away from the null.   

\subsection{Computational Cost}\label{Simu:cost}

One of the major advantage of our proposed method over the permutation test in \cite{zhu2021interpoint} is the computational efficiency. In this section we compare the computational cost per $100$ replications of our method with that of the permutation test under multiple settings; see Table \ref{Tab:Time}. The number of permutations per replication is set to be $300$. As shown in Table \ref{Tab:Time}, it is obvious that our method is much more computationally efficient compared to permutation-based counterpart, which makes up for the slight size distortion of our test under the null.

\begin{table}
    \footnotesize
    \centering
    \begin{tabular}{c|c|c||ccc|cccc}
    \hline\hline
        \multirow{2}{*}{$n$} & \multirow{2}{*}{$m$} & \multirow{2}{*}{$p$} & \multicolumn{3}{c|}{Proposed} & \multicolumn{4}{c}{Permutation} \\ \cline{4-10} 
        & & & $L_2$ & G & L & $L_2$ & G & L & $L_1$ \\ \hline
        \multirow{2}{*}{25} & \multirow{2}{*}{25} &	50 & 0.12 & 0.21 & 0.19 & 20.87 & 22.55 & 22.64 & 22.04 \\ \cline{3-10}
        & & 100	& 0.15 & 0.25 & 0.25 & 33.32 & 35.20 & 35.12 & 36.03 \\ \hline
        \multirow{2}{*}{50} & \multirow{2}{*}{50} & 50 & 0.36 & 0.69 & 0.68 & 82.95 & 94.39 & 95.05 & 86.64 \\ \cline{3-10}
        & & 100	& 0.58 & 1.12 & 1.12 & 148.75 & 158.73 & 159.32 & 152.56 \\ \hline
        \multirow{2}{*}{50} & \multirow{2}{*}{100} & 50 & 0.59 & 1.11 & 1.12 & 155.40 & 171.97 & 173.00 & 168.64 \\ \cline{3-10}
        & & 100 & 0.89 & 1.74 & 1.74 & 263.34 & 278.31 & 276.89 & 283.76 \\ \hline
        \multirow{2}{*}{100} & \multirow{2}{*}{100} & 50 & 1.35 & 2.57 & 2.59 & 328.13 & 369.62 & 367.83 & 332.62 \\ \cline{3-10}
        & & 100	& 2.31 & 4.30 & 4.28 & 576.19 & 625.81 & 626.61 & 601.24 \\ \hline
        \multirow{2}{*}{200} & \multirow{2}{*}{200} & 50 & 5.04 & 9.98 & 10.05 & 1229.59 & 1380.77 & 1387.88 & 1241.33 \\ \cline{3-10}
        & & 100	& 6.69 & 13.43 & 13.24 & 1879.47 & 1989.71 & 1991.56 & 2001.98 \\ \hline\hline
    \end{tabular}
    \caption{Computational cost under multiple settings. All the numerical results are counted in seconds.}\label{Tab:Time}
\end{table}

\section{Discussion}
\label{Sec:Discussion}

In this paper, we have obtained the central limit theorems for studentized sample MMD and derived the explicit rates of convergence
 under the null hypothesis of equal distributions 
when both sample size and dimensionality are diverging. Furthermore, we have also developed a general power theory for the studentized sample MMD  and demonstrated its ability of detecting the difference in distributions. Our proof is built on the argument in \cite{gao2019asymptotic} but we need to develop some new theoretical tools owing to the fact that we are dealing with {\color{black} a general class of kernels} and a two sample U-statistic with high-dimensional observations. In particular, the pooled sample estimate $\cv_{n,m}^{k\ast}(X,Y)$ is proposed and its ratio-consistency as an estimator of $\cv_k^2(Z)$ is shown using a newly developed moment inequality for the multi-sample U-statistics. To deal with {\color{black} a general class of kernels}, we also develop new bounds for the moments of the in-sample and between-sample distances. From a practical viewpoint, our proposed test is simple and easy to implement with much less computational cost compared to permutation test in \cite{zhu2021interpoint}. Finite sample simulations suggest that the size is quite accurate and there is no power loss compared to the permutation-based counterpart. 

As a part of future work, we expect our theory to be useful to the study of asymptotic behavior of sample HSIC in high-dimension, and to test for distributional change in a sequence of high-dimensional data. We leave these topics for future investigation.





\bibliographystyle{agsm}
\bibliography{JMLR-Revision1-arXiv}


\clearpage

\begin{appendices}
\renewcommand\appendixpagename{Supplementary Material}

\appendixpage

{\color{black}
The supplementary material is organized as follows. Appendix \ref{Appendix:examples} contains two additional illustrative examples with the $L_2$ norm and the Laplacian kernel. The statements of the generalized counterparts of Proposition \ref{Prop:PowerAnalysis-1}-\ref{Prop:PowerAnalysis-4} can be found in Appendix \ref{Appendix:generalpower}. In Appendix \ref{Appendix:main}, we present the proofs of all the main results throughout this article, and in Appendix \ref{Appendix:simulation}, we present some additional simulation results. Appendix \ref{Appendix:auxiliary} provides some auxiliary lemmas and their proofs. Finally, Appendix \ref{Appendix:lemma-1}-\ref{Appendix:lemma-3} includes all the lemmas with detailed proofs regarding the proposed test statistic, the rate of convergence, and the power analysis, respectively.}


\section{Illustrative Examples}\label{Appendix:examples}

{\color{black} 
In Section \ref{subsec:example}, we have used the Gaussian kernel as a special example to illustrate the verification of Assumption \ref{Assumpt:uniform-kernel} and condition (\ref{Condition:PowerAnalysis-1}). In this section, we provide two other illustrative examples with the $L_2$ norm and the Laplacian kernel.}

\subsection{Examples to illustrate Assumption \ref{Assumpt:uniform-kernel}}\label{Ex:Assumpt}
\begin{example}\label{Ex:Assumpt-ED}
Consider the $L_2$ norm. For any fixed $p$ and $x,y\in\br^p$, $k^{(p)}(x,y) = |x-y|$ denote the $L_2$ distance between $x$ and $y$. In this case, $f^{(p)}(s) \equiv s$ for any $k^{(p)}\in\cc$. We have $D=[0,\infty)$ and let $D_0 = D$. For simplicity, we use $f$ to represent $f^{(p)}$ and it is trivial that
$f_0(s)=s$, $f_1(s)=1$ and $f_i(s)=0$ for $2\le i\le 7$. By noting the smoothness and continuity of $f$ and its derivatives, Assumption \ref{Assumpt:uniform-kernel}\ref{Assumpt:uniform-kernel-1} is satisfied following the Taylor theorem with the Lagrange form of the remainder. It remains to verify Assumption \ref{Assumpt:uniform-kernel}\ref{Assumpt:uniform-kernel-2}-\ref{Assumpt:uniform-kernel-4}.
\begin{enumerate}[label=(\roman*)]
\setcounter{enumi}{1}
    \item It is trivial that $\sup\limits_{s\in D}|f_i(s)|=1\bone\{i=1\}$ for any $i\geq1$, then for any $s_0\in D_0$, we have $\max\limits_{1\le i\le 7}\sup\limits_{s\in D}|f_i(s)|\cdot|s_0^i|=|s_0|=|f_0(s_0)|$, that is, Assumption \ref{Assumpt:uniform-kernel}\ref{Assumpt:uniform-kernel-2} holds with $\tilde{M}=1$.
    
    \item Note that $\frac{1}{2}|f_1(s)s| = \frac{1}{2}|s| = \frac{1}{2}|f_0(s)|$ and $|-\frac{1}{8}f_1(s)s + \frac{1}{8}f_2(s)s^2| = \frac{1}{8}|s| = \frac{1}{8}|f_0(s)|$, then Assumption \ref{Assumpt:uniform-kernel}\ref{Assumpt:uniform-kernel-3} holds with $\hat{M}=8$.
    
    \item It follows from some simple calculations that $|\frac{1}{16}f_1(s)s - \frac{1}{16}f_2(s)s^2 + \frac{1}{48}f_3(s)s^3| = \frac{1}{16}|s| = \frac{1}{16}|f_0(s)|$ and $|-\frac{5}{128}f_1(s)s + \frac{5}{128}f_2(s)s^2 - \frac{1}{64}f_3(s)s^3 + \frac{1}{384}f_4(s)s^4| = \frac{5}{128}|s| = \frac{5}{128}|f_0(s)|$, which implies that Assumption \ref{Assumpt:uniform-kernel}\ref{Assumpt:uniform-kernel-4} is satisfied with $\hat{M} = \frac{128}{5}$.
\end{enumerate}
\end{example}

\begin{example}\label{Ex:Assumpt-Laplacian}
Now think of the Laplacian kernel. For each fixed $p$ and $x,y\in\br^p$, let $k^{(p)}(x,y) = -\exp\lrp{-\frac{|x-y|}{\gamma}}$ where $\gamma^2$ is a tuning parameter that depends on $p$. Let $f^{(p)}(s) = -\exp\lrp{-\frac{s}{\gamma}}$. Similarly as the previous example, assume that $\cc$ satisfies that for any fixed $p$, there exists some positive constants $\ell,u<\infty$, such that $\blre{|X_1-Y_1|^2}^{1/2}/u \le \gamma \le \blre{|Z_1-Z_2|^2}^{1/2}/\ell$, and again let $D=[0,\infty)$, \\ $D_0=[\blre{|Z_1-Z_2|^2}^{1/2}, \blre{|X_1-Y_1|^2}^{1/2}]$. For simplicity, we use $f$ to represent $f^{(p)}$. The derivatives of $f$ is $f_i(s) = (-1)^{i+1} \frac{1}{\gamma^i}\exp\lrp{-\frac{s}{\gamma}}$ for $0\le i\le 7$. Assumption \ref{Assumpt:uniform-kernel}\ref{Assumpt:uniform-kernel-1} follows from the smoothness, then we verify the remaining conditions.
\begin{enumerate}[label=(\roman*)]
\setcounter{enumi}{1}
     \item It is trivial that $\sup\limits_{s\in D}|f_i(s)| = \frac{1}{\gamma^i}$ for $1\le i\le7$, then with $\blre{|X_1-Y_1|^2}^{1/2}/u \le \gamma \le \blre{|Z_1-Z_2|^2}^{1/2}/\ell$ and $\blre{|Z_1-Z_2|^2}^{1/2} \le s_0 \blre{|X_1-Y_1|^2}^{1/2}$ for any $s_0\in D_0$, then we have $\sup\limits_{s_0\in D_0} \frac{s}{\gamma} \le u$, which implies that Assumption \ref{Assumpt:uniform-kernel}\ref{Assumpt:uniform-kernel-2} holds with $\tilde{M} = \max\{u\exp(u),u^7\exp(u)\}$.
    
    \item It follows some simple calculations that 
    \BEqn
    & & |f_0(s)| 
     =  \exp\lrp{-\frac{s}{\gamma}}, \\
    & & \frac{1}{2}|f_1(s)s| 
     =  \frac{1}{2} \lrp{\frac{s}{\gamma}} \exp\lrp{-\frac{s}{\gamma}}, \\
    & & |-\frac{1}{8}f_1(s)s + \frac{1}{8}f_2(s)s^2| 
     =  \frac{1}{8} \lrp{\lrp{\frac{s}{\gamma}}+\lrp{\frac{s}{\gamma}}^2} \exp\lrp{-\frac{s}{\gamma}},
    \EEqn
    and note that $\inf\limits_{s_0\in D_0} (s_0/\gamma)\geq\ell$
    then Assumption \ref{Assumpt:uniform-kernel}\ref{Assumpt:uniform-kernel-3} holds with $\hat{M}=\max\{2/\ell,8/(\ell+\ell^2)\}$.
    
    \item Note that
    \BEqn
    & & |f_0(s)| 
     =  \exp\lrp{-\frac{s}{\gamma}}, \\
    & & |\frac{1}{16}f_1(s)s - \frac{1}{16}f_2(s)s^2 + \frac{1}{48}f_3(s)s^3| 
     =  \frac{1}{48} \lrp{3\lrp{\frac{s}{\gamma}} + 3\lrp{\frac{s}{\gamma}}^2 + \lrp{\frac{s}{\gamma}}^3} \exp\lrp{-\frac{s}{\gamma}}, \\
    & & |-\frac{5}{128}f_1(s)s + \frac{5}{128}f_2(s)s^2 - \frac{1}{64}f_3(s)s^3 + \frac{1}{384}f_4(s)s^4| \\
    &=&  \frac{1}{384} \lrp{15\lrp{\frac{s}{\gamma}} + 15\lrp{\frac{s}{\gamma}}^2 + 6\lrp{\frac{s}{\gamma}}^3 + \lrp{\frac{s}{\gamma}}^4} \exp\lrp{-\frac{s}{\gamma}},
    \EEqn
    which implies that Assumption \ref{Assumpt:uniform-kernel}\ref{Assumpt:uniform-kernel-4} is satisfied with 
    \begin{equation*}
        \hat{M}
      = \max\lrcp{\frac{48}{3\ell+3\ell^2+\ell^3}, \frac{384}{15\ell+15\ell^2+6\ell^3+\ell^4}}.
    \end{equation*}
\end{enumerate}
\end{example}

\subsection{Examples to illustrate condition (\ref{Condition:PowerAnalysis-1})}\label{Ex:Condition}

\begin{example}\label{Ex:Condition-ED}
We first consider the $L_2$ norm and in this case we have $f(s)=s$, then it is equivalent to verify that $\lrabs{2A_0^{XY}-A_0^{X}-A_0^{Y}} \geq L^{\ast} A_0^{XY}$ for some positive constant $L^{\ast}<\infty$. Note that $A^{XY}=\bbe{|\tilde{X}_1|^2} + \bbe{|\tilde{Y}_1|^2} + |\Delta|^2$, $A^{X} = 2\bbe{|\tilde{X}_1|^2}$ and $A^{Y}=2\bbe{|\tilde{Y}_1|^2}$, it follows that
\begin{equation*}
    f(A_0^{XY}) - f(A_0^{X})
  = \frac{\bbe{|\tilde{Y}_1|^2} - \bbe{|\tilde{X}_1|^2} + |\Delta|^2}{\sqrt{\bbe{|\tilde{X}_1|^2} + \bbe{|\tilde{Y}_1|^2} + |\Delta|^2} + \sqrt{2\bbe{|\tilde{X}_1|^2}}}
\end{equation*}
and
\begin{equation*}
    f(A_0^{XY}) - f(A_0^{Y})
  = \frac{\bbe{|\tilde{X}_1|^2} - \bbe{|\tilde{Y}_1|^2} + |\Delta|^2}{\sqrt{\bbe{|\tilde{X}_1|^2} + \bbe{|\tilde{Y}_1|^2} + |\Delta|^2} + \sqrt{2\bbe{|\tilde{Y}_1|^2}}}.
\end{equation*}
Consequently, we have
\BEqn
& & 2f(A_0^{XY}) - f(A_0^{X}) - f(A_0^{Y}) \\
&=& \frac{2 |\Delta|^2 \sqrt{\bbe{|\tilde{X}_1|^2} + \bbe{|\tilde{Y}_1|^2} + |\Delta|^2}}{\lrp{\sqrt{\bbe{|\tilde{X}_1|^2} + \bbe{|\tilde{Y}_1|^2} + |\Delta|^2} + \sqrt{2\bbe{|\tilde{X}_1|^2}}}\lrp{\sqrt{\bbe{|\tilde{X}_1|^2} + \bbe{|\tilde{Y}_1|^2} + |\Delta|^2} + \sqrt{2\bbe{|\tilde{Y}_1|^2}}}} \\
& & + \frac{|\Delta|^2\lrp{\sqrt{2\bbe{|\tilde{X}_1|^2}} + \sqrt{2\bbe{|\tilde{Y}_1|^2}}}}{\lrp{\sqrt{\bbe{|\tilde{X}_1|^2} + \bbe{|\tilde{Y}_1|^2} + |\Delta|^2} + \sqrt{2\bbe{|\tilde{X}_1|^2}}}\lrp{\sqrt{\bbe{|\tilde{X}_1|^2} + \bbe{|\tilde{Y}_1|^2} + |\Delta|^2} + \sqrt{2\bbe{|\tilde{Y}_1|^2}}}} \\
& & + \frac{\lrp{\sqrt{2\bbe{\tilde{X}_1}^2} - \sqrt{2\bbe{|\tilde{Y}_1|^2}}} \lrp{\bbe{|\tilde{X}_1|^2} - \bbe{|\tilde{Y}_1|^2}}}{\lrp{\sqrt{\bbe{|\tilde{X}_1|^2} + \bbe{|\tilde{Y}_1|^2} + |\Delta|^2} + \sqrt{2\bbe{|\tilde{X}_1|^2}}}\lrp{\sqrt{\bbe{|\tilde{X}_1|^2} + \bbe{|\tilde{Y}_1|^2} + |\Delta|^2} + \sqrt{2\bbe{|\tilde{Y}_1|^2}}}}.
\EEqn
It is trivial that all the three individual terms of $2f(A_0^{XY}) - f(A_0^{X}) - f(A_0^{Y})$ are nonnegative, then 
\begin{equation*}
    \lrabs{2f(A_0^{XY}) - f(A_0^{X}) - f(A_0^{Y})} = 2f(A_0^{XY}) - f(A_0^{X}) - f(A_0^{Y}).
\end{equation*}
Under Assumption \ref{Assumpt:component-dept}\ref{Assumpt:component-dept-3}, we have $C(L_0)\sqrt{p} \le f(A_0^{XY}) \le C(U_0)\sqrt{p}$ where $C(L_0)$ and $C(U_0)$ are both positive constants, and additionally,
\BEqn
& & \lrabs{2f(A_0^{XY}) - f(A_0^{X}) - f(A_0^{Y})} \\
&\geq& \frac{2 |\Delta|^2 \sqrt{\bbe{|\tilde{X}_1|^2} + \bbe{|\tilde{Y}_1|^2} + |\Delta|^2}}{\lrp{\sqrt{\bbe{|\tilde{X}_1|^2} + \bbe{|\tilde{Y}_1|^2} + |\Delta|^2} + \sqrt{2\bbe{|\tilde{X}_1|^2}}}\lrp{\sqrt{\bbe{|\tilde{X}_1|^2} + \bbe{|\tilde{Y}_1|^2} + |\Delta|^2} + \sqrt{2\bbe{|\tilde{Y}_1|^2}}}} \\
& & + \frac{|\Delta|^2\lrp{\sqrt{2\bbe{|\tilde{X}_1|^2}} + \sqrt{2\bbe{|\tilde{Y}_1|^2}}}}{\lrp{\sqrt{\bbe{|\tilde{X}_1|^2} + \bbe{|\tilde{Y}_1|^2} + |\Delta|^2} + \sqrt{2\bbe{|\tilde{X}_1|^2}}}\lrp{\sqrt{\bbe{|\tilde{X}_1|^2} + \bbe{|\tilde{Y}_1|^2} + |\Delta|^2} + \sqrt{2\bbe{|\tilde{Y}_1|^2}}}} \\
&\geq& C(L_0,U_0) \sqrt{p},
\EEqn
with $C(L_0,U_0)$ being some positive constants. Then we may conclude that there exists some positive constant $L^{\ast}=C(L_0,U_0)$, such that $\lrabs{2f(A_0^{XY}) - f(A_0^{X}) - f(A_0^{Y})} \geq L^{\ast} |f(A_0^{XY})|$. This implies that the condition holds naturally for $L_2$ norm.
\end{example}


\begin{example}\label{Ex:Condition-Laplacian}
Finally we consider the Laplacian kernel $f(s)=-\exp\lrp{-\frac{s}{\gamma}}$ with $\bbe{|X_1-Y_1|^2}^{1/2}/u \le \gamma \le \bbe{|Z_1-Z_2|^2}^{1/2}/\ell$. In this case, we have
\begin{equation*}
    \lrp{2f(A_0^{XY}) - f(A_0^{X}) - f(A_0^{Y})}/f(A_0^{XY})
  = 2 - \exp\lrp{-\frac{A_0^{X}-A_0^{XY}}{\gamma}} - \exp\lrp{-\frac{A_0^{Y}-A_0^{XY}}{\gamma}},
\end{equation*}
where 
\begin{equation*}
    A_0^{X}-A_0^{XY}
  = \frac{\bbe{|\tilde{X}_1|^2} - \bbe{|\tilde{Y}_1|^2} - |\Delta|^2}{\sqrt{2\bbe{|\tilde{X}_1|^2}} + \sqrt{\bbe{|\tilde{X}_1|^2} + \bbe{|\tilde{Y}_1|^2} + |\Delta|^2}}
\end{equation*}
and
\begin{equation*}
    A_0^{Y}-A_0^{XY}
  = \frac{\bbe{|\tilde{Y}_1|^2} - \bbe{|\tilde{X}_1|^2} - |\Delta|^2}{\sqrt{2\bbe{|\tilde{Y}_1|^2}} + \sqrt{\bbe{|\tilde{X}_1|^2} + \bbe{|\tilde{Y}_1|^2} + |\Delta|^2}}.
\end{equation*}
Note that at least one of $A_0^{X}-A_0^{XY}$ and $A_0^{Y}-A_0^{XY}$ is negative, then we discuss the following cases to derive the lower bound of $\lrabs{\lrp{2f(A_0^{XY}) - f(A_0^{X}) - f(A_0^{Y})}/f(A_0^{XY})}$.

If $A_0^{X}-A_0^{XY}<0$ and $A_0^{Y}-A_0^{XY}<0$, then it holds for any $\gamma>0$ that 
\begin{equation*}
    \lrp{2f(A_0^{XY}) - f(A_0^{X}) - f(A_0^{Y})}/f(A_0^{XY}) < 0
\end{equation*}
since $\exp(s)>1$ for any $s>0$, and it follows that the condition is naturally satisfied with
\BEqn
& & L^{\ast} \\
&=& \inf\limits_{A_0^{XY}/u \le \gamma \le A_0/\ell} \lrcp{\exp\lrp{-\frac{A_0^{X}-A_0^{XY}}{\gamma}} + \exp\lrp{-\frac{A_0^{Y}-A_0^{XY}}{\gamma}} - 2} \\
&=& \exp\lrp{-\frac{\ell(A_0^{X}-A_0^{XY})}{A_0}} + \exp\lrp{-\frac{\ell(A_0^{Y}-A_0^{XY})}{A_0}} - 2.
\EEqn
    
If $A_0^{X}=A_0^{XY}$ while $A_0^{Y}-A_0^{XY}<0$, it holds for any $\gamma>0$ that
\begin{equation*}
    \lrp{2f(A_0^{XY}) - f(A_0^{X}) - f(A_0^{Y})}/f(A_0^{XY})
  = 1 - \exp\lrp{-\frac{A_0^{Y}-A_0^{XY}}{\gamma}}
  < 0,
\end{equation*}
and the condition is satisfied with
\BEqn
& & L^{\ast} \\
&=& \inf\limits_{A_0^{XY}/u \le \gamma \le A_0/\ell} \lrcp{\exp\lrp{-\frac{A_0^{Y}-A_0^{XY}}{\gamma}} - 1} \\
&=& \exp\lrp{-\frac{\ell(A_0^{Y}-A_0^{XY})}{A_0}} - 1,
\EEqn
and the condition can be verified in a similar way when $A_0^{Y}=A_0^{XY}$ and $A_0^{X}-A_0^{XY}<0$.

Finally, consider the case when $A_0^{Y}-A_0^{XY}>0$ but $A_0^{X}-A_0^{XY}<0$. Let $a = A_0^{X}-A_0^{XY} > 0 $ and $b = A_0^{XY} - A_0^{Y} > 0$, then $\lrp{2f(A_0^{XY}) - f(A_0^{X}) - f(A_0^{Y})}/f(A_0^{XY}) = 2 - h(\gamma)$ where $h(s) = \exp(-\frac{a}{s}) + \exp(\frac{b}{s})$. It is trivial that
\begin{equation*}
    h'(s) = \frac{1}{s^2} \lrp{a\exp\lrp{-\frac{a}{s}} - b\exp\lrp{\frac{b}{s}}},
\end{equation*}
which implies that the minimum of $h(s)$ is attained when $s^{\ast}$ satisfies $a\exp(-\frac{a}{s^{\ast}}) = b\exp(\frac{b}{s^{\ast}})$, and it follows that $h(s^{\ast}) = \frac{a+b}{a} \exp(\frac{b}{s^{\ast}}) > 2$ since in this case we also have 
\begin{equation*}
    b = \bbe{|\tilde{X}_1|^2} - \bbe{|\tilde{Y}_1|^2} + |\Delta|^2 > \bbe{|\tilde{X}_1|^2} - \bbe{|\tilde{Y}_1|^2} - |\Delta|^2 = a.
\end{equation*}
Therefore, $\lrp{2f(A_0^{XY}) - f(A_0^{X}) - f(A_0^{Y})}/f(A_0^{XY}) < 0$ for any $\gamma>0$, then the condition holds with
\BEqn
& & L^{\ast} \\
&=& \inf\limits_{A_0^{XY}/u \le \gamma \le A_0/\ell} \lrcp{\exp\lrp{-\frac{A_0^{X}-A_0^{XY}}{\gamma}} + \exp\lrp{-\frac{A_0^{Y}-A_0^{XY}}{\gamma}} - 2}.
\EEqn
And the same argument is valid when $A_0^{Y}-A_0^{XY}<0$ but $A_0^{X}-A_0^{XY}>0$.

From above, we may conclude that for the Laplacian kernel, there naturally exists some positive constant $L^{\ast}$ such that the condition $\lrabs{2f(A_0^{XY}) - f(A_0^{X}) - f(A_0^{Y})} \geq L^{\ast} |f(A^{XY})|$ holds. 
\end{example}


\section{General Results of Power Analysis}\label{Appendix:generalpower}

\subsection{Generalized version of Proposition \ref{Prop:PowerAnalysis-1}}

\begin{proposition}\label{Prop:GeneralPowerAnalysis-1}
  Suppose that Assumption \ref{Assumpt:uniform-kernel}\ref{Assumpt:uniform-kernel-1}-\ref{Assumpt:uniform-kernel-2}, Assumption \ref{Assumpt:component-dept}\ref{Assumpt:component-dept-1}-\ref{Assumpt:component-dept-3} and Assumption \ref{Assumpt:PowerAnalysis-1} hold, and additionally, for any $k=k^{(p)}\in\cc$, assume that $\blre{k^4(Z_1,Z_2)}<\infty$ and $n/N = \rho + O(1/N^{s})$ as $n,m\rightarrow\infty$ for some $0 < \rho < 1$ and $s>0$. When there exists some positive constant $L^{\ast}<\infty$, such that
  \begin{equation}
     \lrabs{2f(A_0^{XY}) - f(A_0^{X}) - f(A_0^{Y})}
  \geq L^{\ast} |f(A^{XY})|,
  \end{equation}
  it holds that $\blrp{T_{n,m,p}^{k}>C} \rightarrow 1$ as $n,m,p\rightarrow\infty$.
\end{proposition}

\subsection{Generalized version of Proposition \ref{Prop:PowerAnalysis-2}}

\begin{proposition}\label{Prop:GeneralPowerAnalysis-2}
  Suppose that Assumption \ref{Assumpt:uniform-kernel}\ref{Assumpt:uniform-kernel-1}-\ref{Assumpt:uniform-kernel-3} and Assumption \ref{Assumpt:component-dept}\ref{Assumpt:component-dept-1}-\ref{Assumpt:component-dept-4} hold, Assumption \ref{Assumpt:PowerAnalysis-2} holds with $\delta_1\neq\delta_2$ and $\max\{\delta_1,\delta_2\}>0$, and additionally, for any $k=k^{(p)}\in\cc$, assume that $\blre{k^4(Z_1,Z_2)}<\infty$ and $n/N = \rho + O(1/N^{s})$ as $n,m\rightarrow\infty$ for some $0 < \rho < 1$ and $s>0$, then we have
  \begin{enumerate}[label=(\roman*)]
      \item if $0 \le \delta_0 < 2\max\{\delta_1,\delta_2\}-1$, it holds that $\blrp{T_{n,m,p}^{k}>C} \rightarrow 1$ as long as $p = o\lrp{N^{s/(2-2\max\{\delta_1,\delta_2\})}}$.
      \item if $2\max\{\delta_1,\delta_2\}-1 < \delta_0 < \max\{\delta_1,\delta_2\}$, it holds that $\blrp{T_{n,m,p}^{k}>C} \rightarrow 1$ as long as \\ $p = o\lrp{\min\{N^{1/(\delta_0+1-2\max\{\delta_1,\delta_2\})}, N^{s/(1-\delta_0)}\}}$.
  \end{enumerate}
\end{proposition}

\subsection{Generalized version of Proposition \ref{Prop:PowerAnalysis-3}}

\begin{proposition}\label{Prop:GeneralPowerAnalysis-3}
  Suppose that Assumption \ref{Assumpt:uniform-kernel}\ref{Assumpt:uniform-kernel-1}-\ref{Assumpt:uniform-kernel-3} and Assumption \ref{Assumpt:component-dept}\ref{Assumpt:component-dept-1}-\ref{Assumpt:component-dept-4} hold, Assumption \ref{Assumpt:PowerAnalysis-3} holds with $\max\{\delta_3,\delta_4\}>0$, and additionally, for any $k=k^{(p)}\in\cc$, assume that $\blre{k^4(Z_1,Z_2)}<\infty$ and $n/N = \rho + O(1/N^{s})$ as $n,m\rightarrow\infty$ for some $0 < \rho < 1$ and $s>0$, then we have
  \begin{enumerate}[label=(\roman*)]
      \item if  $\delta_0 < \min\{(\delta_3-\delta_4)/\delta_4, \max\{\delta_3/2, \delta_4/(2-\delta_4)\}\}$, it holds that $\blrp{T_{n,m,p}^{k}>C} \rightarrow 1$ as long as $p = o\lrp{\min\{N^{1/(\delta_0+3-2\delta_3)}, N^{s/(1-\delta_0)}\}}$.
      \item if $(\delta_3-\delta_4)/\delta_4 < \delta_0 < \max\{\delta_3/2, \delta_4/(2-\delta_4)\}$, it holds that $\blrp{T_{n,m,p}^{k}>C} \rightarrow 1$ as long as $p = o\lrp{\min\{N^{1/((1-2\delta_4)\delta_0+3-2\delta_4)},N^{s/(1-\delta_0)}\}}$.
  \end{enumerate}
\end{proposition}

\subsection{Generalized version of Proposition \ref{Prop:PowerAnalysis-4}}

\begin{proposition}\label{Prop:GeneralPowerAnalysis-4}
  Suppose that Assumption \ref{Assumpt:uniform-kernel}\ref{Assumpt:uniform-kernel-1}-\ref{Assumpt:uniform-kernel-4} and Assumption \ref{Assumpt:component-dept}\ref{Assumpt:component-dept-1}-\ref{Assumpt:component-dept-4} hold, Assumption \ref{Assumpt:PowerAnalysis-4} holds with $\max\{\delta_5,\delta_6,\delta_7\}>0$, and additionally, for any $k=k^{(p)}\in\cc$, assume that $\blre{k^4(Z_1,Z_2)}<\infty$ and $n/N = \rho + O(1/N^{s})$ as $n,m\rightarrow\infty$ for some $0 < \rho < 1$ and $s>0$, then we have
  \begin{enumerate}[label=(\roman*)]
      \item if $\max\{\delta_5,\delta_6\} > \delta_7$ and $\delta_0 < \max\{\delta_5/(3-\delta_5), \delta_6/(3-\delta_6), (2\delta_7-1)/(3-2\delta_7)\}$, it holds that $\blrp{T_{n,m,p}^{k}>C} \rightarrow 1$ as long as $p = o\lrp{N^{1/((1-4(\delta_5\vee\delta_6))\delta_0+5-2(\delta_5\vee\delta_6))},N^{s/(1-\delta_0)}}$.
      \item if $\delta_7 > \max\{\delta_5,\delta_6\}$ and $\delta_0 < \min\{(1+\delta_5\vee\delta_6-2\delta_7)/(2(\delta_7-\delta_5\vee\delta_6)), \max\{\delta_5/(3-\delta_5), \delta_6/(3-\delta_6), (2\delta_7-1)/(3-2\delta_7)\}$, it holds that $\blrp{T_{n,m,p}^{k}>C} \rightarrow 1$ as long as $p = o\lrp{N^{1/((1-4(\delta_5\vee\delta_6))\delta_0+5-2(\delta_5\vee\delta_6))},N^{s/(1-\delta_0)}}$.
      \item if $\delta_7 > \max\{\delta_5,\delta_6\}$ and $(1+\delta_5\vee\delta_6-2\delta_7)/(2(\delta_7-\delta_5\vee\delta_6)) < \delta_0 < \max\{\delta_5/(3-\delta_5), \delta_6/(3-\delta_6), (2\delta_7-1)/(3-2\delta_7)\}$, it holds that $\blrp{T_{n,m,p}^{k}>C} \rightarrow 1$ as long as $p = o\lrp{\min\{N^{1/((1-4\delta_7)\delta_0+7-4\delta_7)},N^{s/(1-\delta_0)}\}}$.
  \end{enumerate}
\end{proposition}


\section{Proofs of Main Results}
\label{Appendix:main}



\subsection{Proof of Proposition \ref{Prop:est_MMD_simplified}}
It follows from Lemma \ref{Lemma:estMMD_k} that $\ce_{n,m}^{k}(X,Y)$ can be decomposed into two parts, namely, $L_{n,m}^{k}(X,Y)$ and $R_{n,m}^{k}(X,Y)$. By Lemma \ref{Lemma:R=0}, we have that
\BEqn
    R_{n,m}^k(X,Y) 
&=& \frac{4}{n(n-1)m} \sum\limits_{1 \le i_1 < i_2 \le n} \sumjm h^{(2,1)}(X_{i_1},X_{i_2},Y_j) \\
& & + \frac{4}{nm(m-1)} \sumin \sum\limits_{1 \le j_1 < j_2 \le m} h^{(1,2)}(X_i,Y_{j_1},Y_{j_2}) \\
& & + \frac{4}{n(n-1)m(m-1)} \sum\limits_{1 \le i_1 < i_2 \le n} \sum\limits_{1 \le j_1 < j_2 \le m} h^{(2,2)}(X_{i_1},X_{i_2},Y_{j_1},Y_{j_2}) \\
&=& 0,
\EEqn
which further implies that $\ce_{n,m}^{k}(X,Y) = L_{n,m}^{k}(X,Y)$ under both the null and the alternative.

Now we look into $L_{n,m}^{k}(X,Y)$. It follows from Lemma \ref{Lemma:estMMD_k} and Lemma \ref{Lemma:h^(c,d)} that
\BEqn
    L_{n,m}^k(X,Y)
&=& \ce^k(X,Y) 
 + \frac{2}{n}\sum\limits_{i=1}^n h^{(1,0)}(X_i) 
 + \frac{2}{m}\sum\limits_{j=1}^{m} h^{(0,1)}(Y_j) \\
& & + \frac{2}{n(n-1)} \sum\limits_{1 \le i_1 < i_2 \le n} h^{(2,0)}(X_{i_1},X_{i_2})
 + \frac{4}{nm} \sumin \sumjm h^{(1,1)}(X_i,Y_j) \\
& & + \frac{2}{m(m-1)} \sum\limits_{1 \le j_1 < j_2 \le m} h^{(0,2)}(Y_{j_1},Y_{j_2}) \\
&=& \ce^k(X,Y) 
 + \frac{2}{n}\sum\limits_{i=1}^n \lrp{h_{10}(X_i) - \ce^{k}(X,Y)} 
 + \frac{2}{m}\sum\limits_{j=1}^{m} \lrp{h_{01}(Y_j) - \ce^{k}(X,Y)} \\
& & + \frac{2}{n(n-1)} \sum\limits_{1 \le i_1 < i_2 \le n} \lrp{h_{20}(X_{i_1},X_{i_2}) - h_{10}(X_{i_1}) - h_{10}(X_{i_2}) + \ce^{k}(X,Y)} \\
& & + \frac{4}{nm} \sumin \sumjm \lrp{h_{11}(X_i,Y_j) - h_{10}(X_i) - h_{01}(Y_j) + \ce^{k}(X,Y)} \\
& & + \frac{2}{m(m-1)} \sum\limits_{1 \le j_1 < j_2 \le m} \lrp{h_{02}(Y_{j_1},Y_{j_2}) - h_{01}(Y_{j_1}) - h_{01}(Y_{j_2}) + \ce^{k}(X,Y)} \\
&=& 3\ce^{k}(X,Y) - \frac{4}{n}\sum\limits_{i=1}^n h_{10}(X_i) - \frac{4}{m}\sum\limits_{j=1}^m h_{01}(Y_j) + \frac{4}{nm} \sumin \sumjm h_{11}(X_i,Y_j) \\
& & + \frac{2}{n(n-1)} \sum\limits_{1 \le i_1 < i_2 \le n} h_{20}(X_{i_1},X_{i_2}) + \frac{2}{m(m-1)} \sum\limits_{1 \le j_1 < j_2 \le m} h_{02}(Y_{j_1},Y_{j_2}).
\EEqn

If additionally, $X$ and $Y$ are identically distributed, we have $\ce^k(X,Y)=0$. Note that it follows from Lemma \ref{Lemma:h_cd} that
\BEqn
& & h_{10}(X_1) 
 =  \beft{Y}{k(X_1,Y)} + \blre{k(X,Y)} - \beft{X}{k(X_1,X)} - \blre{k(Y,Y')}
 =  0, \\
& & h_{01}(Y_1) 
 =  \beft{X}{k(X,Y_1)} + \blre{k(X,Y)} - \blre{k(X,X')} - \beft{Y}{k(Y_1,Y)}
 =  0,
\EEqn
and consequently, under the null, we additionally have that
\BEqn
& & \ce_{n,m}^{k}(X,Y) = L_{n,m}^{k}(X,Y) \\
&=& \frac{2}{n(n-1)} \sum\limits_{1 \le i_1 < i_2 \le n} h_{20}(X_{i_1},X_{i_2}) + \frac{4}{nm} \sumin \sumjm h_{11}(X_i,Y_j) + \frac{2}{m(m-1)} \sum\limits_{1 \le j_1 < j_2 \le m} h_{02}(Y_{j_1},Y_{j_2}),
\EEqn
which completes the proof.

\subsection{Proof of Proposition \ref{Prop:Var(Lnm)}}
 
For each fixed $p$ and set $X,Y\in\br^p$. If $X$ and $Y$ are identically distributed, then $X,Y,Z$ are all identically distributed. For simplicity, we use $\fZ = \{Z_1,\dots,Z_{n+m}\}$ to represent the pooled sample over $\fX$ and $\fY$, that is, $Z_i = X_i$ for $1\le i\le n$ and $Z_{n+j} = Y_j$ for $1\le j\le m$. Under the null, $\fZ$ is an iid sample. By Proposition \ref{Prop:est_MMD_simplified} and Lemma \ref{Lemma:phi}, we can rewrite $L_{n,m}^{k}(X,Y)$ as
\begin{align*}
   L_{n,m}^k(X,Y)
=& \frac{2}{n(n-1)} \sum\limits_{1 \le i_1 < i_2 \le n} h_{20}(X_{i_1},X_{i_2})
   + \frac{4}{nm} \sumin \sumjm h_{11}(X_i,Y_j) \nonumber \\
 & + \frac{2}{m(m-1)} \sum\limits_{1 \le j_1 < j_2 \le m} h_{02}(Y_{j_1},Y_{j_2}) \\
=& \frac{2}{n(n-1)} \sum\limits_{\ell=1}^{\infty} \lambda_\ell \sum\limits_{1 \le i_1 < i_2 \le n} \phi_\ell(X_{i_1})\phi_\ell(X_{i_2})
   - \frac{2}{nm} \sum\limits_{\ell=1}^{\infty} \lambda_\ell \sumin \sumjm \phi_\ell(X_i)\phi_\ell(Y_j) \nonumber \\
 & + \frac{2}{m(m-1)} \sum\limits_{\ell=1}^{\infty} \lambda_\ell \sum\limits_{1 \le j_1 < j_2 \le m} \phi_\ell(Y_{j_1})\phi_\ell(Y_{j_2}) \\   
=& \sum\limits_{j=1}^{n+m} \sum\limits_{i=1}^{j-1} \sum\limits_{\ell=1}^{\infty} \tau_{i,j,\ell}\phi_\ell(Z_i)\phi_\ell(Z_j) \\
=& \sum\limits_{j=1}^{n+m}\xi_j,
\end{align*}
where
\begin{equation*}
    \xi_j = \sum\limits_{i=1}^{j-1} \sum\limits_{\ell=1}^{\infty} \tau_{i,j,\ell}\phi_\ell(Z_i)\phi_\ell(Z_j),
\end{equation*}
and
\begin{equation*}
    \tau_{i,j,\ell}
  = \left\{
    \begin{array}{ll}
      \displaystyle{\frac{2\lambda_\ell}{n(n-1)}} & 1\le i<j \le n, \vspace{5pt}\\
      \displaystyle{-\frac{2\lambda_\ell}{nm}}   & 1\le i\le n<j\le n+m,  
      \vspace{5pt}\\
      \displaystyle{\frac{2\lambda_\ell}{m(m-1)}} & n+1\le i<j \le n+m.
    \end{array} 
    \right.
\end{equation*}
For $1\le j\le n+m$, define $\cf_j = \sigma(Z_1,\dots,Z_j)$, then $\xi_j$ is $\cf_j$-measurable and $\{\cf_j\}_{j=1}^{n+m}$ forms a filtration. Note that
\begin{align*}
    \blre{\xi_{j+1}|\cf_j}
 =& \blre{\sum\limits_{i=1}^{j} \sum\limits_{\ell=1}^{\infty} \tau_{i,j+1,\ell} \phi_\ell(Z_i) \phi_\ell(Z_{j+1}) | Z_1,\dots,Z_j} \\
 =& \sum\limits_{i=1}^{j} \sum\limits_{\ell=1}^{\infty} \tau_{i,j+1,\ell} \phi_\ell(Z_i) \blre{\phi_\ell(Z_{j+1}) | Z_1,\dots,Z_j} \\
 =& \sum\limits_{i=1}^{j} \sum\limits_{\ell=1}^{\infty} \tau_{i,j+1,\ell} \phi_\ell(Z_i) \blre{\phi_\ell(Z_{j+1})} \\
 =& 0,
\end{align*}
then $\lrcp{\cf_j,\xi_j:1\le j\le n+m}$ is a martingale difference sequence, and $L_{n,m}^{k}(X,Y)$ is thus a martingale with respect to $\cf_{n+m}$. Using the iid property of $Z_i$ and the properties of $\phi_\ell$ stated in Lemma \ref{Lemma:phi}, we have
\begin{align*}
    \blre{\xi_j^2}
 =& \sum\limits_{\ell,r=1}^{\infty} \sum\limits_{1\le i_1,i_2 \le j-1} \tau_{i_1,j,\ell} \tau_{i_2,j,r} \blre{\phi_\ell(Z_{i_1}) \phi_r(Z_{i_2}) \phi_\ell(Z_j) \phi_r(Z_j)} \\
 =& \sum\limits_{\ell,r=1}^{\infty} \sum\limits_{1\le i_1,i_2 \le j-1} \tau_{i_1,j,\ell} \tau_{i_2,j,r} \blre{\phi_\ell(Z_{i_1}) \phi_r(Z_{i_2})} \blre{\phi_\ell(Z_j) \phi_r(Z_j)} \\
 =& \sum\limits_{\ell,r=1}^{\infty} \sum\limits_{1\le i_1,i_2 \le j-1} \tau_{i_1,j,\ell} \tau_{i_2,j,r} \bone\{i_1=i_2,\ell=r\} \\
 =& \sum\limits_{\ell=1}^{\infty} \sum\limits_{i=1}^{j-1} \tau_{i,j,\ell}^2,
\end{align*}
and consequently, by using the property of martingale difference sequence and the results derived in Lemma \ref{Lemma:phi}, we obtain that
\begin{align*}
    \Var(L_{n,m}^{k}(X,Y))
 =& \sum\limits_{j=1}^{n+m}\blre{\xi_j^2}
 =  \sum\limits_{\ell=1}^{\infty} \sum\limits_{j=1}^{n+m} \sum\limits_{i=1}^{j-1} \tau_{i,j,\ell}^2 \\
 =& \sum\limits_{\ell=1}^{\infty} \lrp{\sum\limits_{1\le i<j\le n} + \sum\limits_{i=1}^{n}\sum\limits_{j=n+1}^{n+m} + \sum\limits_{n+1\le i<j\le n+m}}\tau_{i,j,\ell}^2 \\
 =& \sum\limits_{\ell=1}^{\infty} \lrp{\bin{n}{2}\cdot\frac{\lambda_\ell^2}{\bin{n}{2}^2} + nm\cdot\frac{4\lambda_\ell^2}{(nm)^2} + \bin{m}{2}\cdot\frac{\lambda_\ell^2}{\bin{m}{2}^2}} \\
 =& c_{n,m} \sum\limits_{\ell=1}^{\infty} \lambda_\ell^2 
 =  c_{n,m} \cv_k^2(Z),
\end{align*}
which completes the proof.

\subsection{Proof of Proposition \ref{Prop:estDC_k}}
Proposition \ref{Prop:estDC_k} directly follows from Lemma \ref{Lemma:ratio-const-k-H0} and Lemma \ref{Lemma:ratio-const-k-H1}.

\subsection{Proof of Proposition \ref{Prop:ratio-const-k-H0}}
For any $0\le\tau\le1$, it follows from Markov's inequality and the $C_r$ inequality that
\BEqn
& & \blrp{\lrabs{\frac{\cv_{n,m}^{k\ast}(X,Y)}{\cv_k^2(Z)} - 1} > \gamma} \\
&\le& \frac{\blre{\lrabs{\cv_{n,m}^{k\ast}(X,Y) - \cv_k^2(Z)}^{1+\tau}}}{\gamma^{1+\tau}\lrp{\cv_k^2(Z)}^{1+\tau}} \\
&\le& 2^{\tau}
      \lrp{\frac{\blre{\lrabs{\cv_{n,m}^{k\ast}(X,Y) - \blre{\cv_{n,m}^{k\ast}(X,Y)}}^{1+\tau}}}{\gamma^{1+\tau}\lrp{\cv_k^2(Z)}^{1+\tau}} 
      + \frac{\lrabs{\blre{\cv_{n,m}^{k\ast}(X,Y)} - \cv_k^2(Z)}^{1+\tau}}{\gamma^{1+\tau}\lrp{\cv_k^2(Z)}^{1+\tau}}}.
\EEqn
Under the null, it follows from Lemma \ref{Lemma:ratio-const-k-H0} that the second term in the inequality is exactly zero, and it follows from Lemma \ref{Lemma:Consistency_bd_1} that the bound of $\blre{\lrabs{\cv_{n,m}^{k\ast}(X,Y) - \blre{\cv_{n,m}^{k\ast}(X,Y)}}^{1+\tau}}$ reduces to $C(\rho,\tau) N^{-\tau} \blre{\lrabs{\tilde{k}(Z_1,Z_2)}^{2+2\tau}}$ in this case. Therefore, it holds under the null that
\begin{equation*}
    \blrp{\lrabs{\frac{\cv_{n,m}^{k\ast}(X,Y)}{\cv_k^2(Z)} - 1} > \gamma}
\le c(\rho,\tau) \lrp{\frac{\blre{\lrabs{\tilde{k}(Z_1,Z_2)}^{2+2\tau}}}{N^{\tau} \gamma^{1+\tau} \lrp{\cv_k^2(Z)}^{1+\tau}}},
\end{equation*}
and the sample estimate $\cv_{n,m}^{k\ast}(X,Y)$ is ratio-consistent of $\cv_k^2(Z)$ as long as condition (\ref{Equ:consistency_H0}) is satisfied.

\subsection{Proof of Proposition \ref{Prop:ratio-const-k-H1}}
For any $0\le\tau\le1$, it is shown in the proof of Proposition \ref{Prop:ratio-const-k-H0} that
\BEqn
& & \blrp{\lrabs{\frac{\cv_{n,m}^{k\ast}(X,Y)}{\cv_k^2(Z)} - 1} > \gamma} \\
&\le& 2^{\tau} 
      \lrp{\frac{\blre{\lrabs{\cv_{n,m}^{k\ast}(X,Y) - \blre{\cv_{n,m}^{k\ast}(X,Y)}}^{1+\tau}}}{\gamma^{1+\tau}\lrp{\cv_k^2(Z)}^{1+\tau}} 
      + \frac{\lrabs{\blre{\cv_{n,m}^{k\ast}(X,Y)} - \cv_k^2(Z)}^{1+\tau}}{\gamma^{1+\tau}\lrp{\cv_k^2(Z)}^{1+\tau}}}.
\EEqn

Using Lemma \ref{Lemma:Consistency_bd_1} and Lemma \ref{Lemma:Consistency_bd_2}, it holds under the alternative that
\BEqn
& & \blrp{\lrabs{\frac{\cv_{n,m}^{k\ast}(X,Y)}{\cv_k^2(Z)} - 1} > \gamma} \\
&\le& C(\rho,\tau) 
      \left(
      \frac{\blre{\lrabs{\tilde{k}(X_1,X_2)}^{2+2\tau}} + \blre{\lrabs{\tilde{k}(X_1,Y_1)}^{2+2\tau}} + \blre{\lrabs{\tilde{k}(Y_1,Y_2)}^{2+2\tau}}}{N^{\tau} \gamma^{1+\tau} \lrp{\cv_k^2(Z)}^{1+\tau}}
      \right.\\
& & + \left.
      \frac{\lrabs{\ce^k(X,Y)}^{2+2\tau}}{N^{\tau} \gamma^{1+\tau} \lrp{\cv_k^2(Z)}^{1+\tau}}
+ \lrp{\frac{\blre{k^2(Z_1,Z_2)}}{N^s \gamma \cv_k^2(Z)}}^{1+\tau}
      \right).
\EEqn
Consequently, $\cv_{n,m}^{k\ast}(X,Y)$ is a ratio-consistent estimator of $\cv_k^2(Z)$ as long as condition (\ref{Equ:consistency_1}) and condition (\ref{Equ:consistency_2}) are satisfied.

\subsection{Proof of Theorem \ref{Thm:clt}}
If $X$ and $Y$ are identically distributed, when condition (\ref{Equ:clt_1}) is satisfied, it follows from Proposition \ref{Prop:ratio-const-k-H0} that $\cv_{n,m}^{k\ast}(X,Y)$ is a consistent estimator of $\cv_k^2(Z)$. By Slutsky's lemma, it suffices to find the null limiting distribution of
\begin{equation*}
    \tilde{T}_{n,m,p}^{k}
  = \frac{\ce_{n,m}^{k}(X,Y)}{\sqrt{c_{n,m} \cv_{n,m}^{k\ast}(X,Y)}} \cdot \sqrt{\frac{\cv_{n,m}^{k\ast}(X,Y)}{\cv_k^2(Z)}}
  = \frac{\ce_{n,m}^{k}(X,Y)}{\sqrt{c_{n,m} \cv_k^2(Z)}}
\end{equation*}
If follows from Proposition \ref{Prop:est_MMD_simplified} that under the null, $\ce_{n,m}^k(X,Y) = L_{n,m}^{k}(X,Y)$, and consequently,
\begin{equation*}
    \tilde{T}_{n,m,p}^k
  = \frac{L_{n,m}^k(X,Y)}{\sqrt{c_{n,m} \cv_k^2(Z)}}
  = \sum\limits_{j=1}^{n+m} \xi_j.
\end{equation*}
is a martingale with
\begin{equation*}
    \xi_j
  = \sum\limits_{i=1}^{j-1} \tilde{\iota}_{ij} \frac{d^k(Z_i,Z_j)}{\sqrt{\cv_k^2(Z)}},
\end{equation*}
where $\tilde{\iota}_{ij}=\iota_{ij}/\sqrt{c_{n,m}}$ and $\iota_{ij}$ is defined as follows
\begin{equation*}
   \iota_{ij} 
 = \left\{
     \begin{array}{ll}
       -\displaystyle{\frac{2}{n(n-1)}}, & 1 \le i<j \le n  \vspace{3pt} \\
       \displaystyle{\frac{2}{nm}},      & 1 \le i \le n < j \le n+m \vspace{3pt} \\
       -\displaystyle{\frac{2}{m(m-1)}}, & n+1 \le i<j \le n+m
     \end{array}
   \right.    
\end{equation*}

Using Lemma \ref{Lemma:B(1)} and Lemma \ref{Lemma:B(2)}, we have
\BEqn
& & B_{n,m,\tau}^{(1)} \\
&=& \blre{\lrabs{\sum\limits_{j=1}^{n+m} \blre{\xi_j^2 | \cf_{j-1}} - 1}^{1+\tau}} \\
&\le& \frac{2^{\tau}}{\lrp{\cv_k^2(Z)}^{1+\tau}} \lrp{A_{n,m,\tau}^{(1)} \blre{\lrabs{d^k(Z_1,Z_2)}^{2+2\tau}} + A_{n,m,\tau}^{(2)} \lrp{\blre{g^k(Z_1,Z_2,Z_3,Z_4)}}^{(1+\tau)/2}} \\
&\le& \frac{2^{\tau}}{\lrp{\cv_k^2(Z)}^{1+\tau}} \lrp{3^{2+2\tau} A_{n,m,\tau}^{(1)} \blre{\lrabs{\tilde{k}(Z_1,Z_2)}^{2+2\tau}} + A_{n,m,\tau}^{(2)} \lrp{\blre{g^k(Z_1,Z_2,Z_3,Z_4)}}^{(1+\tau)/2}}
\EEqn
and
\begin{equation*}
     B_{n,m,\tau}^{(2)}
 =  \sum\limits_{j=1}^{n+m} \blre{\lrabs{\xi_j}^{2+2\tau}}
 =  A_{n,m,\tau}^{(3)} \frac{\blre{\lrabs{d^k(Z_1,Z_2)}^{2+2\tau}}}{\lrp{\cv_k^2(Z)}^{1+\tau}}
 \le 3^{2+2\tau} A_{n,m,\tau}^{(3)} \frac{\blre{\lrabs{\tilde{k}(Z_1,Z_2)}^{2+2\tau}}}{\lrp{\cv_k^2(Z)}^{1+\tau}}.
\end{equation*}
When condition (\ref{Equ:clt_1}) and condition (\ref{Equ:clt_2}) are satisfied, it follows from Lemma \ref{Lemma:Aorder} that $B_{n,m,\tau}^{(1)} \rightarrow 0$ and $B_{n,m,\tau}^{(2)} \rightarrow 0$, which further implies that
\begin{equation*}
    \sum\limits_{j=1}^{n+m} \blre{\xi_j^2 | \cf_{j-1}} \rightarrow^{p} 1
\end{equation*}
and for any $\varepsilon>0$,
\begin{equation*}
    \sum\limits_{j=1}^{n+m} \blre{\lrabs{\xi_j}^2 \bone\{|\xi_j| > \varepsilon\}} \rightarrow 0.
\end{equation*}
Therefore, $\tilde{T}_{n,m,p}^{k} \rightarrow^d \cn(0,1)$ as $n,m,p\rightarrow\infty$.

\subsection{Proof of Theorem \ref{Thm:BEbd-k}}
Choose the value of $\gamma$ as
\begin{equation*}
    \gamma 
  = \lrp{\frac{\blre{\lrabs{\tilde{k}(Z_1,Z_2)}^{2+2\tau}}}{N^{\tau} \lrp{\cv_k^2(Z)}^{1+\tau}}}^{1/(2+\tau)}.
\end{equation*}
Suppose $\gamma\geq1$, note that it is shown in Lemma \ref{Lemma:cross_eta_2} that $\bbe{g^k(Z_1,Z_2,Z_3,Z_4)}\geq0$, thus with $C(\rho,\tau)=1$, the inequality in Theorem \ref{Thm:BEbd-k} naturally holds since the right-hand side is no smaller than one. It remains to consider the case when $\gamma<1$, and in this case we have that
\begin{equation*}
    \gamma 
  \le \lrp{\frac{\blre{\lrabs{\tilde{k}(Z_1,Z_2)}^{2+2\tau}}}{N^{\tau} \lrp{\cv_k^2(Z)}^{1+\tau}}}^{1/(3+2\tau)}
\end{equation*}
since $\frac{1}{2+\tau} > \frac{1}{3+2\tau}$. It follows from Lemma \ref{Lemma:P2P3} and Lemma \ref{Lemma:P4} that
\begin{equation*}
    P_2 + P_3 + P_4
  \le C(\tau) \gamma
  \le C(\tau) \lrp{\frac{\blre{\lrabs{\tilde{k}(Z_1,Z_2)}^{2+2\tau}}}{N^{\tau} \lrp{\cv_k^2(Z)}^{1+\tau}}}^{1/(3+2\tau)},
\end{equation*}
Also note that it is shown in  Lemma \ref{Lemma:P1} that
\begin{equation*}
    P_1
\le C(\rho,\tau)
    \lrp{\frac{\blre{\lrabs{\tilde{k}(Z_1,Z_2)}^{2+2\tau}}}{N^{\tau}\lrp{\cv_k^2(Z)}^{1+\tau}} + \lrp{\frac{\blre{g^k(Z_1,Z_2,Z_3,Z_4)}}{\lrp{\cv_k^2(Z)}^2}}^{(1+\tau)/2} }^{1/(3+2\tau)},
\end{equation*}
then we may conclude that
\BEqn
& & \sup\limits_{x\in\br} \lrabs{\blrp{T_{n,m,p}^k \le x} - \Phi(x)} \\
&=& 2P_1 + P_2 + P_3 + 2P_4 \\
&\le& C(\rho,\tau)
      \lrp{\frac{\blre{\lrabs{\tilde{k}(Z_1,Z_2)}^{2+2\tau}}}{N^{\tau}\lrp{\cv_k^2(Z)}^{1+\tau}} + \lrp{\frac{\blre{g^k(Z_1,Z_2,Z_3,Z_4)}}{\lrp{\cv_k^2(Z)}^2}}^{(1+\tau)/2} }^{1/(3+2\tau)},
\EEqn
which completes the proof of Theorem \ref{Thm:BEbd-k}.

\subsection{Proof of Proposition \ref{Prop:gk}}

It follows from the definition of $d^k$ and the independence between $Z_1,Z_2$ that
\BEqn
    d^k(Z_1,Z_2)
&=& \tilde{k}(Z_1,Z_2) - \blreft{Z_1}{\tilde{k}(Z_1,Z_2)} - \blreft{Z_2}{\tilde{k}(Z_1,Z_2)} \\
&=& \tilde{k}(Z_1,Z_2) - \blre{\tilde{k}(Z_1,Z_2)|Z_1} - \blre{\tilde{k}(Z_1,Z_2)|Z_2}.
\EEqn
Note that $\{Z_i\}_{i=1}^{N}$ is an iid sample and $\blre{\tilde{k}(Z_1,Z_2)} = 0$, we can expand and organize 
$$d^k(Z_1,Z_2) d^k(Z_1,Z_3) d^k(Z_2,Z_4) d^k(Z_3,Z_4),$$
which leads to the following expression:
\BEqn
& & \blre{g^k(Z_1,Z_2,Z_3,Z_4)} \\
&=& \blre{d^k(Z_1,Z_2) d^k(Z_1,Z_3) d^k(Z_2,Z_4) d^k(Z_3,Z_4)} \\
&=& \blre{\tilde{k}(Z_1,Z_2) \tilde{k}(Z_1,Z_3) \tilde{k}(Z_2,Z_4) \tilde{k}(Z_3,Z_4)} \\
& & - 8\blre{\tilde{k}(Z_1,Z_2) \tilde{k}(Z_1,Z_3) \tilde{k}(Z_2,Z_4) \bbe{\tilde{k}(Z_3,Z_4)|Z_3}} \\
& & + 4\blre{\tilde{k}(Z_1,Z_2) \tilde{k}(Z_1,Z_3) \bbe{\tilde{k}(Z_2,Z_4)|Z_2} \bbe{\tilde{k}(Z_3,Z_4)|Z_3}} \\
& & + 4\blre{\tilde{k}(Z_1,Z_2) \tilde{k}(Z_1,Z_3) \bbe{\tilde{k}(Z_2,Z_4)|Z_4} \bbe{\tilde{k}(Z_3,Z_4)|Z_4}} \\
& & + 4\blre{\tilde{k}(Z_1,Z_2) \bbe{\tilde{k}(Z_1,Z_3)|Z_1} \bbe{\tilde{k}(Z_2,Z_4)|Z_4} \tilde{k}(Z_3,Z_4)} \\
& & - 8\blre{\tilde{k}(Z_1,Z_2) \bbe{\tilde{k}(Z_1,Z_3)|Z_1} \bbe{\tilde{k}(Z_2,Z_4)|Z_4} \bbe{\tilde{k}(Z_3,Z_4)|Z_4}} \\
& & + 2\blre{\bbe{\tilde{k}(Z_1,Z_2)|Z_1} \bbe{\tilde{k}(Z_1,Z_3)|Z_1} \bbe{\tilde{k}(Z_2,Z_4)|Z_4} \bbe{\tilde{k}(Z_3,Z_4)|Z_4}}.
\EEqn
Next we simplify the terms containing the conditional expectations. Again, using the iid property of the $Z_i$'s, we obtain that
\BEqn
& & \blre{\tilde{k}(Z_1,Z_2) \tilde{k}(Z_1,Z_3) \tilde{k}(Z_2,Z_4) \bbe{\tilde{k}(Z_3,Z_4)|Z_3}} \\
&=& \blre{\tilde{k}(Z_1,Z_2) \tilde{k}(Z_1,Z_3) \tilde{k}(Z_2,Z_4) \bbe{\tilde{k}(Z_3,Z_5)|Z_3}} \\
&=& \int \tilde{k}(Z_1,Z_2) \tilde{k}(Z_1,Z_3) \tilde{k}(Z_2,Z_4) \lrp{\int \tilde{k}(Z_3,Z_5) dF_{Z_5}} dF_{Z_1} dF_{Z_2} dF_{Z_3} dF_{Z_4} \\
&=& \int \tilde{k}(Z_1,Z_2) \tilde{k}(Z_1,Z_3) \tilde{k}(Z_2,Z_4) \tilde{k}(Z_3,Z_5) dF_{Z_1} dF_{Z_2} dF_{Z_3} dF_{Z_4} dF_{Z_5} \\
&=& \blre{\tilde{k}(Z_1,Z_2) \tilde{k}(Z_1,Z_3) \tilde{k}(Z_2,Z_4) \tilde{k}(Z_3,Z_5)} \\
&=& \blre{\tilde{k}(Z_1,Z_2) \tilde{k}(Z_1,Z_3) \tilde{k}(Z_2,Z_4) \tilde{k}(Z_4,Z_5)}.
\EEqn
Similarly,
\BEqn
& & \blre{\tilde{k}(Z_1,Z_2) \tilde{k}(Z_1,Z_3) \bbe{\tilde{k}(Z_2,Z_4)|Z_2} \bbe{\tilde{k}(Z_3,Z_4)|Z_3}} \\
&=& \int \tilde{k}(Z_1,Z_2) \tilde{k}(Z_1,Z_3) \lrp{\int \tilde{k}(Z_2,Z_4) dF_{Z_4}} \lrp{\int \tilde{k}(Z_3,Z_5) dF_{Z_5}} dF_{Z_1} dF_{Z_2} dF_{Z_3} \\
&=& \blre{\tilde{k}(Z_1,Z_2) \tilde{k}(Z_1,Z_3) \tilde{k}(Z_2,Z_4) \tilde{k}(Z_4,Z_5)}.
\EEqn
Also, after some computations, we observe that
\BEqn
& & \blre{\tilde{k}(Z_1,Z_2) \tilde{k}(Z_1,Z_3) \bbe{\tilde{k}(Z_2,Z_4)|Z_4} \bbe{\tilde{k}(Z_3,Z_4)|Z_4}} \\
&=& \blre{\tilde{k}(Z_1,Z_2) \tilde{k}(Z_1,Z_3)} \blre{\bbe{\tilde{k}(Z_2,Z_4)|Z_4} \bbe{\tilde{k}(Z_3,Z_4)|Z_4}} \\
&=& \blre{\tilde{k}(Z_1,Z_2) \tilde{k}(Z_1,Z_3)} \blre{\tilde{k}(Z_2,Z_4) \tilde{k}(Z_3,Z_4)} \\
&=& \blre{\tilde{k}(Z_1,Z_2) \tilde{k}(Z_1,Z_3)}^2,
\EEqn
and in a similar manner, we can verify that
\BEqn
& & \blre{\tilde{k}(Z_1,Z_2) \bbe{\tilde{k}(Z_1,Z_3)|Z_1} \bbe{\tilde{k}(Z_2,Z_4)|Z_4} \tilde{k}(Z_3,Z_4)} \\
&=& \blre{\tilde{k}(Z_1,Z_2) \bbe{\tilde{k}(Z_1,Z_3)|Z_1} \bbe{\tilde{k}(Z_2,Z_4)|Z_4} \bbe{\tilde{k}(Z_3,Z_4)|Z_4}} \\
&=& \blre{\bbe{\tilde{k}(Z_1,Z_2)|Z_1} \bbe{\tilde{k}(Z_1,Z_3)|Z_1} \bbe{\tilde{k}(Z_2,Z_4)|Z_4} \bbe{\tilde{k}(Z_3,Z_4)|Z_4}} \\
&=& \blre{\tilde{k}(Z_1,Z_2) \tilde{k}(Z_1,Z_3)}^2.
\EEqn
Consequently, the expression of $\blre{g^k(Z_1,Z_2,Z_3,Z_4)}$ can be simplified as follows,
\begin{eqnarray}
    \blre{g^k(Z_1,Z_2,Z_3,Z_4)}
&=& \blre{\tilde{k}(Z_1,Z_2) \tilde{k}(Z_1,Z_3) \tilde{k}(Z_2,Z_4) \tilde{k}(Z_3,Z_4)} \label{Equ:gk_simple} \\
& & -4 \blre{\tilde{k}(Z_1,Z_2) \tilde{k}(Z_1,Z_3) \tilde{k}(Z_2,Z_4) \tilde{k}(Z_4,Z_5)} \nonumber \\
& & +2 \blre{\tilde{k}(Z_1,Z_2) \tilde{k}(Z_1,Z_3)}^2. \nonumber
\end{eqnarray}
Finally, we expand all the individual terms, and it turns out that with $r=\bbe{k(Z_1,Z_2)}$, it holds that
\BEqn
& & \blre{\tilde{k}(Z_1,Z_2) \tilde{k}(Z_1,Z_3) \tilde{k}(Z_2,Z_4) \tilde{k}(Z_3,Z_4)} \\
&=& \blre{k(Z_1,Z_2) k(Z_1,Z_3) k(Z_2,Z_4) k(Z_3,Z_4)} 
 -  4r \blre{k(Z_1,Z_2) k(Z_1,Z_3) k(Z_2,Z_4)} \\
& & + 4r^2 \blre{k(Z_1,Z_2) k(Z_1,Z_3)} 
 -  r^4, \\[3mm]
& & \blre{\tilde{k}(Z_1,Z_2) \tilde{k}(Z_1,Z_3) \tilde{k}(Z_2,Z_4) \tilde{k}(Z_4,Z_5)} \\
&=& \blre{k(Z_1,Z_2) k(Z_1,Z_3) k(Z_2,Z_4) k(Z_4,Z_5)} 
 -  2r \blre{k(Z_1,Z_2) k(Z_1,Z_3) k(Z_2,Z_4)} \\
& & + r^2 \blre{k(Z_1,Z_2) k(Z_1,Z_3)}, \\[3mm]
& & \blre{\tilde{k}(Z_1,Z_2) \tilde{k}(Z_1,Z_3)}^2 \\
&=& \blre{k(Z_1,Z_2) k(Z_1,Z_3)}^2
 -  2r^2 \blre{k(Z_1,Z_2) k(Z_1,Z_3)}
 +  r^4.
\EEqn
By putting these back to (\ref{Equ:gk_simple}), we obtain the desired result.

\subsection{Proof of Proposition \ref{Prop:BEbd-f}}

For any $k=k^{(p)}\in\cc$, with $\tau=1$ in Theorem \ref{Thm:BEbd-k}, we have
\begin{equation*}
    \sup\limits_{x\in\br} \lrabs{\blrp{T_{n,m,p}^k \le x} - \Phi(x)}
\le C(\rho) \lrcp{\frac{\blre{\tilde{k}^4(Z_1,Z_2)}}{N (\cv_k^2(Z))^2}
      + \frac{\lrabs{\blre{g^k(Z_1,Z_2,Z_3,Z_4)}}}{\lrp{\cv_k^2(Z)}^2}}^{1/5}.
\end{equation*}

From Lemma \ref{Lemma:RateConvg-1}, we have
\begin{equation*}
    |\cv_k^2(Z) - 4c_1^2(A_0) A^{-2} \|\Sigma\|_F^2|
\le C(\tilde{M},U^{\ast},L_0,U_0) c_0^2(A_0) \lrp{\frac{\alpha(p)}{p}}^2,
\end{equation*}
where $c_0(x)$ and $c_1(x)$ are defined in Lemma \ref{Lemma:approx-2}. Then it follows that
\begin{equation*}
    |\bigp{\cv_k^2(Z)}^2 - 16c_1^4(A_0) A^{-4} \|\Sigma\|_F^4|
\le C(\tilde{M},U^{\ast},L_0,U_0) c_0^4(A_0) \lrp{\frac{\alpha(p)}{p}}^3,
\end{equation*}
Using Assumption \ref{Assumpt:uniform-kernel}\ref{Assumpt:uniform-kernel-3} and Assumption \ref{Assumpt:component-dept}\ref{Assumpt:component-dept-4}, there exists some $p_0=p_0(\tilde{M},\hat{M},U^{\ast},L_0,U_0,L_0^{\ast},U_0^{\ast})$, such that for any $p\geq p_0$, it holds that
\begin{equation*}
    \bigp{\cv_k^2(Z)}^2 
\geq C(\tilde{M},\hat{M},U^{\ast},L_0,U_0,L_0^{\ast},U_0^{\ast}) c_0^4(A_0) \lrp{\frac{\alpha(p)}{p}}^2.
\end{equation*}

Then it follows from Lemma \ref{Lemma:RateConvg-2} and Lemma \ref{Lemma:RateConvg-3} that
\BEqn
    \frac{\blre{\tilde{k}^4(Z_1,Z_2)}}{N (\cv_k^2(Z))^2}
&\le& C(\tilde{M},\hat{M},U^{\ast},L_0,U_0,L_0^{\ast},U_0^{\ast})
      \frac{c_1^4(A_0) \lrp{\frac{\alpha(p)}{p}}^2 + c_0^4(A_0) \lrp{\frac{\alpha(p)}{p}}^3}{Nc_0^4(A_0) \lrp{\frac{\alpha(p)}{p}}^2} \\
&\le& C(\tilde{M},\hat{M},U^{\ast},L_0,U_0,L_0^{\ast},U_0^{\ast}) \frac{1}{N}
\EEqn
and
\BEqn
    \frac{\lrabs{\bbe{g^k(Z_1,Z_2,Z_3,Z_4)}}}{(\cv_k^2(Z))^2}
&\le& C(\tilde{M},\hat{M},U^{\ast},L_0,U_0,L_0^{\ast},U_0^{\ast})
      \frac{c_0^4(A_0) \lrp{\frac{\alpha(p)}{p}}^3}{c_0^4(A_0) \lrp{\frac{\alpha(p)}{p}}^2} \\
&\le& C(\tilde{M},\hat{M},U^{\ast},L_0,U_0,L_0^{\ast},U_0^{\ast}) \lrp{\frac{\alpha(p)}{p}},
\EEqn
which jointly imply the desired rate of convergence.

\subsection{Proof of Theorem \ref{Thm:power-k}}

It follows from the definition of $T_{n,m,p}^k$ that
\BEqn
    T_{n,m,p}^k
&=& \frac{\ce_{n,m}^k(X,Y)}{\sqrt{c_{n,m} \cv_{n,m}^{k\ast}(X,Y)}} \\
&=& \frac{\ce^k(X,Y)}{\sqrt{c_{n,m} \cv_k^2(Z)}}
    \lrp{1 + \frac{\ce_{n,m}^k(X,Y) - \ce^k(X,Y)}{\ce^k(X,Y)}}
    \sqrt{\frac{\cv_k^2(Z)}{\cv_{n,m}^{k\ast}(X,Y)}}.
\EEqn
If condition (\ref{Equ:consistency_1}) and condition (\ref{Equ:consistency_2}) are satisfied, then it follows from Proposition \ref{Prop:ratio-const-k-H1} that $\displaystyle{\sqrt{\frac{\cv_k^2(Z)}{\cv_{n,m}^{k\ast}(X,Y)}} \rightarrow^p 1}$.

Under the assumption that $\displaystyle{\frac{N \lrp{\ce^k(X,Y)}^2}{\blre{\lrp{h^k(X_1,X_2,Y_1,Y_2)}^2}} \rightarrow \infty}$, it follows from Lemma \ref{Lemma:diff_E_bound} that
\begin{equation*}
    \blre{\lrabs{\frac{\ce_{n,m}^k(X,Y) - \ce^k(X,Y)}{\ce^k(X,Y)}}^2}
\le C \lrp{\frac{\blre{\lrp{h^k(X_1,X_2,Y_1,Y_2)}^2}}{N \lrp{\ce^k(X,Y)}^2}}
\rightarrow 0,
\end{equation*}
for some universal constant $C>0$, and thus implies that $\displaystyle{\frac{\ce_{n,m}^k(X,Y) - \ce^k(X,Y)}{\ce^k(X,Y)} \rightarrow^p 0}$. Also note that under the assumption $n/N\rightarrow\rho$ it holds that $c_{n,m} = O(N^{-2})$, then by Slutsky's lemma we can obtain that
\begin{equation*}
    \frac{\sqrt{\cv_k^2(Z)}}{N \ce^k(X,Y)} T_{n,m,p}^k \rightarrow^p 1
\end{equation*}
as $n,m,p \rightarrow \infty$. Therefore, under the assumption that $\displaystyle{\frac{N \ce^k(X,Y)}{\sqrt{\cv_k^2(Z)}} \rightarrow \infty}$, we can conclude that $\blrp{T_{n,m,p}^k > C} \rightarrow 1$ as $n,m,p\rightarrow\infty$.

\subsection{Proof of Proposition \ref{Prop:GeneralPowerAnalysis-1}}

It follows from Lemma \ref{Lemma:PowerAnalysis-1}\ref{Lemma:PowerAnalysis-1-1} that, in this case we have
\begin{equation*}
    C_1(\tilde{M},U^{\ast},L^{\ast},L_0,U_0) |c_0(A_0^{XY})|
\le \ce^{k}(X,Y)
\le C_2(\tilde{M},U^{\ast},L^{\ast},L_0,U_0) |c_0(A_0^{XY})|,
\end{equation*}
which by Lemma \ref{Lemma:PowerAnalysis-3}\ref{Lemma:PowerAnalysis-3-3} further implies that
\begin{equation*}
    \bbe{(h^{k}(X_1,X_2,Y_1,Y_2))^2} \le C\lrp{\ce^k(X,Y)}^2.
\end{equation*}
Additionally, from Lemma \ref{Lemma:PowerAnalysis-2}\ref{Lemma:PowerAnalysis-2-1} and Lemma \ref{Lemma:PowerAnalysis-3}\ref{Lemma:PowerAnalysis-3-1}-\ref{Lemma:PowerAnalysis-3-2} that
\begin{equation*}
    C_1(\tilde{M},U^{\ast},L^{\ast},L_0,U_0,\rho) c_0^2(A_0^{XY}) 
\le \cv_k^2(Z) 
\le C_2(\tilde{M},U^{\ast},L^{\ast},L_0,U_0,\rho) c_0^2(A_0^{XY}),
\end{equation*}
and
\BEqn
& & \bbe{k^2(Z_1,Z_2)} \le C(\tilde{M},U^{\ast},L_0,U_0,\rho) c_0^2(A_0^{XY}), \\
& & \max\lrcp{\bbe{\tilde{k}^4(X_1,X_2)}, \bbe{\tilde{k}^4(X_1,Y_1)}, \bbe{\tilde{k}^4(Y_1,Y_2)}}
\le C(\tilde{M},U^{\ast},L_0,U_0,\rho) c_0^4(A_0^{XY}) \lrp{\frac{\alpha(p)}{p}}^2.
\EEqn
It remains to verify the conditions in Theorem \ref{Thm:power-k} with $\tau=1$. Note that
\BEqn
& & \frac{\lrp{\ce^k(X,Y)}^4}{N \lrp{\cv_k^2(Z)}^2}
\le C(\tilde{M},U^{\ast},L^{\ast},L_0,U_0,\rho) \frac{1}{N} 
\rightarrow 0, \\
& & \frac{N^2 \lrp{\ce^k(X,Y)}^2}{\cv_k^2(Z)} 
\geq C(\tilde{M},U^{\ast},L^{\ast},L_0,U_0,\rho) N^2 
\rightarrow \infty, \\
& & \frac{\bbe{k^2(Z_1,Z_2)}}{N^s \cv_k^2(Z)} 
\le C(\tilde{M},U^{\ast},L^{\ast},L_0,U_0,\rho) \frac{1}{N^s} 
\rightarrow 0, \\
& & \frac{N \lrp{\ce^k(X,Y)}^2}{\bbe{(h^k(X_1,X_2,Y_1,Y_2)^2}} 
\geq C \cdot N
\rightarrow \infty,
\EEqn
and
\begin{equation*}
    \frac{\bbe{\tilde{k}^4(X_1,X_2)} + \bbe{\tilde{k}^4(X_1,Y_1)} + \bbe{\tilde{k}^4(Y_1,Y_2)}}{N \lrp{\cv_k^2(Z)}^2} 
\le C(\tilde{M},U^{\ast},L^{\ast},L_0,U_0,\rho) \frac{1}{N} \lrp{\frac{\alpha(p)}{p}}^2
\rightarrow 0,
\end{equation*}
which leads to the proposed result in Proposition \ref{Prop:GeneralPowerAnalysis-1}.

\subsection{Proof of Proposition \ref{Prop:PowerAnalysis-1}}

Proposition \ref{Prop:PowerAnalysis-1} can be viewed as a direct corollary of Proposition \ref{Prop:GeneralPowerAnalysis-1} when $\alpha(p)=O_s(1)$ and $s=1$.

\subsection{Proof of Proposition \ref{Prop:GeneralPowerAnalysis-2}}

It follows from Lemma \ref{Lemma:PowerAnalysis-3}\ref{Lemma:PowerAnalysis-3-1}-\ref{Lemma:PowerAnalysis-3-2} that
\BEqn
& & \bbe{k^2(Z_1,Z_2)} \le C(\tilde{M},U^{\ast},L_0,U_0,\rho) c_0^2(A_0^{XY}), \\
& & \max\lrcp{\bbe{\tilde{k}^4(X_1,X_2)}, \bbe{\tilde{k}^4(X_1,Y_1)}, \bbe{\tilde{k}^4(Y_1,Y_2)}}
\le C(\tilde{M},U^{\ast},L_0,U_0,\rho) c_0^4(A_0^{XY}) \lrp{\frac{\alpha(p)}{p}}^2.
\EEqn
Also note that $p^{2\max\{\delta_1,\delta_2\}-1} \prec p^{\max\{\delta_1,\delta_2\}}$, thus it holds under both cases of Proposition \ref{Prop:GeneralPowerAnalysis-2} that $\alpha(p) = p^{\delta_0} \prec p^{\max\{\delta_1,\delta_2\}}$, and it follows from Lemma \ref{Lemma:PowerAnalysis-1}\ref{Lemma:PowerAnalysis-1-2} that
\BEqn
& & C_1(\tilde{M},\hat{M},U^{\ast},L_0,U_0,L_2,U_2) |c_0(A_0^{XY})| p^{\max\{\delta_1,\delta_2\}-1} \\
&\le& \ce^{k}(X,Y) \\
&\le& C_2(\tilde{M},U^{\ast},L_0,U_0,L_2,U_2) |c_0(A_0^{XY})| p^{\max\{\delta_1,\delta_2\}-1}
\EEqn
when $p\rightarrow\infty$.

\begin{enumerate}[label=(\roman*)]
    \item when $\delta_0 < 2\max\{\delta_1,\delta_2\}-1$, $p^{2\max\{\delta_1,\delta_2\}-1}$ dominates $\alpha(p)$ in this case, and this implies that
    \begin{equation*}
        \bbe{(h^{k}(X_1,X_2,Y_1,Y_2))^2} \le C\lrp{\ce^k(X,Y)}^2.
    \end{equation*}
    Also, from Lemma \ref{Lemma:PowerAnalysis-2}\ref{Lemma:PowerAnalysis-2-2}, it holds when $p\rightarrow\infty$ that 
    \BEqn
    & & C_1(\tilde{M},\hat{M},U^{\ast},L_0,U_0,L_2,U_2,\rho) c_0^2(A_0^{XY}) p^{2\max\{\delta_1,\delta_2\}-2} \\
    &\le& \cv_k^2(Z) \\
    &\le& C_2(\tilde{M},U^{\ast},L_0,U_0,L_2,U_2,\rho) c_0^2(A_0^{XY}) p^{2\max\{\delta_1,\delta_2\}-2}.
    \EEqn
    
    Now we verify the conditions in Theorem \ref{Thm:power-k} with $\tau=1$. It is trivial that
    \BEqn
    & & \frac{\lrp{\ce^k(X,Y)}^4}{N \lrp{\cv_k^2(Z)}^2}
    \le C(\tilde{M},\hat{M},U^{\ast},L_0,U_0,L_2,U_2,\rho) \frac{1}{N} 
    \rightarrow 0, \\
    & & \frac{N^2 \lrp{\ce^k(X,Y)}^2}{\cv_k^2(Z)} 
    \geq C(\tilde{M},\hat{M},U^{\ast},L_0,U_0,L_2,U_2,\rho) N^2 
    \rightarrow \infty, \\
    & & \frac{N \lrp{\ce^k(X,Y)}^2}{\bbe{(h^k(X_1,X_2,Y_1,Y_2)^2}} 
    \geq C \cdot N
    \rightarrow \infty,
    \EEqn
    Also note that
    \BEqn
    & & \frac{\bbe{\tilde{k}^4(X_1,X_2)} + \bbe{\tilde{k}^4(X_1,Y_1)} + \bbe{\tilde{k}^4(Y_1,Y_2)}}{N \lrp{\cv_k^2(Z)}^2} \\
    &\le& C(\tilde{M},\hat{M},U^{\ast},L_0,U_0,L_2,U_2,\rho) \frac{1}{Np^{2(2\max\{\delta_1,\delta_2\}-1-\delta_0)}} 
    \rightarrow 0
    \EEqn
    since $\delta_0 < 2\max\{\delta_1,\delta_2\}-1$. Finally, we have that
    \begin{equation*}
         \frac{\bbe{k^2(Z_1,Z_2)}}{N^s \cv_k^2(Z)} 
    \le C(\tilde{M},\hat{M},U^{\ast},L_0,U_0,L_2,U_2,\rho) \frac{p^{2-2\max\{\delta_1,\delta_2\}}}{N^s}
    \rightarrow 0
    \end{equation*}
    as long as $p = o\lrp{N^{s/(2-2\max\{\delta_1,\delta_2\})}}$, which completes the proof in this case,
    \item when $2\max\{\delta_1,\delta_2\}-1 < \delta_0 < \max\{\delta_1,\delta_2\}$, it holds that $p^{2\max\{\delta_1,\delta_2\}-1} \prec \alpha(p) \prec p^{\max\{\delta_1,\delta_2\}}$. Thus the upper bound of $\bbe{(h^k(X_1,X_2,Y_1,Y_2))^2}$ in this case is
    \begin{equation*}
        \bbe{(h^k(X_1,X_2,Y_1,Y_2))^2} 
    \le C(\tilde{M},U^{\ast},L_0,U_0) c_0^2(A_0^{XY}) \lrp{\frac{\alpha(p)}{p}},
    \end{equation*}
    and it follows from Lemma \ref{Lemma:PowerAnalysis-2}\ref{Lemma:PowerAnalysis-2-2} that
    \BEqn
    & & C_1(\tilde{M},\hat{M},U^{\ast},L_0,U_0,L_2,U_2,\rho) c_0^2(A_0^{XY}) \lrp{\frac{\alpha(p)}{p}} \\
    &\le& \cv_k^2(Z) \\
    &\le& C_2(\tilde{M},U^{\ast},L_0,U_0,L_2,U_2,\rho) c_0^2(A_0^{XY}) \lrp{\frac{\alpha(p)}{p}}.
    \EEqn
    Again, we verify the conditions in Theorem \ref{Thm:power-k} with $\tau=1$, and obtain that
    \begin{equation*}
        \frac{\lrp{\ce^k(X,Y)}^4}{N \lrp{\cv_k^2(Z)}^2}
    \le C(\tilde{M},\hat{M},U^{\ast},L_0,U_0,L_2,U_2,\rho) \frac{1}{N} \lrp{\frac{p^{2\max\{\delta_1,\delta_2\}-1}}{\alpha(p)}}^2
    \rightarrow 0,
    \end{equation*}
    and
    \begin{equation*}
        \frac{\bbe{\tilde{k}^4(X_1,X_2)} + \bbe{\tilde{k}^4(X_1,Y_1)} + \bbe{\tilde{k}^4(Y_1,Y_2)}}{N \lrp{\cv_k^2(Z)}^2} 
    \le C(\tilde{M},\hat{M},U^{\ast},L_0,U_0,L_2,U_2,\rho) \frac{1}{N}
    \rightarrow 0.
    \end{equation*}
    Also note that
    \BEqn
    & & \frac{\bbe{k^2(Z_1,Z_2)}}{N^s \cv_k^2(Z)} 
    \le C(\tilde{M},\hat{M},U^{\ast},L_0,U_0,L_2,U_2,\rho) \frac{1}{N^s} \lrp{\frac{p}{\alpha(p)}}
    \rightarrow 0
    \EEqn
    as long as $\frac{p}{\alpha(p)} = p^{1-\delta_0} =  o\lrp{N^s}$
    and
    \BEqn
    & & \frac{N^2 \lrp{\ce^k(X,Y)}^2}{\cv_k^2(Z)} 
    \geq C(\tilde{M},\hat{M},U^{\ast},L_0,U_0,L_2,U_2,\rho) N^2 \lrp{\frac{p^{2\max\{\delta_1,\delta_2\}-1}}{\alpha(p)}}
    \rightarrow \infty, \\
    & & \frac{N \lrp{\ce^k(X,Y)}^2}{\bbe{(h^k(X_1,X_2,Y_1,Y_2)^2}} 
    \geq C(\tilde{M},\hat{M},U^{\ast},L_0,U_0,L_2,U_2,\rho) N \lrp{\frac{p^{2\max\{\delta_1,\delta_2\}-1}}{\alpha(p)}}
    \rightarrow \infty
    \EEqn
    as long as $p^{1-2\max\{\delta_1,\delta_2\}}\alpha(p) = p^{\delta_0+1-2\max\{\delta_1,\delta_2\}} \prec N$. By putting these conditions together, we obtain the desired result.
\end{enumerate}

\subsection{Proof of Proposition \ref{Prop:PowerAnalysis-2}}

It suffices to consider the case when $\delta_0=0$ and $s=1$ in Proposition \ref{Prop:GeneralPowerAnalysis-2}. Note that $2\max\{\delta_1,\delta_2\}-1>\delta_0$ holds if $\max\{\delta_1,\delta_2\}>1/2$, and in this case the asymptotic power one is obtained as long as $p=o\lrp{N^{1/(2-2\max\{\delta_1,\delta_2\})}}$. On the other hand, when $\max\{\delta_1,\delta_2\}<1/2$, it holds by natural that $2\max\{\delta_1,\delta_2\}-1<\delta_0<\max\{\delta_1,\delta_2\}$, and the constraint that $p=o\lrp{\min\{N^{1/(\delta_0+1-2\max\{\delta_1,\delta_2\})},N^{s/(1-\delta_0)}\}}$ is simplified as $p = o\lrp{N}$. Consequently, the constraints under both cases can be unified as $p=o\lrp{N^{1/(2-2\max\{\delta_1,\delta_2,1/2\})}}$, which completes the proof.

\subsection{Proof of Proposition \ref{Prop:GeneralPowerAnalysis-3}}

From Lemma \ref{Lemma:PowerAnalysis-1}\ref{Lemma:PowerAnalysis-1-3}, when $\max\{\delta_3,\delta_4\}>0$ and $\alpha(p)=o\lrp{p^{\max\{\delta_3/2,\delta_4/(2-\delta_4)\}}}$, we have that
\BEqn
& & C_1(\tilde{M},\hat{M},U^{\ast},L_0,U_0,L_3,U_3) |c_0(A_0^{XY})| \max\{p^{\delta_3-2}, (\alpha(p))^{\delta_4}p^{\delta_4-2}\} \\
&\le& \ce^{k}(X,Y) \\
&\le& C_2(\tilde{M},U^{\ast},L_0,U_0,L_3,U_3) |c_0(A_0^{XY})| \max\{p^{\delta_3-2}, (\alpha(p))^{\delta_4}p^{\delta_4-2}\}.
\EEqn
Since $\max\{\delta_3,\delta_4\}<1$ and $1 \le \alpha(p) \prec p$, then it always holds in this case that $p^{2\delta_3-3} = o\lrp{\alpha(p)}$ and $p^{2\delta_4-3} = o\lrp{(\alpha(p))^{1-2\delta_4}}$, which implies by Lemma \ref{Lemma:PowerAnalysis-3}\ref{Lemma:PowerAnalysis-3-3} that
\begin{equation*}
    \bbe{(h^k(X_1,X_2,Y_1,Y_2))^2} 
\le C(\tilde{M},U^{\ast},L_0,U_0) c_0^2(A_0^{XY}) \lrp{\frac{\alpha(p)}{p}}.
\end{equation*}

Also, from Lemma \ref{Lemma:PowerAnalysis-2}\ref{Lemma:PowerAnalysis-2-3} we have that
\BEqn
& & C_1(\tilde{M},\hat{M},U^{\ast},L_0,U_0,L_3,U_3,\rho) c_0^2(A_0^{XY}) \lrp{\frac{\alpha(p)}{p}} \\
&\le& \cv_k^2(Z) \\
&\le& C_2(\tilde{M},U^{\ast},L_0,U_0,L_3,U_3,\rho) c_0^2(A_0^{XY}) \lrp{\frac{\alpha(p)}{p}}.
\EEqn
and from Lemma \ref{Lemma:PowerAnalysis-3}\ref{Lemma:PowerAnalysis-3-1}-\ref{Lemma:PowerAnalysis-3-2} we have that
\BEqn
& & \bbe{k^2(Z_1,Z_2)} \le C(\tilde{M},U^{\ast},L_0,U_0,\rho) c_0^2(A_0^{XY}), \\
& & \max\lrcp{\bbe{\tilde{k}^4(X_1,X_2)}, \bbe{\tilde{k}^4(X_1,Y_1)}, \bbe{\tilde{k}^4(Y_1,Y_2)}}
\le C(\tilde{M},U^{\ast},L_0,U_0,\rho) c_0^4(A_0^{XY}) \lrp{\frac{\alpha(p)}{p}}^2.
\EEqn
Hence it holds that
\begin{equation*}
    \frac{\bbe{\tilde{k}^4(X_1,X_2)} + \bbe{\tilde{k}^4(X_1,Y_1)} + \bbe{\tilde{k}^4(Y_1,Y_2)}}{N \lrp{\cv_k^2(Z)}^2} 
\le C(\tilde{M},\hat{M},U^{\ast},L_0,U_0,L_3,U_3,\rho) \frac{1}{N} 
\rightarrow 0
\end{equation*}
and
\BEqn
& & \frac{\bbe{k^2(Z_1,Z_2)}}{N^s \cv_k^2(Z)}
\le C(\tilde{M},\hat{M},U^{\ast},L_0,U_0,L_3,U_3,\rho) \frac{1}{N^s} \lrp{\frac{p}{\alpha(p)}}
\rightarrow 0
\EEqn
as long as $\frac{p}{\alpha(p)} = p^{1-\delta_0} = o\lrp{N^s}$. 

Now it remains to verify the remaining conditions case by case.
\begin{enumerate}[label=(\roman*)]
    \item when $\delta_0 < \min\lrcp{(\delta_3-\delta_4)/\delta_4, \max\{\delta_3/2, \delta_4/(2-\delta_4)\}}$, we have $\alpha(p) \prec p^{(\delta_3-\delta_4)/\delta_4}$, and it follows
    \BEqn
    & & C_1(\tilde{M},\hat{M},U^{\ast},L_0,U_0,L_3,U_3) |c_0(A_0^{XY})| p^{\delta_3-2} \\
    &\le& \ce^{k}(X,Y) \\
    &\le& C_2(\tilde{M},U^{\ast},L_0,U_0,L_3,U_3) |c_0(A_0^{XY})| p^{\delta_3-2},
    \EEqn
    and it follows that
    \begin{equation*}
        \frac{\lrp{\ce^k(X,Y)}^4}{N \lrp{\cv_k^2(Z)}^2}
    \le C(\tilde{M},\hat{M},U^{\ast},L_0,U_0,L_3,U_3,\rho) \frac{1}{N} \lrp{\frac{1}{\alpha(p)p^{3-2\delta_3}}}^2
    \rightarrow 0,
    \end{equation*}
    and
    \BEqn
    & & \frac{N^2 \lrp{\ce^k(X,Y)}^2}{\cv_k^2(Z)} 
    \geq C(\tilde{M},\hat{M},U^{\ast},L_0,U_0,L_3,U_3,\rho) N^2 \lrp{\frac{1}{\alpha(p)p^{3-2\delta_3}}}
    \rightarrow \infty, \\
    & & \frac{N \lrp{\ce^k(X,Y)}^2}{\bbe{(h^k(X_1,X_2,Y_1,Y_2)^2}} 
    \geq C(\tilde{M},\hat{M},U^{\ast},L_0,U_0,L_3,U_3,\rho) N \lrp{\frac{1}{\alpha(p)p^{3-2\delta_3}}}
    \rightarrow \infty
    \EEqn
    as long as $\alpha(p)p^{3-2\delta_3} = p^{\delta_0+3-2\delta_3} \prec N$ and thus completes the proof.
    
    \item when $(\delta_3-\delta_4)/\delta_4 < \delta_0 < \max\{\delta_3/2, \delta_4/(2-\delta_4)\}$, we have $p^{(\delta_3-\delta_4)/\delta_4} \prec \alpha(p)$, then we further obtain that
    \BEqn
    & & C_1(\tilde{M},\hat{M},U^{\ast},L_0,U_0,L_3,U_3) |c_0(A_0^{XY})| (\alpha(p))^{\delta_4}p^{\delta_4-2} \\
    &\le& \ce^{k}(X,Y) \\
    &\le& C_2(\tilde{M},U^{\ast},L_0,U_0,L_3,U_3) |c_0(A_0^{XY})| (\alpha(p))^{\delta_4}p^{\delta_4-2},
    \EEqn
    and it follows that
    \BEqn
        \frac{\lrp{\ce^k(X,Y)}^4}{N \lrp{\cv_k^2(Z)}^2}
    &\le& C(\tilde{M},\hat{M},U^{\ast},L_0,U_0,L_3,U_3,\rho) \frac{(\alpha(p))^{4\delta_4-2} p^{4\delta_4-6}}{N} \\
    &\le& C(\tilde{M},\hat{M},U^{\ast},L_0,U_0,L_3,U_3,\rho) \frac{1}{N} \lrp{\frac{1}{\alpha(p) p}}^{4-4\delta_4} \lrp{\frac{\alpha(p)}{p}}^2 \\
    &\rightarrow& 0,
    \EEqn
    and
    \BEqn
    & & \frac{N^2 \lrp{\ce^k(X,Y)}^2}{\cv_k^2(Z)} 
    \geq C(\tilde{M},\hat{M},U^{\ast},L_0,U_0,L_3,U_3,\rho) N^2 \lrp{\frac{1}{(\alpha(p))^{1-2\delta_4} p^{3-2\delta_4}}}
    \rightarrow \infty, \\
    & & \frac{N \lrp{\ce^k(X,Y)}^2}{\bbe{(h^k(X_1,X_2,Y_1,Y_2)^2}} 
    \geq C(\tilde{M},\hat{M},U^{\ast},L_0,U_0,L_3,U_3,\rho) N \lrp{\frac{1}{(\alpha(p))^{1-2\delta_4} p^{3-2\delta_4}}}
    \rightarrow \infty
    \EEqn
    as long as $(\alpha(p))^{1-2\delta_4} p^{3-2\delta_4} = p^{(1-2\delta_4)\delta_0 + 3 - 2\delta_4} \prec N$ and thus completes the proof.
\end{enumerate}

\subsection{Proof of Proposition \ref{Prop:PowerAnalysis-3}}

Based on Proposition \ref{Prop:GeneralPowerAnalysis-3}, we can obtain the simplified results when $s=1$ and $\delta_0=0$. For the first case, with $\max\{\delta_3,\delta_4\}>0$, $\delta_0<\min\{(\delta_3-\delta_4)/\delta_4, \max\{\delta_3/2,\delta_4/(2-\delta_4)\}\}$ holds when $\delta_3>\delta_4$ and the nontrivial power is thus obtained when $p=o\lrp{N^{1/(3-2\delta_3)}}$. As for the second case, $(\delta_3-\delta_4)/\delta_4 < \delta_0 < \max\{\delta_3/2,\delta_4/(2-\delta_4)\}$ holds if and only if $\delta_3<\delta_4$, and with $\delta_0=0$, the constraint for this case can be simplified as $p=o\lrp{N^{1/(3-2\delta_4)}}$. In summary, the unified constraint for both cases is $p=o\lrp{N^{1/(3-2\max\{\delta_3,\delta_4\})}}$.

\subsection{Proof of Proposition \ref{Prop:GeneralPowerAnalysis-4}}

From Lemma \ref{Lemma:PowerAnalysis-1}\ref{Lemma:PowerAnalysis-1-4}, when
\begin{equation*}
    \alpha(p)=o\lrp{p^{\max\{\delta_5/(3-2\delta_5),\delta_6/(3-2\delta_6),(2\delta_7-1)/(3-2\delta_7)\}}},
\end{equation*}
we have that
\begin{equation*}
  \begin{array}{l}
     \hspace{1em}C_1(\tilde{M},\hat{M},U^{\ast},L_0,U_0,L_4,U_4) |c_0(A_0^{XY})| \max\{(\alpha(p))^{2\delta_5}p^{\delta_5-3}, (\alpha(p))^{2\delta_6}p^{\delta_6-3}, (\alpha(p))^{2\delta_7}p^{2\delta_7-4}\} \vspace{2mm}\\
     \le \ce^{k}(X,Y) \vspace{2mm}\\
     \le C_2(\tilde{M},U^{\ast},L_0,U_0,L_4,U_4) |c_0(A_0^{XY})| \max\{(\alpha(p))^{2\delta_5}p^{\delta_5-3}, (\alpha(p))^{2\delta_6}p^{\delta_6-3}, (\alpha(p))^{2\delta_7}p^{2\delta_7-4}\}.
  \end{array}
\end{equation*}
Also, in this case, it follows from Lemma \ref{Lemma:PowerAnalysis-2}\ref{Lemma:PowerAnalysis-2-4} that
\BEqn
& & C_1(\tilde{M},\hat{M},U^{\ast},L_0,U_0,\rho) c_0^2(A_0^{XY}) \lrp{\frac{\alpha(p)}{p}} \\
&\le& \cv_k^2(Z) \\
&\le& C_2(\tilde{M},U^{\ast},L_0,U_0,\rho) c_0^2(A_0^{XY}) \lrp{\frac{\alpha(p)}{p}},
\EEqn
and from Lemma \ref{Lemma:PowerAnalysis-3}\ref{Lemma:PowerAnalysis-3-1}-\ref{Lemma:PowerAnalysis-3-2} we have that
\BEqn
& & \bbe{k^2(Z_1,Z_2)} \le C(\tilde{M},U^{\ast},L_0,U_0,\rho) c_0^2(A_0^{XY}), \\
& & \max\lrcp{\bbe{\tilde{k}^4(X_1,X_2)}, \bbe{\tilde{k}^4(X_1,Y_1)}, \bbe{\tilde{k}^4(Y_1,Y_2)}} \le  C(\tilde{M},U^{\ast},L_0,U_0,\rho) c_0^4(A_0^{XY}) \lrp{\frac{\alpha(p)}{p}}^2,
\EEqn
which jointly imply that
\begin{equation*}
   \frac{\bbe{\tilde{k}^4(X_1,X_2)} + \bbe{\tilde{k}^4(X_1,Y_1)} + \bbe{\tilde{k}^4(Y_1,Y_2)}}{N \lrp{\cv_k^2(Z)}^2}
\le C(\tilde{M},\hat{M},U^{\ast},L_0,U_0,\rho) \frac{1}{N} 
\rightarrow 0 
\end{equation*}
and
\BEqn
& & \frac{\bbe{k^2(Z_1,Z_2)}}{N^s \cv_k^2(Z)}
\le C(\tilde{M},\hat{M},U^{\ast},L_0,U_0,\rho) \frac{1}{N^s} \lrp{\frac{p}{\alpha(p)}}
\rightarrow 0
\EEqn
as long as $\frac{p}{\alpha(p)} = p^{1-\delta_0} = o\lrp{N^s}$. 

Similarly as in the previous proof, we determine the order of $\ce^k(X,Y)$ and verify the remaining conditions case by case.
\begin{enumerate}[label=(\roman*)]
    \item if $\max\{\delta_5,\delta_6\} > \delta_7$, then it naturally holds that $(\alpha(p))^{2\delta_7}p^{2\delta_7-4} \prec (\alpha(p))^{2(\delta_5\vee\delta_6)}p^{\delta_5\vee\delta_6-3}$ and it follows that
    \BEqn
    & & C_1(\tilde{M},\hat{M},U^{\ast},L_0,U_0,L_4,U_4) |c_0(A_0^{XY})| (\alpha(p))^{2(\delta_5\vee\delta_6)} p^{\delta_5\vee\delta_6-3} \\
    &\le& \ce^{k}(X,Y) \\
    &\le& C_2(\tilde{M},U^{\ast},L_0,U_0,L_4,U_4) |c_0(A_0^{XY})| (\alpha(p))^{2(\delta_5\vee\delta_6)} p^{\delta_5\vee\delta_6-3}.
    \EEqn
    Consequently, it holds that
    \BEqn
        \frac{\lrp{\ce^k(X,Y)}^4}{N \lrp{\cv_k^2(Z)}^2}
    &\le& C(\tilde{M},\hat{M},U^{\ast},L_0,U_0,L_4,U_4,\rho) \frac{(\alpha(p))^{8(\delta_5\vee\delta_6)-2} p^{4(\delta_5\vee\delta_6)-10}}{N} \\
    &\le& C(\tilde{M},\hat{M},U^{\ast},L_0,U_0,L_4,U_4,\rho) \frac{1}{N} \lrp{\frac{1}{\alpha(p) p}}^{6-6(\delta_5\vee\delta_6)} \lrp{\frac{\alpha(p)}{p}}^{2(\delta_5\vee\delta_6)+4} \\
    &\rightarrow& 0,
    \EEqn
    and
    \BEqn
    & & \frac{N^2 \lrp{\ce^k(X,Y)}^2}{\cv_k^2(Z)} 
    \geq C(\tilde{M},\hat{M},U^{\ast},L_0,U_0,L_4,U_4,\rho) N^2 
    \lrp{\frac{1}{(\alpha(p))^{1-4(\delta_5\vee\delta_6)} p^{5-2(\delta_5\vee\delta_6)}}}
    \rightarrow \infty, \\
    & & \frac{N \lrp{\ce^k(X,Y)}^2}{\bbe{(h^k(X_1,X_2,Y_1,Y_2)^2}} 
    \geq C(\tilde{M},\hat{M},U^{\ast},L_0,U_0,L_4,U_4,\rho) N 
    \lrp{\frac{1}{(\alpha(p))^{1-4(\delta_5\vee\delta_6)} p^{5-2(\delta_5\vee\delta_6)}}}
    \rightarrow \infty
    \EEqn
    as long as $(\alpha(p))^{1-4(\delta_5\vee\delta_6)} p^{5-2(\delta_5\vee\delta_6)} \prec N$, that is, $p \prec N^{1/((1-4(\delta_5\vee\delta_6))\delta_0+5-2(\delta_5\vee\delta_6))}$.
    
    \item if $\delta_7 > \max\{\delta_5,\delta_6\}$, and $\delta_0 < (1+\delta_5\vee\delta_6-2\delta_7)/(2(\delta_7-\delta_5\vee\delta_6))$, then we still have
    $(\alpha(p))^{2\delta_7}p^{2\delta_7-4} \prec (\alpha(p))^{2(\delta_5\vee\delta_6)}p^{\delta_5\vee\delta_6-3}$, that is,
    \BEqn
    & & C_1(\tilde{M},\hat{M},U^{\ast},L_0,U_0,L_4,U_4)  |c_0(A_0^{XY})| (\alpha(p))^{2(\delta_5\vee\delta_6)} p^{\delta_5\vee\delta_6-3} \\
    &\le& \ce^{k}(X,Y) \\
    &\le& C_2(\tilde{M},U^{\ast},L_0,U_0,L_4,U_4)  |c_0(A_0^{XY})| (\alpha(p))^{2(\delta_5\vee\delta_6)} p^{\delta_5\vee\delta_6-3},
    \EEqn
    thus all the analysis in the previous case still holds.
    
    \item if $\delta_7 > \max\{\delta_5,\delta_6\}$, and $(1+\delta_5\vee\delta_6-2\delta_7)/(2(\delta_7-\delta_5\vee\delta_6)) < \delta_0 < \max\{\delta_5/(3-\delta_5), \delta_6/(3-\delta_6), (2\delta_7-1)/(3-2\delta_7)\}$, then in this case, we have
    \BEqn
    & & C_1(\tilde{M},\hat{M},U^{\ast},L_0,U_0,L_4,U_4) |c_0(A_0^{XY})| (\alpha(p))^{2\delta_7} p^{2\delta_7-4} \\
    &\le& \ce^{k}(X,Y) \\
    &\le& C_2(\tilde{M},U^{\ast},L_0,U_0,L_4,U_4) |c_0(A_0^{XY})| (\alpha(p))^{2\delta_7} p^{2\delta_7-4}
    \EEqn
    implying that
    \BEqn
        \frac{\lrp{\ce^k(X,Y)}^4}{N \lrp{\cv_k^2(Z)}^2}
    &\le& C(\tilde{M},\hat{M},U^{\ast},L_0,U_0,L_4,U_4,\rho) \frac{(\alpha(p))^{8\delta_7-2} p^{8\delta_7-14}}{N} \\
    &\le& C(\tilde{M},\hat{M},U^{\ast},L_0,U_0,L_4,U_4,\rho) \frac{1}{N} \lrp{\frac{1}{\alpha(p) p}}^{8-8\delta_5} \lrp{\frac{\alpha(p)}{p}}^{6} \\
    &\rightarrow& 0,
    \EEqn
    and
    \BEqn
    & & \frac{N^2 \lrp{\ce^k(X,Y)}^2}{\cv_k^2(Z)} 
    \geq C(\tilde{M},\hat{M},U^{\ast},L_0,U_0,L_4,U_4,\rho) N^2 
    \lrp{\frac{1}{(\alpha(p))^{1-4\delta_7} p^{7-4\delta_7}}}
    \rightarrow \infty, \\
    & & \frac{N \lrp{\ce^k(X,Y)}^2}{\bbe{(h^k(X_1,X_2,Y_1,Y_2)^2}} 
    \geq C(\tilde{M},\hat{M},U^{\ast},L_0,U_0,L_4,U_4,\rho) N 
    \lrp{\frac{1}{(\alpha(p))^{1-4\delta_7} p^{7-4\delta_7}}}
    \rightarrow \infty
    \EEqn
    as long as $(\alpha(p))^{1-4\delta_7} p^{7-4\delta_7} \prec N$. By summarizing all these analysis, we arrive at the proposed result.
\end{enumerate}

\subsection{Proof of Proposition \ref{Prop:PowerAnalysis-4}}

When $\delta_0=0$ and $s=1$, we discuss the results in Proposition \ref{Prop:GeneralPowerAnalysis-4} case by case. For the first case, when $\max\{\delta_5,\delta_6\} > \delta_7$, then it follows that $\max\{\delta_5,\delta_6\}>0$ under the assumption that $\max\{\delta_5,\delta_6,\delta_7\}>0$ and it thus holds by natural that $\delta_0<\max\{\delta_5/(3-\delta_5), \delta_6/(3-\delta_6), (2\delta_7-1)/(3-2\delta_7)\}$. Hence, $\blrp{T_{n,m,p}^{k}>C}\rightarrow1$ holds when $p =  o\lrp{N^{1/((1-4(\delta_5\vee\delta_6))\delta_0+5-2(\delta_5\vee\delta_6))},N^{s/(1-\delta_0)}} = o\lrp{N^{1/(5-2\max\{\delta_5,\delta_6\})}}$. 

For the second case when $\delta_7 > \max\{\delta_5,\delta_6\}$, the condition $\delta_0 < \min\{(1+\delta_5\vee\delta_6-2\delta_7)/(2(\delta_7-\delta_5\vee\delta_6)), \max\{\delta_5/(3-\delta_5), \delta_6/(3-\delta_6), (2\delta_7-1)/(3-2\delta_7)\}$ holds only when $\delta_7<(1+\max\{\delta_5,\delta_6\})/2$ and $\max\{\delta_5,\delta_6,2\delta_7-1\}>0$, which is satisfied when $0 < \max\{\delta_5,\delta_6\} < \delta_7 < (1+\max\{\delta_5,\delta_6\})/2$. It follows that the constraint in this case becomes $p=o\lrp{N^{1/(5-2\max\{\delta_5,\delta_6\})}}$. 

As for the last case,  $\delta_7 > \max\{\delta_5,\delta_6\}$ and $(1+\delta_5\vee\delta_6-2\delta_7)/(2(\delta_7-\delta_5\vee\delta_6)) < \delta_0 < \max\{\delta_5/(3-\delta_5), \delta_6/(3-\delta_6), (2\delta_7-1)/(3-2\delta_7)\}$ hold together when $\delta_7>(1+\max\{\delta_5,\delta_6\})/2$, and the corresponding constraint is simplified as $p=o\lrp{N^{1/(7-4\delta_7)}}$. 

Note the the constraints for the first and second case are the same, both of which are $p=o\lrp{N^{1/(5-2\max\{\delta_5,\delta_6\})}}$, and the conditions of these cases can be unified as $\delta_7 < (1+\max\{\delta_5,\delta_6\})/2$. Together with the third case, we can summarize a uniform constraint for all the three cases, that is, $p = o\lrp{N^{1/(7-2\max\{1+\max\{\delta_5,\delta_6\}, 2\delta_7\})}}$.


\section{Additional Simulation Results}\label{Appendix:simulation}

\subsection{Additional Results of Normal Approximation Accuracy}\label{Simu:accuracy}

We follow the same DGP as in Example \ref{Ex:Hist-R1} to investigate the normal approximation accuracy when the sample sizes are unequal with the difference beyond a constant. To be specific, we generate two independent random samples $\fX = \{X_1,\dots,X_n\}$ and $\fY = \{Y_1,\dots,Y_m\}$ from the following data generating process. 

\begin{example}\label{Ex:Hist-R1-supp}
Let 
\begin{equation*}
    X_1,\dots,X_n \stsim{iid} \cn(0,\Sigma), \qquad
    Y_1,\dots,Y_m \stsim{iid} \cn(0,\Sigma),
\end{equation*}
where $\Sigma = \lrp{\sigma_{ij}} \in \br^{p \times p}$ with $\sigma_{ij} = \rho^{|i-j|}$ and $\rho = 0.5$. 
  
We consider the setting that $n \in \{25,50,100,200,400\}$, the sample size ratio $m/n \in \{2,4,8\}$, and the data dimensionality $p \in \{25,50,100,200\}$. As for the kernel $k$, we consider the $L_2$-norm $k_{L_2}(x,y) = |x-y|$, the Gaussian kernel multiplied by -1, that is, $k_G(x,y) = -\exp\lrp{-|x-y|^2/(2\gamma^2)}$ with $\gamma^2 = \text{Median}\{|X_{i_1}-X_{i_2}|^2, |X_i-Y_j|^2, |Y_{j_1}-Y_{j_2}|^2\}$, and the Laplacian kernel multiplied by -1, that is, $k_L(x,y) = -\exp\lrp{-|x-y|/\gamma}$ with $\gamma = \text{Median}\{|X_{i_1}-X_{i_2}|, |X_i-Y_j|, |Y_{j_1}-Y_{j_2}|\}$.
\end{example}

\begin{figure}[h!]
  \begin{center}
    \includegraphics[height=105mm,width=140mm]{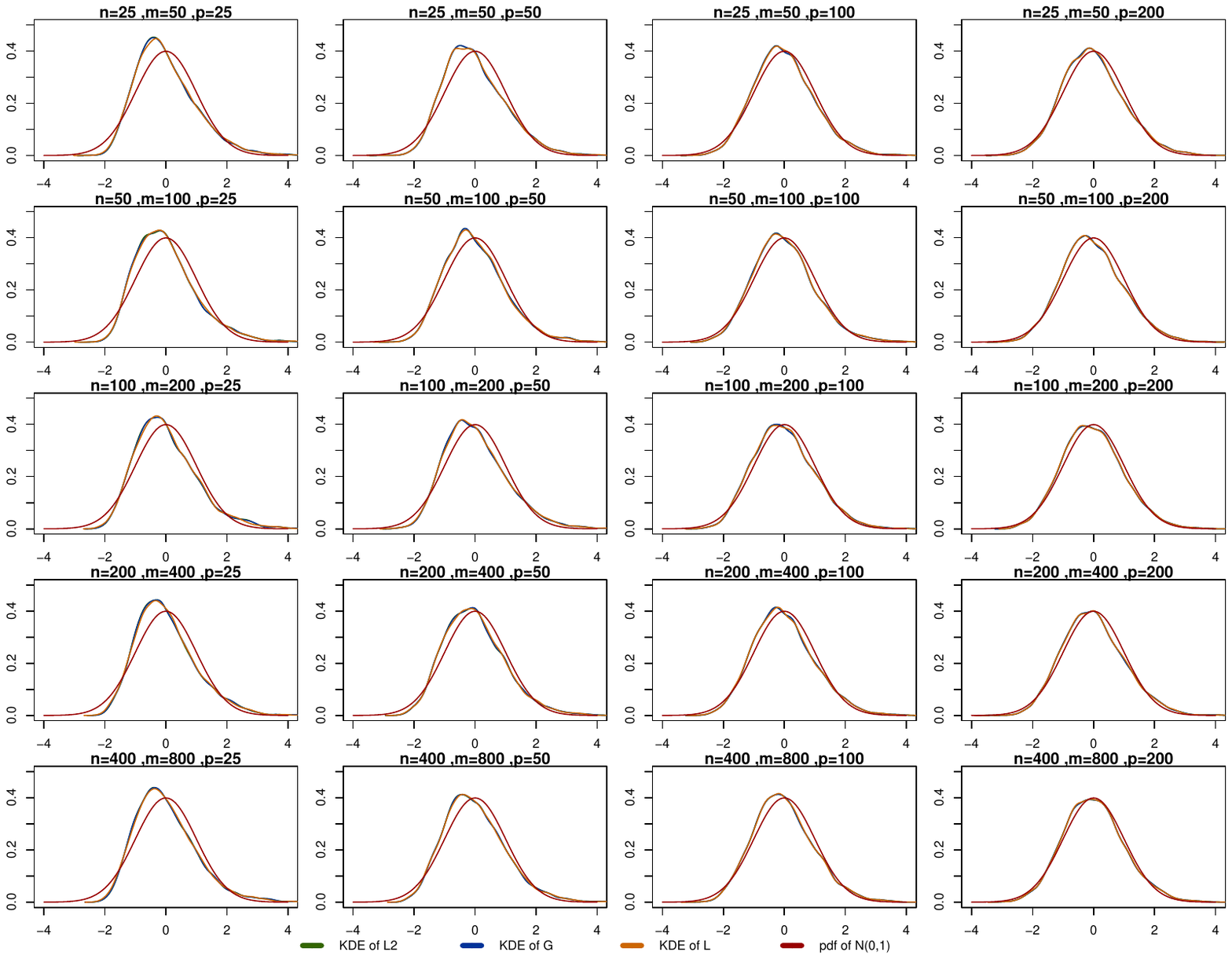}
    \caption{Kernel density estimates of the studentized test statistic $T_{n,m,p}^k$ with different kernels for Example \ref{Ex:Hist-R1-supp} when $m/n=2$. The four columns correspond to different $p$'s and the five rows correspond to different pairs of $(n,m)$.}
    \label{Fig:Hist-R1-r2}
  \end{center}
\end{figure}

\begin{figure}[h!]
  \begin{center}
    \includegraphics[height=105mm,width=140mm]{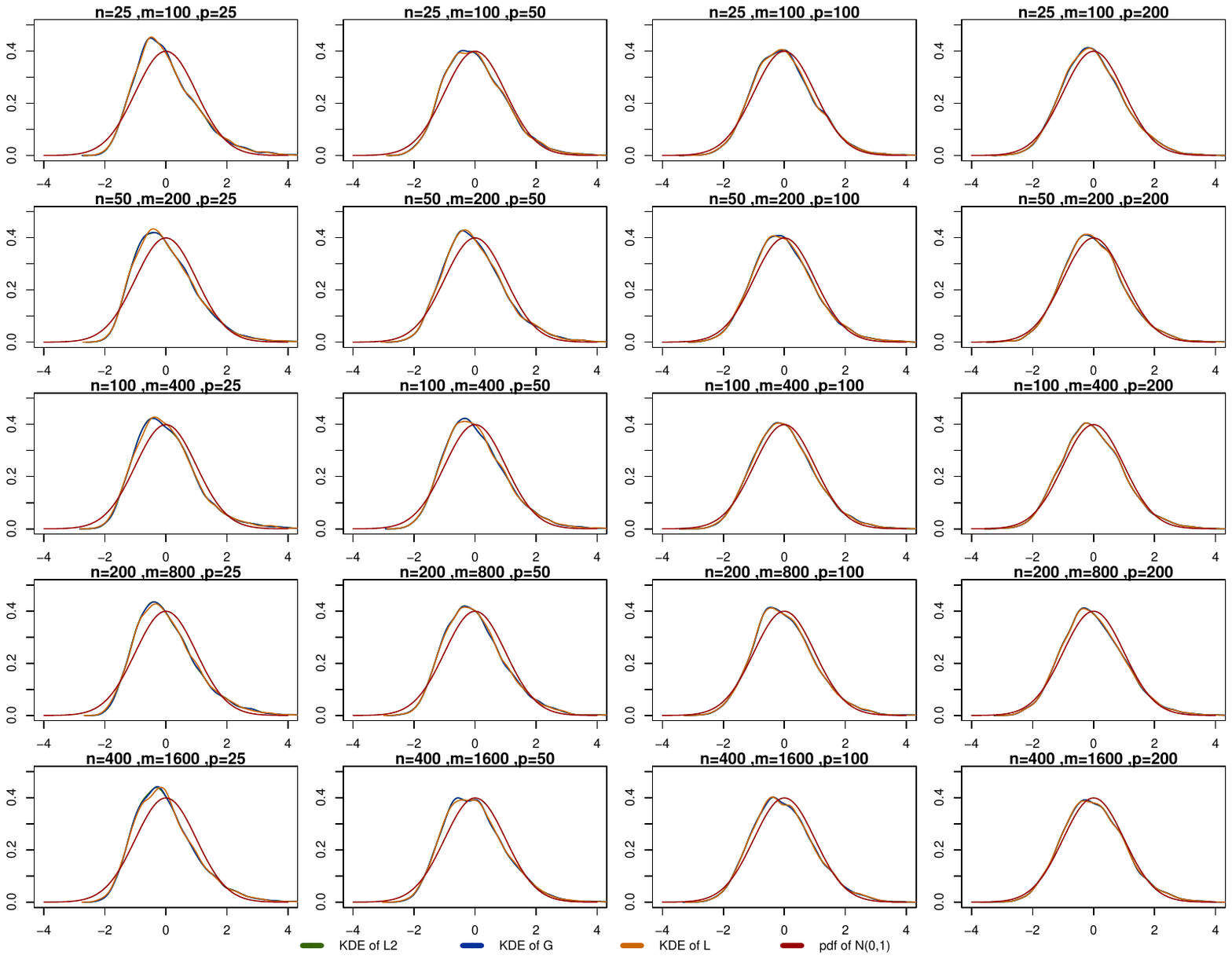}
    \caption{Kernel density estimates of the studentized test statistic $T_{n,m,p}^k$ with different kernels for Example \ref{Ex:Hist-R1-supp} when $m/n=4$. The four columns correspond to different $p$'s and the five rows correspond to different pairs of $(n,m)$.}
    \label{Fig:Hist-R1-r3}
  \end{center}
\end{figure}

\begin{figure}[h!]
  \begin{center}
    \includegraphics[height=105mm,width=140mm]{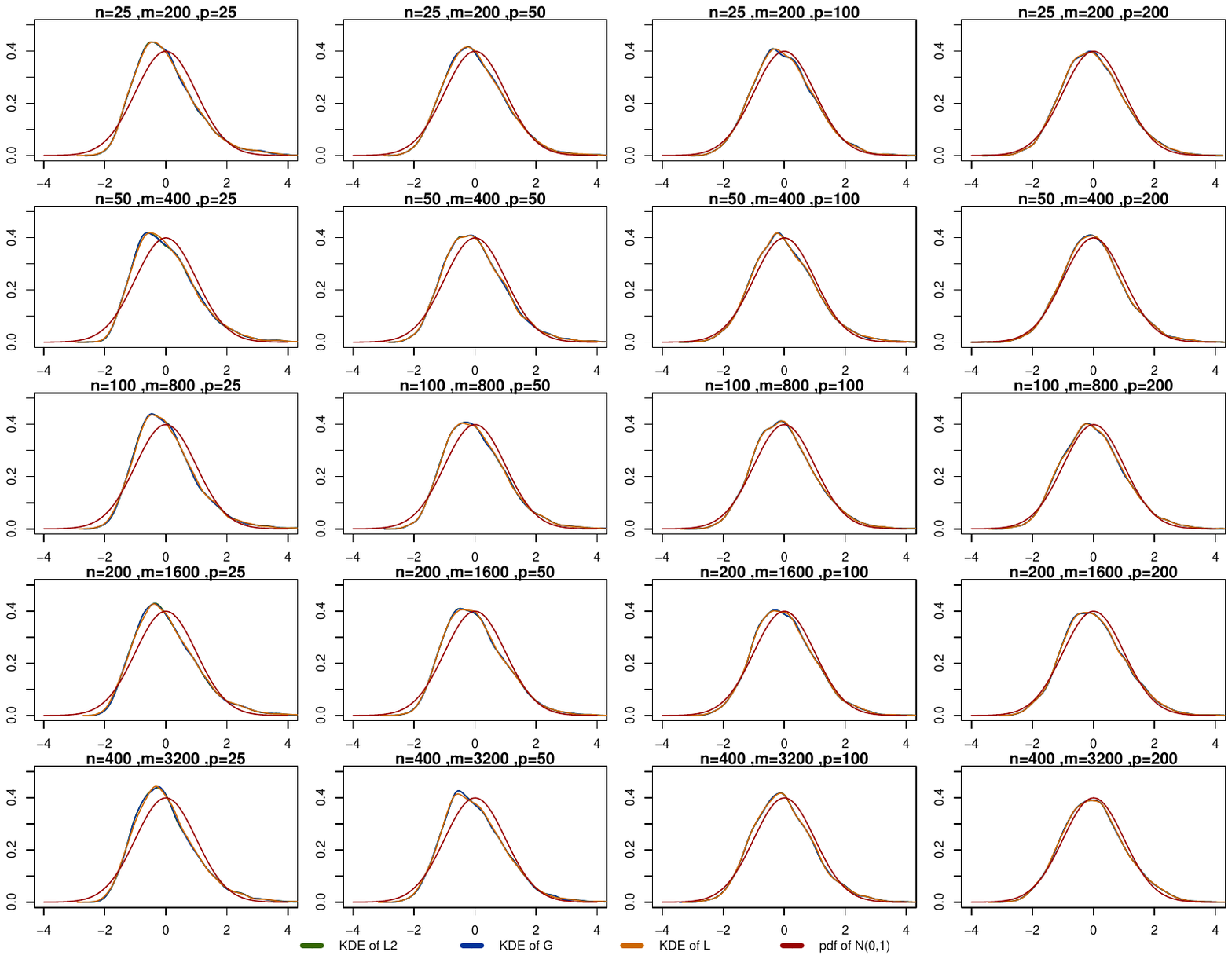}
    \caption{Kernel density estimates of the studentized test statistic $T_{n,m,p}^k$ with different kernels for Example \ref{Ex:Hist-R1-supp} when $m/n=8$. The four columns correspond to different $p$'s and the five rows correspond to different pairs of $(n,m)$.}
    \label{Fig:Hist-R1-r4}
  \end{center}
\end{figure}

The results are based on $5000$ Monte Carlo simulations, and again we plot the kernel density estimates (KDE) for the three kernels and the standard normal density function for each combination of sample size  and dimensionality, see Figure \ref{Fig:Hist-R1-r2}-Figure \ref{Fig:Hist-R1-r4}.

The overall trends shown in Figure \ref{Fig:Hist-R1-r2}-Figure \ref{Fig:Hist-R1-r4} match what we have observed from Figure \ref{Fig:Hist-R1-r1}, that is, the accuracy does not necessarily improve significantly when $N$ increases but $p$ is fixed, whereas we observe higher normal approximation accuracy when $p$ grows but $N$ remains unchanged. In general, the accuracy improves as both $N$ and $p$ increases, though relatively high accuracy is already achieved when $N$ and $p$ are not too larger, say, $N=50$ and $p=100$. By comparing the trends in these figures, we further notice that the sample size ratio has a minimal impact on the normal approximation accuracy, and high accuracy is achievable for an heavily unbalanced data set as long as both $N$ and $p$ are sufficiently large, which matches the theory established in the paper. 

To precisely quantify the discrepancy between the standard normal distribution and the empirical distribution of the proposed test statistic under the null, we further report the Kolmogorov-Smirnov distance and the Wasserstein distance for Example \ref{Ex:Hist-R1} and Example \ref{Ex:Hist-R1-supp}, see Table \ref{Tab:KSstat} and Table \ref{Tab:WSdist}.

\begin{table}[h!]
\centering
\scalebox{0.75}{
\begin{tabular}{c|c|cccc|cccc|cccc}
  \hline\hline
  & & \multicolumn{4}{c|}{$L_2$} & \multicolumn{4}{c|}{$G$} & \multicolumn{4}{c}{$L$} \\ \hline
  $m/n$ & \diagbox{$n$}{$p$} & 25 & 50 & 100 & 200 & 25 & 50 & 100 & 200 & 25 & 50 & 100 & 200  \\ \hline
  \multirow{5}{*}{1} & 25 & 0.055 & 0.053 & 0.035 & 0.026 & 0.055 & 0.053 & 0.035 & 0.026 & 0.051 & 0.053 & 0.034 & 0.028 \\
  & 50 & 0.048 & 0.045 & 0.028 & 0.018 & 0.048 & 0.044 & 0.028 & 0.017 & 0.043 & 0.045 & 0.031 & 0.019 \\
  & 100 & 0.060 & 0.044 & 0.032 & 0.030 & 0.058 & 0.043 & 0.032 & 0.030 & 0.051 & 0.042 & 0.033 & 0.029 \\
  & 200 & 0.067 & 0.044 & 0.044 & 0.028 & 0.066 & 0.044 & 0.043 & 0.028 & 0.059 & 0.041 & 0.043 & 0.029  \\ 
  & 400 & 0.053 & 0.047 & 0.038 & 0.020 & 0.052 & 0.047 & 0.038 & 0.020 & 0.050 & 0.045 & 0.039 & 0.020 \\ \hline
  \multirow{5}{*}{2} & 25 & 0.063 & 0.061 & 0.032 & 0.031 & 0.063 & 0.060 & 0.032 & 0.031 & 0.060 & 0.056 & 0.030 & 0.030 \\
  & 50 & 0.059 & 0.045 & 0.036 & 0.042 & 0.058 & 0.045 & 0.036 & 0.042 & 0.056 & 0.044 & 0.034 & 0.039 \\
  & 100 & 0.055 & 0.040 & 0.032 & 0.033 & 0.055 & 0.041 & 0.032 & 0.032 & 0.049 & 0.038 & 0.031 & 0.033 \\
  & 200 & 0.057 & 0.049 & 0.040 & 0.030 & 0.056 & 0.049 & 0.040 & 0.030 & 0.051 & 0.044 & 0.039 & 0.032 \\
  & 400 & 0.058 & 0.049 & 0.043 & 0.017 & 0.058 & 0.050 & 0.042 & 0.016 & 0.054 & 0.048 & 0.044 & 0.016 \\ \hline
  \multirow{5}{*}{4} & 25 & 0.064 & 0.035 & 0.030 & 0.035 & 0.064 & 0.035 & 0.030 & 0.036 & 0.062 & 0.034 & 0.030 & 0.033  \\
  & 50 & 0.058 & 0.052 & 0.034 & 0.027 & 0.056 & 0.052 & 0.035 & 0.027 & 0.056 & 0.053 & 0.032 & 0.027 \\
  & 100 & 0.045 & 0.044 & 0.029 & 0.032 & 0.045 & 0.046 & 0.029 & 0.032 & 0.044 & 0.038 & 0.027 & 0.033 \\
  & 200 & 0.054 & 0.047 & 0.036 & 0.032 & 0.051 & 0.047 & 0.036 & 0.031 & 0.048 & 0.044 & 0.037 & 0.031 \\
  & 400 & 0.063 & 0.046 & 0.041 & 0.020 & 0.064 & 0.045 & 0.041 & 0.020 & 0.064 & 0.043 & 0.043 & 0.021 \\ \hline 
  \multirow{5}{*}{8} & 25 & 0.061 & 0.041 & 0.032 & 0.027 & 0.060 & 0.042 & 0.032 & 0.027 & 0.060 & 0.038 & 0.030 & 0.026 \\
  & 50 & 0.056 & 0.042 & 0.032 & 0.026 & 0.056 & 0.044 & 0.032 & 0.025 & 0.052 & 0.042 & 0.030 & 0.023 \\
  & 100 & 0.054 & 0.036 & 0.032 & 0.028 & 0.051 & 0.036 & 0.033 & 0.028 & 0.052 & 0.033 & 0.033 & 0.027 \\
  & 200 & 0.054 & 0.048 & 0.030 & 0.023 & 0.053 & 0.047 & 0.030 & 0.023 & 0.047 & 0.048 & 0.028 & 0.022 \\
  & 400 & 0.063 & 0.055 & 0.041 & 0.018 & 0.062 & 0.054 & 0.040 & 0.017 & 0.057 & 0.051 & 0.042 & 0.016 \\ \hline\hline
\end{tabular}}
\caption{Kolmogorov–Smirnov distance between the empirical distribution of $T_{n,m,p}$ and the standard normal distribution under multiple $(n,m,p)$}
\label{Tab:KSstat}
\end{table}

\begin{table}[h!]
\centering
\scalebox{0.75}{
\begin{tabular}{c|c|cccc|cccc|cccc}
  \hline\hline
  & & \multicolumn{4}{c|}{$L_2$} & \multicolumn{4}{c|}{$G$} & \multicolumn{4}{c}{$L$} \\ \hline
  $m/n$ & \diagbox{$n$}{$p$} & 25 & 50 & 100 & 200 & 25 & 50 & 100 & 200 & 25 & 50 & 100 & 200  \\ \hline
  \multirow{5}{*}{1} & 25 & 0.140 & 0.115 & 0.090 & 0.057 & 0.147 & 0.110 & 0.079 & 0.066 & 0.132 & 0.093 & 0.076 & 0.063 \\
  & 50 & 0.131 & 0.108 & 0.069 & 0.043 & 0.125 & 0.110 & 0.066 & 0.046 & 0.119 & 0.103 & 0.055 & 0.066 \\
  & 100 & 0.124 & 0.119 & 0.071 & 0.073 & 0.138 & 0.101 & 0.075 & 0.069 & 0.135 & 0.096 & 0.068 & 0.063 \\
  & 200 & 0.134 & 0.112 & 0.079 & 0.061 & 0.139 & 0.106 & 0.075 & 0.073 & 0.126 & 0.103 & 0.068 & 0.063 \\
  & 400 & 0.131 & 0.097 & 0.076 & 0.034 & 0.125 & 0.104 & 0.088 & 0.042 & 0.109 & 0.086 & 0.101 & 0.041 \\ \hline
  \multirow{5}{*}{2} & 25 & 0.150 & 0.131 & 0.064 & 0.084 & 0.135 & 0.127 & 0.071 & 0.072 & 0.134 & 0.110 & 0.062 & 0.066 \\
  & 50 & 0.161 & 0.110 & 0.083 & 0.046 & 0.136 & 0.097 & 0.082 & 0.061 & 0.140 & 0.092 & 0.077 & 0.079 \\
  & 100 & 0.127 & 0.116 & 0.081 & 0.056 & 0.126 & 0.118 & 0.066 & 0.050 & 0.130 & 0.108 & 0.074 & 0.051 \\
  & 200 & 0.119 & 0.130 & 0.076 & 0.073 & 0.128 & 0.115 & 0.072 & 0.057 & 0.123 & 0.089 & 0.078 & 0.083 \\
  & 400 & 0.132 & 0.097 & 0.081 & 0.060 & 0.125 & 0.116 & 0.074 & 0.057 & 0.128 & 0.105 & 0.099 & 0.057 \\ \hline
  \multirow{5}{*}{4} & 25 & 0.146 & 0.093 & 0.065 & 0.071 & 0.141 & 0.092 & 0.084 & 0.058 & 0.139 & 0.082 & 0.058 & 0.062\\
  & 50 & 0.140 & 0.120 & 0.064 & 0.059 & 0.134 & 0.121 & 0.076 & 0.062 & 0.134 & 0.126 & 0.070 & 0.050 \\
  & 100 & 0.118 & 0.114 & 0.066 & 0.046 & 0.127 & 0.099 & 0.062 & 0.072 & 0.119 & 0.106 & 0.057 & 0.086 \\
  & 200 & 0.149 & 0.109 & 0.081 & 0.065 & 0.138 & 0.102 & 0.069 & 0.072 & 0.126 & 0.106 & 0.085 & 0.078 \\
  & 400 & 0.141 & 0.106 & 0.081 & 0.051 & 0.126 & 0.102 & 0.088 & 0.044 & 0.134 & 0.100 & 0.081 & 0.052 \\ \hline
  \multirow{5}{*}{8} & 25 & 0.154 & 0.095 & 0.052 & 0.055 & 0.155 & 0.108 & 0.069 & 0.053 & 0.143 & 0.084 & 0.057 & 0.053 \\
  & 50 & 0.129 & 0.091 & 0.086 & 0.048 & 0.146 & 0.085 & 0.066 & 0.061 & 0.132 & 0.091 & 0.074 & 0.050 \\
  & 100 & 0.145 & 0.093 & 0.077 & 0.076 & 0.128 & 0.100 & 0.076 & 0.073 & 0.122 & 0.087 & 0.067 & 0.066 \\
  & 200 & 0.146 & 0.117 & 0.076 & 0.063 & 0.131 & 0.120 & 0.064 & 0.064 & 0.129 & 0.110 & 0.085 & 0.051 \\
  & 400 & 0.152 & 0.094 & 0.075 & 0.061 & 0.155 & 0.123 & 0.073 & 0.061 & 0.132 & 0.106 & 0.086 & 0.046 \\ \hline\hline
\end{tabular}}
\caption{Wasserstein Distance between the empirical distribution of $T_{n,m,p}$ and the standard normal distribution under multiple $(n,m,p)$}
\label{Tab:WSdist}
\end{table}

The numerical results presented in Table \ref{Tab:KSstat} and Table \ref{Tab:WSdist} generally match our observations from the figures. To be specific, the discrepancy between the standard normal approximation and the empirical distribution of our proposed test does decrease significantly as both $N$ and $p$ grow. When $p$ is fixed, the normal approximation is not guaranteed to be more accurate when $N$ increases, whereas the improvement of accuracy is observed when $N$ is fixed but $p$ grows. Additionally, we do not observe significant difference in the approximation accuracy among the three kernels considered in Example \ref{Ex:Hist-R1} and Example \ref{Ex:Hist-R1-supp}, which indicates that our proposed method is expected to perform well for a broad range of kernels as far as size is concerned.

\subsection{Additional Results of Empirical Size}\label{Simu:size}

Here we consider a different example from Example \ref{Ex:Size-1} to investigate a bit more into the empirical size. We generate data from two independent multivariate normal distributions as follows.

\begin{example}\label{Ex:Size-2}
  Let
  \begin{equation*}
      X_1,\dots,X_n \stsim{iid} \cn(0,\Sigma), \qquad
      Y_1,\dots,Y_m \stsim{iid} \cn(0,\Sigma),
  \end{equation*}
  where $\Sigma = \lrp{\sigma_{ij}} \in \br^{p \times p}$ with $\sigma = \rho^{|i-j|}$ and $\rho \in \{0.4, 0.7\}$. We consider the setting that $(n,m) \in \{(25,25), (50,50), (50,100), (100,100), (200,200)\}$ and $p \in \{50, 100\}$.
\end{example}

Same as in Example \ref{Ex:Size-1}, we consider the $L_2$ norm, the Gaussian kernel and the Laplacian kernel for our method, and compare our test with the permutation test studied in \cite{zhu2021interpoint} which additionally adopts the $L_1$-norm $k_{L_1}(x,y) = |x-y|_1$. We report the empirical sizes at the significance level $\alpha = 0.05$ in Table \ref{Tab:Size-2-Cr}. 


\begin{table}[h!]
    \centering
    \scalebox{0.75}{
    \begin{tabular}{c|c|c|c||ccc|cccc}
    \hline\hline
        \multirow{2}{*}{$n$} & \multirow{2}{*}{$m$} & \multirow{2}{*}{$p$} & \multirow{2}{*}{$\rho$} & \multicolumn{3}{c}{Proposed} & \multicolumn{4}{c}{Permutation} \\ \cline{5-11} 
        & & & & $L_2$ & G & L & $L_2$ & G & L & $L_1$ \\ \hline
        \multirow{4}{*}{25} & \multirow{4}{*}{25} & \multirow{2}{*}{50} & 0.4 & 6.00 & 6.12 & 5.96 & 5.00 & 4.90 & 4.90 & 5.00 \\
        & & & 0.7 & 6.82 & 6.82 & 6.70 & 5.10 & 5.10 & 5.20 & 6.00 \\ \cline{3-11}
        & & \multirow{2}{*}{100} & 0.4 & 5.78 & 5.78 & 5.74 &  4.10 & 4.30 & 4.70 & 5.00 \\
        & & & 0.7 & 6.38 & 6.42 & 6.38 & 5.20 & 5.20 & 5.10 & 4.90 \\ \hline
        \multirow{4}{*}{50} & \multirow{4}{*}{50} & \multirow{2}{*}{50} & 0.4 & 6.54 & 6.50 & 6.44 & 6.40 & 6.20 & 5.80 & 6.40 \\
        & & & 0.7 & 7.18 & 7.22 & 7.06 & 5.80 & 5.70 & 5.50 & 5.40 \\ \cline{3-11}
        & & \multirow{2}{*}{100} & 0.4 & 5.98 & 5.98 & 5.84 & 5.50 & 5.40 & 5.60 & 5.90 \\
        & & & 0.7 & 6.02 & 5.92 & 6.06 & 5.20 & 5.20 & 5.30 & 5.40 \\ \hline
        \multirow{4}{*}{50} & \multirow{4}{*}{100} & \multirow{2}{*}{50} & 0.4 & 6.30 & 6.30 & 6.28 & 5.40 & 5.40 & 5.50 & 5.50 \\
        & & & 0.7 & 7.28 & 7.28 & 7.06 & 3.90 & 3.70 & 4.10 & 3.60 \\ \cline{3-11}
        & & \multirow{2}{*}{100} & 0.4 & 6.02 & 6.04 & 5.88 & 6.00 & 6.00 & 6.00 & 5.10 \\
        & & & 0.7 & 6.50 & 6.54 & 6.52 & 5.30 & 5.40 & 5.40 & 6.20 \\ \hline
        \multirow{4}{*}{100} & \multirow{4}{*}{100} & \multirow{2}{*}{50} & 0.4 & 6.42 & 6.44 & 6.40 & 5.20 & 5.20 & 5.00 & 5.10 \\
        & & & 0.7 & 6.24 & 6.16 & 6.10 & 2.70 & 2.80 & 3.30 & 3.30 \\ \cline{3-11}
        & & \multirow{2}{*}{100} & 0.4 & 6.08 & 6.08 & 6.08 & 5.20 & 5.20 & 5.10 & 5.20 \\
        & & & 0.7 & 5.98 & 6.04 & 6.02 & 4.50 & 4.70 & 4.90 & 4.80 \\ \hline
        \multirow{4}{*}{200} & \multirow{4}{*}{200} & \multirow{2}{*}{50} & 0.4 & 5.62 & 5.60 & 5.58 & 4.30 & 4.40 & 4.20 & 4.10 \\
        & & & 0.7 & 6.88 & 6.90 & 6.96 & 4.20 & 4.20 & 4.30 & 4.80 \\ \cline{3-11}
        & & \multirow{2}{*}{100} & 0.4 & 6.60 & 6.58 & 6.50 & 5.00 & 5.00 & 5.00 & 5.90 \\
        & & & 0.7 & 6.34 & 6.26 & 6.20 & 4.70 & 4.70 & 4.50 & 5.10 \\ \hline\hline
    \end{tabular}}
    \caption{Size comparison for Example \ref{Ex:Size-2}. All the empirical sizes are reported in percentage.}\label{Tab:Size-2-Cr}
\end{table}

Table \ref{Tab:Size-2-Cr} exhibits similar phenomenon as Table \ref{Tab:Size-1-Cr}, and we observe some mild size distortion in both Example \ref{Ex:Size-1} and Example \ref{Ex:Size-2}, which tends to increase when the componentwise dependence gets stronger. However, the overall distortion level is acceptable even wehn $N$ and $p$ are relatively small.

\subsection{Additional Results of Power Behavior}\label{Simu:power}

We investigate the power behavior under the alternative of mean shift in Example \ref{Ex:Power-1}, and now we consider the scenario where two distributions only differ in their covariance matrices. 

\begin{example}\label{Ex:Power-2}
  Generate samples as
  \begin{equation*}
      X_1,\dots,X_n \stsim{iid} (V^{1/2} \Sigma V^{1/2})^{1/2} Z_X, \qquad
      Y_1,\dots,Y_m \stsim{iid} (V^{\ast 1/2} \Sigma V^{\ast 1/2})^{1/2} Z_Y,    
  \end{equation*}
  where $V$ is an identity matrix and $V^{\ast}$ is a diagonal matrix with $V_{ii}^{\ast 1/2} = 1.2$ for $i=1,2,\dots,\beta p$ and { $V_{ii}^{\ast 1/2} = 1$} for $i = \beta p+1,\dots,p$. All the other settings are the same as in Example \ref{Ex:Power-1}.
\end{example}

\begin{figure}[h!]
  \begin{center}
    \includegraphics[width=120mm, height=90mm]{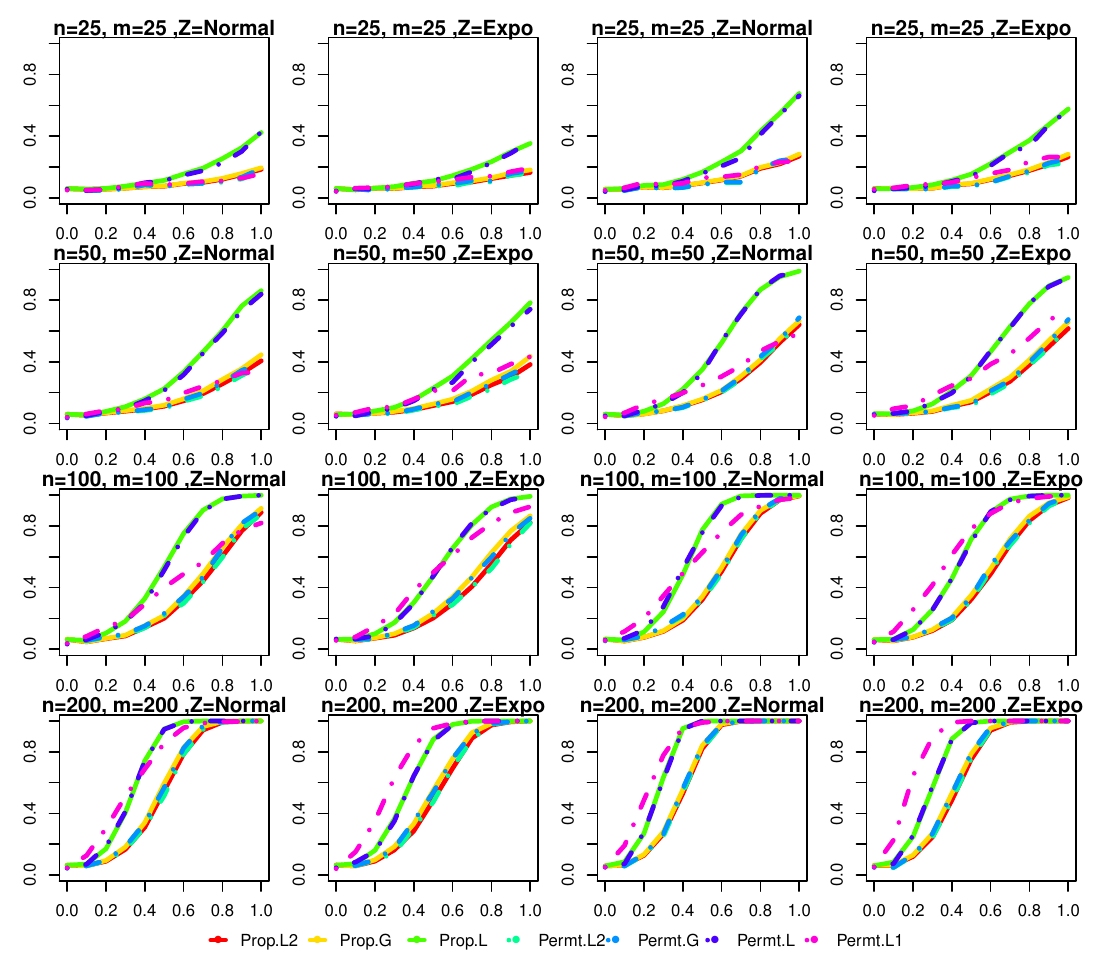}
    \caption{Size-adjusted Power Curves for Example \ref{Ex:Power-2}. The first two columns correspond to $p=50$ while the last two columns correspond to $p=100$.}\label{Fig:Power-2}
  \end{center}
\end{figure}

All the other settings are exactly the same as in Example \ref{Ex:Power-1}, including the comparison methods, the kernels used and the number of replications. Again, we plot the size-adjusted power curves against $\beta$, see Figure \ref{Fig:Power-2}. 

The overall trends of Figure \ref{Fig:Power-2} generally match those of Figure \ref{Fig:Power-1}.  
A major difference between these two scenarios is that, the use of the Laplacian kernel (for both our test and permutation test) appears to 
be advantageous as compared to the Gaussian kernel and $L_2$ norm. The power corresponding to Laplacian kernel can be higher than the one corresponding to $L_1$-norm in some cases, and this power difference is 
very interesting and it would be desirable to provide a theoretical explanation. 


\section{Some Useful Auxiliary Results}\label{Appendix:auxiliary}

\subsection{Part 1: Orders of Some Important Quantities}

Throughout this section, we use $\tilde{j}_1 \le  \dots \le \tilde{j}_s$ to denote the increasing re-ordering of any given $(j_1,\dots,j_s) \in \{1,\dots,p\}^s$.

\begin{lemma}\label{Lemma:component-1}
Under Assumption \ref{Assumpt:component-dept}\ref{Assumpt:component-dept-1}, $\cum\lrp{\tilde{x}_{j_1}, \dots, \tilde{x}_{j_s}} \neq 0$ only when $s\geq2$ and $\max\{\tilde{j}_2-\tilde{j}_1, \dots, \tilde{j}_s-\tilde{j}_{s-1}\} \le \alpha(p)$.
\end{lemma}
\begin{proof}
  When $s=1$ it is trivial that $\cum(\tilde{x}_{j_1}) = \be[\tilde{x}_{j_1}] = 0$. For $s\geq2$, suppose that there exits $2 \le u_0 \le s$, such that $\tilde{j}_{u_0} - \tilde{j}_{u_0-1} > \alpha(p)$, then under Assumption \ref{Assumpt:component-dept}, $\{\tilde{x}_{\tilde{j}_1}, \dots, \tilde{x}_{\tilde{j}_{u_0-1}}\}$ is independent of $\{\tilde{x}_{\tilde{j}_{u_0}}, \dots, \tilde{x}_{\tilde{j}_s}\}$, and it follows from the properties of cumulants that $\cum\lrp{\tilde{x}_{j_1}, \dots, \tilde{x}_{j_s}} = 0$, which completes the proof by contradiction.
\end{proof}

\begin{lemma}\label{Lemma:component-2}
Under Assumption \ref{Assumpt:component-dept}\ref{Assumpt:component-dept-1}, $\cum\lrp{\tilde{x}_{j_1}, \dots, \tilde{x}_{j_s}} \neq 0$ only when $s\geq2$ and $\max\limits_{1 \le u \le s}\{j_u\} - \min\limits_{1 \le u \le s} \{j_u\} \le (s-1)\alpha(p)$.
\end{lemma}
\begin{proof}
  Note that
  \begin{equation*}
      \max\limits_{1 \le u \le s}\{j_u\} - \min\limits_{1 \le u \le s} \{j_u\}
    = \tilde{j}_s - \tilde{j}_1
    = \sum\limits_{u=2}^{s} \lrp{ \tilde{j}_u - \tilde{j}_{u-1} },
  \end{equation*}
  then the statement directly follows from Lemma \ref{Lemma:component-1}.
\end{proof}

\begin{lemma}\label{Lemma:component-3}
Under Assumption \ref{Assumpt:component-dept}\ref{Assumpt:component-dept-1}, it holds for any fixed $s\geq2$ that
\begin{equation*}
    \# \lrcp{(j_1,\dots,j_s):\ \cum\lrp{\tilde{x}_{j_1},\dots,\tilde{x}_{j_s}}\neq0,\ 1 \le j_1,\dots,j_s \le p} \le C(s)(\alpha(p))^{s-1}p.
\end{equation*}

\end{lemma}
\begin{proof}
  It follows from Lemma \ref{Lemma:component-1} that
  \BEqn
  & & \# \lrcp{(j_1,\dots,j_s):\ \cum\bigp{\tilde{x}_{j_1},\dots,\tilde{x}_{j_s}}\neq0,\ 1 \le j_1,\dots,j_s \le p} \\
  &\le& (s!) \# \lrcp{(\tilde{j}_1,\dots,\tilde{j}_s):\ \cum\bigp{\tilde{x}_{\tilde{j}_1},\dots,\tilde{x}_{\tilde{j}_s}}\neq0,\ 1 \le \tilde{j}_1 \le \dots \le \tilde{j}_s \le p} \\
  &\le& (s!) \# \lrcp{(\tilde{j}_1,\dots,\tilde{j}_s):\ \max\limits_{2 \le u \le s}\{\tilde{j}_u - \tilde{j}_{u-1}\} \le \alpha(p),\ 1 \le \tilde{j}_1 \le \dots \le \tilde{j}_s \le p}
  \EEqn
  Let $\tilde{j}_1$ range free over $\{1,\dots,p\}$, and the number of values that $\tilde{j}_1$ can take is at most $p$. Once the value of $\tilde{j}_1$ is fixed, the number of values that $\tilde{j}_2$ can take is at most $\alpha(p)+1$ due to the constraint $\tilde{j}_2-\tilde{j}_1\le\alpha(p)$. Similarly, after the value of $\tilde{j}_2$ is fixed, the number of values that $\tilde{j}_3$ can take is again at most $\alpha(p)+1$, and so on. It follows that
  \BEqn
  & & \# \lrcp{(j_1,\dots,j_s):\ \cum\bigp{\tilde{x}_{j_1},\dots,\tilde{x}_{j_s}}\neq0,\ 1 \le j_1,\dots,j_s \le p} \\
  &\le& (s!) (\alpha(p)+1)^{s-1} p \\
  &\le& 2^{s-1}(s!) (\alpha(p))^{s-1} p
  \EEqn
  and thus completes the derivation.
\end{proof}  

\begin{lemma}\label{Lemma:component-4}
Suppose that $\tilde{j}_u \le \tilde{j}_v$ are fixed, then for any fixed $s\geq2$, it holds under Assumption \ref{Assumpt:component-dept}\ref{Assumpt:component-dept-1} that
\BEqn
& & \# \lrcp{(\tilde{j}_1,\dots,\tilde{j}_{u-1},\tilde{j}_{u+1},\dots,\tilde{j}_{v-1},\tilde{j}_{v+1},\dots,\tilde{j}_s):\ 
    \begin{array}{l}
      \cum\lrp{\tilde{x}_{j_1},\dots,\tilde{x}_{j_s}}\neq0 \vspace{1mm}\\ 
      1 \le \tilde{j}_1 \le \dots \tilde{j}_{s} \le p 
    \end{array}} \\
&\le& C(s)(\alpha(p))^{s-2} \bone\lrcp{\tilde{j}_v - \tilde{j}_u \le (s-1)\alpha(p)}.
\EEqn
\end{lemma}
\begin{proof}
  Let $\ca$ denote the set of interest, then it follows from Lemma \ref{Lemma:component-2} that $\# \ca=0$ if $\tilde{j}_v - \tilde{j}_u > (s-1)\alpha(p)$. Now assume that $\tilde{j}_v - \tilde{j}_u \le (s-1)\alpha(p)$, and we have
  \BEqn
  & & \# \ca \\
  &\le& \# \lrcp{(\tilde{j}_1,\dots,\tilde{j}_{u-1},\tilde{j}_{u+1},\dots,\tilde{j}_{v-1},\tilde{j}_{v+1},\dots,\tilde{j}_s):\ 
    \begin{array}{l}
      \max\limits_{2 \le r \le s} \{\tilde{j}_r - \tilde{j}_{r-1}\} \le \alpha(p) \vspace{1mm}\\ 
      1 \le \tilde{j}_1 \le \dots \tilde{j}_{s} \le p 
    \end{array}} \\
    &=& \# \lrcp{(\tilde{j}_1,\dots,\tilde{j}_{u-1}):\ 
      \max\limits_{2 \le r \le u} \{\tilde{j}_r - \tilde{j}_{r-1}\} \le \alpha(p),\ 
      1 \le \tilde{j}_1 \le \dots \tilde{j}_{u-1} \le \tilde{j}_u} \\
    & & \times \# \lrcp{(\tilde{j}_{u+1},\dots,\tilde{j}_{v-1}):\ 
      \max\limits_{u+1 \le r \le v} \{\tilde{j}_r - \tilde{j}_{r-1}\} \le \alpha(p),\ 
      \tilde{j}_u \le \tilde{j}_{u+1} \le \dots \tilde{j}_{v-1} \le \tilde{j}_{v}} \\
    & & \times \# \lrcp{(\tilde{j}_{v+1},\dots,\tilde{j}_s):\ 
      \max\limits_{v+1 \le r \le s} \{\tilde{j}_r - \tilde{j}_{r-1}\} \le \alpha(p),\
      \tilde{j}_{v} \le \tilde{j}_{v+1} \le \dots \tilde{j}_{s} \le p} \\
    &=:& (\# \ca_1) \times (\# \ca_2) \times (\# \ca_3). 
    \EEqn
    First consider $\# \ca_1$. When $\tilde{j}_u$ is already fixed, then the number of values that $\tilde{j}_{u-1}$ can take is at most $\alpha(p)+1$ under the constraint $\tilde{j}_{u}-\tilde{j}_{u-1} \le \alpha(p)$. After $\tilde{j}_{u-1}$ is set, the number of values that $\tilde{j}_{u-2}$ can take is at most $\alpha(p)+1$ and so on. Thus
    \begin{equation*}
      \# \ca_1 \le (\alpha(p)+1)^{u-1} \le 2^{u-1} (\alpha(p))^{u-1},
    \end{equation*}
    and similarly, we also have
    \begin{equation*}
      \# \ca_3 \le 2^{s-v} (\alpha(p))^{s-v}.
    \end{equation*}
    
    It remains to think of $\# \ca_2$. Note that
    \BEqn
    & & \ca_2 \\
    &=& \lrcp{(\tilde{j}_{u+1},\dots,\tilde{j}_{v-1}):\ 
      \max\limits_{u+1 \le r \le v} \{\tilde{j}_r - \tilde{j}_{r-1}\} \le \alpha(p),\ 
      \tilde{j}_u \le \tilde{j}_{u+1} \le \dots \tilde{j}_{v-1} \le p} \\
    & & \bigcap \lrcp{(\tilde{j}_{u+1},\dots,\tilde{j}_{v-1}):\ 
      \max\limits_{u+1 \le r \le v} \{\tilde{j}_r - \tilde{j}_{r-1}\} \le \alpha(p),\ 
      1 \le \tilde{j}_{u+1} \le \dots \tilde{j}_{v-1} \le \tilde{j}_{v}},
    \EEqn
    it follows that
    \BEqn
    & & \# \ca_2 \\
    &=& \# \lrcp{(\tilde{j}_{u+1},\dots,\tilde{j}_{v-1}):\ 
      \max\limits_{u+1 \le r \le v} \{\tilde{j}_r - \tilde{j}_{r-1}\} \le \alpha(p),\ 
      \tilde{j}_u \le \tilde{j}_{u+1} \le \dots \tilde{j}_{v-1} \le p}  \\
    & & \wedge \# \lrcp{(\tilde{j}_{u+1},\dots,\tilde{j}_{v-1}):\ 
      \max\limits_{u+1 \le r \le v} \{\tilde{j}_r - \tilde{j}_{r-1}\} \le \alpha(p),\ 
      1 \le \tilde{j}_{u+1} \le \dots \tilde{j}_{v-1} \le \tilde{j}_{v}} \\
    &\le& 2^{v-u-1} (\alpha(p))^{v-u-1},
    \EEqn
    where $a \wedge b = \min\{a,b\}$.
    
    In summary, we may conclude that
    \begin{equation*}
      \# \ca \le C(s) (\alpha(p))^{s-2} \bone\{\tilde{j}_v - \tilde{j}_u \le (s-1)\alpha(p)\}
    \end{equation*}
    where $C(s) = 2^{s-2}$.
\end{proof}

\begin{lemma}\label{Lemma:component-5}
Suppose that $j_u,j_v$ are fixed, then for any fixed $s\geq2$, it holds under Assumption \ref{Assumpt:component-dept}\ref{Assumpt:component-dept-1} that
\BEqn
& & \# \lrcp{(j_1,\dots,j_{u-1},j_{u+1},\dots,j_{v-1},j_{v+1},\dots,j_s):\ 
    \begin{array}{l}
      \cum\lrp{\tilde{x}_{j_1},\dots,\tilde{x}_{j_s}}\neq0 \vspace{1mm}\\ 
      1 \le j_1,\dots,j_s \le p 
    \end{array}} \\
&\le& C(s)(\alpha(p))^{s-2} \bone\lrcp{\lrabs{j_v - j_u} \le (s-1)\alpha(p)}. 
\EEqn
\end{lemma}
\begin{proof}
  When $|j_u-j_v|>(s-1)\alpha(p)$, it follows from Lemma \ref{Lemma:component-2} that the quantity of interest is zero. Otherwise, it follows from Lemma \ref{Lemma:component-4} that
  \BEqn
  & & \# \lrcp{(j_1,\dots,j_{u-1},j_{u+1},\dots,j_{v-1},j_{v+1},\dots,j_s):\ 
    \begin{array}{l}
      \cum\lrp{\tilde{x}_{j_1},\dots,\tilde{x}_{j_s}}\neq0 \vspace{1mm}\\ 
      1 \le j_1,\dots,j_s \le p 
    \end{array}} \\[2mm]
  &=& \sum\limits_{1 \le t \neq r \le s} 
  \# \lrcp{(j_1,\dots,j_{u-1},j_{u+1},\dots,j_{v-1},j_{v+1},\dots,j_s):\ 
    \begin{array}{l}
      \cum\lrp{\tilde{x}_{j_1},\dots,\tilde{x}_{j_s}}\neq0 \vspace{1mm}\\ 
      1 \le j_1,\dots,j_s \le p \vspace{1mm}\\
      j_u=\tilde{j}_t,\ j_v=\tilde{j}_r
    \end{array}} \\[2mm]
  &=& \sum\limits_{1 \le t \neq r \le s} (s-2)! 
  \# \lrcp{(\tilde{j}_1,\dots,\tilde{j}_{t \wedge r-1},\tilde{j}_{t \wedge r+1},\dots,\tilde{j}_{t \vee r-1},\tilde{j}_{t \vee r-1},\dots,j_s):\ 
    \begin{array}{l}
      \cum\lrp{\tilde{x}_{\tilde{j}_1},\dots,\tilde{x}_{\tilde{j}_s}}\neq0 \vspace{1mm}\\ 
      1 \le \tilde{j}_1 \le \dots \le \tilde{j}_s \le p
    \end{array}} \\[1mm]
  &\le& \sum\limits_{1 \le t \neq r \le s} 2^{s-2} (s-2)! (\alpha(p))^{s-2} \\
  &\le& C(s) (\alpha(p))^{s-2}
  \EEqn
  with $C(s) = 2^{s-2} s!$.
\end{proof}

A direct generalization of Lemma \ref{Lemma:component-5} is summarized as the following statement, for which we spare the proof.
\begin{lemma}\label{Lemma:component-6}
For any fixed $s\geq2$ and $\cd \subseteq \{1,\dots,s\}$, suppose that the values of $(j_u: u\in\cd)$ are all fixed, then it holds under Assumption \ref{Assumpt:component-dept}\ref{Assumpt:component-dept-1} that
\BEqn
& & \# \lrcp{(j_u:\ u\in\cd^c):\ \cum\lrp{\tilde{x}_{j_1},\dots,\tilde{x}_{j_s}}\neq0,\ 1 \le j_u \le p, \forall u\in \cd^c} \\
&\le& C(s)(\alpha(p))^{s-\#\cd} \bone\lrcp{\max\limits_{u\in\cd}\{j_u\} - \min\limits_{u\in\cd}\{j_u\}\le (s-1)\alpha(p)}. 
\EEqn
\end{lemma}

\begin{lemma}\label{Lemma:component-7}
For any fixed $s\geq2$ and fixed positive constant $c$, it holds under Assumption \ref{Assumpt:component-dept}\ref{Assumpt:component-dept-1} that
\begin{equation*}
    \# \lrcp{(j_1,\dots,j_s):\ \max\limits_{1\le u\le s}\{j_u\} - \min\limits_{1\le u\le s}\{j_u\} \le c\alpha(p),\ 1\le j_1,\dots,j_s\le p} 
    \le C(c,s) (\alpha(p))^{s-1}p.  
\end{equation*}
\end{lemma}
\begin{proof}
  Note that
  \BEqn
  & & \# \lrcp{(j_1,\dots,j_s):\ \max\limits_{1\le u\le s}\{j_u\} - \min\limits_{1\le u\le s}\{j_u\} \le c\alpha(p),\ 1\le j_1,\dots,j_s\le p} \\
  &=& (s!) \# \lrcp{(\tilde{j}_1,\dots,\tilde{j}_s):\ \tilde{j}_s - \tilde{j}_1 \le c\alpha(p),\ 1\le \tilde{j}_1 \le \dots \le \tilde{j}_s\le p} \\
  &=& (s!) 
  \lrp{\# \lrcp{(\tilde{j}_1,\tilde{j}_s):\ \tilde{j}_s - \tilde{j}_1 \le c\alpha(p),\ 1\le \tilde{j}_1 \le \tilde{j}_s\le p}} \\
  & & \times \lrp{\# \lrcp{(\tilde{j}_2,\dots,\tilde{j}_{s-1}):\ \tilde{j}_1 \le \tilde{j}_2 \le \dots \le \tilde{j}_{s-1}\le \tilde{j}_s}}
  \EEqn
  
  Now let $\tilde{j}_1$ range free over the set $\{1,\dots,p\}$, the total number of possible values that $\tilde{j}_1$ can take is at most $s$. Once the value of $\tilde{j}_1$ is fixed, the number of values that $\tilde{j}_s$ can take is at most $c\alpha(p)+1$ due to the constraint $\tilde{j}_s - \tilde{j}_1 \le c\alpha(p)$. In other words, ignoring $\tilde{j}_2,\dots,\tilde{j}_{s-1}$, the total number of possible pairs of $(\tilde{j}_1,\tilde{j}_s)$ is at most $(c\alpha(p)+1)p$.
  
  Assume that $\tilde{j}_1 = a \le b = \tilde{j}_s$ are fixed and $b-a \le c\alpha(p)$, then the values $\tilde{j}_2,\dots,\tilde{j}_{s-1}$ can take are chosen from the set $\{a,a+1,\dots,b-1,b\}$. It follows that the number of all the possible combinations of $(\tilde{j}_2,\dots,\tilde{j}_{s-1})$ is at most $(b-a+1)^{s-2}$, which is upper bounded by $(c\alpha(p)+1)^{s-2}$ regardless of values of $a,b$. 
  
  Finally, we obtain that
  \BEqn
  & & \# \lrcp{(j_1,\dots,j_s):\ \max\limits_{1\le u\le s}\{j_u\} - \min\limits_{1\le u\le s}\{j_u\} \le c\alpha(p),\ 1\le j_1,\dots,j_s\le p} \\
  &\le& (s!)  (c\alpha(p)+1)^{s-1} p
   \le  C(c,s) (\alpha(p))^{s-1}p,
  \EEqn
  where $C(c,s) = (s!) (c+1)^{s-1}$.
\end{proof}

Before presenting the next lemma, we introduce some notations to facilitate the subsequent analysis. For $i\geq1$, define $L_i = \blre{\lrp{|Z_1-Z_2|^2 - A}^i}$ with $A=\bbe{|Z_1-Z_2|^2}$, and
\BEqn
& & L_i^{XY} = \blre{\lrp{|X_1-Y_1|^2 - A^{XY}}^i}, \qquad A^{XY}=\bbe{|X_1-Y_1|^2}. \\
& & L_i^{X} = \blre{\lrp{|X_1-X_2|^2 - A^{X}}^i}, \qquad A^X = \bbe{|X_1-X_2|^2}. \\
& & L_i^{Y} = \blre{\lrp{|Y_1-Y_2|^2 - A^{Y}}^i}, \qquad A^Y = \bbe{|Y_1-Y_2|^2}.
\EEqn

We derive the order of $|L_i^{XY}|$ under Assumption \ref{Assumpt:component-dept}.
\begin{lemma}\label{Lemma:Li-1}
For any $2 \le i \le 64$, it holds under Assumption \ref{Assumpt:component-dept}\ref{Assumpt:component-dept-1}-\ref{Assumpt:component-dept-2} that 
\begin{equation*}
    \lrabs{L_i^{XY}} \le C(i,U^{\ast}) \lrp{\alpha(p)}^{\lrceil{i/2}} p^{\lrfloor{i/2}}
\end{equation*}
where $\lrfloor{x}$ denote the largest integer that is no larger than $x$ and $\lrceil{x}$ denote the smallest integer that is no smaller than $x$.
\end{lemma}
\begin{proof}
  Let $(\cs_1,\cs_2,\cs_3,\cs_4,\cs_5)$ denote a partition of the set $\{1,\dots,i\}$, that is, $\cs_1,\dots,\cs_5$ are pairwise disjoint and $\bigcup\limits_{t=1}^{5}\cs_t = \{1,\dots,s\}$. Define
  \begin{equation*}
      \cp = \lrcp{\pi:\ \pi \mbox{ is a partition of } \cs_1\cup\cs_3\cup\cs_5, \mbox{ and }|\cb|\geq2,\forall \cb\in\pi}
  \end{equation*}
  and similarly, define
  \begin{equation*}
      \cp' = \lrcp{\pi':\ \pi' \mbox{ is a partition of } \cs_2\cup\cs_4\cup\cs_5, \mbox{ and }|\cb'|\geq2,\forall \cb'\in\pi'}
  \end{equation*}
  
  Recall that $L_i^{XY} = \bbe{\bigp{\lrabs{X_1-Y_1}^2 - A^{XY}}^i}$, and note that
  \begin{equation*}
      \lrabs{X_1-Y_1}^2 - A^{XY}
    = \sum\limits_{j=1}^{p} 
    \lrp{  \bigp{\tilde{x}_{1 j}^{2} - \sigma_{X, j}^2}
         + \bigp{\tilde{y}_{1 j}^{2} - \sigma_{Y, j}^2}
         + 2\Delta_j \tilde{x}_{1j}
         - 2\Delta_j \tilde{y}_{1j}
         - 2\tilde{x}_{1 j}\tilde{y}_{1 j}}
  \end{equation*}
  then it follows that
  \BEqn
  & & L_i^{XY} \\
  &=& \sum\limits_{j_1,\dots,j_i=1}^{p}
      \bbe{\prod\limits_{s=1}^{i} 
           \bigp{  (\tilde{x}_{1 j_s}^2 - \sigma_{X, j_s}^2)
                 + (\tilde{y}_{1 j_s}^2 - \sigma_{Y, j_s}^2)
                 + 2\Delta_{j_s}\tilde{x}_{1 j_s}
                 - 2\Delta_{j_s}\tilde{y}_{1 j_s}
                 - 2\tilde{x}_{1 j_s}\tilde{y}_{1 j_s}}} \\
  &=& \sum\limits_{(\cs_1,\dots,\cs_5)} \sum\limits_{j_1,\dots,j_i=1}^{p} 
      \be\big{[}
           \bigp{\prod\limits_{s_1\in\cs_1} (\tilde{x}_{1 j_{s_1}}^2 - \sigma_{X, j_{s_1}}^2)}
           \bigp{\prod\limits_{s_2\in\cs_2} (\tilde{y}_{1 j_{s_2}}^2 - \sigma_{Y, j_{s_2}}^2)}
           \bigp{\prod\limits_{s_3\in\cs_3} (2\Delta_{j_{s_3}}\tilde{x}_{1 j_{s_3}})} \\
  & & \hspace{8em}
           \bigp{\prod\limits_{s_4\in\cs_4} (-2\Delta_{j_{s_4}}\tilde{y}_{1 j_{s_4}})}
           \bigp{\prod\limits_{s_5\in\cs_5} (-2\tilde{x}_{1 j_{s_5}}\tilde{y}_{1 j_{s_5}})} 
      \big{]} \\
  &=& \sum\limits_{(\cs_1,\dots,\cs_5)} \sum\limits_{j_1,\dots,j_i=1}^{p} (-1)^{|\cs_4|+|\cs_5|}(2)^{|\cs_3|+|\cs_4|+|\cs_5|} \\
  & & \hspace{7em}
      \times \bbe{\bigp{\prod\limits_{s_1\in\cs_1} (\tilde{x}_{1 j_{s_1}}^2 - \sigma_{X, j_{s_1}}^2)}
                  \bigp{\prod\limits_{s_3\in\cs_3} (\Delta_{j_{s_3}}\tilde{x}_{1 j_{s_3}})}
                  \bigp{\prod\limits_{s_5\in\cs_5} \tilde{x}_{1 j_{s_5}}}} \\
  & & \hspace{7em}
      \times \bbe{\bigp{\prod\limits_{s_2\in\cs_2} (\tilde{y}_{1 j_{s_2}}^2 - \sigma_{Y, j_{s_2}}^2)}
                  \bigp{\prod\limits_{s_4\in\cs_4} (\Delta_{j_{s_4}}\tilde{y}_{1 j_{s_4}})}
                  \bigp{\prod\limits_{s_5\in\cs_5} \tilde{y}_{1 j_{s_5}}}} \\
  &=& \sum\limits_{(\cs_1,\dots,\cs_5)} \
      \sum\limits_{\substack{1\le j_{s_1}\le p \\ s_1\in\cs_1}} \
      \sum\limits_{\substack{1\le j_{s_2}\le p \\ s_2\in\cs_2}} \ 
      \sum\limits_{\substack{1\le j_{s_3}\le p \\ s_3\in\cs_3}} \ 
      \sum\limits_{\substack{1\le j_{s_4}\le p \\ s_4\in\cs_4}} \ 
      \sum\limits_{\substack{1\le j_{s_5}\le p \\ s_5\in\cs_5}} \ 
      (-1)^{|\cs_4|+|\cs_5|}(2)^{|\cs_3|+|\cs_4|+|\cs_5|} \\
  & & \hspace{2em} 
      \sum\limits_{\pi\in\cp} \prod\limits_{\cb\in\pi} 
      \cum\bigp{\tilde{x}_{1 j_s}^2-\sigma_{X, j_s}^2, \Delta_{j_t}\tilde{x}_{1 j_t}, \tilde{x}_{1 j_u}:\  s\in\cb\cap\cs_1, t\in\cb\cap\cs_3, u\in\cb\cap\cs_5} \\
  & & \hspace{2em} 
      \sum\limits_{\pi'\in\cp'} \prod\limits_{\cb'\in\pi'} 
      \cum\bigp{\tilde{y}_{1 j_{s'}}^2-\sigma_{Y, j_{s'}}^2, \Delta_{j_{t'}}\tilde{y}_{1 j_{t'}}, \tilde{y}_{1 j_{u'}}:\ s'\in\cb'\cap\cs_2, t'\in\cb'\cap\cs_4, u'\in\cb'\cap\cs_5} \\     
  &=& \sum\limits_{(\cs_1,\dots,\cs_5)} \
      \sum\limits_{\pi\in\cp} \
      \sum\limits_{\pi'\in\cp'} \ 
      \sum\limits_{\substack{1\le j_{s_5}\le p \\ s_5\in\cs_5}} \ 
      (-1)^{|\cs_4|+|\cs_5|}(2)^{|\cs_3|+|\cs_4|+|\cs_5|} \\
  & & \hspace{1em} 
      \sum\limits_{\substack{1\le j_{s_1},j_{s_3}\le p \\ s_1\in\cs_1 \\ s_3\in\cs_3}} \
      \prod\limits_{\cb\in\pi}
      \cum\bigp{\tilde{x}_{1 j_s}^2-\sigma_{X, j_s}^2, \Delta_{j_t}\tilde{x}_{1 j_t}, \tilde{x}_{1 j_u}:\  s\in\cb\cap\cs_1, t\in\cb\cap\cs_3, u\in\cb\cap\cs_5} \\
  & & \hspace{1em} 
      \sum\limits_{\substack{1\le j_{s_2},j_{s_4}\le p \\ s_2\in\cs_2 \\ s_4\in\cs_4}} \
      \prod\limits_{\cb'\in\pi'}
      \cum\bigp{\tilde{y}_{1 j_{s'}}^2-\sigma_{Y, j_{s'}}^2, \Delta_{j_{t'}}\tilde{y}_{1 j_{t'}}, \tilde{y}_{1 j_{u'}}:\ s'\in\cb'\cap\cs_2, t'\in\cb'\cap\cs_4, u'\in\cb'\cap\cs_5} \\           
  \EEqn
  
  Note that $\bigcup\limits_{\cb\in\pi} (\cb\cap\cs_1) = \cs_1$ and $\bigcup\limits_{\cb\in\pi} (\cb\cap\cs_3) = \cs_3$, then by using the fact that $\pi$ is a partition of $\cs_1 \cup \cs_3 \cup \cs_5$, we obtain that
  \BEqn
  & & \sum\limits_{\substack{1\le j_{s_1},j_{s_3}\le p \\ s_1\in\cs_1 \\ s_3\in\cs_3}} \
      \prod\limits_{\cb\in\pi}
      \cum\bigp{\tilde{x}_{1 j_s}^2-\sigma_{X, j_s}^2, \Delta_{j_t}\tilde{x}_{1 j_t}, \tilde{x}_{1 j_u}:\  s\in\cb\cap\cs_1, t\in\cb\cap\cs_3, u\in\cb\cap\cs_5} \\
  &=& \prod\limits_{\cb\in\pi}
      \sum\limits_{\substack{1\le j_{s_1}\le p \\ s_1\in\cb\cap\cs_1}} \
      \sum\limits_{\substack{1\le j_{s_3}\le p \\ s_3\in\cb\cap\cs_3}} \
      \cum\bigp{\tilde{x}_{1 j_s}^2-\sigma_{X, j_s}^2, \Delta_{j_t}\tilde{x}_{1 j_t}, \tilde{x}_{1 j_u}:\  s\in\cb\cap\cs_1, t\in\cb\cap\cs_3, u\in\cb\cap\cs_5}
  \EEqn
  and similarly, 
  \BEqn
  & & \sum\limits_{\substack{1\le j_{s_2},j_{s_4}\le p \\ s_2\in\cs_2 \\ s_4\in\cs_4}} \
      \prod\limits_{\cb'\in\pi'}
      \cum\bigp{\tilde{y}_{1 j_{s'}}^2-\sigma_{Y, j_{s'}}^2, \Delta_{j_{t'}}\tilde{y}_{1 j_{t'}}, \tilde{y}_{1 j_{u'}}:\ s'\in\cb'\cap\cs_2, t'\in\cb'\cap\cs_4, u'\in\cb'\cap\cs_5} \\
  &=& \prod\limits_{\cb'\in\pi'}
      \sum\limits_{\substack{1\le j_{s_2}\le p \\ s_2\in\cb'\cap\cs_2}} \
      \sum\limits_{\substack{1\le j_{s_4}\le p \\ s_4\in\cb'\cap\cs_4}} \
      \cum\bigp{\tilde{y}_{1 j_{s'}}^2-\sigma_{Y, j_{s'}}^2, \Delta_{j_{t'}}\tilde{y}_{1 j_{t'}}, \tilde{y}_{1 j_{u'}}:\ s'\in\cb'\cap\cs_2, t'\in\cb'\cap\cs_4, u'\in\cb'\cap\cs_5}
  \EEqn
  
  Then from above, we obtain an upper bound of $\lrabs{L_i^{XY}}$, that is,
  \BEqn
  & & \lrabs{L_i^{XY}} \\[2mm]
  &\le& 2^i
      \sum\limits_{(\cs_1,\dots,\cs_5)} \
      \sum\limits_{\pi\in\cp} \
      \sum\limits_{\pi'\in\cp'} \ 
      \sum\limits_{\substack{1\le j_{s_5}\le p \\ s_5\in\cs_5}} \ 
       \\
  & & \hspace{2em}      
      \prod\limits_{\cb\in\pi}
      \sum\limits_{\substack{1\le j_{s_1},j_{s_3}\le p \\ s_1\in\cb\cap\cs_1 \\ s_3\in\cb\cap\cs_3}} \
      \lrabs{\cum\bigp{\tilde{x}_{1 j_s}^2-\sigma_{X, j_s}^2, \Delta_{j_t}\tilde{x}_{1 j_t}, \tilde{x}_{1 j_u}:\  s\in\cb\cap\cs_1, t\in\cb\cap\cs_3, u\in\cb\cap\cs_5}} \\
  & & \hspace{2em}      
      \prod\limits_{\cb'\in\pi'}
      \sum\limits_{\substack{1\le j_{s_2},j_{s_4}\le p \\ s_2\in\cb'\cap\cs_2 \\ s_4\in\cb'\cap\cs_4}} \
      \lrabs{\cum\bigp{\tilde{y}_{1 j_{s'}}^2-\sigma_{Y, j_{s'}}^2, \Delta_{j_{t'}}\tilde{y}_{1 j_{t'}}, \tilde{y}_{1 j_{u'}}:\ s'\in\cb'\cap\cs_2, t'\in\cb'\cap\cs_4, u'\in\cb'\cap\cs_5}}
  \EEqn
  
  For any $(\cs_1,\dots,\cs_5)$ and any corresponding $\pi\in\cp$, $\cb\in\pi$, and for any fixed $(j_u:\ u\in\cb\cap\cs_5)$, it follows from Lemma \ref{Lemma:component-6} that, under Assumption \ref{Assumpt:component-dept}\ref{Assumpt:component-dept-1}-\ref{Assumpt:component-dept-2},  
  \BEqn
  & & \sum\limits_{\substack{1\le j_{s_1},j_{s_3}\le p \\ s_1\in\cb\cap\cs_1 \\ s_3\in\cb\cap\cs_3}} \
      \lrabs{\cum\bigp{\tilde{x}_{1 j_s}^2-\sigma_{X, j_s}^2, \Delta_{j_t}\tilde{x}_{1 j_t}, \tilde{x}_{1 j_u}:\  s\in\cb\cap\cs_1, t\in\cb\cap\cs_3, u\in\cb\cap\cs_5}} \\
  &\le& C(U^{\ast})
      \# \lrcp{
         \begin{array}{l}
           (j_s,j_t:\ s\in\cb\cap\cs_1, t\in\cb\in\cs_3):  \vspace{1mm}\\
           \cum\bigp{\tilde{x}_{1 j_s}^2-\sigma_{X, j_s}^2, \Delta_{j_t}\tilde{x}_{1 j_t}, \tilde{x}_{1 j_u}:\  s\in\cb\cap\cs_1, t\in\cb\cap\cs_3, u\in\cb\cap\cs_5} \neq 0 
      \end{array}} \\[2mm]
  &\le& C(i,U^{\ast}) (\alpha(p)) ^{|\cb\cap\cs_1|+|\cb\cap\cs_3|} \\
  & & \times \bone\{\max\limits_{u\in\cb\cap\cs_5}\{j_u\} - \min\limits_{u\in\cb\cap\cs_5}\{j_u\} \le \bigp{|\cb\cap\cs_1| + |\cb\cap\cs_3| + |\cb\cap\cs_5| - 1} \alpha(p) \}.
  \EEqn
  
  It follows from the fact that $\bigcup\limits_{\cb\in\pi} \cb\cap\cs_1 = \cs_1$ and $\bigcup\limits_{\cb\in\pi} \cb\cap\cs_3 = \cs_3$ that 
  \begin{equation*}
      \sum\limits_{\cb\in\pi} \lrp{|\cb\cap\cs_1| + |\cb\cap\cs_3|} = |\cs_1| + |\cs_3|
  \end{equation*}
  hence we further have
  \BEqn
  & & \prod\limits_{\cb\in\pi} 
      \sum\limits_{\substack{1\le j_{s_1},j_{s_3}\le p \\ s_1\in\cb\cap\cs_1 \\ s_3\in\cb\cap\cs_3}} \
      \lrabs{\cum\bigp{\tilde{x}_{1 j_s}^2-\sigma_{X, j_s}^2, \Delta_{j_t}\tilde{x}_{1 j_t}, \tilde{x}_{1 j_u}:\  s\in\cb\cap\cs_1, t\in\cb\cap\cs_3, u\in\cb\cap\cs_5}} \\
  &\le& \prod\limits_{\cb\in\pi} C(i,U^{\ast}) (\alpha(p)) ^{|\cb\cap\cs_1|+|\cb\cap\cs_3|} \\
  & & \hspace{2em}
      \times \bone\{\max\limits_{u\in\cb\cap\cs_5}\{j_u\} - \min\limits_{u\in\cb\cap\cs_5}\{j_u\} \le \bigp{|\cb\cap\cs_1| + |\cb\cap\cs_3| + |\cb\cap\cs_5| - 1} \alpha(p)\} \\[2mm]
  &\le& C(i,U^{\ast}) (\alpha(p))^{\sum\limits_{\cb\in\pi}\bigp{|\cb\cap\cs_1|+|\cb\cap\cs_3|}} 
   \prod\limits_{\cb\in\pi} \bone\{\max\limits_{u\in\cb\cap\cs_5}\{j_u\} - \min\limits_{u\in\cb\cap\cs_5}\{j_u\} \le i\alpha(p)\} \\[2mm]
  &=& C(i,U^{\ast}) (\alpha(p))^{|\cs_1|+|\cs_3|} 
   \prod\limits_{\cb\in\pi} \bone\{\max\limits_{u\in\cb\cap\cs_5}\{j_u\} - \min\limits_{u\in\cb\cap\cs_5}\{j_u\} \le i\alpha(p)\}
  \EEqn
  
  Using similar arguments, we also obtain that for any $\pi'\in\cp'$,
  \BEqn
  & & \prod\limits_{\cb'\in\pi'}\sum\limits_{\substack{1\le j_{s_2},j_{s_4}\le p \\ s_2\in\cb'\cap\cs_2 \\ s_4\in\cb'\cap\cs_4}} \
      \lrabs{\cum\bigp{\tilde{y}_{1 j_{s'}}^2-\sigma_{Y, j_{s'}}^2, \Delta_{j_{t'}}\tilde{y}_{1 j_{t'}}, \tilde{y}_{1 j_{u'}}:\ s'\in\cb'\cap\cs_2, t'\in\cb'\cap\cs_4, u'\in\cb'\cap\cs_5}} \\
  &\le& C(i,U^{\ast}) (\alpha(p))^{|\cs_2|+|\cs_4|} 
   \prod\limits_{\cb'\in\pi'} \bone\{\max\limits_{u'\in\cb'\cap\cs_5}\{j_{u'}\} - \min\limits_{u'\in\cb'\cap\cs_5}\{j_{u'}\} \le i\alpha(p)\}
  \EEqn
  
  Now we have
  \BEqn
  & & \lrabs{L_i^{XY}} \\
  &\le& C(i,U^{\ast}) 
        \sum\limits_{(\cs_1,\dots,\cs_5)}
        \sum\limits_{\pi\in\cp} \
        \sum\limits_{\pi'\in\cp'} \ 
        (\alpha(p))^{|\cs_1| + |\cs_2| + |\cs_3| + |\cs_4|} \\
  & & \hspace{4em}
      \times \sum\limits_{\substack{1\le j_{s_5}\le p \\ s_5\in\cs_5}} \ 
      \lrp{\prod\limits_{\cb\in\pi} \bone\{\max\limits_{u\in\cb\cap\cs_5}\{j_u\} - \min\limits_{u\in\cb\cap\cs_5}\{j_u\} \le i\alpha(p)\}} \\
  & & \hspace{8em}
      \times \lrp{\prod\limits_{\cb'\in\pi'} \bone\{\max\limits_{u'\in\cb'\cap\cs_5}\{j_{u'}\} - \min\limits_{u'\in\cb'\cap\cs_5}\{j_{u'}\} \le i\alpha(p)\}} \\
  &=& C(i,U^{\ast}) 
      \sum\limits_{(\cs_1,\dots,\cs_5)}
      \sum\limits_{\pi\in\cp} \
      \sum\limits_{\pi'\in\cp'} \ 
      (\alpha(p))^{|\cs_1| + |\cs_2| + |\cs_3| + |\cs_4|} \\
  & & \hspace{4em}      
      \times \# \lrcp{
             (j_u:\ u\in\cs_5):\ 
             \begin{array}{ll}
                \max\limits_{u\in\cb\cap\cs_5}\{j_u\} - \min\limits_{u\in\cb\cap\cs_5}\{j_u\} \le i\alpha(p), & \forall\cb\in\pi \\
                \max\limits_{u'\in\cb'\cap\cs_5}\{j_{u'}\} - \min\limits_{u'\in\cb'\cap\cs_5}\{j_{u'}\} \le i\alpha(p), & \forall\cb'\in\pi' \\
             \end{array}} \\
  &\le& C(i,U^{\ast})
        \sum\limits_{(\cs_1,\dots,\cs_5)}
        \sum\limits_{\pi\in\cp} \
        \sum\limits_{\pi'\in\cp'} \ 
        (\alpha(p))^{|\cs_1| + |\cs_2| + |\cs_3| + |\cs_4|} \\
  & & \hspace{4em}      
      \times \min\lrcp{ 
             \begin{array}{l}
                \# \bigcp{(j_u:\ u\in\cs_5):\ \max\limits_{u\in\cb\cap\cs_5}\{j_u\} - \min\limits_{u\in\cb\cap\cs_5}\{j_u\} \le i\alpha(p),\ \forall\cb\in\pi}, \vspace{2mm}\\
                \# \bigcp{(j_u:\ u\in\cs_5):\ \max\limits_{u'\in\cb'\cap\cs_5}\{j_{u'}\} - \min\limits_{u'\in\cb'\cap\cs_5}\{j_{u'}\} \le i\alpha(p),\ \forall\cb'\in\pi'}
             \end{array}} \\
  &\le& C(i,U^{\ast}) 
        \sum\limits_{(\cs_1,\dots,\cs_5)}
        \sum\limits_{\pi\in\cp} \
        \sum\limits_{\pi'\in\cp'} \ 
        (\alpha(p))^{|\cs_1| + |\cs_2| + |\cs_3| + |\cs_4|} \\
  & & \hspace{4em}      
      \times \min\lrcp{ 
             \begin{array}{l}
                \prod\limits_{\cb\in\pi} 
                \# \bigcp{(j_u:\ u\in\cb\cap\cs_5):\ \max\limits_{u\in\cb\cap\cs_5}\{j_u\} - \min\limits_{u\in\cb\cap\cs_5}\{j_u\} \le i\alpha(p)}, \vspace{2mm}\\
                \prod\limits_{\cb'\in\pi'}
                \# \bigcp{(j_u:\ u\in\cb'\cap\cs_5):\ \max\limits_{u'\in\cb'\cap\cs_5}\{j_{u'}\} - \min\limits_{u'\in\cb'\cap\cs_5}\{j_{u'}\} \le i\alpha(p)}
             \end{array}}
  \EEqn
  where the last step follows from the fact that $\{\cb\cap\cs_5,\ \cb\in\pi\}$ forms a partition of $\cs_5$ and $\{\cb'\cap\cs_5,\ \cb'\in\pi'\}$ does as well. 
  
  It follows from Lemma \ref{Lemma:component-7} that, 
  \BEqn
  & & \prod\limits_{\cb\in\pi} 
      \# \bigcp{(j_u:\ u\in\cb\cap\cs_5):\ \max\limits_{u\in\cb\cap\cs_5}\{j_u\} - \min\limits_{u\in\cb\cap\cs_5}\{j_u\} \le i\alpha(p)} \\
  &\le& C(i) \prod\limits_{\cb\in\pi} \lrp{(\alpha(p))^{|\cb\cap\cs_5|-1} p} \\
  &=& C(i) (\alpha(p))^{|\cs_5|-|\pi|} p^{|\pi|}
  \EEqn
  Similarly, we obtain that
  \begin{equation*}
      \prod\limits_{\cb'\in\pi'}
      \# \bigcp{(j_u:\ u\in\cb'\cap\cs_5):\ \max\limits_{u'\in\cb'\cap\cs_5}\{j_{u'}\} - \min\limits_{u'\in\cb'\cap\cs_5}\{j_{u'}\} \le i\alpha(p)}
  \le C(i) (\alpha(p))^{|\cs_5|-|\pi'|} p^{|\pi'|}
  \end{equation*}
  which further implies that
  \BEqn
  & & \lrabs{L_i^{XY}} \\
  &\le& C(i,U^{\ast}) 
        \sum\limits_{(\cs_1,\dots,\cs_5)}
        \sum\limits_{\pi\in\cp} \
        \sum\limits_{\pi'\in\cp'} \ 
        (\alpha(p))^{|\cs_1| + |\cs_2| + |\cs_3| + |\cs_4|} 
        \min\bigcp{(\alpha(p))^{|\cs_5|-|\pi|} p^{|\pi|}, (\alpha(p))^{|\cs_5|-|\pi'|} p^{|\pi'|}} \\
  &=& C(i,U^{\ast})
      \sum\limits_{(\cs_1,\dots,\cs_5)} (\alpha(p))^{|\cs_1| + |\cs_2| + |\cs_3| + |\cs_4| + |\cs_5|} 
      \sum\limits_{\pi\in\cp} \
      \sum\limits_{\pi'\in\cp'} \
      \min\bigcp{(\alpha(p))^{-|\pi|} p^{|\pi|}, (\alpha(p))^{-|\pi'|} p^{|\pi'|}} \\
  &=& C(i,U^{\ast}) (\alpha(p))^{i}
      \sum\limits_{(\cs_1,\dots,\cs_5)}    
      \sum\limits_{\pi\in\cp} \
      \sum\limits_{\pi'\in\cp'} \ 
      \lrp{\frac{p}{\alpha(p)}}^{|\pi| \wedge |\pi'|}
  \EEqn
  where the last equality follows from the fact that $(\cs_1,\dots,\cs_5)$ forms a partition of $\{1,\dots,i\}$ and thus $|\cs_1|+\dots+|\cs_5|=i$, and $1\le\alpha(p)\prec p$.
  
  Now it remains to find an upper bound of $|\pi|\wedge|\pi'|$ over all possible $(\cs_1,\dots,\cs_5)$ and $\pi,\pi'$. By the definition of $\pi,\pi'$, we have 
  \begin{equation*}
      |\pi| \le \lrfloor{\bigp{|\cs_1| + |\cs_3| + |\cs_5|}/2}
  \end{equation*}
  since $\pi$ is a partition of $\cs_1\cup\cs_3\cup\cs_5$ and $|\cb|\geq2$ for each $\cb\in\pi$. When $|\cs_1| + |\cs_3| + |\cs_5|$ is even, the upper bound is attained when $|\cb|=2$ for all $\cb\in\pi$, and when $|\cs_1| + |\cs_3| + |\cs_5|$ is odd, it's attained when all but one $\cb\in\pi$ contain 2 elements. Similarly,
  \begin{equation*}
      |\pi'| \le \lrfloor{\bigp{|\cs_2| + |\cs_4| + |\cs_5|}/2}
  \end{equation*}
  and consequently, by noting that the number of all possible combinations of $(\cs_1,\dots,\cs_5)$ and $\pi,\pi'$ is a constant depends only on $i$, we may conclude that
  \BEqn
  & & \lrabs{L_i^{XY}} \\
  &\le& C(i,U^{\ast}) (\alpha(p))^{i}
        \sum\limits_{(\cs_1,\dots,\cs_5)}    
        \sum\limits_{\pi\in\cp} \
        \sum\limits_{\pi'\in\cp'} \ 
        \lrp{\frac{p}{\alpha(p)}}^{\lrfloor{\bigp{|\cs_1| + |\cs_3| + |\cs_5|}/2} \wedge \lrfloor{\bigp{|\cs_2| + |\cs_4| + |\cs_5|}/2}} \\
 &\le&  C(i,U^{\ast}) (\alpha(p))^{i}
        \sum\limits_{(\cs_1,\dots,\cs_5)}    
        \sum\limits_{\pi\in\cp} \
        \sum\limits_{\pi'\in\cp'} \
        \lrp{\frac{p}{\alpha(p)}}^{\lrfloor{i/2}} \\
  &\le& C(i,U^{\ast}) \lrp{\alpha(p)}^{\lrceil{i/2}} p^{\lrfloor{i/2}},
  \EEqn
  which completes the proof.
\end{proof}
\begin{remark}
  For the special case when $\alpha(p)\equiv\alpha$ is a fixed constant independent of $p$, the results is simplified as $|L_i^{XY}| \le C(i,U^{\ast},\alpha) p^{\lrfloor{i/2}}$.
\end{remark}

\begin{lemma}\label{Lemma:Li-2}
For any $2 \le i \le 64$, it holds under Assumption \ref{Assumpt:component-dept}\ref{Assumpt:component-dept-1}-\ref{Assumpt:component-dept-2} that 
\BEqn
& & \lrabs{L_i^{X}} \le C(i,U^{\ast}) \lrp{\alpha(p)}^{\lrceil{i/2}} p^{\lrfloor{i/2}}, \\
& & \lrabs{L_i^{Y}} \le C(i,U^{\ast}) \lrp{\alpha(p)}^{\lrceil{i/2}} p^{\lrfloor{i/2}}
\EEqn
where $\lrfloor{x}$ denote the largest integer that is no larger than $x$ and $\lrceil{x}$ denote the smallest integer that is no smaller than $x$.
\end{lemma}
\begin{proof}
  Recall that $L_i^{X} = \bbe{\bigp{\lrabs{X_1-X_2}^2 - A^X}^i}$ and $X_1,X_2$ are iid copies of $X$, then the proof of \ref{Lemma:Li-1} is still valid for $L_i^{X}$. Similarly, the statement of $L_i^Y$ can be proved in the same way.
\end{proof}

\begin{lemma}\label{Lemma:Li-3}
For any $2 \le i \le 64$, when Assumption \ref{Assumpt:component-dept}\ref{Assumpt:component-dept-1}-\ref{Assumpt:component-dept-2} are satisfied, it holds under the null that 
\begin{equation*}
    \lrabs{L_i} \le C(i,U^{\ast}) \lrp{\alpha(p)}^{\lrceil{i/2}} p^{\lrfloor{i/2}},
\end{equation*}
where $\lrfloor{x}$ denote the largest integer that is no larger than $x$ and $\lrceil{x}$ denote the smallest integer that is no smaller than $x$.
\end{lemma}
\begin{proof}
  By the definition of $Z$, it holds under the null that $X =^d Y =^d Z$, and the result directly follows from Lemma \ref{Lemma:Li-2}.
\end{proof}

Let $\mu = \bbe{Z}$ denote the mean vector of $Z$ and $\tilde{Z} = Z-\mu$ be the centered version of $Z$. We use $\sigma_{j}^2 = \bbe{\tilde{z}_j^2}$ to denote the componentwise variance of $Z$. Additionally, we introduce the following notations:
\BEqn
& & N_{i_1 i_2} = \bbe{\bigp{|Z_1-Z_2|^2 - A}^{i_1} \bigp{|Z_1-Z_3|^2 - A}^{i_2}}, \qquad i_1+i_2\geq2,\\
& & D_{i_1 i_2 i_3} = \bbe{\bigp{|Z_1-Z_2|^2 - A}^{i_1} \bigp{|Z_1-Z_3|^2 - A}^{i_2} \bigp{|Z_2-Z_4^2 - A|}^{i_3}}, \qquad i_1+i_2+i_3\geq3,\\
& & M_{i_1 i_2 i_3 i_4} = \bbe{\bigp{|Z_1-Z_2|^2 - A}^{i_1} \bigp{|Z_1-Z_3|^2 - A}^{i_2} \bigp{|Z_2-Z_4|^2 - A}^{i_3} \bigp{|Z_3-Z_4|^2 - A}^{i_4}}, \\
& & \tilde{M}_{i_1 i_2 i_3 i_4} = \bbe{\bigp{|Z_1-Z_2|^2 - A}^{i_1} \bigp{|Z_1-Z_3|^2 - A}^{i_2} \bigp{|Z_2-Z_4|^2 - A}^{i_3} \bigp{|Z_2-Z_5|^2 - A}^{i_4}}, \\
& & \hspace{30em} i_1+i_2+i_3+i_4\geq4.
\EEqn

\begin{lemma}\label{Lemma:NDM-1}
Suppose that Assumption \ref{Assumpt:component-dept}\ref{Assumpt:component-dept-1}-\ref{Assumpt:component-dept-2} are satisfied, then it holds under the null that 
\begin{equation*}
\begin{array}{ll}
    \lrabs{N_{i_1 i_2}} \le C(i_1,i_2,U^{\ast})\lrp{\alpha(p)}^{\lrceil{(i_1+i_2)/2}} p^{\lrfloor{(i_1+i_2)/2}}, & \forall 2 \le i_1+i_2 \le 64 \vspace{1mm}\\
    \lrabs{D_{i_1 i_2 i_3}} \le C(i_1,i_2,i_3,U^{\ast})\lrp{\alpha(p)}^{\lrceil{(i_1+i_2+i_3)/2}} p^{\lrfloor{(i_1+i_2+i_3)/2}}, & \forall 3 \le i_1+i_2+i_3 \le 64 \vspace{1mm}\\
    \lrabs{M_{i_1 i_2 i_3 i_4}} \le C(i_1,i_2,i_3,i_4,U^{\ast})\lrp{\alpha(p)}^{\lrceil{(i_1+i_2+i_3+i_4)/2}} p^{\lrfloor{(i_1+i_2+i_3+i_4)/2}}, & \forall 4 \le i_1+i_2+i_3+i_4 \le 64 \vspace{1mm}\\
    \lrabs{\tilde{M}_{i_1 i_2 i_3 i_4}} \le C(i_1,i_2,i_3,i_4,U^{\ast})\lrp{\alpha(p)}^{\lrceil{(i_1+i_2+i_3+i_4)/2}} p^{\lrfloor{(i_1+i_2+i_3+i_4)/2}}, & \forall 4 \le i_1+i_2+i_3+i_4 \le 64 \vspace{1mm}\\
\end{array}
\end{equation*}
where $\lrfloor{x}$ denote the largest integer that is no larger than $x$ and $\lrceil{x}$ denote the smallest integer that is no smaller than $x$.
\end{lemma}
\begin{proof}
  We provide the detailed proof for the statement regarding $D_{i_1 i_2 i_3}$. Under the null, it holds that $X =^d Y =^d Z$, then under  Assumption \ref{Assumpt:component-dept}\ref{Assumpt:component-dept-1}-\ref{Assumpt:component-dept-2}, the distribution $Z = (z_1,\dots,z_j)$ has $\alpha(p)$-dependent components and a uniform upper bound over the component moments. Let $(\cs_1,\cs_2,\cs_3)$ be a partition of the set $\{1,\dots,i_1\}$, $(\ct_1,\ct_2,\ct_3)$ be a partition of the set $\{i_1+1,\dots,i_1+i_2\}$, and $(\cv_1,\cv_2,\cv_3)$ be a partition of the set $\{i_1+i_2+1,\dots,i_1+i_2+i_3\}$. Define the notations $\cp_1,\dots,\cp_4$ as follows:
  \begin{equation*}
      \begin{array}{l}
          \cp_1 = \lrcp{\pi_1:\ \pi_1 \mbox{ is a partition of } \bigp{\cs_1\cup\cs_3}\bigcup\bigp{\ct_1\cup\ct_3}, \mbox{ and }|\cb_1|\geq2,\forall \cb_1\in\pi_1} \vspace{1mm}\\
          \cp_2 = \lrcp{\pi_2:\ \pi_2 \mbox{ is a partition of } \bigp{\cs_2\cup\cs_3}\bigcup\bigp{\cv_1\cup\cv_3}, \mbox{ and }|\cb_2|\geq2,\forall \cb_2\in\pi_2} \vspace{1mm}\\
          \cp_3 = \lrcp{\pi_3:\ \pi_3 \mbox{ is a partition of } \ct_2\cup\ct_3, \mbox{ and }|\cb_3|\geq2,\forall \cb_3\in\pi_3} \vspace{1mm}\\
          \cp_4 = \lrcp{\pi_4:\ \pi_4 \mbox{ is a partition of } \cv_2\cup\cv_3, \mbox{ and }|\cb_4|\geq2,\forall \cb_4\in\pi_4} \vspace{1mm}\\
      \end{array}
  \end{equation*}
  
  Recall that $D_{i_1 i_2 i_3} = \bbe{\bigp{|Z_1-Z_2| - A}^{i_1} \bigp{|Z_1-Z_3| - A}^{i_2} \bigp{|Z_2-Z_4|}^{i_3}}$, then using the notations above, we can rewrite $D_{i_1 i_2 i_3}$ as follows:
  \BEqn
  & & D_{i_1 i_2 i_3} \\
  &=& \sum\limits_{j_1,\dots,j_{i_1+i_2+i_3}=1}^{p} 
      \be\big{[}
      \prod\limits_{s=1}^{i_1} \bigp{(\tilde{z}_{1 j_{s}}^2-\sigma_{j_{s}}^2) + (\tilde{z}_{2 j_{s}}^2-\sigma_{j_{s}}^2)) - 2\tilde{z}_{1 j_{s}}\tilde{z}_{2 j_{s}}} \\
  & & \hspace{7em}      
      \times \prod\limits_{t=i_1+1}^{i_1+i_2} \bigp{(\tilde{z}_{1 j_{t}}^2-\sigma_{j_{t}}^2) + (\tilde{z}_{3 j_{t}}^2-\sigma_{j_{t}}^2)) - 2\tilde{z}_{1 j_{t}}\tilde{z}_{3 j_{t}}} \\
  & & \hspace{7em}      
      \times \prod\limits_{v=i_1+i_2+1}^{i_1+i_2+i_3} \bigp{(\tilde{z}_{2 j_{v}}^2-\sigma_{j_{v}}^2) + (\tilde{z}_{4 j_{v}}^2-\sigma_{j_{v}}^2)) - 2\tilde{z}_{2 j_{v}}\tilde{z}_{4 j_{v}}}
      \big{]} \\
  &=& \sum\limits_{(\cs_1,\cs_2,\cs_3,\ct_1,\ct_2,\ct_3,\cv_1,\cv_2,\cv_3)}
      \sum\limits_{j_1,\dots,j_{i_1+i_2+i_3}=1}^{p} 
      (-2)^{|\cs_3|+|\ct_3|+|\cv_3|} \\
  & & \hspace{4em} \times
      \bbe{\bigp{\prod\limits_{s_1\in\cs_1} (\tilde{z}_{1 j_{s_1}}^2-\sigma_{j_{s_1}}^2)}
           \bigp{\prod\limits_{s_3\in\cs_3} \tilde{z}_{1 j_{s_3}}}
           \bigp{\prod\limits_{t_1\in\ct_1} (\tilde{z}_{1 j_{t_1}}^2-\sigma_{j_{t_1}}^2)}
           \bigp{\prod\limits_{t_3\in\ct_3} \tilde{z}_{1 j_{t_3}}}} \\
  & & \hspace{4em} \times
      \bbe{\bigp{\prod\limits_{s_2\in\cs_2} (\tilde{z}_{2 j_{s_2}}^2-\sigma_{j_{s_2}}^2)}
           \bigp{\prod\limits_{s_3\in\cs_3} \tilde{z}_{2 j_{s_2}}}
           \bigp{\prod\limits_{v_1\in\cv_1} (\tilde{z}_{2 j_{v_1}}^2-\sigma_{j_{v_1}}^2)}
           \bigp{\prod\limits_{v_3\in\cv_3} \tilde{z}_{2 j_{v_3}}}} \\
  & & \hspace{4em} \times
      \bbe{\bigp{\prod\limits_{t_2\in\ct_2} (\tilde{z}_{3 j_{t_2}}^2-\sigma_{j_{t_2}}^2)}
           \bigp{\prod\limits_{t_3\in\ct_3} \tilde{z}_{3 j_{t_3}}}} \\
  & & \hspace{4em} \times
      \bbe{\bigp{\prod\limits_{v_2\in\cv_2} (\tilde{z}_{4 j_{v_2}}^2-\sigma_{j_{v_2}}^2)}
           \bigp{\prod\limits_{v_3\in\cv_3} \tilde{z}_{4 j_{v_4}}}} \\
  &=& \sum\limits_{(\cs_1,\cs_2,\cs_3,\ct_1,\ct_2,\ct_3,\cv_1,\cv_2,\cv_3)}
      (-2)^{|\cs_3|+|\ct_3|+|\cv_3|} 
      \sum\limits_{\substack{1\le j_{s_3}\le p \\ s_3\in\cs_3}} \ 
      \sum\limits_{\substack{1\le j_{t_3}\le p \\ t_3\in\ct_3}} \ 
      \sum\limits_{\substack{1\le j_{v_3}\le p \\ v_3\in\cv_3}} \ 
      \\
  & & \hspace{2em} \times
      \sum\limits_{\pi_1\in\cp_1} 
      \sum\limits_{\substack{1\le j_{s_1}\le p \\ s_1\in\cs_1}} \ 
      \sum\limits_{\substack{1\le j_{t_1}\le p \\ t_1\in\ct_1}} \ 
      \prod\limits_{\cb_1\in\pi_1}
      \cum\bigp{\tilde{z}_{1 j_s}^2-\sigma_{j_s}^2, \tilde{z}_{1 j_t}: s\in\cb_1\cap(\cs_1\cup\ct_1), t\in\cb_1\cap(\cs_3\cup\ct_3)} \\
  & & \hspace{2em} \times
      \sum\limits_{\pi_2\in\cp_2}       
      \sum\limits_{\substack{1\le j_{s_2}\le p \\ s_2\in\cs_2}} \ 
      \sum\limits_{\substack{1\le j_{v_1}\le p \\ v_1\in\cv_1}} \ 
      \prod\limits_{\cb_2\in\pi_2} 
      \cum\bigp{\tilde{z}_{2 j_s}^2-\sigma_{j_s}^2, \tilde{z}_{2 j_t}: s\in\cb_2\cap(\cs_2\cup\cv_1), t\in\cb_2\cap(\cs_3\cup\cv_3)} \\
  & & \hspace{2em} \times
      \sum\limits_{\pi_3\in\cp_3}
      \sum\limits_{\substack{1\le j_{t_2}\le p \\ t_2\in\ct_2}} \ 
      \prod\limits_{\cb_3\in\pi_3} 
      \cum\bigp{\tilde{z}_{3 j_s}^2-\sigma_{j_s}^2, \tilde{z}_{3 j_t}: s\in\cb_3\cap\ct_2, t\in\cb_3\cap\ct_3)} \\
  & & \hspace{2em} \times
      \sum\limits_{\pi_4\in\cp_4}
      \sum\limits_{\substack{1\le j_{v_2}\le p \\ v_2\in\cv_2}} \ 
      \prod\limits_{\cb_4\in\pi_4}
      \cum\bigp{\tilde{z}_{4 j_s}^2-\sigma_{j_s}^2, \tilde{z}_{4 j_t}: s\in\cb_4\cap\cv_2, t\in\cb_4\cap\cv_3)} \\
  &=& \sum\limits_{(\cs_1,\cs_2,\cs_3,\ct_1,\ct_2,\ct_3,\cv_1,\cv_2,\cv_3)}
      (-2)^{|\cs_3|+|\ct_3|+|\cv_3|} 
      \sum\limits_{\pi_1\in\cp_1} 
      \sum\limits_{\pi_2\in\cp_2}   
      \sum\limits_{\pi_3\in\cp_3}
      \sum\limits_{\pi_4\in\cp_4}
      \sum\limits_{\substack{1\le j_{s_3}\le p \\ s_3\in\cs_3}} \ 
      \sum\limits_{\substack{1\le j_{t_3}\le p \\ t_3\in\ct_3}} \ 
      \sum\limits_{\substack{1\le j_{v_3}\le p \\ v_3\in\cv_3}} \ 
      \\
  & & \hspace{2em} \times
      \prod\limits_{\cb_1\in\pi_1}
      \sum\limits_{\substack{1\le j_{u_1}\le p \\ u_1\in\cb_1\cap(\cs_1\cup\ct_1)}}
      \cum\bigp{\tilde{z}_{1 j_{u_1}}^2-\sigma_{j_{u_1}}^2, \tilde{z}_{1 j_{w_1}}: u_1\in\cb_1\cap(\cs_1\cup\ct_1), w_1\in\cb_1\cap(\cs_3\cup\ct_3)} \\
  & & \hspace{2em} \times
      \prod\limits_{\cb_2\in\pi_2}
      \sum\limits_{\substack{1\le j_{u_2}\le p \\ u_2\in\cb_2\cap(\cs_2\cup\cv_1)}}
      \cum\bigp{\tilde{z}_{2 j_{u_2}}^2-\sigma_{j_{u_2}}^2, \tilde{z}_{2 j_{w_2}}: u_2\in\cb_2\cap(\cs_2\cup\cv_1), w_2\in\cb_2\cap(\cs_3\cup\cv_3)} \\
  & & \hspace{2em} \times
      \prod\limits_{\cb_3\in\pi_3} 
      \sum\limits_{\substack{1\le j_{u_3}\le p \\ u_3\in\cb_3\cap\ct_2}}
      \cum\bigp{\tilde{z}_{3 j_{u_3}}^2-\sigma_{j_{u_3}}^2, \tilde{z}_{3 j_{w_3}}: u_3\in\cb_3\cap\ct_2, w_3\in\cb_3\cap\ct_3)} \\
  & & \hspace{2em} \times
      \prod\limits_{\cb_4\in\pi_4}
      \sum\limits_{\substack{1\le j_{u_4}\le p \\ u_4\in\cb_4\cap\cv_2}}
      \cum\bigp{\tilde{z}_{4 j_{u_4}}^2-\sigma_{j_{u_4}}^2, \tilde{z}_{4 j_{w_4}}: u_4\in\cb_4\cap\cv_2, w_4\in\cb_4\cap\cv_3)}
  \EEqn
  
  It follows that
  \BEqn
  & & \lrabs{D_{i_1 i_2 i_3}} \\
  &\le& 2^{i_1+i_2+i_3}
        \sum\limits_{(\cs_1,\cs_2,\cs_3,\ct_1,\ct_2,\ct_3,\cv_1,\cv_2,\cv_3)}
        \sum\limits_{\pi_1\in\cp_1} 
        \sum\limits_{\pi_2\in\cp_2}   
        \sum\limits_{\pi_3\in\cp_3}
        \sum\limits_{\pi_4\in\cp_4}
        \sum\limits_{\substack{1\le j_{s_3}\le p \\ s_3\in\cs_3}} \
        \sum\limits_{\substack{1\le j_{t_3}\le p \\ t_3\in\ct_3}} \ 
        \sum\limits_{\substack{1\le j_{v_3}\le p \\ v_3\in\cv_3}} \\ 
  & & \hspace{1em} \times
      \prod\limits_{\cb_1\in\pi_1}
      \sum\limits_{\substack{1\le j_{u_1}\le p \\ u_1\in\cb_1\cap(\cs_1\cup\ct_1)}}
      \lrabs{\cum\bigp{\tilde{z}_{1 j_{u_1}}^2-\sigma_{j_{u_1}}^2, \tilde{z}_{1 j_{w_1}}: u_1\in\cb_1\cap(\cs_1\cup\ct_1), w_1\in\cb_1\cap(\cs_3\cup\ct_3)}} \\
  & & \hspace{1em} \times
      \prod\limits_{\cb_2\in\pi_2}
      \sum\limits_{\substack{1\le j_{u_2}\le p \\ u_2\in\cb_2\cap(\cs_2\cup\cv_1)}}
      \lrabs{\cum\bigp{\tilde{z}_{2 j_{u_2}}^2-\sigma_{j_{u_2}}^2, \tilde{z}_{2 j_{w_2}}: u_2\in\cb_2\cap(\cs_2\cup\cv_1), w_2\in\cb_2\cap(\cs_3\cup\cv_3)}} \\
  & & \hspace{1em} \times
      \prod\limits_{\cb_3\in\pi_3} 
      \sum\limits_{\substack{1\le j_{u_3}\le p \\ u_3\in\cb_3\cap\ct_2}}
      \lrabs{\cum\bigp{\tilde{z}_{3 j_{u_3}}^2-\sigma_{j_{u_3}}^2, \tilde{z}_{3 j_{w_3}}: u_3\in\cb_3\cap\ct_2, w_3\in\cb_3\cap\ct_3)}} \\
  & & \hspace{1em} \times
      \prod\limits_{\cb_4\in\pi_4}
      \sum\limits_{\substack{1\le j_{u_4}\le p \\ u_4\in\cb_4\cap\cv_2}}
      \lrabs{\cum\bigp{\tilde{z}_{4 j_{u_4}}^2-\sigma_{j_{u_4}}^2, \tilde{z}_{4 j_{w_4}}: u_4\in\cb_4\cap\cv_2, w_4\in\cb_4\cap\cv_3)}}
  \EEqn
  
  Under Assumption \ref{Assumpt:component-dept}, it follows from Lemma \ref{Lemma:component-6} that, for any $\pi$, $\cb\in\pi$ and any fixed $(j_{w_1}: w_1\in\cb_1\cap(\cs_3\cup\ct_3))$ that
  \BEqn
  & & \sum\limits_{\substack{1\le j_{u_1}\le p \\ u_1\in\cb_1\cap(\cs_1\cup\ct_1)}}
      \lrabs{\cum\bigp{\tilde{z}_{1 j_{u_1}}^2-\sigma_{j_{u_1}}^2, \tilde{z}_{1 j_{w_1}}: u_1\in\cb_1\cap(\cs_1\cup\ct_1), w_1\in\cb_1\cap(\cs_3\cup\ct_3)}} \\
  &\le& C(U^{\ast}) 
        \# \lrcp{\begin{array}{l}
                    (j_{u_1}: u_1\in\cb_1\cap(\cs_1\cup\ct_1)): \vspace{1mm}\\
                    \qquad \cum\bigp{\tilde{z}_{1 j_{u_1}}^2-\sigma_{j_{u_1}}^2, \tilde{z}_{1 j_{w_1}}: u_1\in\cb_1\cap(\cs_1\cup\ct_1), w_1\in\cb_1\cap(\cs_3\cup\ct_3)} \neq 0 \vspace{2mm}\\
                    \qquad 1 \le j_{u_1} \le p, u_1\in\cb_1\cap(\cs_1\cup\ct_1)
                 \end{array}} \\
  &\le& C(i_1,i_2,i_3,U^{\ast}) (\alpha(p))^{|\cb_1\cap\cs_1|+|\cb_1\cap\ct_1|} \\
  & & \hspace{4em} \times
        \bone\{\max\limits_{u\in\cb_1\cap(\cs_3\cup\ct_3)}\{j_u\}
               - \min\limits_{u\in\cb_1\cap(\cs_3\cup\ct_3)}\{j_u\}
               \le \bigp{\lrabs{\cb_1\cap(\cs_1\cup\cs_3\cup\ct_1\cup\ct_3)} - 1}\alpha(p)\} \\
  &\le& C(i_1,i_2,i_3,U^{\ast}) (\alpha(p))^{|\cb_1\cap\cs_1|+|\cb_1\cap\ct_1|} \\
  & & \hspace{4em} \times
        \bone\{\max\limits_{u\in\cb_1\cap(\cs_3\cup\ct_3)}\{j_u\}
               - \min\limits_{u\in\cb_1\cap(\cs_3\cup\ct_3)}\{j_u\}
               \le (i_1+i_2+i_3)\alpha(p)\}
  \EEqn
  which implies that
  \BEqn
  & & \prod\limits_{\cb_1\in\pi_1} \sum\limits_{\substack{1\le j_{u_1}\le p \\ u_1\in\cb_1\cap(\cs_1\cup\ct_1)}}
      \lrabs{\cum\bigp{\tilde{z}_{1 j_{u_1}}^2-\sigma_{j_{u_1}}^2, \tilde{z}_{1 j_{w_1}}: u_1\in\cb_1\cap(\cs_1\cup\ct_1), w_1\in\cb_1\cap(\cs_3\cup\ct_3)}} \\
  &=& C(i_1,i_2,i_3,U^{\ast}) (\alpha(p))^{|\cs_1| + |\ct_1|} \\[2mm]
  & & \hspace{4em} \times
        \prod\limits_{\cb_1\in\pi_1} 
        \bone\{\max\limits_{u\in\cb_1\cap(\cs_3\cup\ct_3)}\{j_u\}
               - \min\limits_{u\in\cb_1\cap(\cs_3\cup\ct_3)}\{j_u\}
               \le (i_1+i_2+i_3)\alpha(p)\} \\
  \EEqn
  where equality follows from the fact that $\sum\limits_{\cb_1\in\pi_1} \bigp{|\cb_1\cap\cs_1|+|\cb_1\cap\ct_1|} = |\cs_1| + |\ct_1|$. Similarly, we can find the corresponding upper bounds for the other terms, which jointly lead to the following result:
  \BEqn
  & & |D_{i_1 i_2 i_3}| \\
  &\le& C(i_1,i_2,i_3,U^{\ast})
        \sum\limits_{(\cs_1,\cs_2,\cs_3,\ct_1,\ct_2,\ct_3,\cv_1,\cv_2,\cv_3)}
        (\alpha(p))^{|\cs_1| + |\cs_2| + |\ct_1| + |\ct_2| + |\cv_1| + |\cv_2|} \\
  & & \hspace{1em} \times
      \sum\limits_{\pi_1\in\cp_1} 
      \sum\limits_{\pi_2\in\cp_2}   
      \sum\limits_{\pi_3\in\cp_3}
      \sum\limits_{\pi_4\in\cp_4} 
      \sum\limits_{\substack{1\le j_{s_3}\le p \\ s_3\in\cs_3}} \
      \sum\limits_{\substack{1\le j_{t_3}\le p \\ t_3\in\ct_3}} \ 
      \sum\limits_{\substack{1\le j_{v_3}\le p \\ v_3\in\cv_3}} \\
  & & \hspace{2em} \times     
      \left(
      \prod\limits_{\cb_1\in\pi_1} 
      \bone\{\max\limits_{u\in\cb_1\cap(\cs_3\cup\ct_3)}\{j_u\}
               - \min\limits_{u\in\cb_1\cap(\cs_3\cup\ct_3)}\{j_u\}
               \le (i_1+i_2+i_3)\alpha(p)\}
      \right.  \\
  & & \hspace{5em} \times
      \left.
      \prod\limits_{\cb_4\in\pi_4}
      \bone\{\max\limits_{u\in\cb_4\cap\cv_3}\{j_u\}
               - \min\limits_{u\in\cb_4\cap\cv_3}\{j_u\}
               \le (i_1+i_2+i_3)\alpha(p)\} 
      \right) \\
  & & \hspace{2em} \times
      \left(
      \prod\limits_{\cb_2\in\pi_2}
      \bone\{\max\limits_{u\in\cb_2\cap(\cs_3\cup\cv_3)}\{j_u\}
               - \min\limits_{u\in\cb_2\cap(\cs_3\cup\cv_3)}\{j_u\}
               \le (i_1+i_2+i_3)\alpha(p)\}
      \right.  \\
  & & \hspace{5em} \times
      \left.   
      \prod\limits_{\cb_3\in\pi_3} 
      \bone\{\max\limits_{u\in\cb_3\cap\ct_3}\{j_u\}
               - \min\limits_{u\in\cb_3\cap\ct_3}\{j_u\}
               \le (i_1+i_2+i_3)\alpha(p)\} 
      \right) \\
  &\le& C(i_1,i_2,i_3,U^{\ast})
        \sum\limits_{(\cs_1,\cs_2,\cs_3,\ct_1,\ct_2,\ct_3,\cv_1,\cv_2,\cv_3)}
        (\alpha(p))^{|\cs_1| + |\cs_2| + |\ct_1| + |\ct_2| + |\cv_1| + |\cv_2|} \\
  & & \times
      \sum\limits_{\pi_1\in\cp_1} 
      \sum\limits_{\pi_2\in\cp_2}   
      \sum\limits_{\pi_3\in\cp_3}
      \sum\limits_{\pi_4\in\cp_4} \\
  & & \times\min\lrcp{
      \begin{array}{l}
         \# \lrcp{\begin{array}{ll}
                     (j_u: u\in\cs_3\cup\ct_3\cup\cv_3): & \vspace{1mm}\\   
                     \qquad \max\limits_{u\in\cb_1\cap(\cs_3\cup\ct_3)}\{j_u\} - \min\limits_{u\in\cb_1\cap(\cs_3\cup\ct_3)}\{j_u\} \le (i_1+i_2+i_3)\alpha(p) & \forall \cb_1\in\pi_1 \vspace{1mm}\\
                     \qquad \max\limits_{u\in\cb_4\cap\cv_3}\{j_u\} - \min\limits_{u\in\cb_4\cap\cv_3}\{j_u\} \le (i_1+i_2+i_3)\alpha(p) & \forall \cb_4\in\pi_4
                  \end{array}} \vspace{2mm}\\
        \# \lrcp{\begin{array}{ll}
                     (j_u: u\in\cs_3\cup\ct_3\cup\cv_3): & \vspace{1mm}\\   
                     \qquad \max\limits_{u\in\cb_2\cap(\cs_3\cup\cv_3)}\{j_u\} - \min\limits_{u\in\cb_2\cap(\cs_3\cup\cv_3)}\{j_u\} \le (i_1+i_2+i_3)\alpha(p) & \forall \cb_2\in\pi_2 \vspace{1mm}\\
                     \qquad \max\limits_{u\in\cb_3\cap\ct_3}\{j_u\} - \min\limits_{u\in\cb_3\cap\ct_3}\{j_u\} \le (i_1+i_2+i_3)\alpha(p) & \forall \cb_3\in\pi_3
                  \end{array}}                  
      \end{array}}
  \EEqn
  
  Note that
  \BEqn
  & & \# \lrcp{\begin{array}{ll}
                     (j_u: u\in\cs_3\cup\ct_3\cup\cv_3): & \vspace{1mm}\\   
                     \qquad \max\limits_{u\in\cb_1\cap(\cs_3\cup\ct_3)}\{j_u\} - \min\limits_{u\in\cb_1\cap(\cs_3\cup\ct_3)}\{j_u\} \le (i_1+i_2+i_3)\alpha(p) & \forall \cb_1\in\pi_1 \vspace{1mm}\\
                     \qquad \max\limits_{u\in\cb_4\cap\cv_3}\{j_u\} - \min\limits_{u\in\cb_4\cap\cv_3}\{j_u\} \le (i_1+i_2+i_3)\alpha(p) & \forall \cb_4\in\pi_4
                  \end{array}} \\
  &\le& \prod\limits_{\cb_1\in\pi_1} 
        \# \lrcp{(j_u: u\in\cb_1\cap(\cs_3\cup\ct_3)):
                 \max\limits_{u\in\cb_1\cap(\cs_3\cup\ct_3)}\{j_u\} - \min\limits_{u\in\cb_1\cap(\cs_3\cup\ct_3)}\{j_u\} \le (i_1+i_2+i_3)\alpha(p)} \\
  & & \hspace{2em} \times
      \prod\limits_{\cb_4\in\pi_4}
      \# \lrcp{(j_u: u\in\cb_4\cap\cv_3): 
               \max\limits_{u\in\cb_4\cap\cv_3}\{j_u\} - \min\limits_{u\in\cb_4\cap\cv_3}\{j_u\} \le (i_1+i_2+i_3)\alpha(p)} \\
  &\le& C(i_1,i_2,i_3) 
        \lrp{\prod\limits_{\cb_1\in\pi_1} \bigp{(\alpha(p))^{|\cb_1\cap(\cs_3\cup\ct_3)|-1} p}}
        \lrp{\prod\limits_{\cb_4\in\pi_4} \bigp{(\alpha(p))^{|\cb_4\cap\cv_3|-1} p}} \\
  &=& C(i_1,i_2,i_3) 
      (\alpha(p))^{|\cs_3|+|\ct_3|+|\cv_3|} \lrp{\frac{p}{\alpha(p)}}^{|\pi_1|+|\pi_4|}
  \EEqn
  
  Similarly, we also have
  \BEqn
  & & \# \lrcp{\begin{array}{ll}
                     (j_u: u\in\cs_3\cup\ct_3\cup\cv_3): & \vspace{1mm}\\   
                     \qquad \max\limits_{u\in\cb_2\cap(\cs_3\cup\cv`_3)}\{j_u\} - \min\limits_{u\in\cb_2\cap(\cs_3\cup\cv_3)}\{j_u\} \le (i_1+i_2+i_3)\alpha(p) & \forall \cb_2\in\pi_2 \vspace{1mm}\\
                     \qquad \max\limits_{u\in\cb_3\cap\ct_3}\{j_u\} - \min\limits_{u\in\cb_3\cap\ct_3}\{j_u\} \le (i_1+i_2+i_3)\alpha(p) & \forall \cb_3\in\pi_3
                  \end{array}} \\
  &\le& C(i_1,i_2,i_3) 
      (\alpha(p))^{|\cs_3|+|\ct_3|+|\cv_3|} \lrp{\frac{p}{\alpha(p)}}^{|\pi_2|+|\pi_3|}
  \EEqn
  and it follows that
  \BEqn
  & & \lrabs{D_{i_1,i_2,i_3}} \\
  &\le& C(i_1,i_2,i_3,U^{\ast})
        \sum\limits_{(\cs_1,\cs_2,\cs_3,\ct_1,\ct_2,\ct_3,\cv_1,\cv_2,\cv_3)}
        (\alpha(p))^{|\cs_1| + |\cs_2| + |\cs_3| + |\ct_1| + |\ct_2| + |\ct_3| + |\cv_1| + |\cv_2| + |\cv_3|} \\
  & & \times
      \sum\limits_{\pi_1\in\cp_1} 
      \sum\limits_{\pi_2\in\cp_2}   
      \sum\limits_{\pi_3\in\cp_3}
      \sum\limits_{\pi_4\in\cp_4} 
      \lrp{\frac{p}{\alpha(p)}}^{\min\{|\pi_1|+|\pi_4|,|\pi_2|+|\pi_3|\}} \\
  &=& C(i_1,i_2,i_3,U^{\ast}) (\alpha(p))^{i_1+i_2+i_3}
      \sum\limits_{\substack{(\cs_1,\cs_2,\cs_3) \\ (\ct_1,\ct_2,\ct_3) \\ (\cv_1,\cv_2,\cv_3)}}
      \sum\limits_{\pi_1\in\cp_1} 
      \sum\limits_{\pi_2\in\cp_2}   
      \sum\limits_{\pi_3\in\cp_3}
      \sum\limits_{\pi_4\in\cp_4} 
      \lrp{\frac{p}{\alpha(p)}}^{\min\{|\pi_1|+|\pi_4|,|\pi_2|+|\pi_3|\}} \\
  \EEqn
  
  It remains to investigate the upper bound of $\min\{|\pi_1|+|\pi_4|,|\pi_2|+|\pi_3|\}$ over all possible $\cp_1,\cp_2,\cp_3,\cp_4$. Note that
  \begin{equation*}
      \begin{array}{cccc}
         |\pi_1|\le\lrfloor{\bigp{|\cs_1|+|\cs_3|+|\ct_1|+|\ct_3|}/2}, & & &  |\pi_2|\le\lrfloor{\bigp{|\cs_2|+|\cs_3|+|\cv_1|+|\cv_3|}/2}, \\
         |\pi_3|\le\lrfloor{\bigp{|\ct_2|+|\ct_3|}/2}, & & & 
         |\pi_4|\le\lrfloor{\bigp{|\cv_2|+|\cv_3|}/2} 
      \end{array}
  \end{equation*}
  then we have
  \begin{equation*}
      \min\{|\pi_1|+|\pi_4|,|\pi_2|+|\pi_3|\} \le \lrfloor{(i_1+i_2+i_3)/2},
  \end{equation*}
  which completes the proof of the statement regarding $D_{i_1,i_2,i_3}$. Note that all the other statements can be proved using similar techniques, we thus spare the details proofs for the others.
\end{proof}

For $i\geq1$, define $\bar{L}_i = \bbe{\lrabs{|Z_1-Z_2|^2-A}^i}$, and
\BEqn
& & \bar{L}_i^{XY} = \bbe{||X_1-Y_1|-A^{XY}|^i}, \\
& & \bar{L}_i^{X} = \bbe{||X_1-X_2|-A^X|^i}, \\
& & \bar{L}_i^{Y} = \bbe{||Y_1-Y_2|-A^Y|^i}.
\EEqn

\begin{lemma}\label{Lemma:Li_bar-1}
For any $1\le i\le 64$, it holds under Assumption \ref{Assumpt:component-dept}\ref{Assumpt:component-dept-1}-\ref{Assumpt:component-dept-2} that 
\begin{equation*}
\begin{array}{l}
  \bar{L}_i^{XY} \le C(i,U^{\ast}) (\alpha(p))^{i/2} p^{i/2}, \vspace{2mm}\\
  \bar{L}_i^{X} \le C(i,U^{\ast}) (\alpha(p))^{i/2} p^{i/2}, \vspace{2mm}\\
  \bar{L}_i^{Y} \le C(i,U^{\ast}) (\alpha(p))^{i/2} p^{i/2}.
\end{array}
\end{equation*}

Moreover, it holds under Assumption \ref{Assumpt:component-dept}\ref{Assumpt:component-dept-1}-\ref{Assumpt:component-dept-2} and the null hypothesis that
\begin{equation*}
    \bar{L}_i \le C(i,U^{\ast}) (\alpha(p))^{i/2} p^{i/2}
\end{equation*}
\end{lemma}
\begin{proof}
  When $i$ is even, then it directly follows from Lemma \ref{Lemma:Li-1} that
  \begin{equation*}
      \bar{L}_i^{XY} 
    = L_i^{XY} 
    \le C(i,U^{\ast}) (\alpha(p))^{\lrfloor{i/2}} p^{\lrceil{i/2}} 
    = C(i,U^{\ast}) (\alpha(p))^{i/2} p^{i/2}
  \end{equation*}
  When $i$ is odd, then $i-1$ and $i+1$ are both even and by H\"older's inequality, we obtain that
  \BEqn
      \bar{L}_i^{XY} 
  &=& \bbe{\lrabs{|X_1-Y_1|^2-A^{XY}}^{(i-1)/2} \lrabs{|X_1-Y_1|^2-A^{XY}}^{(i+1)/2}} \\
  &\le& (\bar{L}_{i-1}^{XY})^{1/2} (\bar{L}_{i+1}^{XY})^{1/2} \\
  &\le& C(i,U^{\ast}) \bigp{(\alpha(p))^{(i-1)/2} p^{(i-1)/2}} \bigp{(\alpha(p))^{(i+1)/2} p^{(i+1)/2}} \\
  &=& C(i,U^{\ast}) (\alpha(p))^{i/2} p^{i/2}
  \EEqn
  Using similar arguments, we can also verify the statements of $\bar{L}_i^{X},\bar{L}_i^{Y}$ and $\bar{L}_i$.
\end{proof}

\subsection{Part 2: Approximation of the Kernel Function}

\begin{lemma}\label{Lemma:approx-1}
    For $x\geq-1$, it holds that
    \begin{equation*}
    \begin{array}{l}
       \sqrt{1+x} 
      = 1 + c_1(x) x, \vspace{2mm}\\
       \sqrt{1+x} 
      = 1 + \frac{1}{2}x + c_2(x) x^2, \vspace{2mm}\\
       \sqrt{1+x} 
      = 1 + \frac{1}{2}x - \frac{1}{8}x^2 + c_3(x) x^3, \vspace{2mm}\\
       \sqrt{1+x} 
      = 1 + \frac{1}{2}x - \frac{1}{8}x^2 + \frac{1}{16}x^3 - c_4(x) x^4, \vspace{2mm}\\
       \sqrt{1+x} 
      = 1 + \frac{1}{2}x - \frac{1}{8}x^2 + \frac{1}{16}x^3 - \frac{5}{128}x^4 + c_5(x) x^5, \vspace{2mm}\\
       \sqrt{1+x} 
      = 1 + \frac{1}{2}x - \frac{1}{8}x^2 + \frac{1}{16}x^3 - \frac{5}{128}x^4 + \frac{7}{256} x^5 + c_6(x) x^6, \vspace{2mm}\\
       \sqrt{1+x} 
      = 1 + \frac{1}{2}x - \frac{1}{8}x^2 + \frac{1}{16}x^3 - \frac{5}{128}x^4 + \frac{7}{256} x^5 - \frac{21}{1024}x^6 + c_7(x) x^7.    
    \end{array}
    \end{equation*}
    where $c_i(x), i=1,\dots,7$ are some constants that depend only on $x$ and satisfy that $\sup\limits_{x\geq-1}|c_i(x)| \le C$ for some positive constant $C$.
\end{lemma}
\begin{proof}
    We take the fifth order approximation as the example and show that $\sqrt{1+x} 
      = 1 + \frac{1}{2}x - \frac{1}{8}x^2 + \frac{1}{16}x^3 - \frac{5}{128}x^4 + c_5(x) x^5$ holds with $|c_5(x)| \le \frac{128}{35}$. 
    
    For $x\geq-1$ and $x\neq0$, define
    \begin{equation*}
        c_5(x)
      = \frac{\sqrt{1+x} - \lrp{1 + x/2 - x^2/8 + x^3/16 - 5x^4/128}}{x^5}
    \end{equation*}
    and for $x=0$, define $c_5(x)=\frac{7}{256}$, then it's equivalent to calculate $\sup\limits_{x\geq-1}|c_5(x)|$.
    
    By taking the derivative of $c_5(x)$, we obtain that $c_5'(x) = h(x)/x^6$, where
    \begin{equation*}
        h(x)
      = \frac{1}{2}x(1+x)^{-1/2} - 5(x+1)^{1/2} + 5 + 2x - \frac{3}{8}x^2 + \frac{1}{8}x^3 - \frac{5}{128}x^4.
    \end{equation*}
    We further calculate the derivatives of $h$ up to the fifth order, and have
    \begin{align*}
        h'(x)
     =& -\frac{1}{4}x(1+x)^{-3/2} - 2(1+x)^{-1/2} + 2 - \frac{3}{4}x + \frac{3}{8}x^2 - \frac{5}{32}x^3, \\
        h''(x)
     =& \frac{3}{8}x(1+x)^{-5/2} + \frac{3}{4}(1+x)^{-3/2} - \frac{3}{4} + \frac{3}{4}x - \frac{15}{32}x^2, \\
        h'''(x)
     =& -\frac{15}{16}x(1+x)^{-7/2} - \frac{3}{4}(1+x)^{-5/2} + \frac{3}{4} - \frac{15}{16}x, \\
        h^{(4)}(x)
     =& \frac{105}{32}x(1+x)^{-9/2} + \frac{15}{16}(1+x)^{-7/2} - \frac{15}{16}, \\
        h^{(5)}(x)
     =& -\frac{945}{64}x(1+x)^{-11/2}.
    \end{align*}
    Note that $h^{(5)}(x) \geq 0$ for $-1\le x\le 0$ and $h^{(5)}(x) < 0$ for $x>0$, then $h^{(4)}(x)$ is increasing in $x$ for $-1\le x\le 0$ and decreasing for $x>0$, which implies that $h^{(4)}(x) \le 0$ for $x\geq-1$. Repeating this process and we can finally deduce that $h(x)$ is increasing in $x$ for $-1\le x\le 0$ and decreasing for $x>0$. This further indicates that $h(x)\le0$ for $x\geq-1$, and so is $c_5'(X)$, hence $c_5(x)$ is decreasing in $x$ for $x\geq-1$.
    
    Note that $c_5(-1) = \frac{35}{128}$, $c_5(0) = \frac{7}{256}$ and $c_5(\infty)=0$, then we have
    \BEqn
        \frac{7}{256}
    \le &\displaystyle{\frac{\sqrt{1+x} - \lrp{1 + x/2 - x^2/8 + x^3/16 - 5x^4/128}}{x^5}}& 
    \le \frac{35}{128}, 
    \quad\quad -1\le x\le 0, \vspace{3pt} \\
        0
    \le &\displaystyle{\frac{\sqrt{1+x} - \lrp{1 + x/2 - x^2/8 + x^3/16 - 5x^4/128}}{x^5}}& 
    \le \frac{7}{256}, 
    \quad\quad x>0,
    \EEqn
    and consequently,
    \BEqn
        \frac{35}{128}x^5
    \le &\sqrt{1+x} - \lrp{1 + \frac{1}{2}x - \frac{1}{8}x^2 + \frac{1}{16}x^3 - \frac{5}{128}x^4}&    
    \le \frac{7}{256}x^5,
    \quad\quad -1\le x\le 0, \vspace{3pt} \\
        0
    \le &\sqrt{1+x} - \lrp{1 + \frac{1}{2}x - \frac{1}{8}x^2 + \frac{1}{16}x^3 - \frac{5}{128}x^4}&    
    \le \frac{7}{256}x^5,
    \quad\quad x > 0,
    \EEqn
    which implies that $|c_5(x)| \le \frac{35}{128}$. All the other statements can be proved using the same technique, for which we spare the details.
\end{proof}

\begin{lemma}\label{Lemma:approx-2}
Consider $f: [0,\infty) \mapsto \br$ and let $D_0 \subseteq D := [0,\infty)$ be a subset of the domain of $f$. Assume that
\begin{enumerate}[label=(\roman*)]
    \item for any $s\in D$ and $s_0\in D_0$, it holds that
    \begin{equation*}
        f(s) = \sum\limits_{i=0}^{6} \frac{1}{i!} f_i(s_0) (s-s_0)^i + f_7(\xi(s,s_0)) (s-s_0)^7,
    \end{equation*}
    where $f_i$ denotes the $i$-th order derivative of $f$ and $\xi(s,s_0)$ is some value between 
    $s$ and $s_0$.
    
    \item there exists a positive constant $\tilde{M} < \infty$, such that for any $s_0 \in D_0$, it holds that
    \begin{equation*}
        \max\limits_{1\le i\le 7} \sup\limits_{s\in D} |f_i(s)|\cdot |s_0^i| \le \tilde{M} |f_0(s_0)|.
    \end{equation*}
\end{enumerate}
Then for any $s\in D$ and $s_0 \in D_0$, it holds that 
\begin{equation*}
    f(s) = \sum\limits_{i=0}^{6} c_i(s_0) \lrp{\frac{s^2-s_0^2}{s_0^2}}^i + R(s,s_0),
\end{equation*}
where $c_i$'s are defined as follows
\BEqn
& & c_0(x) = f(x) \\
& & c_1(x) = \frac{1}{2} f_1(x) x \\
& & c_2(x) = -\frac{1}{8} f_1(x) x + \frac{1}{8} f_2(x) x^2 \\
& & c_3(x) = \frac{1}{16} f_1(x) x - \frac{1}{16} f_2(x) x^2 + \frac{1}{48} f_3(x) x^3 \\
& & c_4(x) = -\frac{5}{128} f_1(x) x + \frac{5}{128} f_2(x) x^2 - \frac{1}{64} f_3(x) x^3 + \frac{1}{384} f_4(x) x^4 \\
& & c_5(x) = \frac{7}{256} f_1(x) x - \frac{7}{256} f_2(x) x^2 + \frac{3}{256} f_3(x) x^3 - \frac{1}{96} f_4(x) x^4 + \frac{1}{3840} f_5(x) x^5 \\
& & c_6(x) = -\frac{21}{1024} f_1(x) x + \frac{21}{1024} f_2(x) x^2 - \frac{7}{768} f_3(x) x^3 + \frac{7}{3072} f_4(x) x^4 \\
& & \hspace{4em}
    - \frac{1}{3072} f_5(x) x^5 + \frac{1}{46080} f_6(x) x^6
\EEqn
and $|R(s,s_0)| \le C(\tilde{M}) \sum\limits_{i=7}^{16} \lrabs{\frac{s^2-s_0^2}{s_0^2}}^i |f_0(s_0)|$.
\end{lemma}
\begin{proof}
    Note that for any $s,s_0\geq0$, it holds that $s = s_0 \sqrt{1 + \frac{s^2-s_0^2}{s_0^2}}$, then by applying Lemma \ref{Lemma:approx-1}, we obtain a third-order approximation of $s-s_0$, denoted by $\ca_3(s,s_0)$, that is,
    \begin{equation*}
        s - s_0
      = \ca_3(s,s_0)
      = s_0 \lrp{\frac{1}{2}\lrp{\frac{s^2-s_0^2}{s_0^2}} - \frac{1}{8}\lrp{\frac{s^2-s_0^2}{s_0^2}}^2 + e_3(s,s_0)\lrp{\frac{s^2-s_0^2}{s_0^2}}^3}
    \end{equation*}
    where $|e_3(s,s_0)|\le C$ with a universal constant $C$. Let $\delta_{s,s_0} = \frac{s^2-s_0^2}{s_0^2}$, then we have
    \BEqn
    & & \lrabs{\ca_3^5(s,s_0)
        - s_0^5 
          \lrp{  \frac{1}{32} \delta_{s,s_0}^5
               - \frac{5}{128} \delta_{s,s_0}^6}}
        \le C s_0^5 \sum\limits_{i=7}^{15} \lrabs{\delta_{s,s_0}}^i
    \EEqn
    for some constant $C>0$.

    Similarly, for $i=1,\dots,7$, define $\ca_i(s,s_0)$ as the $i$-th order approximation of $s-s_0$. After some simple calculations, with , we have 
    \BEqn
    & & \lrabs{\ca_1^7(s,s_0)}
        \le C s_0^7 \lrabs{\delta_{s,s_0}}^7 \\
    & & \lrabs{\ca_2^6(s,s_0)
        - s_0^6 
          \lrp{  \frac{1}{64} \delta_{s,s_0}^6}}
        \le C s_0^6 \sum\limits_{i=7}^{12} \lrabs{\delta_{s,s_0}}^i \\
    & & \lrabs{\ca_4^4(s,s_0)
        - s_0^4 
          \lrp{  \frac{1}{16} \delta_{s,s_0}^4
               - \frac{1}{4} \delta_{s,s_0}^5
               + \frac{7}{128} \delta_{s,s_0}^6}}
        \le C s_0^4 \sum\limits_{i=7}^{16} \lrabs{\delta_{s,s_0}}^i \\       
    & & \lrabs{\ca_5^3(s,s_0)
        - s_0^3 
          \lrp{  \frac{1}{8} \delta_{s,s_0}^3
               - \frac{3}{32} \delta_{s,s_0}^4
               + \frac{9}{128} \delta_{s,s_0}^5
               - \frac{7}{128} \delta_{s,s_0}^6}}
        \le C s_0^3 \sum\limits_{i=7}^{15} \lrabs{\delta_{s,s_0}}^i \\
    & & \lrabs{\ca_6^2(s,s_0)
        - s_0^2 
          \lrp{  \frac{1}{4} \delta_{s,s_0}^2
               - \frac{1}{8} \delta_{s,s_0}^3
               + \frac{5}{64} \delta_{s,s_0}^4
               - \frac{7}{128} \delta_{s,s_0}^5
               + \frac{21}{512} \delta_{s,s_0}^6}}
        \le C s_0^2 \sum\limits_{i=7}^{12} \lrabs{\delta_{s,s_0}}^i    
    \EEqn
    where $C$ denotes a universal constant that may vary from line to line.
    
    Now it follows from the first assumption that
    \BEqn
    & & f(s) 
     =  f_0(s_0) + f_1(s_0) \ca_7(s,s_0) + \frac{1}{2} f_2(s_0) \ca_6^2(s,s_0) + \frac{1}{6} f_3(s_0) \ca_5^3(s,s_0) \\
    & & \hspace{4em} 
         + \frac{1}{24} f_4(s_0) \ca_4^4(s,s_0) + \frac{1}{120} f_5(s_0) \ca_3^5(s,s_0) + \frac{1}{720} f_6(s_0) \ca_2^6(s,s_0)  \\
    & & \hspace{4em}
        + f_7(\xi(s,s_0)) \ca_1^7(s,s_0)
    \EEqn
    By reorganizing the quantities above, we obtain that
    \BEqn
    & & \lrabs{f(s) - \sum\limits_{i=0}^{6} c_i(s_0) \delta_{s,s_0} ^i} \\
    &\le& C \lrabs{f_7(\xi(s,s_0)) s_0^7 \lrabs{\delta_{s,s_0}}^7} \\
    & & + C \lrabs{f_6(s_0) s_0^6 \sum\limits_{i=7}^{12} \lrabs{\delta_{s,s_0}}^i} 
        + C \lrabs{f_5(s_0) s_0^5 \sum\limits_{i=7}^{15} \lrabs{\delta_{s,s_0}}^i} 
        + C \lrabs{f_4(s_0) s_0^4 \sum\limits_{i=7}^{16} \lrabs{\delta_{s,s_0}}^i} \\
    & & + C \lrabs{f_3(s_0) s_0^3 \sum\limits_{i=7}^{15} \lrabs{\delta_{s,s_0}}^i} 
        + C \lrabs{f_2(s_0) s_0^2 \sum\limits_{i=7}^{12} \lrabs{\delta_{s,s_0}}^i} 
        + C \lrabs{f_1(s_0) s_0 \lrabs{\delta_{s,s_0}}^7} \\
    &\le& C(\tilde{M}) \sum\limits_{i=7}^{16} \lrabs{\delta_{s,s_0}}^i |f_0(s_0)|,
    \EEqn
    where the last step follows from the assumption that for any $s_0\in D_0$ we have $\max\limits_{1\le i\le 7} \sup\limits_{s\in D} |f_i(s)| \cdot |s_0^i| \le \tilde{M} |f_0(s_0)|$.
\end{proof}

\subsection{Part 3: Preliminary Results for the Power Analysis}

\begin{lemma}\label{Lemma:PowerAnalysis-Prep-1}
It holds that
\BEqn
& & \bigp{A^{X}-A^{XY}} + \bigp{A^{Y}-A^{XY}} 
 =  -2|\Delta|^2, \\
& & \bigp{A^{X}-A^{XY}}^2 + \bigp{A^{Y}-A^{XY}}^2 
 =  2\bigp{\bbe{|\tilde{X}_1|^2} - \bbe{|\tilde{Y}_1|^2}}^2 + 2|\Delta|^4, \\
& & \bigp{A^{X}-A^{XY}}^3 + \bigp{A^{Y}-A^{XY}}^3 
 =  -6\bigp{\bbe{|\tilde{X}_1|^2} - \bbe{|\tilde{Y}_1|^2}}^2|\Delta|^2 -2|\Delta|^6, 
\EEqn
and
\begin{equation*}
    \bigp{A^{X}-A^{XY}}^4 + \bigp{A^{Y}-A^{XY}}^4 
  = 2\bigp{\bbe{|\tilde{X}_1|^2} - \bbe{|\tilde{Y}_1|^2}}^4 + 12\bigp{\bbe{|\tilde{X}_1|^2} - \bbe{|\tilde{Y}_1|^2}}^2|\Delta|^4 + 2|\Delta|^8.
\end{equation*}
\end{lemma}
\begin{proof}
  Recall that
  \BEqn
  & & A^{X} 
   =  \bbe{|\tilde{X}_1-\tilde{X}_2|^2}
   =  2\bbe{|\tilde{X}_1|^2}, \\
  & & A^{Y} 
   =  \bbe{|\tilde{Y}_1-\tilde{Y}_2|^2}
   =  2\bbe{|\tilde{Y}_1|^2}, \\
  & & A^{XY} 
   =  \bbe{|\tilde{X}_1-\tilde{Y}_2+\Delta|^2}
   =  \bbe{|\tilde{X}_1|^2} + \bbe{|\tilde{Y}_1|^2} + |\Delta|^2,
  \EEqn
  which implies that
  \BEqn
  & & A^{X} - A^{XY} 
   =  \bigp{\bbe{|\tilde{X}_1|^2} - \bbe{|\tilde{Y}_1|^2}} - |\Delta|^2, \\
  & & A^{Y} - A^{XY} 
   =  -\bigp{\bbe{|\tilde{X}_1|^2} - \bbe{|\tilde{Y}_1|^2}} - |\Delta|^2,
  \EEqn
  and the statements directly follow from some simple calculations.
\end{proof}

\begin{lemma}\label{Lemma:PowerAnalysis-Prep-2}
It hold that 
\begin{equation*}
    \max\{|A^{X}-A^{XY}|, |A^{Y}-A^{XY}|\}
\le \lrabs{\bbe{|\tilde{X}_1|^2}-\bbe{|\tilde{Y}_1|^2}} + |\Delta|^2.
\end{equation*}
\end{lemma}
\begin{proof}
  As shown in the proof of Lemma \ref{Lemma:PowerAnalysis-Prep-1}, we have
  \BEqn
  & & A^{X} - A^{XY} 
   =  \bigp{\bbe{|\tilde{X}_1|^2} - \bbe{|\tilde{Y}_1|^2}} - |\Delta|^2, \\
  & & A^{Y} - A^{XY} 
   =  -\bigp{\bbe{|\tilde{X}_1|^2} - \bbe{|\tilde{Y}_1|^2}} - |\Delta|^2,
  \EEqn
  which implies the desired results.
\end{proof}

\begin{lemma}\label{Lemma:PowerAnalysis-Prep-3}
Suppose that Assumption \ref{Assumpt:component-dept}\ref{Assumpt:component-dept-3} and the assumptions in Lemma \ref{Lemma:approx-2} holds with $A_0^{XY}\in D_0$, then we have
\begin{equation*}
    \max\{|c_0(A_0^{X})|,|c_0(A_0^{Y})|\} \le C(\tilde{M},L_0,U_0) |c_0(A_0^{XY})|.
\end{equation*}
\end{lemma}
\begin{proof}
  Recall that by definition we have $c_0(x)=f(x)$, then with $s=A_0^{X}$ and $s_0=A_0^{XY}\in D_0$ it follows from Lemma \ref{Lemma:approx-2} that
  \begin{equation*}
      c_0(A_0^{X})
   =  c_0(A_0^{XY})
   + \sum\limits_{i=1}^{6} c_i(A_0^{XY}) (A^{XY})^{-i} \bigp{A^{X} - A^{XY}}^i
   + R(A^{X},A^{XY})
  \end{equation*}
  where 
  \begin{equation*}
      \lrabs{R(A^{X},A^{XY})}
  \le C(\tilde{M}) \sum\limits_{i=7}^{16} \lrabs{\frac{A^{X}-A^{XY}}{A^{XY}}}^i |c_0(A_0^{XY})|.
  \end{equation*}
  
  When Assumption \ref{Assumpt:component-dept}\ref{Assumpt:component-dept-3} holds, we have 
  \begin{equation*}
      \lrabs{\frac{A^{X}-A^{XY}}{A^{XY}}}
  \le C(L_0,U_0),
  \end{equation*}
  and under the assumptions in Lemma \ref{Lemma:approx-2} with $A_0^{XY}\in D_0$ we also have $\max\limits_{1\le i\le 7} |c_i(A_0^{XY})| \le C(\tilde{M}) |c_0(A_0^{XY})|$, both of which jointly imply that $|c_0(A_0^X)| \le C(\tilde{M},L_0,U_0) |c_0(A^{XY})|$. Similarly, we can show the counterpart for $c_0(A^{Y})$, which jointly lead to the desired result.
\end{proof}

\begin{lemma}\label{Lemma:PowerAnalysis-Prep-4}
Suppose that Assumption \ref{Assumpt:component-dept}\ref{Assumpt:component-dept-3} and the assumptions in Lemma \ref{Lemma:approx-2} holds with $A_0^{XY}\in D_0$, then
\begin{enumerate}[label=(\roman*)]
    \item \label{Lemma:PowerAnalysis-Prep-4-1}
    when $|\Delta|>0$ we have
    \begin{equation*}
        \begin{array}{ll}
               & \lrabs{\bigp{2c_0(A^{XY}) - c_0(A^{X}) - c_0(A^{Y})} - \lrp{2c_1(A_0^{XY}) \frac{|\Delta|^2}{A^{XY}} - 2c_2(A_0^{XY}) \frac{\lrp{\bbe{|\tilde{X}_1|^2} - \bbe{|\tilde{Y}_1|^2}}^2}{(A^{XY})^2}}} \\
           \le & C(\tilde{M},L_0,U_0) |c_0(A_0^{XY})|
           \lrp{\frac{|\Delta|^4}{p^2} + \frac{\lrabs{\bbe{|\tilde{X}_1|^2} - \bbe{|\tilde{Y}_1|^2}}^3}{p^3}}
        \end{array}
    \end{equation*}
    \item \label{Lemma:PowerAnalysis-Prep-4-2}
    when $|\Delta|=0$ we have
    \begin{equation*}
        \begin{array}{ll}
               & \lrabs{\bigp{2c_0(A^{XY}) - c_0(A^{X}) - c_0(A^{Y})} + 2c_2(A_0^{XY}) \frac{\lrp{\bbe{|\tilde{X}_1|^2} - \bbe{|\tilde{Y}_1|^2}}^2}{(A^{XY})^2}} \\
           \le & C(\tilde{M},L_0,U_0) |c_0(A_0^{XY})| \frac{\lrabs{\bbe{|\tilde{X}_1|^2} - \bbe{|\tilde{Y}_1|^2}}^4}{p^4}
        \end{array}
    \end{equation*}
\end{enumerate}

\end{lemma}
\begin{proof}
\begin{enumerate}[label=(\roman*)]
    \item 
    It follows from Lemma \ref{Lemma:approx-2} that
    \begin{equation*}
        c_0(A_0^{X})
     =  c_0(A_0^{XY})
     + \sum\limits_{i=1}^{6} c_i(A_0^{XY}) (A^{XY})^{-i} \bigp{A^{X} - A^{XY}}^i
     + R(A^{X},A^{XY})
    \end{equation*}
    where 
    \begin{equation*}
        \lrabs{R(A^{X},A^{XY})}
    \le C(\tilde{M}) \sum\limits_{i=7}^{16} \lrabs{\frac{A^{X}-A^{XY}}{A^{XY}}}^i |c_0(A_0^{XY})|
    \end{equation*}
    Similarly, we obtain the approximation of $c_0(A^{Y})$, and it follows that
    \BEqn
        2c_0(A_0^{XY}) - c_0(A_0^{X}) - c_0(A_0^{Y})
    &=& - c_1(A_0^{XY}) (A^{XY})^{-1} \lrp{(A^{X}-A^{XY}) + (A^{Y}-A^{XY})} \\
    & & - c_2(A_0^{XY}) (A^{XY})^{-2} \lrp{(A^{X}-A^{XY})^2 + (A^{Y}-A^{XY})^2} \\
    & & - \sum\limits_{i=3}^{6} c_i(A_0^{XY}) (A^{XY})^{-i} \lrp{(A^{X}-A^{XY})^i + (A^{Y}-A^{XY})^i} \\
    & & - R(A^{X},A^{XY}) - R(A^{Y},A^{XY})
    \EEqn
    Also note that from Lemma \ref{Lemma:PowerAnalysis-Prep-1} that
    \BEqn
    & & \bigp{A^{X}-A^{XY}} + \bigp{A^{Y}-A^{XY}} 
     =  -2|\Delta|^2, \\
    & & \bigp{A^{X}-A^{XY}}^2 + \bigp{A^{Y}-A^{XY}}^2 
     =  2\bigp{\bbe{|\tilde{X}_1|^2} - \bbe{|\tilde{Y}_1|^2}}^2 + 2|\Delta|^4,
    \EEqn
    and it follows from Lemma \ref{Lemma:PowerAnalysis-Prep-2} that
    \begin{equation*}
        \max\{|A^{X}-A^{XY}|, |A^{Y}-A^{XY}|\}
    \le \lrabs{\bbe{|\tilde{X}_1|^2}-\bbe{|\tilde{Y}_1|^2}} + |\Delta|^2,
    \end{equation*}
    which further implies that
    \begin{equation*}
        \max\{|A^{X}-A^{XY}|^i, |A^{Y}-A^{XY}|^i\}
    \le C \lrp{\lrabs{\bbe{|\tilde{X}_1|^2}-\bbe{|\tilde{Y}_1|^2}}^i + |\Delta|^{2i}}
    \end{equation*}
    by the $C_r$ inequality.
    
    Hence we conclude that
    \BEqn
    & & \lrabs{\bigp{2c_0(A_0^{XY}) - c_0(A_0^{X}) - c_0(A_0^{Y})} - \lrp{2c_1(A_0^{XY}) \frac{|\Delta|^2}{A^{XY}} - 2c_2(A_0^{XY}) \frac{\lrp{\bbe{|\tilde{X}_1|^2} - \bbe{|\tilde{Y}_1|^2}}^2}{(A^{XY})^2}}} \\
    &\le& 2|c_2(A_0^{XY})| (A^{XY})^{-2} |\Delta|^4 + \sum\limits_{i=3}^{6} |c_i(A_0^{XY})| (A^{XY})^{-i} \lrp{|A^{X}-A^{XY}|^i + |A^{Y}-A^{XY}|^i} \\
    & & + |R(A^{X},A^{XY})| + |R(A^{Y},A^{XY})| \\
    &\le& C(\tilde{M},L_0,U_0) |c_0(A_0^{XY})| \lrp{\frac{|\Delta|^4}{p^2} + \frac{\lrabs{\bbe{|\tilde{X}_1|^2} - \bbe{|\tilde{Y}_1|^2}}^3}{p^3}}
    \EEqn
    \item
    When $|\Delta|=0$, it follows from Lemma \ref{Lemma:PowerAnalysis-Prep-1} that $(A^{X}-A^{XY})^i + (A^{Y}-A^{XY})^i = 0$ if $i$ is odd, and $(A^{X}-A^{XY})^i + (A^{Y}-A^{XY})^i = 2\bigp{\bbe{|\tilde{X}_1|^2} - \bbe{|\tilde{Y}_1|^2}}^i$ if $i$ is even. Thus in this case, it follows from Lemma \ref{Lemma:approx-2} that
    \BEqn
        \bigp{2c_0(A_0^{XY}) - c_0(A_0^{X}) - c_0(A_0^{Y})}
    &=& - c_2(A_0^{XY}) (A^{XY})^{-2} \lrp{(A^{X}-A^{XY})^2 + (A^{Y}-A^{XY})^2} \\
    & & - c_4(A_0^{XY}) (A^{XY})^{-4} \lrp{(A^{X}-A^{XY})^4 + (A^{Y}-A^{XY})^4} \\
    & & - c_6(A_0^{XY}) (A^{XY})^{-6} \lrp{(A^{X}-A^{XY})^6 + (A^{Y}-A^{XY})^6} \\
    & & - R(A^{X},A^{XY}) - R(A^{Y},A^{XY}) \\
    &=& - 2c_2(A_0^{XY}) (A^{XY})^{-2} \bigp{\bbe{|\tilde{X}_1|^2} - \bbe{|\tilde{Y}_1|^2}}^2 \\
    & & - 2c_4(A_0^{XY}) (A^{XY})^{-4} \bigp{\bbe{|\tilde{X}_1|^2} - \bbe{|\tilde{Y}_1|^2}}^4 \\
    & & - 2c_6(A_0^{XY}) (A^{XY})^{-6} \bigp{\bbe{|\tilde{X}_1|^2} - \bbe{|\tilde{Y}_1|^2}}^6 \\
    & & - R(A^{X},A^{XY}) - R(A^{Y},A^{XY})
    \EEqn
    which implies the proposed result.
\end{enumerate}
\end{proof}

\begin{lemma}\label{Lemma:PowerAnalysis-Prep-5}
Consider $f: [0,\infty) \mapsto \br$ and let $D_0 \subseteq D := [0,\infty)$ be a subset of the domain of $f$ and $A_0^{XY}\in D_0$. Assume that
\begin{enumerate}[label=(\roman*)]
    \item for any $s\in D$ and $s_0\in D_0$, it holds that
    \begin{equation*}
        f_1(s) = f_1(s_0) + f_2(\xi(s,s_0)) (s-s_0),
    \end{equation*}
    where $f_i$ denotes the $i$-th order derivative of $f$ and $\xi(s,s_0)$ is a point between $s$ and $s_0$.
    
    \item there exists a positive constant $\tilde{M} < \infty$, such that for any $s_0 \in D_0$, it holds that
    \begin{equation*}
        \max\limits_{1\le i\le 2} \sup\limits_{s\in D} |f_i(s)|\cdot |s_0^i| \le \tilde{M} |f_0(s_0)|.
    \end{equation*}
\end{enumerate}

Suppose that Assumption \ref{Assumpt:component-dept}\ref{Assumpt:component-dept-3} holds, then it holds that 
\BEqn
& & \lrabs{c_1(A_0^{X}) (A^X)^{-1} - c_1(A_0^{XY}) (A^{XY})^{-1}}
\le C(\tilde{M},L_0,U_0) |c_0(A_0^{XY})| \frac{|A^{X}-A^{XY}|}{p^2}, \\
& & \lrabs{c_1(A_0^{Y}) (A^Y)^{-1} - c_1(A_0^{XY}) (A^{XY})^{-1}}
\le C(\tilde{M},L_0,U_0) |c_0(A_0^{XY})| \frac{|A^{Y}-A^{XY}|}{p^2}.
\EEqn
where $c_i(x)$ is defined as Lemma \ref{Lemma:approx-2}.
\end{lemma}
\begin{proof}
  Let $s = A_0^{X} = A_0^{XY}\sqrt{1+\frac{A^{X}-A^{XY}}{A^{XY}}}$ and $s_0 = A_0^{XY}$, then it follows from Lemma \ref{Lemma:approx-1} that, there exists a constant $c(A^{X},A^{XY})$, such that
  \begin{equation*}
      s - s_0 = c(A^{X},A^{XY}) A_0^{XY} \frac{A^{X}-A^{XY}}{A^{XY}},
  \end{equation*}
  and $c(A^{X},A^{XY}) \le C$ for a universal constant $C<\infty$. Then under the assumptions, we obtain that
  \BEqn
      |f_1(A_0^{X}) - f_1(A_0^{XY})|
  &=& \lrabs{c(A^{X},A^{XY}) f_2(\xi(A_0^{X},A_0^{XY})) A_0^{XY} \frac{A^{X}-A^{XY}}{A^{XY}}} \\
  &=& \lrabs{c(A^{X},A^{XY}) f_2(\xi(A_0^{X},A_0^{XY})) (A_0^{XY})^2} \cdot \frac{|A^{X}-A^{XY}|}{(A_0^{XY})^3} \\
  &\le& C(\tilde{M}) |c_0(A_0^{XY})| \frac{|A^{X}-A^{XY}|}{(A_0^{XY})^3}
  \EEqn
  
  Also note that 
  \begin{equation*}
      (A_0^{X})^{-1} - (A_0^{XY})^{-1}
    = \frac{A_0^{XY}-A_0^{X}}{A_0^{X} A_0^{XY}}
    = \frac{A^{XY} - A^{X}}{A_0^{X} A_0^{XY}(A_0^{X}+A_0^{XY})},
  \end{equation*}
  then it follows from the definition $c_1(s) = \frac{1}{2}f_1(s)s$ that
  \BEqn
  & & \lrabs{c_1(A_0^{X})(A^X)^{-1} - c_1(A_0^{XY})(A^{XY})^{-1}} \\
  &=& \frac{1}{2} \lrabs{f_1(A_0^{X})(A_0^X)^{-1} - f_1(A_0^{XY})(A_0^{XY})^{-1}} \\
  &\le& \frac{1}{2} \lrabs{f_1(A_0^{X})-f_1(A_0^{XY})} (A_0^{X})^{-1} 
      + \frac{1}{2} \lrabs{f_1(A_0^{XY})} \cdot \lrabs{(A_0^{X})^{-1} - (A_0^{XY})^{-1}} \\
  &\le& C(\tilde{M}) 
        \lrp{  |c_0(A_0^{XY})| \frac{|A^{X}-A^{XY}|}{A_0^{X} (A_0^{XY})^3}
             + |f_1(A_0^{XY}) (A_0^{XY})| \frac{|A^{X} - A^{XY}|}{A_0^{X} (A_0^{XY})^2 (A_0^{X}+A_0^{XY})}} \\
  &\le& C(\tilde{M},L_0,U_0) |c_0(A_0^{XY})| \frac{|A^{X}-A^{XY}|}{p^2}.
  \EEqn
  Using the same arguments, we also obtain that
  \begin{equation*}
      \lrabs{c_1(A_0^{Y}) (A^Y)^{-1} - c_1(A_0^{XY}) (A^{XY})^{-1}}
  \le C(\tilde{M},L_0,U_0) |c_0(A_0^{XY})| \frac{|A^{Y}-A^{XY}|}{p^2},
  \end{equation*}
  which completes the proof.
\end{proof}

\begin{lemma}\label{Lemma:PowerAnalysis-Prep-6}
Consider $f: [0,\infty) \mapsto \br$ and let $D_0 \subseteq D := [0,\infty)$ be a subset of the domain of $f$ and $A_0^{XY}\in D_0$. Assume that
\begin{enumerate}[label=(\roman*)]
    \item for any $s\in D$ and $s_0\in D_0$, it holds that
    \begin{equation*}
        f_1(s) = f_1(s_0) + f_2(\xi(s,s_0)) (s-s_0),
    \end{equation*}
    where $f_i$ denotes the $i$-th order derivative of $f$ and $\xi(s,s_0)$ is a point between $s$ and $s_0$.
    
    \item there exists a positive constant $\tilde{M} < \infty$, such that for any $s_0 \in D_0$, it holds that
    \begin{equation*}
        \max\limits_{1\le i\le 2} \sup\limits_{s\in D} |f_i(s)|\cdot |s_0^i| \le \tilde{M} |f_0(s_0)|.
    \end{equation*}
\end{enumerate}

Suppose that Assumption \ref{Assumpt:component-dept}\ref{Assumpt:component-dept-3} holds, then it holds that 
\BEqn
& & \lrabs{c_1^2(A_0^{X}) (A^X)^{-2} - c_1^2(A_0^{XY}) (A^{XY})^{-2}}
\le C(\tilde{M},L_0,U_0) c_0^2(A_0^{XY}) \frac{|A^{X}-A^{XY}|}{p^3}, \\
& & \lrabs{c_1^2(A_0^{Y}) (A^Y)^{-2} - c_1^2(A_0^{XY}) (A^{XY})^{-2}}
\le C(\tilde{M},L_0,U_0) c_0^2(A_0^{XY}) \frac{|A^{Y}-A^{XY}|}{p^3}.
\EEqn
where $c_i(x)$ is defined as Lemma \ref{Lemma:approx-2}.
\end{lemma}
\begin{proof}
  We only provide the proof of the first statement. From Lemma \ref{Lemma:PowerAnalysis-Prep-5}, we have that
  \BEqn
  & & \lrabs{c_1(A_0^{X}) (A^X)^{-1} + c_1(A_0^{XY}) (A^{XY})^{-1}} \\
  &\le& \lrabs{c_1(A_0^{X}) (A^X)^{-1} - c_1(A_0^{XY}) (A^{XY})^{-1}} + 2|c_1(A_0^{XY})| (A^{XY})^{-1} \\
  &\le& C(\tilde{M},L_0,U_0) |c_0(A_0^{XY})| \lrp{\frac{|A^{X}-A^{XY}|}{p^2} + \frac{1}{p}} \\
  &\le& C(\tilde{M},L_0,U_0) |c_0(A_0^{XY})| \frac{1}{p}
  \EEqn
  and consequently, we obtain that
  \BEqn
  & & \lrabs{c_1^2(A_0^{X}) (A^X)^{-2} - c_1^2(A_0^{XY}) (A^{XY})^{-2}} \\
  &=& \lrabs{c_1(A_0^{X}) (A^X)^{-1} - c_1(A_0^{XY}) (A^{XY})^{-1}} 
      \cdot \lrabs{c_1(A_0^{X}) (A^X)^{-1} + c_1(A_0^{XY}) (A^{XY})^{-1}} \\
  &\le& C(\tilde{M},L_0,U_0) c_0^2(A_0^{XY}) \frac{|A^{X}-A^{XY}|}{p^3}.
  \EEqn
\end{proof}

\begin{lemma}\label{Lemma:PowerAnalysis-Prep-7}
Consider $f: [0,\infty) \mapsto \br$ and let $D_0 \subseteq D := [0,\infty)$ be a subset of the domain of $f$ and $A_0^{XY}\in D_0$. Assume that
\begin{enumerate}[label=(\roman*)]
    \item for any $s\in D$ and $s_0\in D_0$, it holds that
    \begin{equation*}
        f_2(s) = f_2(s_0) + f_3(\xi(s,s_0)) (s-s_0),
    \end{equation*}
    where $f_i$ denotes the $i$-th order derivative of $f$ and $\xi(s,s_0)$ is some value between 
    $s$ and $s_0$.
    
    \item there exists a positive constant $\tilde{M} < \infty$, such that for any $s_0 \in D_0$, it holds that
    \begin{equation*}
        \max\limits_{1\le i\le 2} \sup\limits_{s\in D} |f_i(s)|\cdot |s_0^i| \le \tilde{M} |f_0(s_0)|.
    \end{equation*}
\end{enumerate}

Suppose that Assumption \ref{Assumpt:component-dept}\ref{Assumpt:component-dept-3} holds, then it holds that 
\BEqn
& & \lrabs{c_2(A_0^{X}) (A^X)^{-2} - c_2(A_0^{XY}) (A^{XY})^{-2}}
\le C(\tilde{M},L_0,U_0) |c_0(A_0^{XY})| \frac{|A^{X}-A^{XY}|}{p^3}, \\
& & \lrabs{c_2(A_0^{Y}) (A^Y)^{-2} - c_2(A_0^{XY}) (A^{XY})^{-2}}
\le C(\tilde{M},L_0,U_0) |c_0(A_0^{XY})| \frac{|A^{Y}-A^{XY}|}{p^3}.
\EEqn
where $c_i(x)$ is defined as Lemma \ref{Lemma:approx-2}.
\end{lemma}
\begin{proof}
  Following the similar arguments used in the proof of Lemma \ref{Lemma:PowerAnalysis-Prep-5}, we obtain under the approximation assumptions that
  \BEqn
      |f_1(A_0^{X}) - f_1(A_0^{XY})|
  &\le& C(\tilde{M}) |c_0(A_0^{XY})| \frac{|A^{X}-A^{XY}|}{(A_0^{XY})^3} \\
      |f_2(A_0^{X}) - f_2(A_0^{XY})|
  &\le& C(\tilde{M}) |c_0(A_0^{XY})| \frac{|A^{X}-A^{XY}|}{(A_0^{XY})^4}
  \EEqn
  We also have
  \BEqn
  & & (A_0^X)^{-2} - (A_0^{XY})^{-2}
   =  \frac{A^{XY} - A^{X}}{(A_0^X)^2(A_0^{XY})^2} \\[2mm]
  & & (A_0^X)^{-3} - (A_0^{XY})^{-3}
   =  \frac{(A_0^{XY})^3 - (A_0^{X})^3}{(A_0^X)^3(A_0^{XY})^3} 
   =  \frac{\bigp{A^{XY}-A^{X}} \bigp{(A^{XY})^2 + A^{XY}A^{X} + (A^{X})^2}}{(A_0^X)^3 (A_0^{XY})^3 \bigp{(A_0^X)^3+(A_0^{XY})^3}} \\
  \EEqn
  
  Recall that $c_2(x) = -\frac{1}{8}f_1(x)x + \frac{1}{8}f_2(x)x^2$, then we have
  \BEqn
  & & \lrabs{c_2(A_0^{X}) (A^X)^{-2} - c_2(A_0^{XY}) (A^{XY})^{-2}} \\
  &\le& \frac{1}{8} \lrabs{f_1(A_0^{X}) (A_0^{X})^{-3} - f_1(A_0^{XY}) (A_0^{XY})^{-3}}
        + \frac{1}{8} \lrabs{f_2(A_0^{X}) (A_0^{X})^{-2} - f_2(A_0^{XY}) (A_0^{XY})^{-2}} \\
  &\le& \frac{1}{8} \lrabs{f_1(A_0^{X}) - f_1(A_0^{XY})} (A_0^{X})^{-3}
        + \frac{1}{8} \lrabs{f_1(A_0^{XY})} \cdot \lrabs{(A_0^{X})^{-3} - (A_0^{XY})^{-3}} \\
  & &   + \frac{1}{8} \lrabs{f_2(A_0^{X}) - f_2(A_0^{XY})} (A_0^{X})^{-2}
        + \frac{1}{8} \lrabs{f_2(A_0^{XY})} \cdot \lrabs{(A_0^{X})^{-2} - (A_0^{XY})^{-2}} \\[2mm]
  &\le& C(\tilde{M},L_0,U_0) \lrabs{c_0(A_0^{XY})} \frac{\lrabs{A^{X} - A^{XY}}}{p^3},
  \EEqn
  and using the same technique, we can obtain the other statement.
\end{proof}

\begin{lemma}\label{Lemma:PowerAnalysis-Prep-8}

It holds that
\begin{equation*}
    2L_2^{XY} - L_2^{X} - L_2^{Y}
  = -4 \|\Sigma_X - \Sigma_Y\|_F^2 + 8 \Delta^{\top} (\Sigma_X + \Sigma_Y) \Delta + 8 \fone^{\top} (\Gamma_X - \Gamma_Y) \Delta,
\end{equation*}
where $\Gamma_X = (\gamma_{X,j_1 j_2})_{p\times p}$ and $\Gamma_Y = (\gamma_{Y,j_1 j_2})_{p\times p}$ are both $p\times p$ matrices with
  \begin{equation*}
      \gamma_{X, j_1 j_2} = \cum(\tilde{x}_{1 j_1}^2 - \sigma_{X, j_1}^2, \tilde{x}_{2 j_2}), \qquad
      \gamma_{Y, j_1 j_2} = \cum(\tilde{y}_{1 j_1}^2 - \sigma_{Y, j_1}^2, \tilde{y}_{2 j_2}).
  \end{equation*}
\end{lemma}
\begin{proof}
  By the definition, it follows some simple calculations that
  \BEqn
  & & 2L_2^{XY} - L_2^{X} - L_2^{Y} \\ 
  &=& 8\sum\limits_{j_1,j_2=1}^{p} \cum(\tilde{x}_{1 j_1}, \tilde{x}_{1 j_2}) \cum(\tilde{y}_{1 j_1}, \tilde{y}_{1 j_2})
      - 4\sum\limits_{j_1,j_2=1}^{p} \cum^2(\tilde{x}_{1 j_1}, \tilde{x}_{1 j_2})
      - 4\sum\limits_{j_1,j_2=1}^{p} \cum^2(\tilde{y}_{1 j_1}, \tilde{y}_{1 j_2}) \\
  & & + 8\sum\limits_{j_1,j_2=1}^{p} \bigp{\cum(\tilde{x}_{1 j_1}, \tilde{x}_{1 j_2}) + \cum(\tilde{y}_{1 j_1}, \tilde{y}_{1 j_2})} \Delta_{j_1} \Delta_{j_2} \\
  & & + 8\sum\limits_{j_1,j_2=1}^{p} \bigp{\cum(\tilde{x}_{1 j_1}^2-\sigma_{j_1}^2, \tilde{x}_{1 j_2}) - \cum(\tilde{y}_{1 j_1}^2-\sigma_{j_1}^2, \tilde{y}_{1 j_2})} \Delta_{j_2} \\
  &=& -4 \|\Sigma_X - \Sigma_Y\|_F^2 + 8 \Delta^{\top} (\Sigma_X + \Sigma_Y) \Delta + 8 \fone^{\top} (\Gamma_X - \Gamma_Y) \Delta,
  \EEqn
  which completes the proof.
\end{proof}

\begin{lemma}\label{Lemma:PowerAnalysis-Prep-9}

Suppose that Assumption \ref{Assumpt:component-dept}\ref{Assumpt:component-dept-1}-\ref{Assumpt:component-dept-2} holds and
\begin{equation*}
    |\Delta|^2 
= \|\Sigma_X - \Sigma_Y\|_F = 0,
\end{equation*}
then it holds that
\BEqn
& & 2L_2^{XY} - L_2^{X} - L_2^{Y} = 0, \\
& & 2L_3^{XY} - L_3^{X} - L_3^{Y} 
 =  8\sum\limits_{j_1,j_2,j_3=1}^{p}\lrp{\cum(\tilde{x}_{1 j_1}, \tilde{x}_{1 j_2}, \tilde{x}_{1 j_3}) - \cum(\tilde{y}_{1 j_1}, \tilde{y}_{1 j_2}, \tilde{y}_{1 j_3})}^2 \\
& & \hspace{9em}
    + 12\sum\limits_{j_2=1}^{p}\lrp{\sum\limits_{j_1=1}^{p} \bigp{\cum(\tilde{x}_{1 j_1}^2 - \sigma_{X, j_1}^2, \tilde{x}_{1 j_1}) - \cum(\tilde{y}_{1 j_1}^2 - \sigma_{Y, j_1}^2, \tilde{y}_{1 j_1})}}^2
\EEqn
and
\BEqn
& & \lrabs{2L_4^{XY} - L_4^{X} - L_4^{Y} + 6\lrp{\sum\limits_{j_1,j_2=1}^{p} \bigp{\cum(\tilde{x}_{1 j_1}^2 - \sigma_{X, j_1}^2, \tilde{x}_{1 j_2}^2 - \sigma_{X, j_2}^2) - \cum(\tilde{y}_{1 j_1}^2 - \sigma_{Y, j_1}^2, \tilde{y}_{1 j_2}^2 - \sigma_{Y, j_2}^2)}}^2} \\
&\le& C(U^{\ast}) (\alpha(p))^3 p
\EEqn
\end{lemma}
\begin{proof}
  It follows from Lemma \ref{Lemma:PowerAnalysis-Prep-8} that $2L_2^{XY} - L_2^{X} - L_2^{Y} = 0$ when $|\Delta|^2=\|\Sigma_x - \Sigma_Y\|_F=0$. The second statement follows from direct computation and we hence spare the details. As for the last statement, note that
  \BEqn
  & & 2L_4^{XY} - L_4^{X} - L_4^{Y} \\
  &=& - 6 \lrp{\sum\limits_{j_1,j_2=1}^{p} \bigp{\cum(\tilde{x}_{1 j_1}^2 - \sigma_{X, j_1}^2, \tilde{x}_{1 j_2}^2 - \sigma_{X, j_2}^2) - \cum(\tilde{y}_{1 j_1}^2 - \sigma_{Y, j_1}^2, \tilde{y}_{1 j_2}^2 - \sigma_{Y, j_2}^2)}}^2 \\
  & & - 16 \sum\limits_{j_1,j_2,j_3,j_4=1}^{p} \lrp{\cum(\tilde{x}_{1 j_1}, \tilde{x}_{1 j_2}, \tilde{x}_{1 j_3}, \tilde{x}_{1 j_4}) - \cum(\tilde{y}_{1 j_1}, \tilde{y}_{1 j_2}, \tilde{y}_{1 j_3}, \tilde{y}_{1 j_4})}^2 \\
  & & + 64 \sum\limits_{j_1,j_2,j_3,j_4=1}^{p} \lrp{\cum(\tilde{x}_{1 j_1}^2 - \sigma_{X, j_1}^2, \tilde{x}_{1 j_2}, \tilde{x}_{1 j_3}, \tilde{x}_{1 j_4})) - \cum(\tilde{y}_{1 j_1}^2 - \sigma_{Y, j_1}^2, \tilde{y}_{1 j_2}, \tilde{y}_{1 j_3}, \tilde{y}_{1 j_4})} \\
  & & \hspace{6em}
      \times \lrp{\cum(\tilde{x}_{1 j_2}, \tilde{x}_{1 j_3}, \tilde{x}_{1 j_4})) - \cum(\tilde{y}_{1 j_2}, \tilde{y}_{1 j_3}, \tilde{y}_{1 j_4})} \\
  & & + 192 \sum\limits_{j_1,j_2,j_3,j_4=1}^{p} \lrp{\cum(\tilde{x}_{1 j_1}^2 - \sigma_{X, j_1}^2, \tilde{x}_{1 j_2})\cum(\tilde{x}_{1 j_3}, \tilde{x}_{1 j_4})) - \cum(\tilde{y}_{1 j_1}^2 - \sigma_{Y, j_1}^2, \tilde{y}_{1 j_2})\cum(\tilde{y}_{1 j_3}, \tilde{y}_{1 j_4})} \\
  & & \hspace{6em}
      \times \lrp{\cum(\tilde{x}_{1 j_2}, \tilde{x}_{1 j_3}, \tilde{x}_{1 j_4})) - \cum(\tilde{y}_{1 j_2}, \tilde{y}_{1 j_3}, \tilde{y}_{1 j_4})} \\
  & & + 48 \sum\limits_{j_1,j_2,j_3,j_4=1}^{p} \lrp{\cum(\tilde{x}_{1 j_1}^2 - \sigma_{X, j_1}^2, \tilde{x}_{1 j_2}^2 - \sigma_{X, j_2}^2, \tilde{x}_{1 j_4}) - \cum(\tilde{y}_{1 j_1}^2 - \sigma_{Y, j_1}^2, \tilde{y}_{1 j_2}^2 - \sigma_{Y, j_2}^2, \tilde{y}_{1 j_4})} \\
  & & \hspace{6em}
      \times \lrp{\cum(\tilde{x}_{1 j_3}^2 - \sigma_{X, j_3}^2, \tilde{x}_{1 j_4}) - \cum(\tilde{y}_{1 j_3}^2 - \sigma_{Y, j_3}^2, \tilde{y}_{1 j_4})} \\
  & & - 48 \sum\limits_{j_1,j_2,j_3,j_4=1}^{p} \lrp{\cum(\tilde{x}_{1 j_1}^2 - \sigma_{X, j_1}^2, \tilde{x}_{1 j_3}, \tilde{x}_{1 j_4}) - \cum(\tilde{y}_{1 j_1}^2 - \sigma_{Y, j_1}^2, \tilde{y}_{1 j_3}, \tilde{y}_{1 j_4})} \\
  & & \hspace{6em}
      \times \lrp{\cum(\tilde{x}_{1 j_2}^2 - \sigma_{X, j_2}^2, \tilde{x}_{1 j_3}, \tilde{x}_{1 j_4}) - \cum(\tilde{y}_{1 j_2}^2 - \sigma_{Y, j_2}^2, \tilde{y}_{1 j_3}, \tilde{y}_{1 j_4})}
  \EEqn
  
  Under Assumption \ref{Assumpt:component-dept}\ref{Assumpt:component-dept-1}, the total numbers of nonzero individual terms in all but the first summations are uniformly bounded by $(\alpha(p))^3 p$, then the result follows from Assumption \ref{Assumpt:component-dept}\ref{Assumpt:component-dept-2}.
\end{proof}


\section{Lemmas Regarding the Test Statistic}\label{Appendix:lemma-1}

\subsection{Lemmas for Proposition \ref{Prop:est_MMD_simplified}}

The main theoretical tool used to prove Proposition \ref{Prop:est_MMD_simplified} is the Hoeffding decomposition of the U-statistics. We list some relevant results in \cite{lee1990u} to facilitate the subsequent analysis. 

\begin{definition}[Section 2.2 in \cite{lee1990u}]\label{Def:psi^(cd)}
  Assume that $X_1,\dots,X_{n_x}$ are independently sampled from the distribution $F_X$ and $Y_1,\dots,Y_{n_y}$ are independently sampled from the distribution $F_Y$, and $\{X_{i}\}_{i=1}^{n_x}$ is independent of $\{Y_{j}\}_{j=1}^{n_y}$. Let $\psi$ be a kernel with $k_x + k_y$ arguments
  \begin{equation*}
    \psi(x_1,\dots,x_{k_x};y_1,\dots,y_{k_y}),
  \end{equation*}
  which is symmetric in $x_1,\dots,x_{k_x}$ and $y_1,\dots,y_{k_y}$. Let $G_x$ denote the density function of a single point mass at $x$, then for $1 \le c \le k_x$ and $1 \le d \le k_y$, we can define
  \begin{align*}
   & \psi^{(c,d)}(x_1,\dots,x_c;y_1,\dots,y_d) \\
  =& \int \dots \int \psi(u_1,\dots,u_{k_x}; v_1,\dots,v_{k_y}) \prod\limits_{i=1}^{c}(dG_{x_i}(u_i)-dF_X(u_i)) \prod\limits_{i=c+1}^{k_x}dF_X(u_i) \\
   & \quad \times \prod\limits_{j=1}^{d}(dG_{y_j}(v_j)-dF_Y(v_j)) \prod\limits_{j=d+1}^{k_y}dF_Y(v_j).
  \end{align*}
\end{definition}

\begin{definition}[Section 2.2 in \cite{lee1990u}]\label{Def:H^(cd)}
  Assume that $X_1,\dots,X_{n_x}$ are independently sampled from the distribution $F_X$ and $Y_1,\dots,Y_{n_y}$ are independently sampled from the distribution $F_Y$, and $\{X_{i}\}_{i=1}^{n_x}$ is independent of $\{Y_{j}\}_{j=1}^{n_y}$. Following the notations in Definition \ref{Def:psi^(cd)}, we can define
  \begin{equation*}
    H_{n_x,n_y}^{(c,d)} = \bin{n_x}{c}^{-1} \bin{n_y}{d}^{-1} \sum\limits_{(n_x,c)} \sum\limits_{(n_y,d)} \psi^{(c,d)}(X_{i_1},\dots,X_{i_c}; Y_{j_1},\dots,Y_{j_d}).
  \end{equation*}
\end{definition}

\begin{lemma}[Theorem 3 in Section 2.2, \cite{lee1990u}]\label{Lemma:H-decomp}
  The U-statistic based on the kernel $\psi$ is a statistic of the form
  \begin{equation*}
    U_{n_x,n_y} = \bin{n_x}{k_x}^{-1} \bin{n_y}{k_y}^{-1} \sum\limits_{(n_x,k_x)} \sum\limits_{(n_y,k_y)} \psi(S_x,S_y),
  \end{equation*}
  where the sum is summed over all $k_x$-subsets $S_x$ of $(x_1,\dots,x_{n_x})$ and all $k_y$-subsets $S_y$ of $(y_1,\dots,y_{n_y})$. Then $U_{n_x,n_y}$ is an unbiased estimator of $\be[\psi(X_1,\dots,X_{k_x};Y_1,\dots,Y_{k_y})]$, and it admits the representation
  \begin{equation*}
    U_{n_x,n_y} = \sum\limits_{c=0}^{k_x} \sum\limits_{d=0}^{k_y} \bin{k_x}{c} \bin{k_y}{d} H_{n_x,n_y}^{(c,d)},
  \end{equation*}
  where $H_{n_x,n_y}^{(c,d)}$ is the U-statistic based on the kernel $\psi^{(c,d)}$ and is given by
  \begin{equation*}
    H_{n_x,n_y}^{(c,d)} = \bin{n_x}{c}^{-1} \bin{n_y}{d}^{-1} \sum\limits_{(n_x,c)} \sum\limits_{(n_y,d)} \psi^{(c,d)}(X_{i_1},\dots,X_{i_c}; Y_{j_1},\dots,Y_{j_d}).
  \end{equation*}
  Furthermore, the functions $\psi^{(c,d)}$ satisfy
  \begin{enumerate}[label=(\roman*)]
    \item $\bbe{\psi^{(c,d)}(X_1,\dots,X_c; Y_1,\dots,Y_d)} = 0$;
    \item $\Cov(\psi^{(c,d)}(S_x,S_y), \psi^{(c',d')}(S_x',S_y')) = 0$ for all integers $c,d,c',d'$ and sets $S_x,S_y,S_x',S_y'$ unless $c=c', d=d'$ and $S_x=S_x', S_y=S_y'$.
  \end{enumerate}
\end{lemma}

Using the Hoeffding decomposition, we can break the two-sample U-statistic $\ce_{n,m}^{k}(X,Y)$ into two parts, namely, $L_{n,m}^{K}(X,Y)$ and $R_{n,m}^{k}(X,Y)$; see Lemma \ref{Lemma:estMMD_k}.

\begin{lemma}\label{Lemma:estMMD_k}
  Let $F_X$ and $F_Y$ denote the distribution functions of $X$ and $Y$, respectively, and $G_x$ denote the distribution function of a single point mass at $x$. For $0\le c,d \le 2$, define
  \begin{align*}
  &  h^{(c,d)}(X_1,\dots,X_c;Y_1,\dots,Y_d)  \\
  =& \int\dots\int h^{k}(u_1,u_2,v_1,v_2) \prod\limits_{i=1}^{c}\lrp{dG_{X_i}(u_i) - dF_X(u_i)} \prod\limits_{i=c+1}^{2}dF_X(u_i) \\
   &\quad \times \prod\limits_{j=1}^{d}\lrp{dG_{Y_j}(v_j) - dF_Y(v_j)} \prod\limits_{j=d+1}^{2}dF_Y(v_j).
  \end{align*}
  The two-sample U-statistic $\ce_{n,m}^k(X,Y)$ admits the following decomposition
  \begin{equation*}
    \ce_{n,m}^k(X,Y) = L_{n,m}^k(X,Y) + R_{n,m}^k(X,Y),
  \end{equation*}
  where
  \begin{align*}
    L_{n,m}^k(X,Y)
  =& \ce^k(X,Y) 
     + \frac{2}{n}\sum\limits_{i=1}^n h^{(1,0)}(X_i) 
     + \frac{2}{m}\sum\limits_{j=1}^{m} h^{(0,1)}(Y_j) \\
   & + \frac{2}{n(n-1)} \sum\limits_{1 \le i_1 < i_2 \le n} h^{(2,0)}(X_{i_1},X_{i_2})
     + \frac{4}{nm} \sumin \sumjm h^{(1,1)}(X_i,Y_j) \nonumber \\
   & + \frac{2}{m(m-1)} \sum\limits_{1 \le j_1 < j_2 \le m} h^{(0,2)}(Y_{j_1},Y_{j_2})
  \end{align*}
  is the dominant term of $\ce_{n,m}^k(X,Y)$, and
  \begin{align*}\label{Equ:Rnm}
    R_{n,m}^k(X,Y)
  =& \frac{4}{n(n-1)m} \sum\limits_{1 \le i_1 < i_2 \le n} \sumjm h^{(2,1)}(X_{i_1},X_{i_2},Y_j) \\
  & + \frac{4}{nm(m-1)} \sumin \sum\limits_{1 \le j_1 < j_2 \le m} h^{(1,2)}(X_i,Y_{j_1},Y_{j_2}) \\
  & + \frac{4}{n(n-1)m(m-1)} \sum\limits_{1 \le i_1 < i_2 \le n} \sum\limits_{1 \le j_1 < j_2 \le m} h^{(2,2)}(X_{i_1},X_{i_2},Y_{j_1},Y_{j_2})
  \end{align*}
  is the remainder term.
\end{lemma}  
\begin{proof}
  As stated in Proposition \ref{Prop:Enm-k}, $\ce_{n,m}^{k}(X,Y)$ is a U-statistic with a two-sample kernel $h^k(X_1,X_2,Y_1,Y_2)$, then the expressions of $L_{n,m}^{k}(X,Y)$ and $R_{n,m}^{k}(X,Y)$ directly follow from Lemma \ref{Lemma:H-decomp}.
\end{proof}

For $0 \le c,d \le 2$, define
\begin{equation}\label{Def:h_cd}
    h_{c d}^{k}
  = h_{c d}^{k}(X_1,\dots,X_c,Y_1,\dots,Y_d)
  = \blreft{X_{c+1},\dots,X_2,Y_{d+1},\dots,Y_2}{h^{k}(X_1,X_2,Y_1,Y_2)}.
\end{equation}
For simplicity, we  use $h_{cd}$ as the shortened notation of $h_{cd}^{k}$ in the subsequent analysis. In the following lemma, we express $h^{(c,d)}$ in terms of $h_{cd}$.

\begin{lemma}\label{Lemma:h^(c,d)}
Let $X_1, X_2$ be independent copies of $X$ and $Y_1, Y_2$ be independent copies of $Y$, then it holds that
\BEqn
& & h^{(1,0)}(X_1) = h_{10}(X_1) - \ce^{k}(X,Y), \\
& & h^{(0,1)}(Y_1) = h_{01}(Y_1) - \ce^{k}(X,Y), \\
& & h^{(2,0)}(X_1,X_2) = h_{20}(X_1,X_2) - h_{10}(X_1) - h_{10}(X_2) + \ce^{k}(X,Y), \\
& & h^{(1,1)}(X_1,Y_1) = h_{11}(X_1,Y_1) - h_{10}(X_1) - h_{01}(Y_1) + \ce^{k}(X,Y), \\
& & h^{(0,2)}(Y_1,Y_2) = h_{02}(Y_1,Y_2) - h_{01}(Y_1) - h_{01}(Y_2) + \ce^{k}(X,Y),
\EEqn
and
\BEqn
    h^{(2,1)}(X_1,X_2,Y_1) 
&=& h_{21}(X_1,X_2,Y_1) - h_{20}(X_1,X_2) - h_{11}(X_1,Y_1) - h_{11}(X_2,Y_1) \\
& & + h_{10}(X_1) + h_{10}(X_2) + h_{01}(Y_1) - \ce^k(X,Y), \\
    h^{(1,2)}(X_1,Y_1,Y_2) 
&=& h_{12}(X_1,Y_1,Y_2) - h_{02}(Y_1,Y_2) - h_{11}(X_1,Y_1) - h_{11}(X_1,Y_2) \\
& & + h_{10}(X_1) + h_{01}(Y_1) + h_{01}(Y_2) - \ce^k(X,Y), \\
    h^{(2,2)}(X_1,X_2,Y_1,Y_2) 
&=& h^{k}(X_1,X_2,Y_1,Y_2) - h_{21}(X_1,X_2,Y_1) - h_{21}(X_1,X_2,Y_2) \\
& & - h_{12}(X_1,Y_1,Y_2) - h_{12}(X_2,Y_1,Y_2) + h_{20}(X_1,X_2) + h_{02}(Y_1,Y_2) \\
& & + h_{11}(X_1,Y_1) + h_{11}(X_1,Y_2) + h_{11}(X_2,Y_1) + h_{11}(X_2,Y_2) \\
& & - h_{10}(X_1) - h_{10}(X_2) - h_{01}(Y_1) - h_{01}(Y_2) + \ce^k(X,Y).
\EEqn
\end{lemma}
\begin{proof}
The results directly follow from some simple calculations, and we use $h^{(1,1)}(X_1,Y_1)$ as the example to illustrate the steps. Recall the definition of $h^{(c,d)}$ stated in Lemma \ref{Lemma:estMMD_k} and the notation (\ref{Def:h_cd}), then we have that
\BEqn
& & h^{(1,1)}(X_1,Y_1) \\
&=& \int\dots\int h^{k}(u_1,u_2,v_1,v_2) \lrp{dG_{X_1}(u_1) - dF_X(u_1)} dF_X(u_2) \lrp{dG_{Y_1}(v_1) - dF_Y(v_1)} dF_Y(v_2) \\
&=& \int\dots\int h^{k}(X_1,u_2,Y_1,v_2) dF_X(u_2) dF_Y(v_2) \\
& & - \int\dots\int h^{k}(X_1,u_2,v_1,v_2)  dF_X(u_2) dF_Y(v_1) dF_Y(v_2) \\
& & - \int\dots\int h^{k}(u_1,u_2,Y_1,v_2) dF_X(u_1) dF_X(u_2) dF_Y(v_2) \\
& & + \int\dots\int h^{k}(u_1,u_2,v_1,v_2) dF_X(u_1) dF_X(u_2) dF_Y(v_1) dF_Y(v_2) \\
&=& h_{11}(X_1,Y_1) - h_{10}(X_1) - h_{01}(Y_1) + \ce^{k}(X,Y), 
\EEqn
where the last step uses the fact that $\bbe{h^k(X_1,X_2,Y_1,Y_2)} = \ce^k(X,Y)$. All the other statements can be shown in a similar way, and we spare the details.
\end{proof}

To facilitate the following analysis, we further present the explicit expressions of each $h_{cd}$ in terms of the original kernel $k$ in the following lemma. Similar results can be found at Lemma 4.3, Lemma 4.6 and Lemma B.2 in \cite{huang2017efficient}, although only the special case $k(x,y)=|x-y|$ is considered in that paper.

\begin{lemma}\label{Lemma:h_cd}
Let $X_1, X_2, X'$ be independent copies of $X$ and $Y_1, Y_2, Y'$ be independent copies of $Y$. If $\blre{\lrabs{k(X,X')}} + \blre{\lrabs{k(X,Y)}} + \blre{\lrabs{k(Y,Y')}}<\infty$, then for independent random vectors $X$ and $Y$, we have
\BEqn
& & h_{10}(X_1) 
 =  \beft{Y}{k(X_1,Y)} + \blre{k(X,Y)} - \beft{X}{k(X_1,X)} - \blre{k(Y,Y')}, \\
& & h_{01}(Y_1) 
 =  \beft{X}{k(X,Y_1)} + \blre{k(X,Y)} - \blre{k(X,X')} - \beft{Y}{k(Y_1,Y)}, \\
& & h_{20}(X_1,X_2) 
 =  \beft{Y}{k(X_1,Y)} + \beft{Y}{k(X_2,Y)} - k(X_1,X_2) - \blre{k(Y,Y')}, \\
& & h_{02}(Y_1,Y_2) 
 =  \beft{X}{k(X,Y_1)} + \beft{X}{k(X,Y_2)} - k(Y_1,Y_2) - \blre{k(X,X')}, \\
& & h_{11}(X_1,Y_1)
 =  \frac{1}{2} k(X_1,Y_1) + \frac{1}{2} \beft{X}{k(X,Y_1)} + \frac{1}{2} \beft{Y}{k(X_1,Y)} + \frac{1}{2} \blre{k(X,Y)} \\
& & \hspace{6em} - \beft{X}{k(X_1,X)} - \beft{Y}{k(Y_1,Y)},
\EEqn
and
\BEqn
    h_{21}(X_1,X_2,Y_1)
&=& \frac{1}{2}k(X_1,Y_1) + \frac{1}{2}k(X_2,Y_1) + \frac{1}{2}\blreft{Y}{k(X_1,Y)} + \frac{1}{2}\blreft{Y}{k(X_2,Y)} \\
& & - k(X_1,X_2) - \blreft{Y}{k(Y_1,Y)}, \\
    h_{12}(X_1,Y_1,Y_2)
&=& \frac{1}{2}k(X_1,Y_1) + \frac{1}{2}k(X_1,Y_2) + \frac{1}{2}\blreft{X}{k(X,Y_1)} + \frac{1}{2}\blreft{X}{k(X,Y_2)} \\
& & - k(Y_1,Y_2) - \blreft{X}{k(X_1,X)}.
\EEqn
\end{lemma}
\begin{proof}
Recall that $h^k(X_1,X_2,Y_1,Y_2) = \frac{1}{2} \sum_{i=1}^2 \sum_{j=1}^2 k(X_i, Y_j) - k(X_1, X_2) - k(Y_1, Y_2)$ is defined in Proposition \ref{Prop:Enm-k}, then the statements just follow from some simple calculations.
\end{proof}

Next we present some interesting properties of $h^{(2,1)}(X_1,X_2,Y_1)$, $h^{(1,2)}(X_1,Y_1,Y_2)$ and $h^{(2,2)}(X_1,X_2,Y_1,Y_2)$.
\begin{lemma}\label{Lemma:R=0}
Let $X_1, X_2$ be independent copies of $X$ and $Y_1, Y_2$ be independent copies of $Y$. If $\blre{\lrabs{k(X,X')}} + \blre{\lrabs{k(X,Y)}} + \blre{\lrabs{k(Y,Y')}}<\infty$, then it holds that \begin{equation*}
    h^{(2,1)}(X_1,X_2,Y_1) = h^{(1,2)}(X_1,Y_1,Y_2) = h^{(2,2)}(X_1,X_2,Y_1,Y_2) = 0.
\end{equation*}
\end{lemma}
\begin{proof}
We have shown in Lemma \ref{Lemma:h^(c,d)} that
\BEqn
    h^{(2,1)}(X_1,X_2,Y_1) 
&=& h_{21}(X_1,X_2,Y_1) - h_{20}(X_1,X_2) - h_{11}(X_1,Y_1) - h_{11}(X_2,Y_1) \\
& & + h_{10}(X_1) + h_{10}(X_2) + h_{01}(Y_1) - \ce^k(X,Y), 
\EEqn
and the expression of each individual term is presented in Lemma \ref{Lemma:h_cd}, then it follows from some simple calculations that
\BEqn
& & h^{(2,1)}(X_1,X_2,Y_1) \\
&=& \frac{1}{2} \lrp{k(X_1,Y_1) + k(X_2,Y_1) + \blreft{Y}{k(X_1,Y)} + \blreft{Y}{k(X_2,Y)} - 2k(X_1,X_2) - 2\blreft{Y}{k(Y_1,Y)}} \\
& & - \lrp{\blreft{Y}{k(X_1,Y)} + \blreft{Y}{k(X_2,Y)} - k(X_1,X_2) - \blre{k(Y,Y')}} \\
& & - \frac{1}{2} \lrp{k(X_1,Y_1) + \blre{k(X_1,Y_1)} + \blreft{Y}{k(X_1,Y)} + \blreft{X}{k(X,Y_1)} - 2\blreft{X}{k(X_1,X)} - 2\blreft{Y}{k(Y_1,Y)}} \\
& & - \frac{1}{2} \lrp{k(X_2,Y_1) + \blre{k(X_1,Y_1)} + \blreft{Y}{k(X_2,Y)} + \blreft{X}{k(X,Y_1)} - 2\blreft{X}{k(X_2,X)} - 2\blreft{Y}{k(Y_1,Y)}} \\
& & + \lrp{\blreft{Y}{k(X_1,Y)} + \blre{k(X_1,Y_1)} - \blreft{X}{k(X_1,X)} - \blre{k(Y_1,Y_2)}} \\
& & + \lrp{\blreft{Y}{k(X_2,Y)} + \blre{k(X_1,Y_1)} - \blreft{X}{k(X_2,X)} - \blre{k(Y_1,Y_2)}} \\
& & + \lrp{\blreft{X}{k(X,Y_1)} + \blre{k(X_1,Y_1)} - \blre{k(X_1,X_2)} - \blreft{Y}{k(Y_1,Y)}} \\
& & - \lrp{2\blre{k(X_1,Y_1)} - \blre{k(X_1,X_2)} - \blre{k(Y_1,Y_2)}} \\
&=& 0,
\EEqn
and similarly, we obtain that $h^{(1,2)}(X_1,Y_1,Y_2) = 0$. Also, it further implies that 
\BEqn
& & h_{21}(X_1,X_2,Y_1) \\
&=& h_{20}(X_1,X_2) + h_{11}(X_1,Y_1) + h_{11}(X_2,Y_1) - h_{10}(X_1) - h_{10}(X_2) - h_{01}(Y_1) + \ce^k(X,Y), \\[2mm]
& & h_{21}(X_1,X_2,Y_2) \\
&=& h_{20}(X_1,X_2) + h_{11}(X_1,Y_2) + h_{11}(X_2,Y_2) - h_{10}(X_1) - h_{10}(X_2) - h_{01}(Y_2) + \ce^k(X,Y), \\[2mm]
& & h_{12}(X_1,Y_1,Y_2) \\
&=& h_{02}(Y_1,Y_2) + h_{11}(X_1,Y_1) + h_{11}(X_1,Y_2) - h_{10}(X_1) - h_{01}(Y_1) - h_{01}(Y_2) + \ce^k(X,Y), \\[2mm]
& & h_{12}(X_2,Y_1,Y_2) \\ 
&=& h_{02}(Y_1,Y_2) + h_{11}(X_2,Y_1) + h_{11}(X_2,Y_2) - h_{10}(X_2) - h_{01}(Y_1) - h_{01}(Y_2) + \ce^k(X,Y).
\EEqn
Using the expression of $h^{(2,2)}(X_1,X_2,Y_1,Y_2)$ derived in Lemma \ref{Lemma:h^(c,d)}, we have
\BEqn
& & h^{(2,2)}(X_1,X_2,Y_1,Y_2) \\
&=& h^{k}(X_1,X_2,Y_1,Y_2) - h_{21}(X_1,X_2,Y_1) - h_{21}(X_1,X_2,Y_2) - h_{12}(X_1,Y_1,Y_2) - h_{12}(X_2,Y_1,Y_2) \\
& & + h_{20}(X_1,X_2) + h_{02}(Y_1,Y_2) + h_{11}(X_1,Y_1) + h_{11}(X_1,Y_2) + h_{11}(X_2,Y_1) + h_{11}(X_2,Y_2) \\
& & - h_{10}(X_1) - h_{10}(X_2) - h_{01}(Y_1) - h_{01}(Y_2) + \ce^k(X,Y) \\
&=& h^{k}(X_1,X_2,Y_1,Y_2) \\
& & - \lrp{h_{20}(X_1,X_2) + h_{11}(X_1,Y_1) + h_{11}(X_2,Y_1) - h_{10}(X_1) - h_{10}(X_2) - h_{01}(Y_1) + \ce^k(X,Y)} \\
& & - \lrp{h_{20}(X_1,X_2) + h_{11}(X_1,Y_2) + h_{11}(X_2,Y_2) - h_{10}(X_1) - h_{10}(X_2) - h_{01}(Y_2) + \ce^k(X,Y)} \\
& & - \lrp{h_{02}(Y_1,Y_2) + h_{11}(X_1,Y_1) + h_{11}(X_1,Y_2) - h_{10}(X_1) - h_{01}(Y_1) - h_{01}(Y_2) + \ce^k(X,Y)} \\
& & - \lrp{h_{02}(Y_1,Y_2) + h_{11}(X_2,Y_1) + h_{11}(X_2,Y_2) - h_{10}(X_2) - h_{01}(Y_1) - h_{01}(Y_2) + \ce^k(X,Y)} \\
& & + h_{20}(X_1,X_2) + h_{02}(Y_1,Y_2) + h_{11}(X_1,Y_1) + h_{11}(X_1,Y_2) + h_{11}(X_2,Y_1) + h_{11}(X_2,Y_2) \\
& & - h_{10}(X_1) - h_{10}(X_2) - h_{01}(Y_1) - h_{01}(Y_2) + \ce^k(X,Y) \\
&=& h^{k}(X_1,X_2,Y_1,Y_2) - h_{11}(X_1,Y_1) - h_{11}(X_1,Y_2) - h_{11}(X_2,Y_1) - h_{11}(X_2,Y_2) \\
& & - h_{20}(X_1,X_2) - h_{02}(Y_1,Y_2) + 2\lrp{h_{10}(X_1) + h_{10}(X_2) + h_{01}(Y_1) + h_{01}(Y_2)} - 3\ce^k(X,Y) \\
&=& 0,
\EEqn
where the last steps follows from the expressions derived in Lemma \ref{Lemma:h_cd} and thus completes the proof.
\end{proof}

\subsection{Lemmas for Proposition \ref{Prop:Var(Lnm)}}

\begin{lemma}\label{Lemma:psd}
  If $\blre{\lrabs{k(X_1,X_2)}} < \infty$, then $-d^k(X_1,X_2)$ is a positive definite kernel.
\end{lemma}
\begin{proof}
  If $k$ is a semimetric of strong negative type, then it is natural that $k$ is also a semimetric of negative type, and the statement can thus be proved using the arguments for Lemma 4.13 of \cite{huang2017statistically}. Then it suffices to consider the case when $(\cf,-k)$ is characteristic on $(\br^{p},\cp)$. By \cite{sejdinovic2013equivalence} Corollary 16, 
  \begin{equation*}
      \rho(x_1,x_2) := - k(x_1,x_1) - k(x_2,x_2) + 2k(x_1,x_2)
  \end{equation*}
  defines a valid semimetric of negative type. Using the results above, $-\rho(X_1,X_2) + \bbeft{X_1}{\rho(X_1,X_2)} + \bbeft{X_2}{\rho(X_1,X_2)} - \bbe{\rho(X_1,X_2)}$ is a positive definite kernel, and the proof is completed by noting that
  \BEqn
 & & -\rho(X_1,X_2) + \bbeft{X_1}{\rho(X_1,X_2)} + \bbeft{X_2}{\rho(X_1,X_2)} - \bbe{\rho(X_1,X_2)} \\
  &=& k(X_1,X_1) + k(X_2,X_2) - 2k(X_1,X_2) \\
 & & - \bbe{k(X_1,X_1)} - k(X_2,X_2) + 2\bbeft{X_1}{k(X_1,X_2)} \\
 & & - k(X_1,X_1) - \bbe{k(X_2,X_2)} + 2\bbeft{X_2}{k(X_1,X_2)} \\
 & & + \bbe{k(X_1,X_1)}  + \bbe{k(X_2,X_2)} - 2\bbe{k(X_1,X_2)} \\
  &=& -2\lrp{k(X_1,X_2) - \bbeft{X_1}{k(X_1,X_2)} - \bbeft{X_2}{k(X_1,X_2)} + \bbe{k(X_1,X_2)}} \\
  &=& -2d^k(X_1,X_2).
  \EEqn
\end{proof}

\begin{lemma}\label{Lemma:phi}
  Assume that $X$ and $Y$ are identically distributed and $\bbe{k^2(X_1,X_2)} < \infty$, then there exist functions $\phi_1(\cdot), \phi_2(\cdot), \dots$ such that
  \begin{align*}
    h_{20}(X_1,X_2) &= \sum\limits_{\ell=1}^{\infty}\lambda_{\ell}\phi_{\ell}(X_1)\phi_{\ell}(X_2), \\
    h_{02}(Y_1,Y_2) &= \sum\limits_{\ell=1}^{\infty}\lambda_{\ell}\phi_{\ell}(Y_1)\phi_{\ell}(Y_2), \\
    h_{11}(X_1,Y_1) &= -\frac{1}{2}\sum\limits_{\ell=1}^{\infty}\lambda_{\ell}\phi_{\ell}(X_1)\phi_{\ell}(Y_1),
  \end{align*}
  where $\be[\phi_{\ell}(X)]=0$, $\be[\phi_{\ell}(X)^2]=1$ and $\be[\phi_{\ell}(X)\phi_{k}(X)]=0, \ell=1,2,\dots,\infty, \ell \neq k$, and $\lambda_1 \geq \lambda_2 \geq \dots \geq0$, such that
  \begin{align*}
    \sum\limits_{\ell=1}^{\infty}\lambda_{\ell} = \bbe{k(X,X')} - \bbe{k(X,X)}, \quad\quad
    \sum\limits_{\ell=1}^{\infty}\lambda_{\ell}^2 = \cv_k^2(X),
  \end{align*}
  where $\cv_k^2(X)$ is the HSIC of $X$ with itself.
\end{lemma}
\begin{proof}
   If $X$ and $Y$ are identically distributed and $\bbe{k^2(X_1,X_2)} < \infty$, then it follows from Lemma \ref{Lemma:h_cd} and Lemma \ref{Lemma:psd} that $h_{20}(\cdot,\cdot)$, $h_{02}(\cdot,\cdot)$ and $-h_{11}(\cdot,\cdot)$ are all positive definite kernels, hence the existence of $\phi_\ell$ and their basic properties directly follow Mercer's Theorem. Consequently,
   \begin{equation*}
       \sum\limits_{\ell=1}^{\infty} \lambda_{\ell} 
     = \blreft{X}{h_{20}(X,X)} 
     = \blre{k(X,X')} - \blre{k(X,X)},
   \end{equation*}
   and
   \begin{equation*}
       \sum\limits_{\ell=1}^{\infty} \lambda_{\ell}^2
     = \sum\limits_{\ell,k=1}^{\infty} \lambda_{\ell}\lambda_k \blre{\phi_\ell(X_1) \phi_k(X_1)} \blre{\phi_\ell(X_2) \phi_k(X_2)}
     = \blre{h_{20}^2(X_1,X_2)}
     = \cv_k^2(X).
   \end{equation*}
\end{proof}

 \subsection{Lemmas for Proposition \ref{Prop:estDC_k}}

\begin{lemma}\label{Lemma:Combn_Kernel}
  For any fixed $k=k^{(p)}\in\cc$ and for any $a_1,a_2,a_3,a_4\in\br^p$, define
  \begin{eqnarray}\label{Equ:psi}
  & & \psi(a_1,a_2,a_3,a_4) \\
  &=& \frac{1}{2}\sum\limits_{1\le i<j \le4} k^2(a_i,a_j)
    - \frac{1}{4}\sum\limits_{i=1}^{4}\lrp{\sum\limits_{j \neq i} k(a_i,a_j)}^2
    + \frac{1}{6}\lrp{\sum\limits_{1\le i<j \le4} k(a_i,a_j)}^2. \nonumber
  \end{eqnarray}
  Assume that $a_0^k = k(X,X)$ is a constant independent of $X$, then $\cv_{n,m}^{k\ast}(X,Y)$ admits the following representation:
  \begin{align*}
    & \cv_{n,m}^{k\ast}(X,Y) \\
   =& \bin{N}{4}^{-1} 
      \left( 
        \sum\limits_{1\le i_1<i_2<i_3<i_4 \le n} \psi(X_{i_1},X_{i_2},X_{i_3},X_{i_4})
      + \sum\limits_{1\le i_1<i_2<i_3 \le n}\sum\limits_{j=1}^{m} \psi(X_{i_1},X_{i_2},X_{i_3},Y_{j})
      \right. \\
   & \hspace{5em} 
      + \sum\limits_{1\le i_1<i_2 \le n}\sum\limits_{1\le j_1<j_2 \le m} \psi(X_{i_1},X_{i_2},Y_{j_1},Y_{j_2})
      + \sum\limits_{i=1}^{n}\sum\limits_{1\le j_1<j_2<j_3 \le m} \psi(X_{i},Y_{j_1},Y_{j_2},Y_{j_3}) \\
   & \hspace{5em}
      \left.
      + \sum\limits_{1\le j_1<j_2<j_3<j_4 \le m} \psi(Y_{j_1},Y_{j_2},Y_{j_3},Y_{j_4}) 
      \right).
  \end{align*}
\end{lemma}
\begin{proof}
  Define 
  \BEqn
  & & n(r) = \lrcp{(i_1,\dots,i_r):1\le i_1,\dots,i_r\le n \text{ pairwise distinct}}, \\
  & & m(r) = \lrcp{(j_1,\dots,j_r):1\le j_1,\dots,j_r\le m \text{ pairwise distinct}}.
  \EEqn
  Recall the notations introduced in Proposition \ref{Prop:estDC_k}, we obtain the following equations,
  \begin{equation*}
      \sum\limits_{s=1}^{N} a_{s,t}^{k} = (N-2) \tilde{a}_{\cdot t}^{k}, \quad
      \sum\limits_{t=1}^{N} a_{s,t}^{k} = (N-2) \tilde{a}_{s \cdot}^{k}, \quad
      \sum\limits_{s,t=1}^{N} a_{s,t}^{k} = (N-1)(N-2) \tilde{a}_{\cdot \cdot}^{k},
  \end{equation*}
  and $\sum\limits_{s=1}^{N} \tilde{a}_{s \cdot}^{k} = \sum\limits_{t=1}^{N} \tilde{a}_{\cdot t}^{k} = (N-1) \tilde{a}_{\cdot \cdot}^{k}$. Also by noting the symmetry of $a_{s,t}^{k}$, we further have $\tilde{a}_{s \cdot}^{k} = \tilde{a}_{\cdot s}^{k}$ and thus $(\tilde{a}_{s \cdot}^{k})^2 = \tilde{a}_{s \cdot}^{k} \tilde{a}_{\cdot s}^{k} = (\tilde{a}_{\cdot s}^{k})^2$. Then by expanding $\lrp{A_{st}^{k\ast}}^2$ in terms of $k(X_{i_1},X_{i_2})$, $k(X_i,Y_j)$ and $k(Y_{j_1},Y_{j_2})$, we have
  \begin{align*}
    & \sum\limits_{1\le s\neq t\le N}\lrp{A_{st}^{k\ast}}^2 \\
   =& \hspace{1em}
        \sum\limits_{1\le s\neq t\le N} (a_{s,t}^{k})^2 
      - 2(N-2) \sum\limits_{s=1}^{N}(\tilde{a}_{s \cdot}^{k})^2
      + (N-1)(N-2) (\tilde{a}_{\cdot \cdot}^{k})^2
      + 2(N-2) a_0^{k} \tilde{a}_{\cdot \cdot}^{k} \\
   =& \hspace{1em}
      \frac{N-3}{N-1} 
      \big(
        \sum\limits_{(i_1,i_2) \in n(2)} k^2(X_{i_1},X_{i_2})
        + \sum\limits_{(j_1,j_2) \in m(2)} k^2(Y_{j_1},Y_{j_2}) 
        + 2\sum\limits_{i=1}^{n}\sum\limits_{j=1}^{m} k^2(X_{i},Y_{j})
      \big) \\
   & -\frac{2(N-3)}{(N-1)(N-2)}
      \big(
        \sum\limits_{(i_1,i_2,i_3) \in n(3)} k(X_{i_1},X_{i_2})  k(X_{i_1},X_{i_3})
        + \sum\limits_{(j_1,j_2,j_3) \in m(3)} k(Y_{j_1},Y_{j_2})  k(Y_{j_1},Y_{j_3})
      \big) \\
   & + \frac{1}{(N-1)(N-2)} 
      \big(
        \sum\limits_{(i_1,i_2,i_3,i_4) \in n(4)} k(X_{i_1},X_{i_2}) k(X_{i_3},X_{i_4})
        + 4\sum\limits_{(i_1,i_2,i_3) \in n(3)} \sum\limits_{j=1}^{m} k(X_{i_1},X_{i_2})  k(X_{i_3},Y_{j})  \\
   & \hspace{8em}
        + \sum\limits_{(j_1,j_2,j_3,j_4) \in m(4)} k(Y_{j_1},Y_{j_2}) k(Y_{j_3},Y_{j_4})
        + 4\sum\limits_{(j_1,j_2,j_3) \in m(3)} \sum\limits_{i=1}^{n} k(X_{i},Y_{j_1})  k(Y_{j_2},Y_{j_3})
      \big) \\
   & - \frac{2(N-3)}{(N-1)(N-2)} 
      \big(
        2\sum\limits_{(i_1,i_2) \in n(2)} \sum\limits_{j=1}^{m} k(X_{i_1},X_{i_2}) k(X_{i_1},Y_{j}) 
        + \sum\limits_{(i_1,i_2) \in n(2)} \sum\limits_{j=1}^{m} k(X_{i_1},Y_{j}) k(X_{i_2},Y_{j}) \\
   & \hspace{9em}    
        + 2\sum\limits_{(j_1,j_2) \in m(2)} \sum\limits_{i=1}^{n} k(X_{i},Y_{j_1}) k(Y_{j_1},Y_{j_2}) 
        + \sum\limits_{(j_1,j_2) \in m(2)} \sum\limits_{i=1}^{n} k(X_{i},Y_{j_1}) k(X_{i},Y_{j_2})
      \big) \\
   & + \frac{2}{(N-1)(N-2)}
      \big(
        2 \sum\limits_{(i_1,i_2) \in n(2)} \sum\limits_{(j_1,j_2) \in m(2)}k(X_{i_1},Y_{j_1}) k(X_{i_2},Y_{j_2}) \\
   & \hspace{9em}        
        + \sum\limits_{(i_1,i_2) \in n(2)} \sum\limits_{(j_1,j_2) \in m(2)}k(X_{i_1},X_{i_2}) k(Y_{j_1},Y_{j_2}) 
      \big) \\
   & + \frac{N}{N-1}\lrp{a_0^k}^2
  \end{align*}
  On the other hand, it follows from the definition of $\psi$ that
  \BEqn
 & & \sum\limits_{1\le i_1<i_2<i_3<i_4 \le n}  \psi(X_{i_1},X_{i_2},X_{i_3},X_{i_4}) \\
  &=& \frac{(n-2)(n-3)}{24} \sum\limits_{(i_1,i_2) \in n(2)} k^2(X_{i_1},X_{i_2}) 
      - \frac{n-3}{12} \sum\limits_{(i_1,i_2,i_3) \in n(3)} k(X_{i_1},X_{i_2}) k(X_{i_1},X_{i_3}) \\
 & & + \frac{1}{24} \sum\limits_{(i_1,i_2,i_3,i_4) \in n(4)} k(X_{i_1},X_{i_2}) k(X_{i_3},X_{i_4}) 
  \\ \\ \\
 & & \sum\limits_{1\le i_1<i_2<i_3 \le n}\sum\limits_{j=1}^{m} \psi(X_{i_1},X_{i_2},X_{i_3},Y_{j}) \\
  &=& \frac{(n-2)m}{12} \sum\limits_{(i_1,i_2) \in n(2)} k^2(X_{i_1},X_{i_2}) 
      + \frac{(n-1)(n-2)}{12} \sum\limits_{i=1}^{n}\sum\limits_{j=1}^{m} k^2(X_i,Y_j) \\
 & & - \frac{m}{12} \sum\limits_{(i_1,i_2,i_3) \in n(3)} k(X_{i_1},X_{i_2}) k(X_{i_1},X_{i_3}) 
      -\frac{n-2}{6}\sum\limits_{(i_1,i_2) \in n(2)}\sum\limits_{j=1}^{m} k(X_{i_1},X_{i_2}) k(X_{i_1},Y_{j}) \\
 & & - \frac{n-2}{12} \sum\limits_{(i_1,i_2) \in n(2)}\sum\limits_{j=1}^{m} k(X_{i_1},Y_{j}) k(X_{i_2},Y_{j})
      + \frac{1}{6} \sum\limits_{(i_1,i_2,i_3) \in n(3)}\sum\limits_{j=1}^{m} k(X_{i_1},X_{i_2}) k(X_{i_3},Y_{j})
  \\ \\ \\
 & & \sum\limits_{1\le i_1<i_2 \le n}\sum\limits_{1\le j_1<j_2 \le m} \psi(X_{i_1},X_{i_2},Y_{j_1},Y_{j_2}) \\
  &=& \frac{m(m-1)}{24} \sum\limits_{(i_1,i_2) \in n(2)} k^2(X_{i_1},X_{i_2})
      + \frac{(n-1)(m-1)}{6} \sum\limits_{i=1}^{n}\sum\limits_{j=1}^{m} k^2(X_{i},Y_{j}) \\
 & & + \frac{n(n-1)}{24} \sum\limits_{(j_1,j_2) \in m(2)} k^2(Y_{j_1},Y_{j_2}) 
      - \frac{m-1}{6} \sum\limits_{(i_1,i_2) \in n(2)}\sum\limits_{j=1}^{m} k(X_{i_1},X_{i_2}) k(X_{i_1},Y_{j}) \\
 & & - \frac{n-1}{6} \sum\limits_{(j_1,j_2) \in m(2)}\sum\limits_{i=1}^{n} k(X_{i},Y_{j_1}) k(Y_{j_1},Y_{j_2}) 
      - \frac{m-1}{12} \sum\limits_{(i_1,i_2) \in n(2)}\sum\limits_{j=1}^{m} k(X_{i_1},Y_{j}) k(X_{i_2},Y_{j}) \\
 & & - \frac{n-1}{12} \sum\limits_{(j_1,j_2) \in m(2)}\sum\limits_{i=1}^{n} k(X_{i},Y_{j_1}) k(X_{i},Y_{j_2})
      + \frac{1}{12} \sum\limits_{(i_1,i_2) \in n(2)} \sum\limits_{(j_1,j_2) \in m(2)} k(X_{i_1},X_{i_2}) k(Y_{j_1},Y_{j_2}) \\
 & & + \frac{1}{6} \sum\limits_{(i_1,i_2) \in n(2)} \sum\limits_{(j_1,j_2) \in m(2)} k(X_{i_1},Y_{j_1}) k(X_{i_2},Y_{j_2})
  \\ \\ \\
 & & \sum\limits_{i=1}^{n}\sum\limits_{1\le j_1<j_2<j_3 \le m} \psi(X_{i},Y_{j_1},Y_{j_2},Y_{j_3}) \\ 
  &=& \frac{n(m-2)}{12} \sum\limits_{(j_1,j_@) \in m(2)} k^2(Y_{j_1},Y_{j_2}) 
      + \frac{(m-1)(m-2)}{12} \sum\limits_{1}^{n} \sum\limits_{j=1}^{m} k^2(X_{i},Y_{j}) \\
 & & - \frac{n}{12} \sum\limits_{(j_1,j_2,j_3) \in m(3)} k(Y_{j_1},Y_{j_2}) k(Y_{j_1},Y_{j_3})
      - \frac{m-2}{6} \sum\limits_{(j_1,j_2) \in m(2)}\sum\limits_{i=1}^{n} k(X_{i},Y_{j_1}) k(Y_{j_1},Y_{j_2}) \\
 & & - \frac{m-2}{12} \sum\limits_{(j_1,j_2) \in m(2)}\sum\limits_{i=1}^{n} k(X_{i},Y_{j_1}) k(X_{i},Y_{j_2})
      + \frac{1}{6} \sum\limits_{(j_1,j_2,j_3) \in m(3)} \sum\limits_{i=1}^{n} k(X_{i},Y_{j_1}) k(Y_{j_2},Y_{j_3}) 
  \\ \\ \\
 & & \sum\limits_{1\le j_1<j_2<j_3<j_4 \le m} \psi(Y_{j_1},Y_{j_2},Y_{j_3},Y_{j_4}) \\
  &=& \frac{(m-2)(m-3)}{24} \sum\limits_{(j_1,j_2) \in m(2)} k^2(Y_{j_1},Y_{j_2}) 
      - \frac{m-3}{12} \sum\limits_{(j_1,j_2,j_3) \in m(3)} k(Y_{j_1},Y_{j_2}) k(Y_{j_1},Y_{j_3}) \\
 & & + \frac{1}{24} \sum\limits_{(j_1,j_2,j_3,j_4) \in m(4)} k(Y_{j_1},Y_{j_2}) k(Y_{j_3},Y_{j_4}) 
  \EEqn
  By putting these terms together and noting that $N=n+m$, we thus have
  \BEqn
 & & \cv_{n,m}^{k\ast}(X,Y) \\
  &=& \frac{1}{N(N-3)} \sum\limits_{1\le s\neq t \le N} \lrp{A_{st}^{k\ast}}^2 
   - \frac{1}{(N-1)(N-3)}\lrp{a_0^k}^2 \\ 
  &=& \bin{N}{4}^{-1} 
      \left( 
        \sum\limits_{1\le i_1<i_2<i_3<i_4 \le n} \psi(X_{i_1},X_{i_2},X_{i_3},X_{i_4})
      + \sum\limits_{1\le i_1<i_2<i_3 \le n}\sum\limits_{j=1}^{m} \psi(X_{i_1},X_{i_2},X_{i_3},Y_{j})
      \right. \\
 & & \hspace{5em} 
      + \sum\limits_{1\le i_1<i_2 \le n}\sum\limits_{1\le j_1<j_2 \le m} \psi(X_{i_1},X_{i_2},Y_{j_1},Y_{j_2})
      + \sum\limits_{i=1}^{n}\sum\limits_{1\le j_1<j_2<j_3 \le m} \psi(X_{i},Y_{j_1},Y_{j_2},Y_{j_3}) \\
 & & \hspace{5em}
      \left.
      + \sum\limits_{1\le j_1<j_2<j_3<j_4 \le m} \psi(Y_{j_1},Y_{j_2},Y_{j_3},Y_{j_4}) 
      \right).
  \EEqn
\end{proof}

\begin{lemma}\label{Lemma:ratio-const-k-H0}
  Under the null, for any fixed $k=k^{(p)}\in\cc$ it holds that $\blre{\cv_{n,m}^{k\ast}(X,Y)} = \cv_k^2(Z)$.
\end{lemma}
\begin{proof}
  Under the null, we have $X =^d Y =^d Z$, and it follows from Lemma \ref{Lemma:Combn_Kernel} that
  \BEqn
  & & \bin{N}{4} \blre{\cv_{n,m}^{k\ast}(X,Y)} \\
  &=& \sum\limits_{1\le i_1<i_2<i_3<i_4 \le n} \blre{\psi(X_{i_1},X_{i_2},X_{i_3},X_{i_4})}
      + \sum\limits_{1\le i_1<i_2<i_3 \le n}\sum\limits_{j=1}^{m} \blre{\psi(X_{i_1},X_{i_2},X_{i_3},Y_{j})} \\
  & & + \sum\limits_{1\le i_1<i_2 \le n}\sum\limits_{1\le j_1<j_2 \le m} \blre{\psi(X_{i_1},X_{i_2},Y_{j_1},Y_{j_2})}
      + \sum\limits_{i=1}^{n}\sum\limits_{1\le j_1<j_2<j_3 \le m} \blre{\psi(X_{i},Y_{j_1},Y_{j_2},Y_{j_3})} \\
  & & + \sum\limits_{1\le j_1<j_2<j_3<j_4 \le m} \blre{\psi(Y_{j_1},Y_{j_2},Y_{j_3},Y_{j_4})} \\
  &=& \sum\limits_{1\le i_1<i_2<i_3<i_4 \le n} \blre{\psi(Z_{i_1},Z_{i_2},Z_{i_3},Z_{i_4})}
      + \sum\limits_{1\le i_1<i_2<i_3 \le n}\sum\limits_{j=1}^{m} \blre{\psi(Z_{i_1},Z_{i_2},Z_{i_3},Z_{j})} \\
  & & + \sum\limits_{1\le i_1<i_2 \le n}\sum\limits_{1\le j_1<j_2 \le m} \blre{\psi(Z_{i_1},Z_{i_2},Z_{j_1},Z_{j_2})}
      + \sum\limits_{i=1}^{n}\sum\limits_{1\le j_1<j_2<j_3 \le m} \blre{\psi(Z_{i},Z_{j_1},Z_{j_2},Z_{j_3})} \\
  & & + \sum\limits_{1\le j_1<j_2<j_3<j_4 \le m} \blre{\psi(Z_{j_1},Z_{j_2},Z_{j_3},Z_{j_4})} \\
  &=& \lrp{ \bin{n}{4} + \bin{n}{3}\bin{m}{1} + \bin{n}{2}\bin{m}{2} + \bin{n}{1}\bin{m}{3} + \bin{m}{4} } \blre{\psi(Z_1,Z_2,Z_3,Z_4)} \\
  &=& \bin{N}{4}\blre{\psi(Z_1,Z_2,Z_3,Z_4)}.
  \EEqn
  Finally, by noting that
  \BEqn
      \blre{\psi(Z_1,Z_2,Z_3,Z_4)}
  &=& \blre{k^2(Z_1,Z_2)} - 2\blre{k(Z_1,Z_2) k(Z_1,Z_3)} + \blre{k(Z_1,Z_2)}^2 \\
  &=& \cv_k^2(Z),
  \EEqn
  we can complete the proof. 
\end{proof}

\begin{lemma}\label{Lemma:ratio-const-k-H1}
  If $n/N \rightarrow \rho$ with $0 < \rho < 1$, then for any fixed $k=k^{(p)}\in\cc$ it holds under the alternative that  $\blre{\cv_{n,m}^{k\ast}(X,Y)} \rightarrow \cv_k^2(Z)$ as $n,m \rightarrow \infty$.
\end{lemma}
\begin{proof}
  We first compute $\blre{\psi(X_1,X_2,X_3,X_4)}$, $\blre{\psi(X_1,X_2,X_3,Y_1)}$, $\blre{\psi(X_1,X_2,Y_1,Y_2)}$, $\blre{\psi(X_1,Y_1,Y_2,Y_3)}$ and $\blre{\psi(Y_1,Y_2,Y_3,Y_4)}$ respectively. It simply follows from (\ref{Equ:psi}) that
  \BEqn
      \blre{\psi(X_1,X_2,X_3,X_4)}
  &=& \blre{k^2(X_1,X_2)} - 2\blre{k(X_1,X_2) k(X_1,X_3)} + \blre{k(X_1,X_2)}^2, \\
      \blre{\psi(Y_1,Y_2,Y_3,Y_4)}
  &=& \blre{k^2(Y_1,Y_2)} - 2\blre{k(Y_1,Y_2) k(Y_1,Y_3)} + \blre{k(Y_1,Y_2)}^2.
  \EEqn
  Also, it can be computed that
  \BEqn
  & & \blre{\psi(X_1,X_2,X_3,Y_1)} \\
  &=& \frac{1}{2}\blre{k^2(X_1,X_2)} + \frac{1}{2}\blre{k^2(X_1,Y_1)} + \blre{k(X_1,X_2)}\blre{k(X_1,Y_1)}  \\
  & & - \frac{1}{2}\blre{k(X_1,X_2)k(X_1,X_3)} - \blre{k(X_1,X_2)k(X_1,Y_1)} - \frac{1}{2}\blre{k(X_1,Y_1)k(X_2,Y_1)}, \\[3mm]
  & & \blre{\psi(X_1,Y_1,Y_2,Y_3)} \\
  &=& \frac{1}{2}\blre{k^2(X_1,Y_1)} + \frac{1}{2}\blre{k^2(Y_1,Y_2)} + \blre{k(X_1,Y_1)}\blre{k(Y_1,Y_2)}  \\
  & & - \frac{1}{2}\blre{k(X_1,Y_1)k(X_1,Y_2)} - \blre{k(X_1,Y_1)k(Y_1,Y_2)} - \frac{1}{2}\blre{k(Y_1,Y_2)k(Y_1,Y_3)}.
  \EEqn
  Finally, we can also obtain that
  \BEqn
      \blre{\psi(X_1,X_2,Y_1,Y_2)} 
  &=& \frac{1}{6}\blre{k^2(X_1,X_2)} + \frac{2}{3}\blre{k^2(X_1,Y_1)} + \frac{1}{6}\blre{k^2(Y_1,Y_2)} \\
  & & - \frac{2}{3}\blre{k(X_1,X_2)k(X_1,Y_1)} - \frac{2}{3}\blre{k(X_1,Y_1)k(Y_1,Y_2)} \\
  & & - \frac{1}{3}\blre{k(X_1,Y_1)k(X_1,Y_2)} - \frac{1}{3}\blre{k(X_1,Y_1)k(X_2,Y_1)} \\
  & & + \frac{1}{3}\blre{k(X_1,X_2)}\blre{k(Y_1,Y_2)} + \frac{2}{3}\blre{k(X_1,Y_1)}^2.
  \EEqn
  Note that
  \BEqn
      \blre{\cv_{n,m}^{k\ast}(X,Y)} 
  &=& \frac{n(n-1)(n-2)(n-3)}{N(N-1)(N-2)(N-3)} \blre{\psi(X_1,X_2,X_3,X_4)} \\
  & & + \frac{4n(n-1)(n-2)m}{N(N-1)(N-2)(N-3)} \blre{\psi(X_1,X_2,X_3,Y_1)} \\
  & & + \frac{6n(n-1)m(m-1)}{N(N-1)(N-2)(N-3)} \blre{\psi(X_1,X_2,Y_1,Y_2)} \\
  & & + \frac{4nm(m-1)(m-2)}{N(N-1)(N-2)(N-3)} \blre{\psi(X_1,Y_1,Y_2,Y_3)} \\
  & & + \frac{m(m-1)(m-2)(m-3)}{N(N-1)(N-2)(N-3)} \blre{\psi(Y_1,Y_2,Y_3,Y_4)},
  \EEqn
  then after some computations we can obtain that
  \BEqn
  & & \blre{\cv_{n,m}^{k\ast}(X,Y)} \\
  &=& \frac{n(n-1)}{N(N-1)} \blre{k^2(X_1,X_2)} 
   +  \frac{2nm}{N(N-1)} \blre{k^2(X_1,Y_1)} 
   +  \frac{m(m-1)}{N(N-1)} \blre{k^2(Y_1,Y_2)} \\
  & & - \frac{2n(n-1)(n-2)}{N(N-1)(N-2)} \blre{k(X_1,X_2)k(X_1,X_3)}
   -  \frac{4n(n-1)m}{N(N-1)(N-2)} \blre{k(X_1,X_2)k(X_1,Y_1)} \\
  & & - \frac{2n(n-1)m}{N(N-1)(N-2)} \blre{k(X_1,Y_1)k(X_2,Y_1)}
   -  \frac{2nm(m-1)}{N(N-1)(N-2)} \blre{k(X_1,Y_1)k(X_1,Y_2)} \\
  & & - \frac{4nm(m-1)}{N(N-1)(N-2)} \blre{k(X_1,Y_1)k(Y_1,Y_2)}
   -  \frac{2m(m-1)(m-2)}{N(N-1)(N-2)} \blre{k(Y_1,Y_2)k(Y_1,Y_3)} \\
  & & + \frac{n(n-1)(n-2)(n-3)}{N(N-1)(N-2)(N-3)} \blre{k(X_1,X_2)}^2 
   +  \frac{m(m-1)(m-2)(m-3)}{N(N-1)(N-2)(N-3)} \blre{k(Y_1,Y_2)}^2 \\
  & & + \frac{4n(n-1)m(m-1)}{N(N-1)(N-2)(N-3)} \blre{k(X_1,Y_1)}^2 \\
  & & + \frac{2n(n-1)m(m-1)}{N(N-1)(N-2)(N-3)} \blre{k(X_1,X_2)}\blre{k(Y_1,Y_2)} \\
  & & + \frac{4n(n-1)(n-2)m}{N(N-1)(N-2)(N-3)} \blre{k(X_1,X_2)}\blre{k(X_1,Y_1)} \\
  & &  + \frac{4nm(m-1)(m-2)}{N(N-1)(N-2)(N-3)} \blre{k(X_1,Y_1)}\blre{k(Y_1,Y_2)}
  \EEqn
  Under the assumption that $n/N \rightarrow \rho$, it follows that
  \BEqn
  & & \blre{\cv_{n,m}^{k\ast}(X,Y)} \\
  &\rightarrow& \rho^2 \blre{k^2(X_1,X_2)} 
   +  2\rho(1-\rho) \blre{k^2(X_1,Y_1)} 
   +  (1-\rho)^2 \blre{k^2(Y_1,Y_2)} \\
  & & - 2\rho^3 \blre{k(X_1,X_2)k(X_1,X_3)}
   -  4\rho^2(1-\rho) \blre{k(X_1,X_2)k(X_1,Y_1)} \\
  & & - 2\rho^2(1-\rho) \blre{k(X_1,Y_1)k(X_2,Y_1)}
   -  2\rho(1-\rho)^2 \blre{k(X_1,Y_1)k(X_1,Y_2)} \\
  & & - 4\rho(1-\rho)^2 \blre{k(X_1,Y_1)k(Y_1,Y_2)}
   -  2(1-\rho)^3 \blre{k(Y_1,Y_2)k(Y_1,Y_3)} \\
  & & + \rho^4 \blre{k(X_1,X_2)}^2 
   +  4\rho^2(1-\rho)^2 \blre{k(X_1,Y_1)}^2 
   +  (1-\rho)^4 \blre{k(Y_1,Y_2)}^2 \\
  & & + 4\rho^3(1-\rho) \blre{k(X_1,X_2)}\blre{k(X_1,Y_1)}
   +  4\rho(1-\rho)^3 \blre{k(X_1,Y_1)}\blre{k(Y_1,Y_2)} \\
  & & + 2\rho^2(1-\rho)^2 \blre{k(X_1,X_2)}\blre{k(Y_1,Y_2)} \\
  &=& \blre{k^2(Z_1,Z_2)} - 2\blre{k(Z_1,Z_2)k(Z_1,Z_3)} + \blre{k(Z_1,Z_2)}^2 \\
  &=& \cv_k^2(Z),
  \EEqn
  which thus completes the proof.
\end{proof}

\subsection{Lemmas for Proposition \ref{Prop:ratio-const-k-H0} and Proposition \ref{Prop:ratio-const-k-H1}}

Next we provide some necessary lemmas to prove the ratio-consistency. 
\begin{lemma}\label{Lemma:psi_property_k}
  For any fixed $k=k^{(p)}\in\cc$ and for any $a_1,a_2,a_3,a_4\in\br^p$, it holds that
  \BEqn
  & & \psi(a_1,a_2,a_3,a_4) \\
  &=& \frac{1}{24}\lrp{2k(a_1,a_2) + 2k(a_3,a_4) - k(a_1,a_3) - k(a_1,a_4) - k(a_2,a_3) - k(a_2,a_4)}^2 \\
  & & + \frac{1}{8}\lrp{k(a_1,a_3) + k(a_2,a_4) - k(a_1,a_4) - k(a_2,a_3)}^2.
  \EEqn
\end{lemma}
\begin{proof}
  It follows from direct computations that
  \BEqn
      \text{RHS}
  &=& \frac{1}{6}
      \big(
        k^2(a_1,a_2) + k^2(a_1,a_3) + k^2(a_1,a_4)
        + k^2(a_2,a_3) + k^2(a_2,a_4) + k^2(a_3,a_4)
      \big) \\
 & & - \frac{1}{6}
      \big(
        k(a_1,a_2)k(a_1,a_3) + k(a_1,a_2)k(a_1,a_4) + k(a_1,a_3)k(a_1,a_4) \\
 & & \hspace{2em}        
        + k(a_1,a_2)k(a_2,a_3) + k(a_1,a_2)k(a_2,a_4) + k(a_2,a_3)k(a_2,a_4) \\
 & & \hspace{2em}        
        + k(a_1,a_3)k(a_2,a_3) + k(a_1,a_3)k(a_3,a_4) + k(a_2,a_3)k(a_3,a_4) \\    
 & & \hspace{2em}        
        + k(a_1,a_4)k(a_2,a_4) + k(a_1,a_4)k(a_3,a_4) + k(a_2,a_4)k(a_3,a_4)      
      \big) \\
 & & + \frac{1}{3}    
      \big(
        k(a_1,a_2)k(a_3,a_4) + k(a_1,a_3)k(a_2,a_4) + k(a_1,a_4)k(a_2,a_3)
      \big) \\
  &=& \psi(a_1,a_2,a_3,a_4),
  \EEqn
  which completes the proof.
\end{proof}

Define $r_X = \blre{k(X_1,X_2)}, r_{XY} = \blre{k(X_1,Y_1)}, r_Y = \blre{k(Y_1,Y_2)}$, then the centered version of $k$ can be expressed as follows:
\begin{equation*}
    \tilde{k}(X_1,X_2) = k(X_1,X_2) - r_X, \qquad
    \tilde{k}(Y_1,Y_2) = k(Y_1,Y_2) - r_Y,
\end{equation*}
and $\tilde{k}(X_1,Y_1) = k(X_1,Y_1) - r_{XY}$. With these notations, we can apply Lemma \ref{Lemma:psi_property_k} to rewrite the kernel functions in Lemma \ref{Lemma:Combn_Kernel}, as summarized in the following lemma.

\begin{lemma}\label{Lemma:psi_property_d^k}
  For $i=1,2,3,4$, let $X_i$ and $Y_i$ be independent copies of $X$ and $Y$, respectively, then for any fixed $k=k^{(p)}\in\cc$ it holds that
  \BEqn
 & & \psi(X_1,X_2,X_3,X_4) \\
  &=& \frac{1}{24}
      \big{(}
        2\tilde{k}(X_1,X_2) + 2\tilde{k}(X_3,X_4)
        - \tilde{k}(X_1,X_3) - \tilde{k}(X_1,X_4) 
        - \tilde{k}(X_2,X_3) - \tilde{k}(X_2,X_4)
      \big{)}^2 \\
 & & + \frac{1}{8}
      \big{(}
        \tilde{k}(X_1,X_3) + \tilde{k}(X_2,X_4) 
        - \tilde{k}(X_1,X_4) - \tilde{k}(X_2,X_3)
      \big{)}^2, \\
      \\
 & & \psi(X_1,X_2,X_3,Y_1) \\
  &=& \frac{1}{24} 
      \big{(}
        2\tilde{k}(X_1,X_2) + 2\tilde{k}(X_3,Y_1) 
        - \tilde{k}(X_1,X_3) - \tilde{k}(X_1,Y_1) 
        - \tilde{k}(X_2,X_3) - \tilde{k}(X_2,Y_1)
      \big{)}^2 \\
 & & + \frac{1}{8}
      \big{(}
        \tilde{k}(X_1,X_3) + \tilde{k}(X_2,Y_1) 
        - \tilde{k}(X_1,Y_1) - \tilde{k}(X_2,X_3) 
      \big{)}^2, \\
      \\
 & & \psi(X_1,X_2,Y_1,Y_2) \\
  &=& \frac{1}{24}
      \big{(}
        2\tilde{k}(X_1,X_2) + 2\tilde{k}(Y_1,Y_2)
        - \tilde{k}(X_1,Y_1) - \tilde{k}(X_1,Y_2)
        - \tilde{k}(X_2,Y_1) - \tilde{k}(X_2,Y_2)
        - 2\ce^k(X,Y)
      \big{)}^2 \\
 & & + \frac{1}{8}
      \big{(}
        \tilde{k}(X_1,Y_1) + \tilde{k}(X_2,Y_2) 
        - \tilde{k}(X_1,Y_2) - \tilde{k}(X_2,Y_1)
      \big{)}^2, \\
      \\
 & & \psi(X_1,Y_1,Y_2,Y_3) \\
  &=& \frac{1}{24}
      \big{(}
        2\tilde{k}(X_1,Y_1) + 2\tilde{k}(Y_2,Y_3)
        - \tilde{k}(X_1,Y_2) - \tilde{k}(X_1,Y_3)
        - \tilde{k}(Y_1,Y_2) - \tilde{k}(Y_1,Y_3)
      \big{)}^2 \\
 & & + \frac{1}{8}
      \big{(}
        \tilde{k}(X_1,Y_2) + \tilde{k}(Y_1,Y_3)
        - \tilde{k}(X_1,Y_3) - \tilde{k}(Y_1,Y_2)
      \big{)}^2, \\
      \\
 & & \psi(Y_1,Y_2,Y_3,Y_4) \\
  &=& \frac{1}{24}
      \big{(}
        2\tilde{k}(Y_1,Y_2) + 2\tilde{k}(Y_3,Y_4) 
        - \tilde{k}(Y_1,Y_3) - \tilde{k}(Y_1,Y_4) 
        - \tilde{k}(Y_2,Y_3) - \tilde{k}(Y_2,Y_4)
      \big{)}^2 \\
 & & + \frac{1}{8}
      \big{(}
        \tilde{k}(Y_1,Y_3) + \tilde{k}(Y_2,Y_4) 
        - \tilde{k}(Y_1,Y_4) - \tilde{k}(Y_2,Y_3)
      \big{)}^2, 
  \EEqn
\end{lemma}

Note that we've connected $\cv_{n,m}^{k\ast}(X,Y)$ with a combination of U-statistics, then we can study the behavior of the estimator via the moment inequalities for the two-sample U-statistics developed in the following lemma.
\begin{lemma}\label{Lemma:moment_ineq}
  Suppose that the two-sample U-statistic $U_{n,m}$ has mean $\theta$ and admits the following representation:
  \begin{equation*}
      U_{n,m} - \theta
    = \bin{n}{k_x}^{-1} \bin{m}{k_y}^{-1} \sum\limits_{1\le i_1<\dots<i_{k_x}\le n} \sum\limits_{1\le j_1<\dots<j_{k_y}\le m} \psi(X_{i_1},\dots,X_{i_{k_x}},Y_{j_1},\dots,Y_{j_{k_y}}),
  \end{equation*}
  then for any $1 \le r \le 2$, it holds that
  \begin{equation*}
      \blre{\lrabs{U_{n,m}-\theta}^r}
  \le C(k_x,k_y) \sum\limits_{c=0}^{k_x} \sum\limits_{d=0}^{k_y} \bin{n}{c}^{-r+1} \bin{m}{d}^{-r+1} \blre{\lrabs{\psi(X_1,\dots,X_{k_x},Y_1,\dots,Y_{k_y})}^r}.
  \end{equation*}
\end{lemma}
\begin{proof}
  It follows from Lemma \ref{Lemma:H-decomp} and the Hoeffding decomposition that the two-sample U-statistic $U_{n,m}$ admits the following representation
  \begin{equation*}
      U_{n,m} - \theta
    = \sum\limits_{c=d}^{k_x} \sum\limits_{d=0}^{k_y} \bin{k_x}{c} \bin{k_y}{d} \bin{n}{c}^{-1} \bin{m}{d}^{-1} S_{n,m}^{(c,d)},
  \end{equation*}
  where
  \begin{equation*}
      S_{n,m}^{(c,d)}
    = \sum\limits_{1 \le i_1 < \dots < i_c \le n} \sum\limits_{1 \le j_1 < \dots < j_d \le m} \psi^{(c,d)}(X_{i_1},\dots,X_{i_c},Y_{j_1},\dots,Y_{j_d})
  \end{equation*}
  and $\psi^{(c,d)}$ is defined by Definition \ref{Def:psi^(cd)}. Define $\cf_0 = \emptyset$, for $1\le i\le n$, define $\cf_i = \sigma(X_1,\dots,X_n)$ and for $1\le j\le m$, define $\cf_{n+j} = \sigma(X_1,\dots,X_n,Y_1,\dots,Y_j)$. It's trivial that $\{\cf_i: 0\le i\le n+m\}$ forms a filtration, and it follows from the definition of $\psi^{(c,d)}$ that
  \begin{equation*}
      \left\{
        \begin{array}{ll}
            \blre{S_{n,m}^{(c,d)} | \cf_k} = S_{k,0}^{(c,0)} & \text{if } c\le k\le n, \vspace{2mm}\\
            \blre{S_{n,m}^{(c,d)} | \cf_k} = S_{n,0}^{(c,0)} & \text{if } n+1\le k< n+d, \vspace{2mm}\\
            \blre{S_{n,m}^{(c,d)} | \cf_{n+k}} = S_{n,k}^{(c,d)} & \text{if } d\le k\le m, \\
        \end{array}
      \right.
  \end{equation*}
  which further implies that, if $c\le k\le n$, it holds that
  \BEqn
  & & \blre{S_{n,m}^{(c,d)} | \cf_k} - \blre{S_{n,m}^{(c,d)} | \cf_{k-1}} \\
  &=& \sum\limits_{1 \le i_1 < \dots < i_{c-1} \le {k-1}} \psi^{(c,0)}(X_{i_1},\dots,X_{i_{c-1}},X_k) \\
  &=& \sum\limits_{i_{c-1}=c-1}^{k-1} \sum\limits_{i_{c-2}=c-2}^{i_{c-1}-1} \cdots \sum\limits_{i_2=2}^{i_3-1} \sum\limits_{i_1=1}^{i_2-1} \psi^{(c,0)}(X_{i_1},\dots,X_{i_{c-1}},X_k),
  \EEqn
  and if $d\le k\le m$, we have
  \BEqn
  & & \blre{S_{n,m}^{(c,d)} | \cf_{n+k}} - \blre{S_{n,m}^{(c,d)} | \cf_{n+k-1}} \\
  &=& \sum\limits_{1 \le i_1 < \dots < i_c \le n} \sum\limits_{1 \le j_1 < \dots < j_{d-1} \le k-1} \psi^{(c,d)}(X_{i_1},\dots,X_{i_c},Y_{j_1},\dots,Y_{j_{d-1}},Y_k) \\
  &=& \sum\limits_{i_c=c}^{n} \sum\limits_{i_{c-1}=c-1}^{i_c-1} \cdots \sum\limits_{i_2=2}^{i_3-1} \sum\limits_{i_1=1}^{i_2-1} \sum\limits_{j_{d-1}=d-1}^{k-1} \sum\limits_{j_{d-2}=d-2}^{j_{d-1}-1} \cdots \sum\limits_{j_2=2}^{j_3-1} \sum\limits_{j_1=1}^{j_2-1} \psi^{(c,d)}(X_{i_1},\dots,X_{i_c},Y_{j_1},\dots,Y_{j_{d-1}},Y_k).
  \EEqn
  Consequently, we have
  \begin{equation*}
      U_{n,m} - \theta
    = \blre{U_{n,m} | \cf_{n+m}} - \blre{U_{n,m} | \cf_{0}}
    = \sum\limits_{t=1}^{n+m} \xi_{t},
  \end{equation*}
  where
  \BEqn
      \xi_t
  &=& \blre{U_{n,m} | \cf_{t}} - \blre{U_{n,m} | \cf_{t-1}} \\
  &=& \sum\limits_{c=0}^{k_x} \sum\limits_{d=0}^{k_y} \bin{k_x}{c} \bin{k_y}{d} \bin{n}{c}^{-1} \bin{m}{d}^{-1} \lrp{\blre{S_{n,m}^{(c,d)} | \cf_{t}} - \blre{S_{n,m}^{(c,d)} | \cf_{t-1}}}.
  \EEqn
  Finally, by Theorem 2.1.1 of \cite{korolyuk2013theory}, with 
  \begin{equation*}
      \alpha_r=\sup\limits_{x\in\br}\lrp{|x|^{-r}(|1+x|^r-1-rx)}\le2^{2-r},
  \end{equation*}
  it holds for any $1\le r\le 2$ that
  \begin{equation*}
      \blre{\lrabs{U_{n,m}-\theta}^{r}}  
  \le \alpha_r \sum\limits_{t=1}^{n+m} \blre{\lrabs{\xi_{t}}^{r}}
  = \alpha_r 
    \lrp{\sum\limits_{t=c}^{n} \blre{\lrabs{\xi_{t}}^{r}}
         + \sum\limits_{t=n+d}^{n+m} \blre{\lrabs{\xi_{t}}^{r}}}
  \end{equation*}
  
  It remains to consider $\blre{\lrabs{\xi_{t}}^{r}}$. If $c \le t \le n$, it follows from the $c_r$ inequality that
  \BEqn
  & & \blre{\lrabs{\xi_{t}}^{r}} \\
  &\le& (k_x+1)^{r-1}(k_y+1)^{r-1}
        \sum\limits_{c=0}^{k_x} \sum\limits_{d=0}^{k_y} \bin{k_x}{c}^r \bin{k_y}{d}^r \bin{n}{c}^{-r} \bin{m}{d}^{-r} 
        \blre{\lrabs{\blre{S_{n,m}^{(c,d)} | \cf_{t}} - \blre{S_{n,m}^{(c,d)} | \cf_{t-1}}}^r}
  \EEqn
  Note that the quantity within each summation of $\blre{S_{n,m}^{(c,d)} | \cf_{t}} - \blre{S_{n,m}^{(c,d)} | \cf_{t-1}}$ forms a martingale, then by repeatedly using Theorem 2.1.1 in \cite{korolyuk2013theory}, we obtain that
  \BEqn
  & & \blre{\lrabs{\xi_{t}}^{r}} \\
  &\le& (k_x+1)^{r-1}(k_y+1)^{r-1}
        \sum\limits_{c=0}^{k_x} \sum\limits_{d=0}^{k_y} \bin{k_x}{c}^r \bin{k_y}{d}^r \bin{n}{c}^{-r} \bin{m}{d}^{-r} \\
  & & \hspace{2em} \times
        \alpha_r^{c-1}
        \sum\limits_{1 \le i_1 < \dots < i_{c-1} \le {t-1}} \blre{\lrabs{\psi^{(c,0)}(X_{i_1},\dots,X_{i_{c-1}},X_k)}^r} \\
  &\le& C(k_x,k_y) \alpha_r^{c-1} \sum\limits_{c=0}^{k_x} \sum\limits_{d=0}^{k_y} \bin{n}{c}^{-r} \bin{m}{d}^{-r} \bin{t-1}{c-1} \blre{\lrabs{\psi(X_1,\dots,X_{k_x},Y_1,\dots,Y_{k_y})}^{r}},
  \EEqn
  then by noting that $\sum\limits_{t=c}^{n}\bin{t-1}{c-1}=\bin{n}{c}$, we further have that
  \begin{equation*}
     \sum\limits_{t=c}^{n} \blre{\lrabs{\xi_{t}}^r}
  \le C(k_x,k_y) \sum\limits_{c=0}^{k_x} \sum\limits_{d=0}^{k_y} \alpha_r^{c-1} \bin{n}{c}^{-r+1} \bin{m}{d}^{-r} \blre{\lrabs{\psi(X_1,\dots,X_{k_x},Y_1,\dots,Y_{k_y})}^{r}}
  \end{equation*}
  Similarly, we have that
  \begin{equation*}
      \sum\limits_{t=n+d}^{n+m} \blre{\lrabs{\xi_{t}}^{r}}
  \le C(k_x,k_y) \sum\limits_{c=0}^{k_x} \sum\limits_{d=0}^{k_y} \alpha_r^{c+d-1} \bin{n}{c}^{-r+1} \bin{m}{d}^{-r+1} \blre{\lrabs{\psi(X_1,\dots,X_{k_x},Y_1,\dots,Y_{k_y})}^{r}},
  \end{equation*}
  which leads to the proposed result.
\end{proof}


\begin{lemma}\label{Lemma:Consistency_bd_1}
  Assume that $n/N \rightarrow \rho$ for some constant $0<\rho<1$ as $n,m\rightarrow\infty$, then for any $0\le\tau\le1$ and any fixed $k=k^{(p)}\in\cc$ it holds that
  \BEqn
 & & \blre{\lrabs{\cv_{n,m}^{k\ast}(X,Y) - \blre{\cv_{n,m}^{k\ast}(X,Y)}}^{1+\tau}} \\
 &\le& C(\rho,\tau) N^{-\tau}
   \lrp{
     \blre{\lrabs{\tilde{k}(X_1,X_2)}^{2+2\tau}} 
   + \blre{\lrabs{\tilde{k}(X_1,Y_1)}^{2+2\tau}} 
   + \blre{\lrabs{\tilde{k}(Y_1,Y_2)}^{2+2\tau}}
   + \lrabs{\ce^k(X,Y)}^{2+2\tau}}
  \EEqn
\end{lemma}
\begin{proof}
  By Lemma \ref{Lemma:Combn_Kernel}, we have
  \BEqn
  & & \blre{\lrabs{\cv_{n,m}^{k\ast}(X,Y) - \blre{\cv_{n,m}^{k\ast}(X,Y)}}^{1+\tau}} \\
  &\le& C(\tau) \bin{N}{4}^{-(1+\tau)} 
    \left(
      \blre{\lrabs{\sum\limits_{1\le i_1<i_2<i_3<i_4\le n} \lrp{\psi(X_{i_1},X_{i_2},X_{i_3},X_{i_4}) - \blre{\psi(X_{i_1},X_{i_2},X_{i_3},X_{i_4})}}}^{1+\tau}}
    \right. \\
  & & \hspace{6em} 
      + \blre{\lrabs{\sum\limits_{1\le i_1<i_2<i_3\le n} \sum\limits_{j=1}^{m} \lrp{\psi(X_{i_1},X_{i_2},X_{i_3},Y_{j}) - \blre{\psi(X_{i_1},X_{i_2},X_{i_3},Y_{j})}}}^{1+\tau}} \\
  & & \hspace{6em} 
      + \blre{\lrabs{\sum\limits_{1\le i_1<i_2\le n} \sum\limits_{1\le j_1<j_2\le m} \lrp{\psi(X_{i_1},X_{i_2},Y_{j_1},Y_{j_2}) - \blre{\psi(X_{i_1},X_{i_2},Y_{j_1},Y_{j_2})}}}^{1+\tau}} \\
  & & \hspace{6em} 
      + \blre{\lrabs{\sum\limits_{1\le j_1<j_2<j_3\le m} \sum\limits_{i=1}^{n} \lrp{\psi(X_{i},Y_{j_1},Y_{j_2},Y_{j_3}) - \blre{\psi(X_{i},Y_{j_1},Y_{j_2},Y_{j_3})}}}^{1+\tau}} \\
  & & \hspace{6em} 
    \left.
      + \blre{\lrabs{\sum\limits_{1\le j_1<j_2<j_3<j_4\le n} \lrp{\psi(Y_{j_1},Y_{j_2},Y_{j_3},Y_{j_4}) - \blre{\psi(Y_{j_1},Y_{j_2},Y_{j_3},Y_{j_4})}}}^{1+\tau}}
    \right).
  \EEqn
  Also note that it follows from Lemma \ref{Lemma:psi_property_d^k} and Lemma \ref{Lemma:moment_ineq} that
  \BEqn
  & & \blre{\lrabs{\sum\limits_{1\le i_1<i_2<i_3<i_4\le n} \lrp{\psi(X_{i_1},X_{i_2},X_{i_3},X_{i_4}) -  \blre{\psi(X_{i_1},X_{i_2},X_{i_3},X_{i_4})}}}^{1+\tau}} \\
  &=& \bin{n}{4}^{1+\tau} \blre{\lrabs{\bin{n}{4}^{-1} \sum\limits_{1\le i_1<i_2<i_3<i_4\le n} \lrp{\psi(X_{i_1},X_{i_2},X_{i_3},X_{i_4}) - \blre{\psi(X_{i_1},X_{i_2},X_{i_3},X_{i_4})}}}^{1+\tau}} \\
  &\le& C(\tau) \bin{n}{4}^{1+\tau} n^{-\tau} \blre{\lrabs{\psi(X_1,X_2,X_3,X_4)}^{1+\tau}} \\
  &\le& C(\tau) \bin{n}{4}^{1+\tau} n^{-\tau} \blre{\lrabs{\tilde{k}(X_1,X_2)}^{2+2\tau}}.
  \EEqn
  Similarly, it can also be derived that
  \BEqn
  & & \blre{\lrabs{\sum\limits_{1\le i_1<i_2<i_3\le n} \sum\limits_{j=1}^{m} \lrp{\psi(X_{i_1},X_{i_2},X_{i_3},Y_{j}) - \blre{\psi(X_{i_1},X_{i_2},X_{i_3},Y_{j})}}}^{1+\tau}} \\
  &\le& C(\tau) \bin{n}{3}^{1+\tau} \bin{m}{1}^{1+\tau} \lrp{n^{-\tau} + m^{-\tau}} 
  \big(
      \blre{\lrabs{\tilde{k}(X_1,X_2)}^{2+2\tau}}
    + \blre{\lrabs{\tilde{k}(X_1,Y_1)}^{2+2\tau}}
  \big) 
  \\ \\
  & & \blre{\lrabs{\sum\limits_{1\le i_1<i_2\le n} \sum\limits_{1\le j_1<j_2\le m} \lrp{\psi(X_{i_1},X_{i_2},Y_{j_1},Y_{j_2}) - \blre{\psi(X_{i_1},X_{i_2},Y_{j_1},Y_{j_2})}}}^{1+\tau}} \\
  &\le& C(\tau) \bin{n}{2}^{1+\tau} \bin{m}{2}^{1+\tau} \lrp{n^{-\tau} + m^{-\tau}} 
  \big(
      \blre{\lrabs{\tilde{k}(X_1,X_2)}^{2+2\tau}}
    + \blre{\lrabs{\tilde{k}(Y_1,Y_2)}^{2+2\tau}} \\
  & & \hspace{16em}       
    + \blre{\lrabs{\tilde{k}(X_1,Y_1)}^{2+2\tau}}
    + \lrabs{\ce^k(X,Y)}^{2+2\tau}
  \big) 
  \\ \\
  & & \blre{\lrabs{\sum\limits_{1\le j_1<j_2<j_3\le m} \sum\limits_{i=1}^{n} \lrp{\psi(X_{i},Y_{j_1},Y_{j_2},Y_{j_3}) - \blre{\psi(X_{i},Y_{j_1},Y_{j_2},Y_{j_3})}}}^{1+\tau}} \\
  &\le& C(\tau) \bin{n}{1}^{1+\tau} \bin{m}{3}^{1+\tau} \lrp{n^{-\tau} + m^{-\tau}} 
  \big(
      \blre{\lrabs{\tilde{k}(X_1,Y_1)}^{2+2\tau}}
    + \blre{\lrabs{\tilde{k}(Y_1,Y_2)}^{2+2\tau}}
  \big) 
  \\ \\
  & & \blre{\lrabs{\sum\limits_{1\le j_1<j_2<j_3<j_4\le n} \lrp{\psi(Y_{j_1},Y_{j_2},Y_{j_3},Y_{j_4}) - \blre{\psi(Y_{j_1},Y_{j_2},Y_{j_3},Y_{j_4})}}}^{1+\tau}} \\
  &\le& C(\tau) \bin{m}{4}^{1+\tau} m^{-\tau} \blre{\lrabs{\tilde{k}(Y_1,Y_2)}^{2+2\tau}},
  \EEqn
  which implies the desired result under the assumption that $n/N \rightarrow \rho$.
\end{proof}

\begin{lemma}\label{Lemma:Consistency_bd_2}
  Assume that $\bbe{k^2(Z_1,Z_2)}<\infty$ and $n/N = \rho + O(1/N^s)$ for some constant $0<\rho<1$ and $s>0$ as $n,m\rightarrow\infty$, then for any fixed $k=k^{(p)}\in\cc$ it holds that
  \begin{equation}
      \lrabs{\blre{\cv_{n,m}^{k\ast}(X,Y)} - \cv_k^2(Z)}
  \le C(\rho) \lrp{\frac{1}{N^s}} \blre{k^2(Z_1,Z_2)}
  \end{equation}
  as $n,m\rightarrow\infty$.
\end{lemma}
\begin{proof}
  We've derived the explicit expression of $\blre{\cv_{n,m}^{k\ast}(X,Y)}$ under the alternative in the proof of Lemma \ref{Lemma:ratio-const-k-H1}. Under the assumption that $n/N = \rho + O(1/N^s)$, we have
  \begin{equation*}
      \frac{n(n-1)}{N(N-1)} = \rho^2 + O(\frac{1}{N^s}), \quad
      \frac{nm}{N(N-1)} = \rho(1-\rho) + O(\frac{1}{N^s}), \quad
      \frac{m(m-1)}{N(N-1)} = (1-\rho)^2 + O(\frac{1}{N^s}).
  \end{equation*}
  Similarly, we can approximate the other coefficients in terms of $\rho$. Also note that 
  \begin{equation*}
      \rho^2 \blre{k^2(X_1,X_2)} \le \blre{k^2(Z_1,Z_2)}, \qquad
      2\rho(1-\rho) \blre{k^2(X_1,Y_1)} \le \blre{k^2(Z_1,Z_2)}, 
  \end{equation*}
  and $(1-\rho)^2 \blre{k^2(Y_1,Y_2)} \le \blre{k^2(Z_1,Z_2)}$. Using H\"older's inequality, we further have that all the quantities in form of $\lrabs{\bbe{k(X_1,X_2)k(Y_1,Y_2)}}$ is upper bounded by $C(\rho)\bbe{k^2(Z_1,Z_2)}$, then it follows that
  \BEqn
  & & \lrabs{\blre{\cv_{n,m}^{k\ast}(X,Y)} - \cv_k^2(Z)} \\
  &\le& \frac{C}{N^s}
      \big|
         \blre{k^2(X_1,X_2)}
       + \blre{k^2(X_1,Y_1)} 
       + \blre{k^2(Y_1,Y_2)} \\
  & & \hspace{4em}
       + \blre{k(X_1,X_2)k(X_1,X_3)}
       + \blre{k(X_1,X_2)k(X_1,Y_1)} 
       + \blre{k(X_1,Y_1)k(X_2,Y_1)} \\
  & & \hspace{4em}
       + \blre{k(X_1,Y_1)k(X_1,Y_2)}
       + \blre{k(X_1,Y_1)k(Y_1,Y_2)} 
       + \blre{k(Y_1,Y_2)k(Y_1,Y_3)} \\       
  & & \hspace{4em}
       + \blre{k(X_1,X_2)}^2
       + \blre{k(Y_1,Y_2)}^2 
       + \blre{k(X_1,X_2)}\blre{k(Y_1,Y_2)} \\
  & & \hspace{4em}
       + \blre{k(X_1,Y_1)}^2
       + \blre{k(X_1,X_2)}\blre{k(X_1,Y_1)}
       + \blre{k(X_1,Y_1)}\blre{k(Y_1,Y_2)}
      \big| \\
  &\le& C(\rho) \lrp{\frac{1}{N^s}} \blre{k^2(Z_1,Z_2)},
  \EEqn
  which completes the proof.
\end{proof}


\section{Lemmas Regarding the Rate of Convergence}\label{Appendix:lemma-2}

\subsection{Lemmas for Theorem \ref{Thm:clt} and Theorem \ref{Thm:BEbd-k}}

To derive the rate of convergence under the null, we start from $P_1$ and attempt to find the upper bound of $P_1$ using the Berry-Esseen inequality.

Under the null hypothesis, it is trivial that $X =^d Y =^d Z$, then without the loss of generality, we assume that $Z$ has mean zero and use $\fZ = \lrp{\fX,\fY}$ to denote the pooled sample over $\fX$ and $\fY$.
Suppose that $\blre{k^2(X_1,X_2)} = \blre{k^2(Y_1,Y_2)} < \infty$, we have shown in Lemma \ref{Lemma:h_cd} that it holds under the null that $h_{20}(X_1,X_2) = -d^k(X_1,X_2)$, $h_{11}(X_1,Y_1) = \frac{1}{2}d^k(X_1,Y_1)$ and $h_{02}(Y_1,Y_2) = -d^{k}(Y_1,Y_2)$. Consequently, we can rewrite $L_{n,m}^k(X,Y)$ introduced in Proposition \ref{Prop:est_MMD_simplified} into the following form
\begin{equation*}
    L_{n,m}^{k}(X,Y) 
  = \sum\limits_{j=1}^{n+m} \sum\limits_{i=1}^{j-1} \iota_{ij} d^k(Z_i,Z_j),
\end{equation*}
where
\begin{equation*}
   \iota_{ij} 
 = \left\{
     \begin{array}{ll}
       \displaystyle{-\frac{2}{n(n-1)}}, & 1 \le i<j \le n  \vspace{3pt} \\
       \displaystyle{\frac{2}{nm}},      & 1 \le i \le n < j \le n+m \vspace{3pt} \\
       \displaystyle{-\frac{2}{m(n-m)}}, & n+1 \le i<j \le n+m
     \end{array}
   \right.    
\end{equation*}
and $d^{k}$ is defined as (\ref{Equ:d^k}). Under the null, it is shown in the proof to Proposition \ref{Prop:Var(Lnm)} that
\begin{equation*}
    \Var(L_{n,m}^{k}(X,Y)) = c_{n,m} \cv_k^2(Z).
\end{equation*}
Define the scaled statistic
\begin{equation*}
    \tilde{L}_{n,m}^{k}(X,Y) = \frac{L_{n,m}^{k}(X,Y)}{\sqrt{\Var(L_{n,m}^{k}(X,Y))}},
\end{equation*}
then under $H_0$, it holds that
\begin{equation*}
    \tilde{L}_{n,m}^{k}(X,Y) 
  = \sum\limits_{j=1}^{n+m} \sum\limits_{i=1}^{j-1} \tilde{\iota}_{ij} \frac{d^k(Z_i,Z_j)}{\sqrt{\cv_K^2(Z)}},
\end{equation*}
where $\displaystyle{\tilde{\iota}_{ij} = \frac{\iota_{ij}}{\sqrt{c_{n,m}}}}$.

By applying the Berry-Esseen inequality established in Theorem 1 of \cite{haeusler1988rate}, we can obtain an upper bound of $P_1$, as stated in the following proposition.

\begin{lemma}\label{Lemma:P1-inter}
For each $1 \le i \le N$, define the $\sigma$-algebra $\cf_i = \sigma(Z_1,\dots,Z_i)$. For $0<\tau\le1$, define
\begin{equation*}
    B_{n,m,\tau}^{(1)}
= \blre{\lrabs{\sum\limits_{j=1}^{n+m} \bbe{\bigp{\sum\limits_{i=1}^{j-1}\tilde{\iota}_{ij}\frac{d^k(Z_i,Z_j)}{\sqrt{\cv_k^2(Z)}}}^2 | \cf_{j-1}} - 1}^{1+\tau}},
\end{equation*}
and 
\begin{equation*}
    B_{n,m,\tau}^{(2)}
= \sum\limits_{j=1}^{n+m} \blre{\lrabs{\sum\limits_{i=1}^{j-1} \tilde{\iota}_{ij}\frac{d^k(Z_i,Z_j)}{\sqrt{\cv_k^2(Z)}}}^{2+2\tau}}.
\end{equation*}
Under the null, there exists a strictly positive constant $C_{\tau}$ depending only on $\tau$, such that
\begin{equation*}
    P_1 
  = \sup\limits_{x\in\br} \lrabs{\blrp{\tilde{L}_{n,m}^{k}(X,Y) \le x} - \Phi(x)}
  \le C_{\tau}\lrp{B_{n,m,\tau}^{(1)} + B_{n,m,\tau}^{(2)}}^{1/(3+2\tau)}.
\end{equation*}
\end{lemma}

Lemma \ref{Lemma:P1-inter} implies that, to derive an explicit upper bound of $P_1$, it suffices to bound $B_{n,m,\tau}^{(1)}$ and $ B_{n,m,\tau}^{(2)}$, respectively. We first look into $B_{n,m,\tau}^{(1)}$. For $1\le j\le N$, define
\begin{equation*}
    \zeta_{j} 
  = \blre{\lrp{\sum\limits_{i=1}^{j-1} \tilde{\iota}_{ij} \frac{d^{k}(Z_i,Z_j)}{\sqrt{\cv_k^2(Z)}}}^2 | \cf_{j-1}},
\end{equation*}
then under the null, it holds that
\begin{equation*}
    B_{n,m,\tau}^{(1)} 
 = \blre{\lrabs{\sum\limits_{j=1}^{N} \zeta_j - 1}^{1+\tau}}.
\end{equation*}
By the independence of the sample $Z$, we have
\BEqn
& & \zeta_j \\
&=& \frac{1}{\cv_k^2(Z)} \sum\limits_{i_1,i_2=1}^{j-1} \tilde{\iota}_{i_1,j} \tilde{\iota}_{i_2,j} \blre{d^k(Z_{i_1},Z_j) d^k(Z_{i_2},Z_j) | \cf_{j-1}} \\
&=& \sum\limits_{i=1}^{j-1} \frac{\tilde{\iota}_{ij}^2}{\cv_k^2(Z)} \blre{\lrp{d^k(Z_i,Z_j)}^2 | Z_i} 
 +  \sum\limits_{1\le i_1<i_2\le j-1} \frac{2 \tilde{\iota}_{i_1,j} \tilde{\iota}_{i_2,j}}{\cv_k^2(Z)} \blre{d^k(Z_{i_1},Z_j) d^k(Z_{i_2},Z_j) | Z_{i_1},Z_{i_2}} \\
 &=:& \zeta_j^{(1)} + \zeta_j^{(2)}.
\EEqn
It follows that
\BEqn
    B_{n,m,\tau}^{(1)}
&=& \blre{\lrabs{ \lrp{\sum\limits_{j=1}^{N} \zeta_j^{(1)} - 1} + \sum\limits_{j=1}^{N} \zeta_j^{(2)} }^{1+\tau}} \\
&\le& 2^{\tau} \blre{\lrabs{\sum\limits_{j=1}^{N} \zeta_j^{(1)}-1}^{1+\tau}} + 2^{\tau} \blre{\lrabs{\sum\limits_{j=1}^{N} \zeta_j^{(2)}}^{1+\tau}} \\
&=:& 2^{\tau}\lrp{\fzeta^{(1)} + \fzeta^{(2)}}.
\EEqn
The bounds of $\fzeta^{(1)}$ and $\fzeta^{(2)}$ are given by the following lemmas.

For each $1 \le i_1 < i_2 \le N$, define
\begin{equation*}
    \eta(Z_{i_1},Z_{i_2})
  = \blre{d^k(Z_{i_1},Z) d^k(Z_{i_2},Z) | Z_{i_1},Z_{i_2}},
\end{equation*}
where $Z$ is independent of $Z_{i_1}$ and $Z_{i_2}$.

\begin{lemma}\label{Lemma:cross_eta_1}
For any $1\le i_1<i_2 \le N$ and $1\le i_3<i_4 \le N$, if $(i_1,i_2)\neq(i_3,i_4)$, then it holds that
\begin{equation*}
    \blre{\eta(Z_{i_1},Z_{i_2}) \eta(Z_{i_3},Z_{i_4})} = 0.
\end{equation*}
\end{lemma}
\begin{proof}
  It follows from the definition of $d^k(Z_1,Z_2)$ that $\blre{d^k(Z_1,Z_2)|Z_1} = 0$, and consequently,
  \begin{equation*}
      \blre{\eta(Z_1,Z_2)}
    = \blre{d^k(Z_1,Z_3) d^k(Z_2,Z_3)}
    = \blre{\blre{d^k(Z_1,Z_3)|Z_3} \blre{d^k(Z_2,Z_3)|Z_3}}
    = 0.
  \end{equation*}
  If $i_1,i_2,i_3$ and $i_4$ are pairwise distinct, then by the independence, it is trivial that
  \begin{equation*}
      \blre{\eta(Z_{i_1},Z_{i_2}) \eta(Z_{i_3},Z_{i_4})}
    = \lrp{\blre{\eta(Z_1,Z_2)}}^2
    = 0.
  \end{equation*}
  If $\{i_1,i_2\}$ and $\{i_3,i_4\}$ have and only have one index in common, without the loss of generality, we assume that $i_1 = i_3$ but $i_2 \neq i_4$. Let $Z,Z'$ denote an independent copy of $Z_j$ and be independent of $Z_{i_1},Z_{i_2},Z_{i_4}$, we thus have
  \BEqn
 & & \blre{\eta(Z_{i_1},Z_{i_2}) \eta(Z_{i_1},Z_{i_4})} \\
  &=& \blre{\blre{d^k(Z_{i_1},Z_j) d^k(Z_{i_2},Z_j) | Z_{i_1},Z_{i_2}} \cdot \blre{d^k(Z_{i_1},Z_j) d^k(Z_{i_4},Z_j) | Z_{i_1},Z_{i_4}}} \\
  &=& \blre{\lrp{\int d^k(Z_{i_1},Z) d^k(Z_{i_2},Z) dF_Z} \cdot \lrp{\int d^k(Z_{i_1},Z') d^k(Z_{i_4},Z') dF_{Z'}}} \\
  &=& \blre{\int\int d^k(Z_{i_1},Z) d^k(Z_{i_1},Z') d^k(Z_{i_2},Z) d^k(Z_{i_4},Z') dF_Z dF_{Z'}} \\
  &=& \blre{d^k(Z_{i_1},Z) d^k(Z_{i_1},Z') d^k(Z_{i_2},Z) d^k(Z_{i_4},Z')} \\
  &=& \blre{\blre{d^k(Z_{i_1},Z) d^k(Z_{i_1},Z') | Z,Z'} \blre{d^k(Z_{i_2},Z) | Z} \blre{d^k(Z_{i_4},Z') | Z'}} \\
  &=& 0,
  \EEqn
  which completes the proof.
\end{proof}

\begin{lemma}\label{Lemma:cross_eta_2}
For any $1\le i_1<i_2 \le N$, it holds that
\begin{equation*}
    \blre{\eta^2(Z_{i_1},Z_{i_2})}
  = \blre{g^k(Z_1,Z_2,Z_3,Z_4)},
\end{equation*}
where the functional $g^k$ is defined as (\ref{Def:gk}).
\end{lemma}
\begin{proof}
  In fact,
  \BEqn
 & & \blre{\lrp{\eta(Z_1,Z_2)}^2} \\
  &=& \blre{\blre{d^k(Z_1,Z_3) d^k(Z_2,Z_3) | Z_1,Z_2} \cdot \blre{d^k(Z_1,Z_4) d^k(Z_2,Z_4) | Z_1,Z_2}} \\
  &=& \blre{\lrp{\int d^k(Z_1,Z_3) d^k(Z_2,Z_3) dF_{Z_3}} \lrp{\int d^k(Z_1,Z_4) d^k(Z_2,Z_4) dF_{Z_4}}} \\
  &=& \blre{\int\int d^k(Z_1,Z_3) d^k(Z_2,Z_3) d^k(Z_1,Z_4) d^k(Z_2,Z_4) dF_{Z_3}dF_{Z_4}} \\
  &=& \blre{g^k(Z_1,Z_2,Z_3,Z_4)},
  \EEqn
  which arrives at the desired result.
\end{proof}

\begin{lemma}\label{Lemma:fzeta_2}
    Under the null, it holds that
    \begin{equation*}
        \displaystyle{\fzeta^{(2)} \le A_{n,m,\tau}^{(2)} \frac{\lrp{\blre{g^k(Z_1,Z_2,Z_3,Z_4)}}^{(1+\tau)/2}}{\lrp{\cv_k^2(Z)}^{1+\tau}}}.
    \end{equation*}
    where $A_{n,m,\tau}^{(2)} = \lrp{\sum\limits_{1\le i_1<i_2\le N} a_{i_1,i_2}^2}^{(1+\tau)/2}$, and $a_{i_1,i_2} = \frac{2}{c_{n,m}} \sum\limits_{j=i_2+1}^{N} \iota_{i_1,j} \iota_{i_2,j}$.
\end{lemma}
\begin{proof}
  Note that with $0<\tau\le1$, it follows from Jensen's inequality that
  \begin{equation*}
      \fzeta^{(2)}
    = \blre{\lrabs{\sum\limits_{j=1}^{N} \zeta_j^{(2)}}^{1+\tau}} 
    \le \lrp{\blre{\lrabs{\sum\limits_{j=1}^{N} \zeta_j^{(2)}}^2}}^{(1+\tau)/2}
    =: \lrp{\tilde{\fzeta}^{(2)}}^{(1+\tau)/2},
  \end{equation*}
  then it follows from Lemma \ref{Lemma:cross_eta_1} and \ref{Lemma:cross_eta_2} that
  \BEqn
      \tilde{\fzeta}^{(2)}
  &=& \blre{\lrabs{\sum\limits_{j=1}^{N} \sum\limits_{1\le i_1<i_2\le j-1} \frac{2\tilde{\iota}_{i_1,j}\tilde{\iota}_{i_2,j}}{\cv_k^2(Z)} \eta(Z_{i_1},Z_{i_2})}^2} \\
  &=& \blre{\lrabs{\sum\limits_{1\le i_1<i_2\le N} \sum\limits_{j=i_2+1}^{N}  \frac{2\tilde{\iota}_{i_1,j}\tilde{\iota}_{i_2,j}}{\cv_k^2(Z)} \eta(Z_{i_1},Z_{i_2})}^2} \\
  &=& \blre{\lrp{\sum\limits_{1\le i_1<i_2\le N} a_{i_1,i_2} \frac{\eta(Z_{i_1},Z_{i_2})}{\cv_k^2(Z)}}^2} \\
  &=& \sum\limits_{1\le i_1<i_2 \le N} \sum\limits_{1\le i_3<i_4 \le N} a_{i_1,i_2} a_{i_3,i_4} \frac{\blre{\eta(Z_{i_1},Z_{i_2})\eta(Z_{i_3},Z_{i_4})}}{\lrp{\cv_k^2(Z)}^2} \\
  &=& \sum\limits_{1\le i_1<i_2 \le N} a_{i_1,i_2}^2 \frac{\blre{g^k(Z_1,Z_2,Z_3,Z_4)}}{\lrp{\cv_k^2(Z)}^2},
  \EEqn
  and consequently,
  \begin{equation*}
      \fzeta^{(2)} 
  \le \lrp{\tilde{\fzeta}^{(2)}}^{(1+\tau)/2}
   =  A_{n,m,\tau}^{(2)} \frac{\lrp{\blre{g^k(Z_1,Z_2,Z_3,Z_4)}}^{(1+\tau)/2}}{\lrp{\cv_k^2(Z)}^{1+\tau}}.
  \end{equation*}
\end{proof}

For each $j$ and $i<j$, define the functional
\begin{equation*}
    \tilde{\eta}(Z_i)
  = \eta(Z_i,Z_i)
  = \blre{\lrp{d^k(Z_i,Z_j)}^2 | Z_i}.
\end{equation*}
\begin{lemma}\label{Lemma:fzeta_1}
    Under the null, it holds that
    \begin{equation*}
        \fzeta^{(1)}
    \le A_{n,m,\tau}^{(1)} \frac{\blre{\lrabs{d^k(Z_1.Z_2)}^{2+2\tau}}}{\lrp{\cv_k^2(Z)}^{1+\tau}},
    \end{equation*}
    where  
    \begin{equation*}
        A_{n,m,\tau}^{(1)} = \sum\limits_{i=1}^{N}\frac{\lrabs{a_{ii}}^{1+\tau}}{1-D_{\tau}}
        \quad\text{with}\quad
        D_{\tau} = \frac{13.52}{(2.6(1+\tau))^{1+\tau}} \Gamma(1+\tau) \sin(\tau\pi/2) < 1,
    \end{equation*}
    and $a_{ii}$ is define as in Lemma \ref{Lemma:fzeta_2}.
\end{lemma}
\begin{proof}
  Note that $\sum\limits_{1\le i<j\le N} \iota_{ij}^2 = c_{n,m}$ and it holds under the null that $\blre{\lrp{d^k(Z_1,Z_2)}^2} = \cv_k^2(Z)$, then
  \BEqn
      \blre{\sum\limits_{j=1}^{N} \zeta_j^{(1)}}
  &=& \sum\limits_{j=1}^{N} \sum\limits_{i=1}^{j-1} \frac{\tilde{\iota}_{ij}^2}{\cv_k^2(Z)} \blre{\bbe{\lrp{d^k(Z_i,Z_j)}^2 | Z_i}} \\
  &=& \sum\limits_{1\le i<j\le N} \frac{\iota_{ij}^2}{c_{n,m}} \frac{\bbe{\lrp{d^k(Z_1,Z_2)}^2}}{\cv_k^2(Z)} \\
  &=& 1,
  \EEqn
  and consequently,
  \begin{equation*}
      \fzeta^{(1)}
    = \blre{\lrabs{\sum\limits_{j=1}^{N} \zeta_j^{(1)} - 1}^{1+\tau}}
    = \blre{\lrabs{\sum\limits_{j=1}^{N} \lrp{\zeta_j^{(1)} - \bbe{\zeta_j^{(1)}}}}^{1+\tau}}
    = \blre{\lrabs{\sum\limits_{j=1}^{N} \tilde{\zeta}_j^{(1)}}^{1+\tau}},
  \end{equation*}
  where $\tilde{\zeta}_j^{(1)} = \zeta_j^{(1)} - \bbe{\zeta_j^{(1)}}$.
  
  Recall that $a_{i_1,i_2} = \frac{2}{c_{n,m}} \sum\limits_{j=i_2+1}^{N} \iota_{i_1,j} \iota_{i_2,j}$ as defined in Lemma \ref{Lemma:fzeta_2}, then it is trivial that $\sum\limits_{j=i+1}^{N} \tilde{\iota}_{ij}^2 = \frac{1}{2}a_{ii}$, and it follows that
  \BEqn
      \sum\limits_{j=1}^{N} \tilde{\zeta}_j^{(1)}
  &=& \sum\limits_{j=1}^{N} \sum\limits_{i=1}^{j-1} \frac{\tilde{\iota}_{ij}^2}{\cv_k^2(Z)} \lrp{\tilde{\eta}(Z_i) - \blre{\tilde{\eta}(Z_i)}} \\
  &=& \sum\limits_{i=1}^{N} \sum\limits_{j=i+1}^{N} \tilde{\iota}_{ij}^2 \frac{\tilde{\eta}(Z_i) - \blre{\tilde{\eta}(Z_i)}}{\cv_k^2(Z)} \\
  &=& \sum\limits_{i=1}^{N} a_{ii} \frac{\tilde{\eta}(Z_i) - \blre{\tilde{\eta}(Z_i)}}{2\cv_k^2(Z)}.
  \EEqn
  Also note that 
  \begin{equation*}
      \blre{\lrabs{\tilde{\eta}(Z_i) - \blre{\tilde{\eta}(Z_i)}}^{1+\tau}}
  \le 2^{\tau} \lrp{\blre{\lrabs{\tilde{\eta}(Z_i)}^{1+\tau}} + \lrabs{\blre{\tilde{\eta}(Z_i)}}^{1+\tau}}
  \le 2^{1+\tau} \blre{\lrabs{\tilde{\eta}(Z_i)}^{1+\tau}},
  \end{equation*}
  and furthermore, it follows from Jensen's inequality that
  \begin{equation*}
      \blre{\lrabs{\tilde{\eta}(Z_i)}^{1+\tau}}
    = \blre{\lrabs{\blre{\lrp{d^k(Z_1,Z_2)}^2 | Z_1}}^{1+\tau}}
    \le \blre{\lrabs{d^k(Z_1,Z_2)}^{2+2\tau}},
  \end{equation*}
  hence
  \begin{equation*}
      \blre{\lrabs{\tilde{\eta}(Z_i) - \blre{\tilde{\eta}(Z_i)}}^{1+\tau}}
  \le 2^{1+\tau} \blre{\lrabs{d^k(Z_1,Z_2)}^{2+2\tau}}.
  \end{equation*}
  
  For each $i$, the term $\displaystyle{a_{ii} \frac{\tilde{\eta}(Z_i) - \blre{\tilde{\eta}(Z_i)}}{2\cv_k^2(Z)}}$ is a function depending only on $Z_i$, thus for any $i\neq j$, they are independent with mean zero. Apply the von Bahr-Esseen inequality derived in Theorem 9.3.a of \cite{lin2011probability}, we can obtain the following bound of $\fzeta^{(1)}$.
  \BEqn
      \fzeta^{(1)}
  &=& \blre{\lrabs{\sum\limits_{i=1}^{N} a_{ii} \frac{\tilde{\eta}(Z_i) - \blre{\tilde{\eta}(Z_i)}}{2\cv_k^2(Z)}}^{1+\tau}} \\
  &\le& \frac{1}{1-D_{\tau}} \sum\limits_{i=1}^{N} \blre{\lrabs{a_{ii} \frac{\tilde{\eta}(Z_i) - \blre{\tilde{\eta}(Z_i)}}{2\cv_k^2(Z)}}^{1+\tau}} \\
  &=& \frac{1}{1-D_{\tau}} \sum\limits_{i=1}^{N} \lrp{\frac{\lrabs{a_{ii}}}{2}}^{1+\tau} \frac{\blre{\lrabs{\tilde{\eta}(Z_i) - \blre{\tilde{\eta}(Z_i)}}^{1+\tau}}}{\lrp{\cv_k^2(Z)}^{1+\tau}} \\
  &\le& \frac{1}{1-D_{\tau}} \sum\limits_{i=1}^{N} \lrabs{a_{ii}}^{1+\tau} \frac{\blre{\lrabs{d^k(Z_1,Z_2)}^{2+2\tau}}}{\lrp{\cv_k^2(Z)}^{1+\tau}} \\
  &=& A_{n,m,\tau}^{(1)} \frac{\blre{\lrabs{d^k(Z_1.Z_2)}^{2+2\tau}}}{\lrp{\cv_k^2(Z)}^{1+\tau}}.
  \EEqn
\end{proof}

Combining the result of Lemma \ref{Lemma:fzeta_2} and Lemma \ref{Lemma:fzeta_1}, we obtain a bound of $B_{n,m,\tau}^{(1)}$, which is summarized in the following lemma.
\begin{lemma}\label{Lemma:B(1)}
    Under the null, with $0<\tau\le1$, it holds that
    \begin{equation*}
        B_{n,m,\tau}^{(1)} 
    \le \frac{2^{\tau}}{\lrp{\cv_k^2(Z)}^{1+\tau}} \lrp{A_{n,m,\tau}^{(1)} \blre{\lrabs{d^k(Z_1,Z_2)}^{2+2\tau}} + A_{n,m,\tau}^{(2)} \lrp{\blre{g^k(Z_1,Z_2,Z_3,Z_4)}}^{(1+\tau)/2}}.
\end{equation*}
\end{lemma}

Next we derive the upper bound of $B_{n,m,\tau}^{(2)}$.
\begin{lemma}\label{Lemma:B(2)}
    Under the null, for any $0<\tau\le1$, it holds that
    \begin{equation*}
        B_{n,m,\tau}^{(2)}
      = A_{n,m,\tau}^{(3)} \frac{\blre{\lrabs{d^k(Z_1,Z_2)}^{2+2\tau}}}{\lrp{\cv_k^2(Z)}^{1+\tau}},
    \end{equation*}
    where $A_{n,m,\tau}^{(3)} = B_{\tau} \sum\limits_{j=1}^{n+m}\lrp{1 + (j-1)^{\tau}} \sum\limits_{i=1}^{j-1} \lrabs{\tilde{\iota}_{ij}}^{2+2\tau}$ and $B_{\tau}$ is positive constant depending only on $\tau$.
\end{lemma}
\begin{proof}
    For any $Z_i,Z_j$, let $d^k(Z_i,Z_j)|Z_j$ denote the random variable conditioning on $Z_j$, then it's a random variable only of $Z_i$ with mean zero. With $0<\tau\le1$, it follows from Rosenthal's inequality that
    \BEqn
   & & \blre{\lrabs{\sum\limits_{i=1}^{j-1} \tilde{\iota}_{ij}  \frac{d^k(Z_i,Z_j)}{\sqrt{\cv_k^2(Z)}}}^{2+2\tau} | Z_j} \\
    &\le& B_{\tau} \lrp{\sum\limits_{i=1}^{j-1} \blre{\lrabs{\tilde{\iota}_{ij}     \frac{d^k(Z_i,Z_j)}{\sqrt{\cv_k^2(Z)}}}^{2+2\tau} | Z_j} + \lrp{\sum\limits_{i=1}^{j-1}\blre{\lrabs{\tilde{\iota}_{ij} \frac{d^k(Z_i,Z_j)}{\sqrt{\cv_k^2(Z)}}}^2 | Z_j}}^{1+\tau} },
    \EEqn
    where $B_{\tau}$ is a positive constant depending only on $\tau$. Taking expectation of the individual terms, we have
    \begin{equation*}
        \blre{\blre{\lrabs{\tilde{\iota}_{ij}     \frac{d^k(Z_i,Z_j)}{\sqrt{\cv_k^2(Z)}}}^{2+2\tau} | Z_j}}
    =  \lrabs{\tilde{\iota}_{ij}}^{2+2\tau}     \frac{\blre{\lrabs{d^k(Z_i,Z_j)}^{2+2\tau}}}{\lrp{\cv_k^2(Z)}^{1+\tau}},
    \end{equation*}
    and 
    \BEqn
        \blre{\lrp{\sum\limits_{i=1}^{j-1}\blre{\lrabs{\tilde{\iota}_{ij} \frac{d^k(Z_i,Z_j)}{\sqrt{\cv_k^2(Z)}}}^2 | Z_j}}^{1+\tau} } 
    &=& \blre{\lrp{e_j \frac{\blre{\lrp{d^k(Z_1,Z_j)}^2 | Z_j}}{\cv_k^2(Z)}}^{1+\tau}} \\
    &\le& \blre{e_j^{1+\tau} \frac{\blre{\lrabs{d^k(Z_1,Z_j)}^{2+2\tau} | Z_j}}{\lrp{\cv_k^2(Z)}^{1+\tau}}} \\
    &=& e_j^{1+\tau} \frac{\blre{\lrabs{d^k(Z_1,Z_2)}^{2+2\tau}}}{\lrp{\cv_k^2(Z)}^{1+\tau}},
    \EEqn
    where $e_j = \sum\limits_{i=1}^{j-1}\tilde{\iota}_{ij}^2$. It's trivial that $e_j^{1+\tau} \le (j-1)^{\tau}\sum\limits_{i=1}^{j-1}\lrabs{\tilde{\iota}_{ij}}^{2+2\tau}$.
    We combine the aforementioned results and finally obtain that
    \BEqn
        B_{n,m,\tau}^{(2)}
    &=& \sum\limits_{j=1}^{n+m} \blre{\blre{\lrabs{\sum\limits_{i=1}^{j-1} \tilde{\iota}_{ij} \frac{d^k(Z_i,Z_j)}{\sqrt{\cv_k^2(Z)}}}^{2+2\tau} | Z_j}} \\
    &\le& B_{\tau} \sum\limits_{j=1}^{n+m} \lrp{\sum\limits_{i=1}^{j-1} \lrabs{\tilde{\iota}_{ij}}^{2+2\tau} + e_j^{1+\tau}} \frac{\blre{\lrabs{d^k(Z_1,Z_2)}^{2+2\tau}}}{\lrp{\cv_k^2(Z)}^{1+\tau}} \\
    &\le& A_{n,m,\tau}^{(3)} \frac{\blre{\lrabs{d^k(Z_1,Z_2)}^{2+2\tau}}}{\lrp{\cv_k^2(Z)}^{1+\tau}},
    \EEqn 
    which completes the proof.
\end{proof}

To derive the explicit expressions of $B_{n,m,\tau}^{(1)}$ and $B_{n,m,\tau}^{(2)}$, it  remains to find the explicit orders of $A_{n,m,\tau}^{(1)}$, $A_{n,m,\tau}^{(2)}$ and $A_{n,m,\tau}^{(3)}$ as introduced in previous lemmas.

\begin{lemma}\label{Lemma:Aorder}
    Assume that there exists $0<\rho<1$ such that $n/N \rightarrow \rho$ as $n,m\rightarrow\infty$ for some constant $s>0$, then with $0<\tau\le1$ and sufficiently large $N$, it holds under the null that
    \begin{equation*}
        A_{n,m,\tau}^{(1)} \le C(\rho,\tau) N^{-\tau},
        \quad\quad
        A_{n,m,\tau}^{(2)} \le C(\rho,\tau),
        \quad\quad
        A_{n,m,\tau}^{(3)} \le C(\rho,\tau) N^{-\tau},
    \end{equation*}
    where $A_{n,m,\tau}^{(i)},i=1,2,3$ are defined as in Lemma \ref{Lemma:fzeta_2}, Lemma \ref{Lemma:fzeta_1} and Lemma \ref{Lemma:B(2)}.
\end{lemma}
\begin{proof}
    Under the assumption that $n/N \rightarrow \rho$ as $n,m\rightarrow\infty$, there exists some positive constant $C(\rho)$ that depends only on $\rho$, such that
    \begin{equation*}
        c_{n,m} 
     = \frac{2}{n(n-1)} + \frac{4}{nm} + \frac{2}{m(m-1)}
     \rightarrow C(\rho) N^{-2}
    \end{equation*}
    W hen $N$ is sufficiently large, there exists some constants, such that
    \begin{equation*}
        C_1(\rho)N^{-2} \le c_{n,m} \le C_2(\rho)N^{-2}
    \end{equation*}
    
    We start from $A_{n,m,\tau}^{(1)}$. Note that
    \begin{equation*}
        \sum\limits_{j=i+1}^{N} \iota_{ij}^2
      = \left\{
        \begin{array}{ll}
            \displaystyle{\frac{4(n-i)}{n^2(n-1)^2} + \frac{4}{n^2m}}, & 1 \le i \le n \vspace{3pt} \\
            \displaystyle{\frac{4(N-i)}{m^2(m-1)^2}}, & n+1 \le i \le N
        \end{array}
        \right.
    \end{equation*}
    it follows that
    \BEqn
   & & \sum\limits_{i=1}^{N} \lrabs{a_{ii}}^{1+\tau} \\
    &=& \frac{2^{1+\tau}}{c_{n,m}^{1+\tau}} \sum\limits_{i=1}^{N}\lrp{\sum\limits_{j=i+1}^{N} \iota_{ij}^2}^{1+\tau} \\
    &=& \frac{2^{1+\tau}}{c_{n,m}^{1+\tau}} \lrp{\sum\limits_{i=1}^{n}\lrp{\frac{4(n-i)}{n^2(n-1)^2} + \frac{4}{n^2m}}^{1+\tau} + \sum\limits_{i=n+1}^{N}\lrp{\frac{4(N-i)}{m^2(m-1)^2}}^{1+\tau}} \\
    &\le& \frac{2^{1+\tau}}{c_{n,m}^{1+\tau}} \lrp{2^{\tau}\sum\limits_{i=1}^{n}\lrp{\frac{4(n-i)}{n^2(n-1)^2}}^{1+\tau} + 2^{\tau}\sum\limits_{i=1}^{n}\lrp{\frac{4}{n^2m}}^{1+\tau} + \sum\limits_{i=n+1}^{N}\lrp{\frac{4(N-i)}{m^2(m-1)^2}}^{1+\tau}} \\
    &\le& \frac{2^{1+\tau}}{c_{n,m}^{1+\tau}} \lrp{\frac{2^{2+3\tau}}{n^{\tau}(n-1)^{2+2\tau}} + \frac{2^{2+3\tau}}{n^{1+2\tau}m^{1+\tau}} + \frac{2^{2+2\tau}}{m^{\tau}(m-1)^{2+2\tau}}} \\
    &\le& C(\rho,\tau) N^{-\tau},
    \EEqn
    and consequently, 
    \begin{equation*}
        A_{n,m,\tau}^{(1)}
      = \frac{1}{1-D_{\tau}}\sum\limits_{i=1}^{N} \lrabs{a_{ii}}^{1+\tau}
      \le C(\rho,\tau) N^{-\tau}
    \end{equation*}
    as $n,m\rightarrow\infty$.
    
    Then we consider $A_{n,m,\tau}^{(2)}$. Observe that
    \begin{equation*}
        \sum\limits_{j=i_2+1}^{N} \iota_{i_1,j} \iota_{i_2,j} 
      = \left\{
        \begin{array}{ll}
            \displaystyle{\frac{4(n-i_2)}{n^2(n-1)^2} + \frac{4}{n^2m}},  & 1 \le i_1<i_2 \le n \vspace{3pt} \\
            \displaystyle{-\frac{4(N-i_2)}{nm^2(m-1)}}, & 1 \le i_1 \le n < i_2 \le N \vspace{3pt} \\
            \displaystyle{\frac{4(N-i_2)}{m^2(m-1)^2}}, & n+1 \le i_1<i_2 \le N
        \end{array}
        \right.
    \end{equation*}
    then we have
    \BEqn
   & & \sum\limits_{1 \le i_1<i_2 \le N} a_{i_1,i_2}^2 \\
    &=& \frac{4}{c_{n,m}^2} \left(\sum\limits_{1\le i_1<i_2\le n}\lrp{\frac{4(n-i_2)}{n^2(n-1)^2} + \frac{4}{n^2m}}^2 + \sum\limits_{i_1=1}^{n}\sum\limits_{i_2=n+1}^{N}\lrp{-\frac{4(N-i_2)}{nm^2(m-1)}}^2 \right. \\
   & & \quad\quad\quad \left. + \sum\limits_{n+1\le i_1<i_2 \le N}\lrp{\frac{4(N-i_2)}{m^2(m-1)^2}}^2 \right) \\
    &\le& \frac{4}{c_{n,m}^2} \left( \sum\limits_{i_2=1}^{n}\sum\limits_{i_1=1}^{i_2-1} \frac{32(n-i_2)^2}{n^4(n-1)^4} + \sum\limits_{i_2=1}^{n}\sum\limits_{i_1=1}^{i_2-1} \frac{32}{n^4 m^2} + \sum\limits_{i_1=1}^{n}\sum\limits_{i_2=n+1}^{N} \frac{16(N-i_2)^2}{n^2 m^4 (m-1)^2} \right. \\
   & & \quad\quad\quad \left. + \sum\limits_{i_2=n+1}^{N}\sum\limits_{i_1=n+1}^{i_2-1} \frac{16(N-i_2)^2}{m^4 (m-1)^4} \right) \\
    &\le& \frac{4}{c_{n,m}^2} \lrp{\frac{32}{(n-1)^4} + \frac{16}{n^2 m^2} + \frac{16}{nm(m-1)^2} + \frac{16}{(m-1)^4}} \\
    &\le& C(\rho),
    \EEqn
    which implies that
    \begin{equation*}
        A_{n,m,\tau}^{(2)}
      = \lrp{\sum\limits_{1\le i_1<i_2 \le N} a_{i_1,i_2}^2}^{(1+\tau)/2}
      \le C(\rho,\tau).
    \end{equation*}
    
    Lastly, we calculate the order of $A_{n,m,\tau}^{(3)}$. Note that
    \begin{equation*}
        \sum\limits_{i=1}^{j-1} \lrabs{\tilde{\iota}_{ij}}^{2+2\tau}
      = \frac{1}{c_{n,m}^{1+\tau}} \sum\limits_{i=1}^{j-1} \lrabs{\iota_{ij}}^{2+2\tau},
    \end{equation*}
    then we have
    \BEqn
        \sum\limits_{j=1}^{N}\sum\limits_{i=1}^{j-1} \lrabs{\tilde{\iota}_{ij}}^{2+2\tau}
    &=& \frac{1}{c_{n,m}^{1+\tau}} \sum\limits_{1\le i<j \le N} \lrabs{\iota_{ij}}^{2+2\tau} \\
    &=& \frac{1}{c_{n,m}^{1+\tau}} \lrp{\frac{2^{1+2\tau}}{n^{1+2\tau}(n-1)^{1+2\tau}} + \frac{2^{2+2\tau}}{n^{1+2\tau}m^{1+2\tau}} + \frac{2^{1+2\tau}}{m^{1+2\tau}(m-1)^{1+2\tau}}} \\
    &\le& C(\rho,\tau) N^{-2\tau},
    \EEqn
    and consequently,
    \begin{equation*}
        A_{n,m,\tau}^{(3)} 
      = B_{\tau} \sum\limits_{j=1}^{N}\lrp{1 + (j-1)^{\tau}} \sum\limits_{i=1}^{j-1} \lrabs{\tilde{\iota}_{ij}}^{2+2\tau}
    \le 2 B_{\tau} N^{\tau} \sum\limits_{j=1}^{N} \sum\limits_{i=1}^{j-1} \lrabs{\tilde{\iota}_{ij}}^{2+2\tau} 
    = C(\rho,\tau) N^{-\tau}.
    \end{equation*}
\end{proof}

Based on the previous results, we obtain an upper bound of $P_1$.
\begin{lemma}\label{Lemma:P1}
  Assume that $n/N \rightarrow \rho$ for some $0<\rho<1$ as $n,m\rightarrow\infty$, then with $0<\tau\le1$ it holds under the null that 
\begin{equation*}
    P_1
\le C(\rho,\tau)
    \lrp{\frac{\blre{\lrabs{\tilde{k}(Z_1,Z_2)}^{2+2\tau}}}{N^{\tau}\lrp{\cv_k^2(Z)}^{1+\tau}} + \lrp{\frac{\blre{g^k(Z_1,Z_2,Z_3,Z_4)}}{\lrp{\cv_k^2(Z)}^2}}^{(1+\tau)/2} }^{1/(3+2\tau)},
\end{equation*}
where $C(\rho,\tau)$ is a positive constant depending only on $\rho$ and $\tau$.
\end{lemma}
\begin{proof}
  Recall that $\tilde{k}(Z_1,Z_2) = k(Z_1,Z_2) - \bbe{k(Z_1,Z_2)}$, then by noting that
  \begin{equation*}
      d^k(Z_1,Z_2)
      = \tilde{k}(Z_1,Z_2)
      - \blre{\tilde{k}(Z_1,Z_2)|Z_1}
      - \blre{\tilde{k}(Z_1,Z_2)|Z_2},
  \end{equation*}
  we have $\blre{\lrabs{d^k(Z_1,Z_2)}^{2+2\tau}} \le 3^{2+2\tau}\blre{\lrabs{\tilde{k}(Z_1,Z_2)}^{2+2\tau}}$, and the bound directly follows from Lemma \ref{Lemma:P1-inter}, Lemma \ref{Lemma:B(1)}, Lemma \ref{Lemma:B(2)} and Lemma \ref{Lemma:Aorder}.
\end{proof}

To prove Theorem \ref{Thm:BEbd-k}, it remains to find the bounds of $P_2,P_3$ and $P_4$, which are summarized in the following lemmas.
\begin{lemma}\label{Lemma:P2P3}
For any $0<\gamma<1$, there exists some positive constant $C<\infty$ such that $P_2 + P_3 \le C \gamma$.
\end{lemma}
\begin{proof}
    Note that $\displaystyle{\sqrt{1+\gamma} = 1 + \frac{\gamma}{2} - \frac{\gamma^2}{8} + o(\gamma^2)}$ for $\gamma$ close to zero, then for any $0<\gamma<1$, let $\delta = \sqrt{1+\gamma}$, and we have $|\delta-1| \le C\gamma$ for some positive constant $C<\infty$ when $\gamma$ is close to zero. There exists some $x^{\ast}$ between $x$ and $x\sqrt{1+\gamma}$, such that
    \begin{equation*}
        \lrabs{\Phi(x) - \Phi(x\sqrt{1+\gamma})}
      = \phi(x^{\ast}) \lrabs{x\sqrt{1+\gamma} - x}
      = \phi(x^{\ast}) \lrabs{x} \cdot \lrabs{\delta-1}
      \le C \gamma \phi(x^{\ast})|x|,
    \end{equation*}
    If $|x|\le2$, then there exists $C>0$, such that $\phi(x) \le C\phi(0)$, and thus $\lrabs{\Phi(x) - \Phi(x\sqrt{1+\gamma})} \le C\gamma$. If $|x|>2$, note that $|x|\phi(x)$ is decreasing in $|x|$, then $|x|\phi(x) \le 2\phi(2)$, and thus $\lrabs{\Phi(x) - \Phi(x\sqrt{1+\gamma})} \le C\gamma$. This implies that
    \begin{equation*}
        P_2
      = \sup\limits_{x\in\br} \lrabs{\Phi(x) - \Phi(x\sqrt{1+\gamma})} 
    \le C\gamma.
    \end{equation*}
    Also note that a similar argument may apply to $\lrabs{\Phi(x) - \Phi(x\sqrt{1-\gamma})}$, then we may conclude that $P_2+P_3 \le C\gamma$ for some constant $C>0$.
\end{proof}

\begin{lemma}\label{Lemma:P4}
For any $\gamma>0$ and $0<\tau\le1$, assume that $\blre{\lrabs{k(X_1,X_2)}^{2+2\tau}} + \blre{\lrabs{k(Y_1,Y_2)}^{2+2\tau}} < \infty$, then under the null, it holds that
\begin{equation*}
    P_4 
  \le 2^{\tau}
      \lrp{\frac{\blre{\lrabs{\tilde{k}(Z_1,Z_2)}^{2+2\tau}}}{N^{\tau} \gamma^{1+\tau} \lrp{\cv_k^2(Z)}^{1+\tau}}}.
\end{equation*}
\end{lemma}
\begin{proof}
    Under the null, we have $\Var(L_{n,m}^k) = c_{n,m} \cv_k^2(Z)$, and it directly follows from the proof of Proposition \ref{Prop:ratio-const-k-H0} that for any $0<\tau\le1$, we have
    \begin{equation*}
        P_4
      = \blrp{\lrabs{\frac{c_{n,m} \cv_{n,m}^{k\ast}(X,Y)}{\Var(L_{n,m}^k)} - 1} > \gamma} \\
      \le 2^{\tau}
          \lrp{\frac{\blre{\lrabs{\tilde{k}(Z_1,Z_2)}^{2+2\tau}}}{N^{\tau} \gamma^{1+\tau} \lrp{\cv_k^2(Z)}^{1+\tau}}},
    \end{equation*}
    which completes the proof.
\end{proof}

\subsection{Lemmas for Proposition \ref{Prop:BEbd-f}}

\begin{lemma}\label{Lemma:RateConvg-1}
Suppose Assumption \ref{Assumpt:uniform-kernel}\ref{Assumpt:uniform-kernel-1}-\ref{Assumpt:uniform-kernel-2} and Assumption \ref{Assumpt:component-dept}\ref{Assumpt:component-dept-1}-\ref{Assumpt:component-dept-2} hold, then under the null, it holds for any $k = k^p \in \cc$ that
\begin{equation*}
    |\cv_k^2(Z) - 4c_1^2(A_0) A^{-2} \|\Sigma\|_F^2|
\le C(\tilde{M},U^{\ast},L_0,U_0) c_0^2(A_0) \lrp{\frac{\alpha(p)}{p}}^2,
\end{equation*}
where $c_1(x)$ is defined as Lemma \ref{Lemma:approx-2} and $\Sigma_Z = \Cov(Z)$ denotes the covariance matrix of $Z$.
\end{lemma}
\begin{proof}
  Under the null, it holds that $X =^d Y =^d Z$. With $s=|Z_1-Z_2|$ and $s_0=A_0$, it follows from Lemma \ref{Lemma:approx-1} that
  \begin{equation*}
      k(Z_1,Z_2) = f(|Z_1-Z_2|) = \sum\limits_{i=0}^{6} c_i(A_0) A^{-i} \bigp{|Z_1-Z_2|^2-A}^i + R(Z_1,Z_2),
  \end{equation*}
  where 
  \begin{equation*}
      |R(Z_1,Z_2)|
  \le C(\tilde{M}) \sum\limits_{i=7}^{16} A^{-i} \lrabs{|Z_1-Z_2|^2-A}^i |c_0(A_0)|.
  \end{equation*}
  
  Under Assumption \ref{Assumpt:uniform-kernel}\ref{Assumpt:uniform-kernel-2}, it holds that $\max\limits_{1\le i\le 7}|c_i(A_0)| \le C(\tilde{M}) |c_0(A_0)|$. Also, it follows from Lemma \ref{Lemma:Li-3}, Lemma \ref{Lemma:Li_bar-1} and Assumption \ref{Assumpt:component-dept}\ref{Assumpt:component-dept-2} that 
  \begin{equation*}
      A^{-i} |L_i| \le C(i,U^{\ast},L_0,U_0) \lrp{\frac{\alpha(p)}{p}}^{\lrceil{i/2}},
  \end{equation*}
  and 
  \begin{equation*}
      A^{-i} \bar{L}_i \le C(i,U^{\ast},L_0,U_0) \lrp{\frac{\alpha(p)}{p}}^{i/2}.
  \end{equation*}
  Under the assumption that $1\le\alpha\le p$, then we have 
  \begin{equation*}
      \max\lrcp{\max\limits_{3\le i\le 6} \{A^{-i}| L_i|\},
            \max\limits_{7\le i\le 16} \{A^{-i} \bar{L}_i\}}
  \le C(U^{\ast},L_0,U_0) \lrp{\frac{\alpha(p)}{p}}^{2},
  \end{equation*}
  By definition, we have $A = \bbe{|Z_1-Z_2|^2}$, then it follows that
  \BEqn
  & & \lrabs{\bbe{k(Z_1,Z_2)} - \bigp{c_0(A_0) + c_2(A_0)A^{-2}L_2}} \\
  &\le& \sum\limits_{i=3}^{6} |c_i(A_0)| A^{-i}|L_i| 
        + C(\tilde{M}) \sum\limits_{i=7}^{16} |c_0(A_0)| A^{-i} \bar{L}_i \\
  &\le& C(\tilde{M},U^{\ast},L_0,U_0) |c_0(A_0)| \lrp{\frac{\alpha(p)}{p}}^2,
  \EEqn
  
  Using similar techniques, we obtain that
  \BEqn
  & & \lrabs{\bbe{k^2(Z_1,Z_2)} - \bigp{c_0^2(A_0) + (2c_0(A_0)c_2(A_0) + c_1^2(A_0))A^{-2}L_2}} \\
  &\le& C(\tilde{M},U^{\ast},L_0,U_0) c_0^2(A_0) \lrp{\frac{\alpha(p)}{p}}^2, \\[2mm]
  & & \lrabs{\bbe{k(Z_1,Z_2)k(Z_1,Z_3)} - \bigp{c_0^2(A_0) + 2c_0(A_0)c_2(A_0)A^{-2}L_2 + c_1^2(A_0)A^{-2}N_{11}}} \\
  &\le& C(\tilde{M},U^{\ast},L_0,U_0) c_0^2(A_0) \lrp{\frac{\alpha(p)}{p}}^2, \\[2mm]
  & & \lrabs{\bbe{k(Z_1,Z_2)}^2 - \bigp{c_0^2(A_0) + 2c_0(A_0)c_2(A_0)A^{-2}L_2}} \\
  &\le& C(\tilde{M},U^{\ast},L_0,U_0) c_0^2(A_0) \lrp{\frac{\alpha(p)}{p}}^2.
  \EEqn
  
  Recall that $\cv_k^2(Z) = \bbe{k^2(Z_1,Z_2)} - 2\bbe{k(Z_1,Z_2)k(Z_1,Z_3)} + \bbe{k(Z_1,Z_2)}^2$, then from above we have
  \begin{equation*}
      \lrabs{\cv_k^2(Z) - c_1^2(A_0)A^{-2}(L_2 - 2N_{11})}
  \le C(\tilde{M},U^{\ast},L_0,U_0) c_0^2(A_0) \lrp{\frac{\alpha(p)}{p}}^2.
  \end{equation*}
  It remains to investigate $L_2 - 2N_{11}$. By definition, it holds that
  \BEqn
      L_2 - 2N_{11}
  &=& \bbe{\bigp{|Z_1-Z_2|^2-A}^2} - 2\bbe{\bigp{|Z_1-Z_2|^2-A} \bigp{|Z_1-Z_3|^2-A}} \\
  &=& 4\sum\limits_{j_1,j_2=1}^{p} \bbe{\tilde{z}_{1 j_1} \tilde{z}_{1 j_2}}^2 \\
  &=& 4 \|\Sigma_Z\|_F^2,
  \EEqn
  which completes the proof.
\end{proof}

\begin{lemma}\label{Lemma:RateConvg-2}
Suppose Assumption \ref{Assumpt:uniform-kernel}\ref{Assumpt:uniform-kernel-1}-\ref{Assumpt:uniform-kernel-2} and Assumption \ref{Assumpt:component-dept}\ref{Assumpt:component-dept-1}-\ref{Assumpt:component-dept-2} hold, then under the null, it holds for any $k = k^p \in \cc$ that
\begin{equation*}
    \bbe{\tilde{k}^4(Z_1,Z_2)}
\le C(\tilde{M},U^{\ast},L_0,U_0) 
    \lrp{c_1^4(A_0) \lrp{\frac{\alpha(p)}{p}}^2 + c_0^4(A_0) \lrp{\frac{\alpha(p)}{p}}^3},
\end{equation*}
where $c_0(x)$ and $c_1(x)$ are defined as Lemma \ref{Lemma:approx-2}.
\end{lemma}
\begin{proof}
  Recall that $\tilde{k}(Z_1,Z_2) = k(Z_1,Z_2) - \bbe{k(Z_1,Z_2)}$, and it follows from Lemma \ref{Lemma:approx-2} that
  \BEqn
  & & k(Z_1,Z_2) 
   =  \sum\limits_{i=0}^{6} c_i(A_0) A^{-i} \bigp{|Z_1-Z_2|^2-A}^i + R(Z_1,Z_2), \\
  & & \bbe{k(Z_1,Z_2)}
   =  \sum\limits_{i=0}^{6} c_i(A_0) A^{-i} L_i + \bbe{R(Z_1,Z_2)},
  \EEqn
  where
  \begin{equation*}
      |R(Z_1,Z_2)|
  \le C(\tilde{M}) \sum\limits_{i=7}^{16} A^{-i} \lrabs{|Z_1-Z_2|^2-A}^i |c_0(A_0)|.
  \end{equation*}
  Then it follows that
  \begin{equation*}
      \tilde{k}(Z_1,Z_2)
    = \sum\limits_{i=1}^{6} c_i(A_0) A^{-i} \lrp{\bigp{|Z_1-Z_2|^2-A}^i - L_i} + \bigp{R(Z_1,Z_2) - \bbe{R(Z_1,Z_2)}},
  \end{equation*}
  and consequently,
  \BEqn
  & & \tilde{k}^4(Z_1,Z_2) \\
  &=& \sum\limits_{\substack{s_1+\dots+s_7=4 \\ s_1,\dots,s_7\geq0}} \bin{4}{s_1 \cdots s_7} 
      \prod\limits_{i=1}^{6} \lrp{c_i(A_0) A^{-i} \bigp{\bigp{|Z_1-Z_2|^2-A}^i - L_i}}^{s_i} \bigp{R(Z_1,Z_2) - \bbe{R(Z_1,Z_2)}}^{s_7} \\
  &=& \sum\limits_{\substack{s_1+\dots+s_6=4 \\ s_1,\dots,s_6\geq0, s_7=0}} \bin{4}{s_1 \cdots s_6\ 0} 
      \prod\limits_{i=1}^{6} \lrp{c_i(A_0) A^{-i} \bigp{\bigp{|Z_1-Z_2|^2-A}^i - L_i}}^{s_i} \\
  &+& \sum\limits_{\substack{s_1+\dots+s_7=4 \\ s_1,\dots,s_6\geq0, s_7>0}} \bin{4}{s_1 \cdots s_7} 
      \prod\limits_{i=1}^{6} \lrp{c_i(A_0) A^{-i} \bigp{\bigp{|Z_1-Z_2|^2-A}^i - L_i}}^{s_i} \bigp{R(Z_1,Z_2) - \bbe{R(Z_1,Z_2)}}^{s_7} \\      
  &=& \ck_1 + \ck_2,
  \EEqn
  where $\ck_1,\ck_2$ denote the two summations in the second to the last step respectively. 
  
  To find an upper bound of $\bbe{\tilde{k}^4(Z_1,Z_2)}$, now it suffices to find the upper bounds of $\bbe{\ck_1}$ and $\bbe{\ck_2}$.
  
  We first look into $\bbe{\ck_1}$. By expanding each $\bigp{\bigp{|Z_1-Z_2|^2-A}^i - L_i}^{i_s}$ using the binomial equation, we obtain that
  \BEqn
  & & \bbe{\ck_1} \\
  &=& \sum\limits_{\substack{s_1+\dots+s_6=4 \\ s_1,\dots,s_6\geq0, s_7=0}} \bin{4}{s_1 \cdots s_6\ 0} 
      \bigp{\prod\limits_{i=1}^{6} c_i^{s_i}(A_0)} A^{-(s_1+2s_2+3s_3+4s_4+5s_5+6s_6)} \\
  & & \hspace{2em}
      \times \sum\limits_{t_2=0}^{s_2} \cdots \sum\limits_{t_6=0}^{s_6} \bigp{\prod\limits_{i=2}^{6} \bin{s_i}{t_i}}
      L_{s_1+2t_2+3t_3+4t_4+5t_5+6t_6} \bigp{\prod\limits_{i=2}^{6}(-L_i)^{s_i-t_i}} \\
  &=& c_1^4(A_0)A^{-4} L_4 \\
  & & + \sum\limits_{\substack{s_1+\dots+s_6=4 \\ s_1,\dots,s_6\geq0, s_1<4, s_7=0}} \bin{4}{s_1 \cdots s_6\ 0} 
      \bigp{\prod\limits_{i=1}^{6} c_i^{s_i}(A_0)} \sum\limits_{t_2=0}^{s_2} \cdots \sum\limits_{t_6=0}^{s_6} \bigp{\prod\limits_{i=2}^{6} \bin{s_i}{t_i}} \\
  & & \hspace{2em}
      \times
      A^{-(s_1+2t_2+3t_3+4t_4+5t_5+6t_6)}L_{s_1+2t_2+3t_3+4t_4+5t_5+6t_6} \prod\limits_{i=2}^{6}(-A^{-i}L_i)^{s_i-t_i},
  \EEqn
  thus by Lemma \ref{Lemma:Li-3} we have
  \BEqn
  & & |\bbe{\ck_1}| \\
  &\le& c_1^4(A_0)A^{-4} L_4 \\
  & &   + C(\tilde{M},U^{\ast},L_0,U_0) c_0^4(A_0)
        \sum\limits_{\substack{s_1+\dots+s_6=4 \\ s_1,\dots,s_6\geq0 \\ s_1<4, s_7=0}}  \sum\limits_{i=2}^{6}
        \sum\limits_{t_i=0}^{s_i} \lrp{\frac{\alpha(p)}{p}}^{\lrceil{(s_1+2t_2+3t_3+4t_4+5t_5+6t_6)/2}+\sum\limits_{i=2}^{6}(s_i-t_i)\lrceil{i/2}} \\
  &\le& C(\tilde{M},U^{\ast},L_0,U_0)
        \lrp{c_1^4(A_0) \lrp{\frac{\alpha(p)}{p}}^2 + c_0^4(A_0) \lrp{\frac{\alpha(p)}{p}}^3}.
  \EEqn
  
  Similarly,
  \BEqn
  & & \bbe{\ck_2} \\
  &=& \sum\limits_{\substack{s_1+\dots+s_7=4 \\ s_1,\dots,s_6\geq0, s_7>0}} \bin{4}{s_1 \cdots s_7}
      \bigp{\prod\limits_{i=1}^{6} c_i^{s_i}(A_0)} A^{-(s_1+2s_2+3s_3+4s_4+5s_5+6s_6)} \\
  & & \hspace{2em}
      \times \sum\limits_{t_2=0}^{s_2} \cdots \sum\limits_{t_7=0}^{s_7} \bigp{\prod\limits_{i=2}^{7} \bin{s_i}{t_i}}
      \bbe{\bigp{|Z_1-Z_2|^2-A}^{s_1+2t_2+3t_3+4t_4+5t_5+6t_6} R^{t_7}(Z_1,Z_2)} \\
  & & \hspace{2em}
      \times \bigp{\prod\limits_{i=2}^{6}(-L_i)^{s_i-t_i} (-\bbe{R(Z_1,Z_2)})^{s_7-t_7}},
  \EEqn
  then by using the fact that 
  \begin{equation*}
      |R(Z_1,Z_2)|
  \le C(\tilde{M}) \sum\limits_{i=7}^{16} A^{-i} \lrabs{|Z_1-Z_2|^2-A}^i |c_0(A_0)|,
  \end{equation*}
  we can obtain an upper bound of $\bbe{\ck_2}$, that is,
  \begin{equation*}
      |\bbe{\ck_2}| 
  \le C(\tilde{M},U^{\ast},L_0,U_0) c_0^4(A_0) \lrp{\frac{\alpha(p)}{p}}^3,
  \end{equation*}
  
  In summary, we conclude that
  \BEqn
      \bbe{\tilde{k}^4(Z_1,Z_2)} 
  &\le& |\bbe{\ck_1}| + |\bbe{\ck_2}| \\
  &\le& C(\tilde{M},U^{\ast},L_0,U_0) \lrp{c_1^4(A_0) \lrp{\frac{\alpha(p)}{p}}^2 + c_0^4(A_0) \lrp{\frac{\alpha(p)}{p}}^3}.
  \EEqn
\end{proof}

\begin{lemma}\label{Lemma:RateConvg-3}
Suppose Assumption \ref{Assumpt:uniform-kernel}\ref{Assumpt:uniform-kernel-1}-\ref{Assumpt:uniform-kernel-2} and Assumption \ref{Assumpt:component-dept}\ref{Assumpt:component-dept-1}-\ref{Assumpt:component-dept-2} hold, then under the null, it holds for any $k = k^p \in \cc$ that
\begin{equation*}
    |\bbe{g^k(Z_1,Z_2,Z_3,Z_4)}|
\le C(\tilde{M},U^{\ast},L_0,U_0) c_0^4(A_0) \lrp{\frac{\alpha(p)}{p}}^3,
\end{equation*}
where $c_0(x)$ is defined as Lemma \ref{Lemma:approx-2}.
\end{lemma}
\begin{proof}
  To obtain the upper bound of $|\bbe{g^k(Z_1,Z_2,Z_3,Z_4)}|$, it suffices to find the upper bound of $G_1,G_2,G_3$ and $G_4$ proposed in Proposition \ref{Prop:gk}. We first look into $G_2$. By using the approximation derived in Lemma \ref{Lemma:approx-2} as well as the results obtained in Lemma \ref{Lemma:Li-3} and Lemma \ref{Lemma:NDM-1}, we have
  \BEqn
  & & \lrabs{\bbe{k(Z_1,Z_2)} - \sum\limits_{i_1=0}^{6} c_{i_1}(A_0) A^{-i_1} L_{i_1}} \\
  &\le& C(\tilde{M},U^{\ast},L_0,U_0) |c_0(A_0)| \lrp{\frac{\alpha(p)}{p}}^3, \\[2mm]
  & & \lrabs{\bbe{k(Z_1,Z_2) k(Z_1,Z_3) k(Z_2,Z_4)} - \sum\limits_{i_2,i_3,i_4=0}^{6} c_{i_2}(A_0) c_{i_3}(A_0) c_{i_4}(A_0) A^{-(i_2+i_3+i_4)} D_{i_2 i_3 i_4}} \\
  &\le& C(\tilde{M},U^{\ast},L_0,U_0) |c_0^3(A_0)| \lrp{\frac{\alpha(p)}{p}}^3,
  \EEqn
  and it follows from the definition of $G_2$ that
  \BEqn
      \lrabs{G_2 - \tilde{G}_2}
  \le C(\tilde{M},U^{\ast},L_0,U_0) |c_0^4(A_0)| \lrp{\frac{\alpha(p)}{p}}^3,
  \EEqn
  where 
  \begin{equation*}
      \tilde{G}_2 
    = 4\sum\limits_{i_1,i_2,i_3,i_4=0}^{6} c_{i_1}(A_0) c_{i_2}(A_0) c_{i_3}(A_0) c_{i_4}(A_0) A^{-(i_1+i_2+i_3+i_4)} L_{i_1}D_{i_2 i_3 i_4}.
  \end{equation*}
  
  Similarly, we obtain the leading term of $G_1,G_3$ and $G_4$ after some calculations:
  \BEqn
  & & \tilde{G}_1
   =  \sum\limits_{i_1,i_2,i_3,i_4=0}^{6} c_{i_1}(A_0) c_{i_2}(A_0) c_{i_3}(A_0) c_{i_4}(A_0) A^{-(i_1+i_2+i_3+i_4)} \bigp{M_{i_1 i_2 i_3 i_4} -4\tilde{M}_{i_1 i_2 i_3 i_4} + 2N_{i_1 i_2} N_{i_3 i_4}}, \\
  & & \tilde{G}_3
   =  -4\sum\limits_{i_1,i_2,i_3,i_4=0}^{6} c_{i_1}(A_0) c_{i_2}(A_0) c_{i_3}(A_0) c_{i_4}(A_0)  A^{-(i_1+i_2+i_3+i_4)} L_{i_1} L_{i_2} N_{i_3 i_4}, \\
  & & \tilde{G}_4
   =  \sum\limits_{i_1,i_2,i_3,i_4=0}^{6} c_{i_1}(A_0) c_{i_2}(A_0) c_{i_3}(A_0) c_{i_4}(A_0) A^{-(i_1+i_2+i_3+i_4)} L_{i_1} L_{i_2} L_{i_3} L_{i_4},
  \EEqn
  and
  \begin{equation*}
      \max\{|G_1-\tilde{G}_1|, |G_3-\tilde{G}_3|, |G_4-\tilde{G}_4|\} 
  \le C(\tilde{M},U^{\ast},L_0,U_0) |c_0^4(A_0)| \lrp{\frac{\alpha(p)}{p}}^3.
  \end{equation*}
  
  Let $\tilde{G} = \sum\limits_{i=1}^{4} \tilde{G}_i$, then it follows that 
  \begin{equation*}
      \lrabs{\bbe{g^k(Z_1,Z_2,Z_3,Z_4)}-\tilde{G}} \le C(\tilde{M},U^{\ast},L_0,U_0) |c_0^4(A_0)| \lrp{\frac{\alpha(p)}{p}}^3.
  \end{equation*}
  To obtain the upper bound of $|\bbe{g^k(Z_1,Z_2,Z_3,Z_4)}|$, now it remains to investigate the order of $\tilde{G}$. 
  
  Note that
  \BEqn
  & & \tilde{G} \\
  &=& \sum\limits_{i_1,i_2,i_3,i_4=0}^{6} c_{i_1}(A_0) c_{i_2}(A_0) c_{i_3}(A_0) c_{i_4}(A_0)  A^{-(i_1+i_2+i_3+i_4)} \\
  & & \hspace{4em}
      \times \bigp{M_{i_1 i_2 i_3 i_4} -4\tilde{M}_{i_1 i_2 i_3 i_4} + 2N_{i_1 i_2} N_{i_3 i_4} + 4L_{i_1}D_{i_2 i_3 i_4} -4L_{i_1} L_{i_2} N_{i_3 i_4} + L_{i_1} L_{i_2} L_{i_3} L_{i_4}} \\
  &=& \sum\limits_{s=0}^{4}
      \sum\limits_{\substack{i_1+\dots+i_4=s \\ 0\le i_1,\dots,i_4\le6}}
      c_{i_1}(A_0) c_{i_2}(A_0) c_{i_3}(A_0) c_{i_4}(A_0) A^{-s} \\
  & & \hspace{4em}
      \times \bigp{M_{i_1 i_2 i_3 i_4} -4\tilde{M}_{i_1 i_2 i_3 i_4} + 2N_{i_1 i_2} N_{i_3 i_4} + 4L_{i_1}D_{i_2 i_3 i_4} -4L_{i_1} L_{i_2} N_{i_3 i_4} + L_{i_1} L_{i_2} L_{i_3} L_{i_4}} \\
  & & + \sum\limits_{\substack{i_1+\dots+i_4\geq5 \\ 0\le i_1,\dots,i_4\le6}}
      c_{i_1}(A_0) c_{i_2}(A_0) c_{i_3}(A_0) c_{i_4}(A_0) A^{-(i_1+i_2+i_3+i_4)} \\
  & & \hspace{4em}
      \times \bigp{M_{i_1 i_2 i_3 i_4} -4\tilde{M}_{i_1 i_2 i_3 i_4} + 2N_{i_1 i_2} N_{i_3 i_4} + 4L_{i_1}D_{i_2 i_3 i_4} -4L_{i_1} L_{i_2} N_{i_3 i_4} + L_{i_1} L_{i_2} L_{i_3} L_{i_4}}.
  \EEqn
  
  After some tedious calculations, it follows from the definition of $L_{i_1},N_{i_1 i_2},D_{i_1 i_2 i_3}, M_{i_1 i_2 i_3 i_4}$ and $\tilde{M}_{i_1 i_2 i_3 i_4}$ that
  \begin{equation*}
      \sum\limits_{\substack{i_1+\dots+i_4=s \\ 0\le i_1,\dots,i_4\le6}}
      \bigp{M_{i_1 i_2 i_3 i_4} -4\tilde{M}_{i_1 i_2 i_3 i_4} + 2N_{i_1 i_2} N_{i_3 i_4} + 4L_{i_1}D_{i_2 i_3 i_4} -4L_{i_1} L_{i_2} N_{i_3 i_4} + L_{i_1} L_{i_2} L_{i_3} L_{i_4}} = 0
  \end{equation*} 
  for $s=0,1,2,3$. Additionally, when $i_1+i_2+i_3+i_4\geq5$, it follows from Lemma \ref{Lemma:Li-3} and Lemma \ref{Lemma:NDM-1} that
  \BEqn
      \lrabs{A^{-(i_1+i_2+i_3+i_4)}
      \bigp{M_{i_1 i_2 i_3 i_4} -4\tilde{M}_{i_1 i_2 i_3 i_4}}}
  &\le& C(U^{\ast},L_0,U_0) \lrp{\frac{\alpha(p)}{p}}^{\lrceil{(i_1+i_2+i_3+i_4)/2}} \\
  &\le& C(U^{\ast},L_0,U_0) \lrp{\frac{\alpha(p)}{p}}^3, \\[2mm]
      \lrabs{A^{-(i_1+i_2+i_3+i_4)} 2N_{i_1 i_2} N_{i_3 i_4}}
  &\le& C(U^{\ast},L_0,U_0) \lrp{\frac{\alpha(p)}{p}}^{\lrceil{(i_1+i_2)/2}+\lrceil{(i_3+i_4)/2}} \\
  &\le& C(U^{\ast},L_0,U_0) \lrp{\frac{\alpha(p)}{p}}^3, \\[2mm]
      \lrabs{A^{-(i_1+i_2+i_3+i_4)} 4L_{i_1}D_{i_2 i_3 i_4}}
  &\le& C(U^{\ast},L_0,U_0) \lrp{\frac{\alpha(p)}{p}}^{\lrceil{i_1/2}+\lrceil{(i_2+i_3+i_4)/2}} \\
  &\le& C(U^{\ast},L_0,U_0) \lrp{\frac{\alpha(p)}{p}}^3, \\[2mm]      
      \lrabs{A^{-(i_1+i_2+i_3+i_4)} \bigp{-4L_{i_1} L_{i_2} N_{i_3 i_4}}}
  &\le& C(U^{\ast},L_0,U_0) \lrp{\frac{\alpha(p)}{p}}^{\lrceil{i_1/2}+\lrceil{i_2/2}+\lrceil{(i_3+i_4)/2}} \\
  &\le& C(U^{\ast},L_0,U_0) \lrp{\frac{\alpha(p)}{p}}^3, \\[2mm]      
      \lrabs{A^{-(i_1+i_2+i_3+i_4)} L_{i_1} L_{i_2} L_{i_3} L_{i_4}}
  &\le& C(U^{\ast},L_0,U_0) \lrp{\frac{\alpha(p)}{p}}^{\lrceil{i_1/2}+\lrceil{i_2/2}+\lrceil{i_3/2}+\lrceil{i_4/2}} \\
  &\le& C(U^{\ast},L_0,U_0) \lrp{\frac{\alpha(p)}{p}}^3, \\[2mm]      
  \EEqn
  which implies that
  \BEqn
  & & \lrabs{
      \begin{array}{l}
         \sum\limits_{\substack{i_1+\dots+i_4\geq5 \\ 0\le i_1,\dots,i_4\le6}} c_{i_1}(A_0) c_{i_2}(A_0) c_{i_3}(A_0) c_{i_4}(A_0) A^{-(i_1+i_2+i_3+i_4)} \\
         \hspace{4em}\times \bigp{M_{i_1 i_2 i_3 i_4} -4\tilde{M}_{i_1 i_2 i_3 i_4} + 2N_{i_1 i_2} N_{i_3 i_4} + 4L_{i_1}D_{i_2 i_3 i_4} -4L_{i_1} L_{i_2} N_{i_3 i_4} + L_{i_1} L_{i_2} L_{i_3} L_{i_4}}
  \end{array}} \\[2mm]
  &\le& C(\tilde{M},U^{\ast},L_0,U_0) c_0^4(A_0) \lrp{\frac{\alpha(p)}{p}}^3.
  \EEqn
  
  Consequently, it remains to consider the case when $i_1+i_2+i_3+i_4=4$. Under Assumption \ref{Assumpt:component-dept}\ref{Assumpt:component-dept-1}-\ref{Assumpt:component-dept-2}, we obtain that
  \BEqn
  & & \sum\limits_{\substack{i_1+\dots+i_4=4 \\ 0\le i_1,\dots,i_4\le6}}
      c_{i_1}(A_0) c_{i_2}(A_0) c_{i_3}(A_0) c_{i_4}(A_0) A^{-4} \\
  & & \hspace{4em}
      \times \bigp{M_{i_1 i_2 i_3 i_4} -4\tilde{M}_{i_1 i_2 i_3 i_4} + 2N_{i_1 i_2} N_{i_3 i_4} + 4L_{i_1}D_{i_2 i_3 i_4} -4L_{i_1} L_{i_2} N_{i_3 i_4} + L_{i_1} L_{i_2} L_{i_3} L_{i_4}} \\
  &=& c_1^4(A_0) A^{-4} \bigp{M_{1111}-4\tilde{M}_{1111}+2(N_{11}^2)} \\
  &=& 16c_1^4(A_0) A^{-4} 
      \sum\limits_{j_1,\dots,j_4=1}^{p} \bbe{\tilde{z}_{1 j_1}\tilde{z}_{1 j_2}} 
                                       \bbe{\tilde{z}_{1 j_1}\tilde{z}_{1 j_3}}     
                                       \bbe{\tilde{z}_{1 j_2}\tilde{z}_{1 j_4}}
                                       \bbe{\tilde{z}_{1 j_3}\tilde{z}_{1 j_4}} \\
  &\le& C(U^{\ast}) c_1^4(A_0) A^{-4}
        \times
        \# \lrcp{(j_1,\dots,j_4): 
                 \begin{array}{l}
                    \cum(\tilde{z}_{1 j_1},\tilde{z}_{1 j_2}) \neq 0, \cum(\tilde{z}_{1 j_1},\tilde{z}_{1 j_3}) \neq 0, \\
                    \cum(\tilde{z}_{1 j_2},\tilde{z}_{1 j_4}) \neq 0,
                                  \cum(\tilde{z}_{1 j_3},\tilde{z}_{1 j_4}) \neq 0. \\ 
                    1\le j_1,\dots,j_4\le p
                 \end{array}} \\
  &\le& C(U^{\ast},L_0,U_0) c_1^4(A_0) \lrp{\frac{\alpha(p)}{p}}^3 \\
  &\le& C(\tilde{M},U^{\ast},L_0,U_0) c_0^4(A_0) \lrp{\frac{\alpha(p)}{p}}^3,
  \EEqn
  implying that $|\tilde{G}| \le C(\tilde{M},U^{\ast},L_0,U_0) c_0^4(A_0) \lrp{\frac{\alpha(p)}{p}}^3$. Thus in summary, we conclude that
  \BEqn
        \lrabs{\bbe{g^k(Z_1,Z_2,Z_3,Z_4)}}
  &\le& \lrabs{\tilde{G}} + C(\tilde{M},U^{\ast},L_0,U_0) c_0^4(A_0) \lrp{\frac{\alpha(p)}{p}}^3 \\
  &\le& C(\tilde{M},U^{\ast},L_0,U_0) c_0^4(A_0) \lrp{\frac{\alpha(p)}{p}}^3.
  \EEqn
\end{proof}


\section{Lemmas Regarding the Power Analysis}\label{Appendix:lemma-3}

\subsection{Lemmas for Theorem \ref{Thm:power-k}}

\begin{lemma}\label{Lemma:diff_E_bound}
Assume that $n/N \rightarrow \rho$ for $0<\rho<1$ as $n,m \rightarrow \infty$, then it holds when $N$ is sufficiently large that
\begin{equation*}
    \blre{\lrabs{ \ce_{n,m}^k(X,Y) - \ce^k(X,Y) }^2} 
\le C \cdot N^{-1} \blre{\lrp{h^k(X_1,X_2,Y_1,Y_2)}^2},
\end{equation*}
for some universal constant $0<C<\infty$ and $h^k$ is the two-sample kernel defined in Proposition \ref{Prop:Enm-k}.
\end{lemma}
\begin{proof}
    Note that $h^k(X_1,X_2,Y_1,Y_2)$ is a two-sample kernel function of order $(2,2)$, then it follows from Lemma \ref{Lemma:moment_ineq} that
    \begin{equation*}
        \blre{\lrabs{\ce_{n,m}^k(X,Y) - \ce^k(X,Y)}^2}
    \le C \cdot N^{-1} \blre{\lrp{h^k(X_1,X_2,Y_1,Y_2)}^2}.
    \end{equation*}
\end{proof}

\subsection{Lemmas for Proposition \ref{Prop:PowerAnalysis-1} - Proposition \ref{Prop:PowerAnalysis-4}}


\begin{lemma}\label{Lemma:PowerAnalysis-1}
Suppose that $f=f^{(p)}$ satisfies Assumption \ref{Assumpt:uniform-kernel}\ref{Assumpt:uniform-kernel-1}-\ref{Assumpt:uniform-kernel-2}, and Assumption \ref{Assumpt:component-dept}\ref{Assumpt:component-dept-1}-\ref{Assumpt:component-dept-3} hold, then 
\begin{enumerate}[label=(\roman*)]
    \item \label{Lemma:PowerAnalysis-1-1}
    if additionally, Assumption \ref{Assumpt:PowerAnalysis-1} holds and there exists some positive constant $L^{\ast}<\infty$, such that
    \begin{equation*}
       \lrabs{2c_0(A_0^{XY}) - c_0(A_0^{X}) - c_0(A_0^{Y})}
    \geq L^{\ast} |c_0(A^{XY})|,
    \end{equation*}
    then there exists some $p_0=p_0(\tilde{M},U^{\ast},L^{\ast},L_0,U_0)$, such that it holds for any $p \geq p_0$ that
    \BEqn
        C_1(\tilde{M},U^{\ast},L^{\ast},L_0,U_0) |c_0(A_0^{XY})|
    \le \ce^{k}(X,Y)
    \le C_2(\tilde{M},U^{\ast},L^{\ast},L_0,U_0) |c_0(A_0^{XY})|.
    \EEqn

    \item \label{Lemma:PowerAnalysis-1-2}
    if additionally, Assumption \ref{Assumpt:uniform-kernel}\ref{Assumpt:uniform-kernel-3} and Assumption \ref{Assumpt:PowerAnalysis-2} hold, and $\alpha(p)=o\lrp{p^{\max\{\delta_1,\delta_2\}}}$ with $\max\{\delta_1,\delta_2\}>0$, then there exists some $p_0=p_0(\tilde{M},\hat{M},U^{\ast},L_0,U_0,L_2,U_2)$, such that it holds for any $p \geq p_0$ that
    \BEqn
    & & C_1(\tilde{M},\hat{M},U^{\ast},L_0,U_0,L_2,U_2) |c_0(A_0^{XY})| p^{\max\{\delta_1,\delta_2\}-1} \\
    &\le& \ce^{k}(X,Y) \\
    &\le& C_2(\tilde{M},U^{\ast},L_0,U_0,L_2,U_2) |c_0(A_0^{XY})| p^{\max\{\delta_1,\delta_2\}-1}.
    \EEqn

    \item \label{Lemma:PowerAnalysis-1-3}
    if additionally, Assumption \ref{Assumpt:uniform-kernel}\ref{Assumpt:uniform-kernel-3} and Assumption \ref{Assumpt:PowerAnalysis-3} hold, and $\alpha(p)=o\lrp{p^{\max\{\delta_3/2,\delta_4/(2-\delta_4)\}}}$ with $\max\{\delta_3,\delta_4\}>0$, then there exists some $p_0=p_0(\tilde{M},\hat{M},U^{\ast},L_0,U_0,L_3,U_3)$, such that it holds for any $p \geq p_0$ that
    \BEqn
    & & C_1(\tilde{M},\hat{M},U^{\ast},L_0,U_0,L_3,U_3) |c_0(A_0^{XY})| \max\{p^{\delta_3-2}, (\alpha(p))^{\delta_4}p^{\delta_4-2}\} \\
    &\le& \ce^{k}(X,Y) \\
    &\le& C_2(\tilde{M},U^{\ast},L_0,U_0,L_3,U_3) |c_0(A_0^{XY})| \max\{p^{\delta_3-2}, (\alpha(p))^{\delta_4}p^{\delta_4-2}\}.
    \EEqn
    
    \item \label{Lemma:PowerAnalysis-1-4}
    if additionally, Assumption \ref{Assumpt:uniform-kernel}\ref{Assumpt:uniform-kernel-4} and Assumption \ref{Assumpt:PowerAnalysis-4} hold, and  
    \begin{equation*}
       \alpha(p)=o\lrp{p^{\max\{\delta_5/(3-2\delta_5),\delta_6/(3-2\delta_6),(2\delta_7-1)/(3-2\delta_7)\}}},
    \end{equation*}
    holds with $\max\{\delta_5,\delta_6,\delta_7\}>0$,  then there exists some $p_0=p_0(\tilde{M},\hat{M},U^{\ast},L_0,U_0,L_4,U_4)$, such that it holds for any $p \geq p_0$ that
    \begin{equation*}
      \begin{array}{l}
         \hspace{1em}C_1(\tilde{M},\hat{M},U^{\ast},L_0,U_0,L_4,U_4) |c_0(A_0^{XY})| \max\{(\alpha(p))^{2\delta_5}p^{\delta_5-3}, (\alpha(p))^{2\delta_6}p^{\delta_6-3}, (\alpha(p))^{2\delta_7}p^{2\delta_7-4}\} \vspace{2mm}\\
         \le \ce^{k}(X,Y) \vspace{2mm}\\
         \le C_2(\tilde{M},U^{\ast},L_0,U_0,L_4,U_4) |c_0(A_0^{XY})| \max\{(\alpha(p))^{2\delta_5}p^{\delta_5-3}, (\alpha(p))^{2\delta_6}p^{\delta_6-3}, (\alpha(p))^{2\delta_7}p^{2\delta_7-4}\}.
      \end{array}
    \end{equation*}
\end{enumerate}
\end{lemma}
\begin{proof}
  For any fixed $k=k^{(p)}\in\cc$, it follows from Lemma \ref{Lemma:approx-2} that 
  \BEqn
  & & \ce^k(X,Y) \\
  &=& 2\bbe{k(X_1,Y_1)} - \bbe{k(X_1,X_2)} - \bbe{k(Y_1,Y_2)} \\
  &=& \lrp{2c_0(A_0^{XY}) - c_0(A^{X}) - c_0(A^{Y})} \\
  & & + \sum\limits_{i=2}^{6} 
        \lrp{  2c_i(A_0^{XY}) (A^{XY})^{-i} L_i^{XY} 
             - c_i(A_0^{X}) (A^{X})^{-i} L_i^{X}
             - c_i(A_0^{Y}) (A^{Y})^{-i} L_i^{Y}} \\
  & & + \bigp{2\bbe{R(X_1,Y_1) - \bbe{R(X_1,X_2)} - \bbe{R(Y_1,Y_2)}}},
  \EEqn
  where
  \BEqn
  & & \lrabs{2\bbe{R(X_1,Y_1) - \bbe{R(X_1,X_2)} - \bbe{R(Y_1,Y_2)}}} \\
  &\le& C(\tilde{M},U^{\ast},L_0,U_0) 
        \lrp{|c_0(A_0^{XY})| + |c_0(A^{XY})| + |c_0(A^{Y})|}
        \lrp{\frac{\alpha(p)}{p}}^{7/2} \\
  &\le& C(\tilde{M},U^{\ast},L_0,U_0) |c_0(A^{XY})| \lrp{\frac{\alpha(p)}{p}}^{7/2},
  \EEqn
  where the last inequality follows from Lemma \ref{Lemma:PowerAnalysis-Prep-3}. 
  
  Next we prove the statement case by case.
  \begin{enumerate}[label=(\roman*)]
    \item Recall the expressions of $c_i(x)$ introduced in Lemma \ref{Lemma:approx-2}, we have
    \begin{equation*}
        |c_i(A_0^{X})| 
    \le C \sum\limits_{s=1}^{i} |f_s(A_0^{X})| (A_0^{X})^s
    \le C \sum\limits_{s=1}^{i} |f_s(A_0^{X})| (A_0^{XY})^s \lrp{\frac{A_0^{X}}{A_0^{XY}}}^s
    \le C(\tilde{M},L_0,U_0) |c_0(A_0^{XY})|,
    \end{equation*}
    where the last inequality follows from Assumption \ref{Assumpt:uniform-kernel}\ref{Assumpt:uniform-kernel-2} and Assumption \ref{Assumpt:component-dept}\ref{Assumpt:component-dept-3}. Together with the results derived in Lemma \ref{Lemma:Li-1} and Lemma \ref{Lemma:Li-2}, we obtain that
      \BEqn
      & & \lrabs{\sum\limits_{i=2}^{6} 
                 \lrp{  2c_i(A_0^{XY}) (A^{XY})^{-i} L_i^{XY} 
                      - c_i(A_0^{X}) (A^{X})^{-i} L_i^{X}
                      - c_i(A_0^{Y}) (A^{Y})^{-i} L_i^{Y}}} \\
      &\le& C(\tilde{M},U^{\ast},L_0,U_0) |c_0(A^{XY})| \lrp{\frac{\alpha(p)}{p}}.
      \EEqn
      Additionally, we have that
      \BEqn
      & & \lrabs{2\bbe{R(X_1,Y_1) - \bbe{R(X_1,X_2)} - \bbe{R(Y_1,Y_2)}}} \\
      &\le& C(\tilde{M},U^{\ast},L_0,U_0) |c_0(A^{XY})| \lrp{\frac{\alpha(p)}{p}}^{7/2} \\
      &\le& C(\tilde{M},U^{\ast},L_0,U_0) |c_0(A^{XY})| \lrp{\frac{\alpha(p)}{p}}.
      \EEqn
      By using Lemma \ref{Lemma:PowerAnalysis-Prep-3} and the condition in Assumption\ref{Assumpt:PowerAnalysis-1}, we have that
      \begin{equation*}
          L^{\ast} |c_0(A_0^{XY})|
      \le \lrabs{2c_0(A_0^{XY}) - c_0(A_0^{X}) - c_0(A_0^{Y})} 
      \le C(\tilde{M},L_0,U_0) |c_0(A_0^{XY})|,
      \end{equation*}
      that is, $\lrabs{2c_0(A_0^{XY}) - c_0(A_0^{X}) - c_0(A_0^{Y})}$ is of order $|c_0(A_0^{XY})|$, implying that $2c_0(A_0^{XY}) - c_0(A_0^{X}) - c_0(A_0^{Y})$ becomes the leading term of $\ce^{k}(X,Y)$. Consequently, the desired lower and upper bounds of $\ce^k(X,Y)$ is obtained when $p \geq p_0$ for some $p_0=p_0(\tilde{M},U^{\ast},L^{\ast},L_0,U_0)$.

    \item From Lemma \ref{Lemma:PowerAnalysis-Prep-4}\ref{Lemma:PowerAnalysis-Prep-4-1} we have
      \begin{equation*}
          \begin{array}{l}
             \hspace{1em}
             \lrabs{\bigp{2c_0(A^{XY}) - c_0(A^{X}) - c_0(A^{Y})} 
           - \lrp{2c_1(A_0^{XY}) \frac{|\Delta|^2}{A^{XY}} - 2c_2(A_0^{XY}) \frac{\lrp{\bbe{|\tilde{X}_1|^2} - \bbe{|\tilde{Y}_1|^2}}^2}{(A^{XY})^2}}} \\
             \le C(\tilde{M},L_0,U_0) |c_0(A_0^{XY})| \lrp{\frac{|\Delta|^2}{p^2} + \frac{\lrabs{\bbe{|\tilde{X}_1|^2} - \bbe{|\tilde{Y}_1|^2}}^3}{p^3}}.
          \end{array}
      \end{equation*}
      Under Assumption \ref{Assumpt:uniform-kernel}\ref{Assumpt:uniform-kernel-3}, $|c_1(A_0^{XY})|$ and $|c_2(A_0^{XY})|$ are of the same order. When Assumption \ref{Assumpt:component-dept}\ref{Assumpt:component-dept-2} is satisfied and Assumption \ref{Assumpt:PowerAnalysis-2} holds with $\delta_1\neq\delta_2$, either $2c_1(A_0^{XY}) \frac{|\Delta|^2}{A^{XY}}$ or $2c_2(A_0^{XY}) \frac{\lrp{\bbe{|\tilde{X}_1|^2} - \bbe{|\tilde{Y}_1|^2}}^2}{(A^{XY})^2}$ solely dominates $2c_1(A_0^{XY}) \frac{|\Delta|^2}{A^{XY}} - 2c_2(A_0^{XY}) \frac{\lrp{\bbe{|\tilde{X}_1|^2} - \bbe{|\tilde{Y}_1|^2}}^2}{(A^{XY})^2}$, then by using Assumption \ref{Assumpt:uniform-kernel}\ref{Assumpt:uniform-kernel-2}-\ref{Assumpt:uniform-kernel-3} we conclude that
      \BEqn
      & & C_1(\tilde{M},\hat{M},L_0,U_0,L_2,U_2) |c_0(A_0^{XY})| p^{\max\{\delta_1,\delta_2\}-1} \\
      &\le& \lrabs{2c_1(A_0^{XY}) \frac{|\Delta|^2}{A^{XY}} - 2c_2(A_0^{XY}) \frac{\lrp{\bbe{|\tilde{X}_1|^2} - \bbe{|\tilde{Y}_1|^2}}^2}{(A^{XY})^2}} \\
      &\le& C_2(\tilde{M},L_0,U_0,L_2,U_2) |c_0(A_0^{XY})| p^{\max\{\delta_1,\delta_2\}-1},
      \EEqn
      where the lower bound holds when $p \geq p_0$ for some $p_0=p_0(\tilde{M},\hat{M},L_0,U_0,L_2,U_2)$.
      
      Also, we have
      \begin{equation*}
          \frac{|\Delta|^4}{p^2} \le C(L_2,U_2) p^{2(\delta_1-1)},
      \end{equation*}
      and
      \begin{equation*}
          \frac{\lrabs{\bbe{|\tilde{X}_1|^2} - \bbe{|\tilde{Y}_1|^2}}^3}{p^3}
      \le C(L_2,U_2) p^{3(\delta_2-1)/2}.
      \end{equation*}
      
      Recall that
      \BEqn
      & & \lrabs{\sum\limits_{i=2}^{6} 
                 \lrp{  2c_i(A_0^{XY}) (A^{XY})^{-i} L_i^{XY} 
                      - c_i(A_0^{X}) (A^{X})^{-i} L_i^{X}
                      - c_i(A_0^{Y}) (A^{Y})^{-i} L_i^{Y}}} \\
      & & + \lrabs{\bigp{2\bbe{R(X_1,Y_1) - \bbe{R(X_1,X_2)} - \bbe{R(Y_1,Y_2)}}}} \\
      &\le& C(\tilde{M},U^{\ast},L_0,U_0) |c_0(A^{XY})| \lrp{\frac{\alpha(p)}{p}},
      \EEqn
      we further obtain that
      \BEqn
      & & \lrabs{\ce^k(X,Y) - \lrp{2c_1(A_0^{XY}) \frac{|\Delta|^2}{A^{XY}} - 2c_2(A_0^{XY}) \frac{\lrp{\bbe{|\tilde{X}_1|^2} - \bbe{|\tilde{Y}_1|^2}}^2}{(A^{XY})^2}}} \\
      &\le& C(\tilde{M},L_0,U_0) |c_0(A_0^{XY})| \lrp{\frac{|\Delta|^2}{p^2} + \frac{\lrabs{\bbe{|\tilde{X}_1|^2} - \bbe{|\tilde{Y}_1|^2}}^3}{p^3}} \\
      & & + \lrabs{\sum\limits_{i=2}^{6} 
                 \lrp{  2c_i(A_0^{XY}) (A^{XY})^{-i} L_i^{XY} 
                      - c_i(A_0^{X}) (A^{X})^{-i} L_i^{X}
                      - c_i(A_0^{Y}) (A^{Y})^{-i} L_i^{Y}}} \\
      & & + \lrabs{\bigp{2\bbe{R(X_1,Y_1) - \bbe{R(X_1,X_2)} - \bbe{R(Y_1,Y_2)}}}} \\
      &\le& C(\tilde{M},U^{\ast},L_0,U_0,L_2,U_2) |c_0(A_0^{XY})| \lrp{p^{2(\delta_1-1)} + p^{3(\delta_2-1)/2} + \lrp{\frac{\alpha(p)}{p}}}.
      \EEqn
      Therefore, under the assumption that $\max\{\delta_1,\delta_2\}<1$ and $\alpha(p) = o\lrp{p^{\max\{\delta_1,\delta_2\}}}$, it holds when $p \geq p_0$ with some $p_0=p_0(\tilde{M},\hat{M},U^{\ast},L_0,U_0,L_2,U_2)$ that
      \BEqn
      & & C_1(\tilde{M},\hat{M},U^{\ast},L_0,U_0,L_2,U_2) |c_0(A_0^{XY})| p^{\max\{\delta_1,\delta_2\}-1} \\
      &\le& \ce^k(X,Y) \\
      &\le& C_2(\tilde{M},U^{\ast},L_0,U_0,L_2,U_2) |c_0(A_0^{XY})| p^{\max\{\delta_1,\delta_2\}-1},
      \EEqn
      which completes the proof.

    \item Note that
      \BEqn
          \ce^k(X,Y)
      &=& \lrp{2c_0(A_0^{XY}) - c_0(A^{X}) - c_0(A^{Y})} \\
      & & + \lrp{2c_2(A_0^{XY}) (A^{XY})^{-2} L_2^{XY} 
                 - c_2(A_0^{X}) (A^{X})^{-2} L_2^{X}
                 - c_2(A_0^{Y}) (A^{Y})^{-2} L_2^{Y}} \\
      & & + \sum\limits_{i=3}^{6} 
            \lrp{  2c_i(A_0^{XY}) (A^{XY})^{-i} L_i^{XY} 
                 - c_i(A_0^{X}) (A^{X})^{-i} L_i^{X}
                 - c_i(A_0^{Y}) (A^{Y})^{-i} L_i^{Y}} \\
      & & + \bigp{2\bbe{R(X_1,Y_1) - \bbe{R(X_1,X_2)} - \bbe{R(Y_1,Y_2)}}},
      \EEqn
      and we have
      \BEqn
      & & \lrabs{\sum\limits_{i=3}^{6} 
                 \lrp{  2c_i(A_0^{XY}) (A^{XY})^{-i} L_i^{XY} 
                      - c_i(A_0^{X}) (A^{X})^{-i} L_i^{X}
                      - c_i(A_0^{Y}) (A^{Y})^{-i} L_i^{Y}}} \\
      & & + \lrabs{\bigp{2\bbe{R(X_1,Y_1) - \bbe{R(X_1,X_2)} - \bbe{R(Y_1,Y_2)}}}} \\
      &\le& C(\tilde{M},U^{\ast},L_0,U_0) |c_0(A^{XY})| \lrp{\frac{\alpha(p)}{p}}^2,
      \EEqn
      then to determine the leading term of $\ce^{k}(X,Y)$ in this case, it remains to calculate the orders of $2c_0(A_0^{XY}) - c_0(A^{X}) - c_0(A^{Y})$ and $2c_2(A_0^{XY}) (A^{XY})^{-2} L_2^{XY} - c_2(A_0^{X}) (A^{X})^{-2} L_2^{X} - c_2(A_0^{Y}) (A^{Y})^{-2} L_2^{Y}$ and compare them with $|c_0(A_0^{XY})|\lrp{\frac{\alpha(p)}{p}}^2$.
      
      From Lemma \ref{Lemma:PowerAnalysis-Prep-4}\ref{Lemma:PowerAnalysis-Prep-4-2}, it holds when $|\Delta|=0$ that
      \BEqn
      & & \lrabs{\bigp{2c_0(A^{XY}) - c_0(A^{X}) - c_0(A^{Y})} + 2c_2(A_0^{XY}) \frac{\lrp{\bbe{|\tilde{X}_1|^2} - \bbe{|\tilde{Y}_1|^2}}^2}{(A^{XY})^2}} \\
      &\le& C(\tilde{M},L_0,U_0) |c_0(A_0^{XY})|
      \frac{\lrabs{\bbe{|\tilde{X}_1|^2} - \bbe{|\tilde{Y}_1|^2}}^4}{p^4},
      \EEqn
      which implies that under Assumption \ref{Assumpt:uniform-kernel}\ref{Assumpt:uniform-kernel-3}, $2c_2(A_0^{XY}) \frac{\lrp{\bbe{|\tilde{X}_1|^2} - \bbe{|\tilde{Y}_1|^2}}^2}{(A^{XY})^2}$ determines the order of $2c_0(A^{XY}) - c_0(A^{X}) - c_0(A^{Y})$ when $p \geq p_0$ for some $p_0=p_0(\tilde{M},\hat{M},L_0,U_0,L_3,U_3)$.
      
      Now consider the leading term of $2c_2(A_0^{XY}) (A^{XY})^{-2} L_2^{XY} - c_2(A_0^{X}) (A^{X})^{-2} L_2^{X} - c_2(A_0^{Y}) (A^{Y})^{-2} L_2^{Y}$. Note that
      \BEqn
      & & 2c_2(A_0^{XY}) (A^{XY})^{-2} L_2^{XY} - c_2(A_0^{X}) (A^{X})^{-2} L_2^{X} - c_2(A_0^{Y}) (A^{Y})^{-2} L_2^{Y} \\
      &=& c_2(A_0^{XY}) (A^{XY})^{-2} \bigp{2L_2^{XY}-L_2^{X}-L_2^{Y}} \\
      & & - \lrp{c_2(A_0^{X}) (A^X)^{-2} - c_2(A_0^{XY}) (A^{XY})^{-2}} L_2^{X} \\
      & & - \lrp{c_2(A_0^{Y}) (A^Y)^{-2} - c_2(A_0^{XY}) (A^{XY})^{-2}} L_2^{Y}.
      \EEqn
      The smoothness condition proposed in Lemma \ref{Lemma:PowerAnalysis-Prep-7} is naturally satisfied under Assumption \ref{Assumpt:uniform-kernel}\ref{Assumpt:uniform-kernel-1}. Hence it follows from Lemma \ref{Lemma:PowerAnalysis-Prep-2}, Lemma \ref{Lemma:PowerAnalysis-Prep-7} and Lemma \ref{Lemma:PowerAnalysis-Prep-8} that
      \BEqn
      & & \lrabs{
          \begin{array}{ll}
               & \bigp{2c_2(A_0^{XY}) (A^{XY})^{-2} L_2^{XY} - c_2(A_0^{X}) (A^{X})^{-2} L_2^{X} - c_2(A_0^{Y}) (A^{Y})^{-2} L_2^{Y}} \vspace{2mm}\\
             + & 4c_2(A_0^{XY}) (A^{XY})^{-2} \|\Sigma_X-\Sigma_Y\|_F^2
      \end{array}} \\
      &\le& C(\tilde{M},U^{\ast},L_0,U_0) |c_0(A_0^{XY})| \frac{\bigp{|A^{X}-A^{XY}|+|A^{Y}-A^{XY}|}\alpha(p)}{p^2} \\
      &=& C(\tilde{M},U^{\ast},L_0,U_0) |c_0(A_0^{XY})| \frac{\lrabs{\bbe{|\tilde{X}_1|^2}-\bbe{|\tilde{Y}_1|^2}} \alpha(p)}{p^2}.
      \EEqn
      
      In summary, we have
      \BEqn
      & & \lrabs{\ce^{k}(X,Y) + 2c_2(A_0^{XY}) (A^{XY})^{-2} \lrp{\lrp{\bbe{|\tilde{X}_1|^2} - \bbe{|\tilde{Y}_1|^2}}^2 + 2\|\Sigma_X-\Sigma_Y\|_F^2}} \\
      &\le& C(\tilde{M},U^{\ast},L_0,U_0) |c_0(A_0^{XY})| \\
      & & \hspace{2em} \times
          \lrp{  \frac{\lrabs{\bbe{|\tilde{X}_1|^2} - \bbe{|\tilde{Y}_1|^2}}^4}{p^4}
               + \frac{\lrabs{\bbe{|\tilde{X}_1|^2}-\bbe{|\tilde{Y}_1|^2}} \alpha(p)}{p^2}
               + \lrp{\frac{\alpha(p)}{p}}^2} \\
      &\le& C(\tilde{M},U^{\ast},L_0,U_0,L_3,U_3) |c_0(A_0^{XY})|
            \lrp{p^{2(\delta_3-2)} + (\alpha(p)) p^{\delta_3/2-2} + (\alpha(p))^2 p^{-2}}.
      \EEqn
      Under the assumption $\alpha(p)=o\lrp{p^{\max\{\delta_3/2,\delta_4/(2-\delta_4)\}}}$, we have that \\ $2c_2(A_0^{XY}) (A^{XY})^{-2} \lrp{\lrp{\bbe{|\tilde{X}_1|^2} - \bbe{|\tilde{Y}_1|^2}}^2 + 2\|\Sigma_X-\Sigma_Y\|_F^2}$ is the leading term of $\ce^{k}(X,Y)$ when $p \geq p_0$ for some $p_0=p_0(\tilde{M},\hat{M},U^{\ast},L_0,U_0,L_3,U_3)$, which further leads to the desired result.

      \item Note that when $|\Delta|=\|\Sigma_X-\Sigma_Y\|_F=0$, it is trivial that $A^{X} = A^{XY} = A^{Y}$ and consequently, we have $c_i(A_0^{X}) = c_i(A_0^{XY}) = c_i(A_0^{Y})$. Then the approximation of $\ce^{k}(X,Y)$ can be simplified as the following expression:
      \BEqn
          \ce^{k}(X,Y)
      &=&   c_2(A_0^{XY}) \bigp{2L_2^{XY} - L_2^{X} - L_2^{Y}} \\
      & & + c_3(A_0^{XY}) \bigp{2L_3^{XY} - L_3^{X} - L_3^{Y}} \\
      & & + c_4(A_0^{XY}) \bigp{2L_4^{XY} - L_4^{X} - L_4^{Y}} \\
      & & + \sum\limits_{i=5}^{6} c_i(A_0^{XY}) \bigp{2L_i^{XY} - L_i^{X} - L_i^{Y}} \\
      & & + \bigp{2\bbe{R(X_1,Y_1)} - \bbe{R(X_1,X_2)} - \bbe{R(Y_1,Y_2)}}.
      \EEqn
      Note that
      \begin{equation*}
          \lrabs{\sum\limits_{i=5}^{6} c_i(A_0^{XY}) \bigp{2L_i^{XY} - L_i^{X} - L_i^{Y}}}
      \le C(M,L_0,U_0) |c_0(A_0^{XY})| \lrp{\frac{\alpha(p)}{p}}^3,
      \end{equation*}
      and
      \BEqn
      & & \lrabs{2\bbe{R(X_1,Y_1)} - \bbe{R(X_1,X_2)} - \bbe{R(Y_1,Y_2)}} \\
      &\le& C(\tilde{M},U^{\ast},L_0,U_0) |c_0(A_0^{XY})| \lrp{\frac{\alpha(p)}{p}}^{7/2} \\
      &\le& C(\tilde{M},U^{\ast},L_0,U_0) |c_0(A_0^{XY})| \lrp{\frac{\alpha(p)}{p}}^3.
      \EEqn
      
      Define
      \begin{equation*}
          \begin{array}{l}
              \hspace{1em}
              E \\
            = c_3(A_0^{XY}) (A^{XY})^{-3}
              \left(
                8\sum\limits_{j_1,j_2,j_3=1}^{p}\lrp{\cum(\tilde{x}_{1 j_1}, \tilde{x}_{1 j_2}, \tilde{x}_{1 j_3}) - \cum(\tilde{y}_{1 j_1}, \tilde{y}_{1 j_2}, \tilde{y}_{1 j_3})}^2 
              \right. \\  
              \hspace{9em}
              \left.
                + 12\sum\limits_{j_2=1}^{p}\lrp{\sum\limits_{j_1=1}^{p} \cum(\tilde{x}_{1 j_1}^2 - \sigma_{X, j_1}^2, \tilde{x}_{1 j_1}) - \cum(\tilde{y}_{1 j_1}^2 - \sigma_{Y, j_1}^2, \tilde{y}_{1 j_1})}^2
                 \right) \\
              - 6c_4(A_0^{XY}) (A^{XY})^{-4}
                \lrp{\sum\limits_{j_1,j_2=1}^{p} \cum(\tilde{x}_{1 j_1}^2 - \sigma_{X, j_1}^2, \tilde{x}_{1 j_2}^2 - \sigma_{X, j_2}^2) - \cum(\tilde{y}_{1 j_1}^2 - \sigma_{Y, j_1}^2, \tilde{y}_{1 j_2}^2 - \sigma_{Y, j_2}^2)}^2.
          \end{array}
      \end{equation*}
      
      By using the results derived in Lemma \ref{Lemma:PowerAnalysis-Prep-9}, we obtain that
      \begin{equation*}
         \lrabs{\ce^{k}(X,Y) - E} 
      \le C(\tilde{M},U^{\ast},L_0,U_0) |c_0(A_0^{XY})| \lrp{\frac{\alpha(p)}{p}}^3.
      \end{equation*}
      Also note that under Assumption \ref{Assumpt:PowerAnalysis-4}, the order of $E$ is determined by 
      \begin{equation*}
          |c_0(A_0^{XY})| \max\{(\alpha(p))^{2\delta_5} p^{\delta_5-3},
                                (\alpha(p))^{2\delta_6} p^{\delta_6-3},
                                (\alpha(p))^{2\delta_7} p^{2\delta_7-4}\},
      \end{equation*}
      thus $E$ is the leading term of $\ce^{k}(X,Y)$ when
      \begin{equation*}
        \alpha(p)
      = o\lrp{p^{\max\{\delta_5/(3-2\delta_5),\delta_6/(3-2\delta_6),(2\delta_7-1)/(3-2\delta_7)\}}},
      \end{equation*}
      and $p \geq p_0$ for some $p_0=p_0(\tilde{M},\hat{M},U^{\ast},L_0,U_0,L_4,U_4)$, which jointly lead to the claimed result.
  \end{enumerate}
\end{proof}

\begin{lemma}\label{Lemma:PowerAnalysis-2}
Suppose that $f=f^{(p)}$ satisfies Assumption \ref{Assumpt:uniform-kernel}\ref{Assumpt:uniform-kernel-1}-\ref{Assumpt:uniform-kernel-2} and Assumption \ref{Assumpt:component-dept}\ref{Assumpt:component-dept-1}-\ref{Assumpt:component-dept-3}, then 
\begin{enumerate}[label=(\roman*)]
    \item \label{Lemma:PowerAnalysis-2-1}
    if additionally, Assumption \ref{Assumpt:PowerAnalysis-1} holds and there exists some positive constant $L^{\ast}<\infty$, such that
    \begin{equation*}
       \lrabs{2c_0(A_0^{XY}) - c_0(A_0^{X}) - c_0(A_0^{Y})}
    \geq L^{\ast} |c_0(A^{XY})|,
    \end{equation*}
    then there exists some $p_0 = p_0(\tilde{M},U^{\ast},L^{\ast},L_0,U_0,L_1,U_1,\rho)$ such that it holds when $p \geq p_0$ that
    \BEqn
    & & C_1(\tilde{M},U^{\ast},L^{\ast},L_0,U_0,\rho) c_0^2(A_0^{XY}) \\
    &\le& \cv_k^2(Z) \\
    &\le& C_2(\tilde{M},U^{\ast},L^{\ast},L_0,U_0,\rho) c_0^2(A_0^{XY}).
    \EEqn
    
    \item \label{Lemma:PowerAnalysis-2-2}
    if additionally, Assumption \ref{Assumpt:uniform-kernel}\ref{Assumpt:uniform-kernel-3} and Assumption \ref{Assumpt:component-dept}\ref{Assumpt:component-dept-4} hold, Assumption \ref{Assumpt:PowerAnalysis-2} holds with $\delta_1\neq\delta_2$, and $\alpha(p)=o\lrp{p^{\max\{\delta_1,\delta_2\}}}$ with $\max\{\delta_1,\delta_2\}>0$, then there exists some $p_0 = p_0(\tilde{M},\hat{M},U^{\ast},L_0,U_0,L_2,U_2,\rho)$ such that it holds when $p \geq p_0$ that
    \BEqn
    & & C_1(\tilde{M},\hat{M},U^{\ast},L_0,U_0,L_2,U_2,\rho) c_0^2(A_0^{XY}) p^{2\max\{\delta_1,\delta_2\}-2} \\
    &\le& \cv_k^2(Z) \\
    &\le& C_2(\tilde{M},U^{\ast},L_0,U_0,L_2,U_2,\rho) c_0^2(A_0^{XY}) p^{2\max\{\delta_1,\delta_2\}-2}
    \EEqn
    if $\alpha(p)=o\lrp{p^{2\max\{\delta_1,\delta_2\}-1}}$, and
    \BEqn
    & & C_1(\tilde{M},\hat{M},U^{\ast},L_0,U_0,L_2,U_2,\rho) c_0^2(A_0^{XY}) \lrp{\frac{\alpha(p)}{p}} \\
    &\le& \cv_k^2(Z) \\
    &\le& C_2(\tilde{M},U^{\ast},L_0,U_0,L_2,U_2,\rho) c_0^2(A_0^{XY}) \lrp{\frac{\alpha(p)}{p}}
    \EEqn
    otherwise.
    
    \item \label{Lemma:PowerAnalysis-2-3}
    if additionally, Assumption \ref{Assumpt:uniform-kernel}\ref{Assumpt:uniform-kernel-3}, Assumption \ref{Assumpt:component-dept}\ref{Assumpt:component-dept-4} and Assumption \ref{Assumpt:PowerAnalysis-3} hold, then there exists some $p_0 = p_0(\tilde{M},\hat{M},U^{\ast},L_0,U_0,L_3,U_3,\rho)$ such that it holds when $p \geq p_0$ that
    \BEqn
    & & C_1(\tilde{M},\hat{M},U^{\ast},L_0,U_0,L_3,U_3,\rho) c_0^2(A_0^{XY}) \lrp{\frac{\alpha(p)}{p}} \\
    &\le& \cv_k^2(Z) \\
    &\le& C_2(\tilde{M},U^{\ast},L_0,U_0,L_3,U_3,\rho) c_0^2(A_0^{XY}) \lrp{\frac{\alpha(p)}{p}}.
    \EEqn
    
    \item \label{Lemma:PowerAnalysis-2-4}
    if additionally, Assumption \ref{Assumpt:uniform-kernel}\ref{Assumpt:uniform-kernel-3}, Assumption \ref{Assumpt:component-dept}\ref{Assumpt:component-dept-4} and Assumption \ref{Assumpt:PowerAnalysis-4} hold, then there exists some $p_0 = p_0(\tilde{M},\hat{M},U^{\ast},L_0,U_0,\rho)$ such that it holds when $p \geq p_0$ that
    \BEqn
    & & C_1(\tilde{M},\hat{M},U^{\ast},L_0,U_0,\rho) c_0^2(A_0^{XY}) \lrp{\frac{\alpha(p)}{p}} \\
    &\le& \cv_k^2(Z) \\
    &\le& C_2(\tilde{M},U^{\ast},L_0,U_0,\rho) c_0^2(A_0^{XY}) \lrp{\frac{\alpha(p)}{p}}.
    \EEqn
    
\end{enumerate}
\end{lemma}
\begin{proof}
  Recall that $\cv_k^2(Z) = \bbe{k^2(Z_1,Z_2)} - 2\bbe{k(Z_1,Z_2)k(Z_1,Z_3)} + \bbe{k(Z_1,Z_2)}^2$ and $Z$ is a mixture distribution of $X$ and $Y$, then it follows that
  \begin{equation*}
      \cv_k^2(Z) = V_1 + V_2 + V_3 + V_4 + V_5 + V_6,
  \end{equation*}
  where
  \BEqn
  V_1 &=& \rho^2 \bbe{k^2(X_1,X_2)} - 2\rho^3 \bbe{k(X_1,X_2)k(X_1,X_3)} + \rho^4\bbe{k(X_1,X_2)}^4, \\
  V_2 &=& 2\rho(1-\rho) \bbe{k^2(X_1,Y_1)}  - 2\rho(1-\rho)^2 \bbe{k(X_1,Y_1)k(X_1,Y_2)} \\
      & & - 2\rho^2(1-\rho) \bbe{k(X_1,Y_1)k(X_2,Y_1)} + 4\rho^2(1-\rho)^2 \bbe{k(X_1,Y_1)}^2, \\
  V_3 &=& (1-\rho)^2 \bbe{k^2(Y_1,Y_2)} - 2(1-\rho)^3 \bbe{k(Y_1,Y_2)k(Y_1,Y_3)} + (1-\rho)^4 \bbe{k(Y_1,Y_2)}^2, \\
  V_4 &=& -4\rho^2(1-\rho) \bbe{k(X_1,X_2)k(X_1,Y_1)} + 4\rho^3(1-\rho) \bbe{k(X_1,X_2)}\bbe{k(X_1,Y_1)}, \\
  V_5 &=& -4\rho(1-\rho)^2 \bbe{k(X_1,Y_1)k(Y_1,Y_2)} + 4\rho(1-\rho)^3 \bbe{k(X_1,Y_1)}\bbe{k(Y_1,Y_2)}, \\
  V_6 &=& 2\rho^2(1-\rho)^2 \bbe{k(X_1,X_2)}\bbe{k(Y_1,Y_2)}.
  \EEqn
  
  By repeatedly applying Lemma \ref{Lemma:approx-2}, we obtain the approximation of each $V_i$. Then by reorganizing and combining all the individual terms, we claim without showing all the detailed calculations that
  \BEqn
      \lrabs{\cv_k^2(Z) - \cw_1 - \cw_2 - \cw_3} 
  &\le& C(\tilde{M},U^{\ast},L_0,U_0,\rho) \bigp{c_0^2(A_0^{X}) + c_0^2(A_0^{XY}) + c_0(A_0^{Y})} \lrp{\frac{\alpha(p)}{p}}^2 \\
  &\le& C(\tilde{M},U^{\ast},L_0,U_0,\rho) c_0^2(A_0^{XY}) \lrp{\frac{\alpha(p)}{p}}^2,
  \EEqn
  where the last step follows from Lemma \ref{Lemma:PowerAnalysis-Prep-3} and 
  \BEqn
  \cw_1 &=& \rho^2(1-\rho)^2 \lrp{2c_0(A_0^{XY}) - c_0(A_0^{X}) - c_0(A_0^{Y})}^2, \\
  \cw_2 &=& 2\rho^2(1-\rho)^2 \lrp{2c_0(A_0^{XY}) - c_0(A_0^{X}) - c_0(A_0^{Y})} \\
        & & \times \lrp{  2c_2(A_0^{XY}) (A_0^{XY})^{-2} L_2^{XY} 
                        - c_2(A_0^{X}) (A^{X})^{-2} L_2^{X}
                        - c_2(A_0^{Y}) (A^{Y})^{-2} L_2^{Y}},
  \EEqn
  and
  \BEqn
  & & \cw_3 \\
  &=& \rho^2 c_1^2(A_0^{X}) (A^{X})^{-2} 
      \lrp{\bbe{(|X_1-X_2|^2-A^{X})^2} - 2\bbe{(|X_1-X_2|^2-A^{X})(|X_1-X_3|^2-A^{X})}} \\
  & & + 2\rho(1-\rho) c_1^2(A_0^{XY}) (A^{XY})^{-2} 
        \big{(}\bbe{(|X_1-Y_1|^2-A^{XY})^2} \\
  & & \hspace{13em} - \bbe{(|X_1-Y_1|^2-A^{XY})(|X_1-Y_2|^2-A^{XY})} \\
  & & \hspace{13em} - \bbe{(|X_1-Y_1|^2-A^{XY})(|X_2-Y_1|^2-A^{XY})}\big{)} \\
  & & + (1-\rho)^2 c_1^2(A_0^{Y}) (A^{Y})^{-2} 
        \lrp{\bbe{(|Y_1-Y_2|^2-A^{Y})^2} - 2\bbe{(|Y_1-Y_2|^2-A^{Y})(|Y_1-Y_3|^2-A^{Y})}} \\ 
  & & + 2\rho^2(1-\rho) 
        \big{(} c_1^2(A_0^{X}) (A^{X})^{-2} \bbe{(|X_1-X_2|^2-A^{X})(|X_1-X_3|^2-A^{X})} \\
  & & \hspace{6em} -2c_1(A_0^{X}) c_1(A_0^{XY}) (A^{X})^{-1} (A^{XY})^{-1} \bbe{(|X_1-X_2|^2-A^{X})(|X_1-Y_1|^2-A^{XY})} \\
  & & \hspace{6em} + c_1^2(A_0^{XY}) (A^{XY})^{-2} \bbe{(|X_1-Y_1|^2-A^{XY})(|X_1-Y_2|^2-A^{XY})} \big{)} \\
  & & + 2\rho(1-\rho)^2
        \big{(} c_1^2(A_0^{Y}) (A^{Y})^{-2} \bbe{(|Y_1-Y_2|^2-A^{Y})(|Y_1-Y_3|^2-A^{Y})} \\
  & & \hspace{6em} -2c_1(A_0^{Y}) c_1(A_0^{XY}) (A^{Y})^{-1} (A^{XY})^{-1} \bbe{(|X_1-Y_1|^2-A^{XY})(|Y_1-Y_2|^2-A^{Y})} \\
  & & \hspace{6em} + c_1^2(A_0^{XY}) (A^{XY})^{-2} \bbe{(|X_1-Y_1|^2-A^{XY})(|X_2-Y_1 |^2-A^{XY})} \big{)}.
  \EEqn
  
  \begin{enumerate}[label=(\roman*)]
    \item
      It follows from Lemma \ref{Lemma:PowerAnalysis-Prep-3} that
      \BEqn
      & & \lrabs{\cw_2+\cw_3} \\
      &\le& C(\tilde{M},U^{\ast},L_0,U_0,\rho) 
            \bigp{c_0^2(A_0^{X}) + c_0^2(A_0^{XY}) + c_0^2(A_0^{Y})} \lrp{\frac{\alpha(p)}{p}} \\
      &\le& C(\tilde{M},U^{\ast},L_0,U_0,\rho) c_0^2(A_0^{XY}) \lrp{\frac{\alpha(p)}{p}},
      \EEqn
      which further implies that
      \BEqn
      & & \lrabs{\cv_k^2(Z) - \cw_1} \\
      &\le& \lrabs{\cw_2+\cw_3} 
      + C(\tilde{M},U^{\ast},L_0,U_0,\rho) c_0^2(A_0^{XY}) \lrp{\frac{\alpha(p)}{p}}^2 \\
      &\le& C(\tilde{M},U^{\ast},L_0,U_0,\rho) c_0^2(A_0^{XY}) \lrp{\frac{\alpha(p)}{p}}.
      \EEqn
      Thus under Assumption \ref{Assumpt:PowerAnalysis-1} and the condition
      \begin{equation*}
          \lrabs{2c_0(A_0^{XY}) - c_0(A_0^{X}) - c_0(A_0^{Y})}
      \geq L^{\ast} |c_0(A^{XY})|,
      \end{equation*}
      we obtain that $\cw_1 = \rho^2(1-\rho)^2 \lrp{2c_0(A_0^{XY}) - c_0(A_0^{X}) - c_0(A_0^{Y})}^2$ becomes the leading term of $\cv_k^2(Z)$, which further implies the proposed result when $p \geq p_0$ for some $p_0 = p_0(\tilde{M},U^{\ast},L^{\ast},L_0,U_0,\rho)$.
    \item By Lemma \ref{Lemma:PowerAnalysis-Prep-4}\ref{Lemma:PowerAnalysis-Prep-4-1}, we have that
      \begin{equation*}
        \begin{array}{ll}
               & \lrabs{\bigp{2c_0(A^{XY}) - c_0(A^{X}) - c_0(A^{Y})} - \lrp{2c_1(A_0^{XY}) \frac{|\Delta|^2}{A^{XY}} - 2c_2(A_0^{XY}) \frac{\lrp{\bbe{|\tilde{X}_1|^2} - \bbe{|\tilde{Y}_1|^2}}^2}{(A^{XY})^2}}} \\
           \le & C(\tilde{M},L_0,U_0) |c_0(A_0^{XY})|
           \lrp{\frac{|\Delta|^4}{p^2} + \frac{\lrabs{\bbe{|\tilde{X}_1|^2} - \bbe{|\tilde{Y}_1|^2}}^3}{p^3}}.
        \end{array}
      \end{equation*}
      Under Assumption \ref{Assumpt:uniform-kernel}\ref{Assumpt:uniform-kernel-3}, when $\max\lrcp{|\Delta|^2,\lrabs{\bbe{|\tilde{X}_1|^2}-\bbe{|\tilde{Y}_1|^2}}} = o(p)$ and $\delta_1\neq\delta_2$, we obtain that $2c_1(A_0^{XY}) \frac{|\Delta|^2}{A^{XY}} - 2c_2(A_0^{XY}) \frac{\lrp{\bbe{|\tilde{X}_1|^2} - \bbe{|\tilde{Y}_1|^2}}^2}{(A^{XY})^2}$ is the leading term of $2c_0(A^{XY}) - c_0(A^{X}) - c_0(A^{Y})$, whose order is $|c_0(A_0^{XY})| p^{\max\{\delta_1,\delta_2\}-1}$.
      
      Also, it holds that
      \BEqn
          \lrabs{\cw_2} 
      &\le& C(\tilde{M},U^{\ast},L_0,U_1,\rho) |2c_0(A_0^{XY}) - c_0(A_0^{X}) - c_0(A_0^{Y})| \cdot |c_0(A_0^{XY})| \lrp{\frac{\alpha(p)}{p}} \\
      &=& o(\cw_1),
      \EEqn
      since $|c_0(A_0^{XY})| \lrp{\frac{\alpha(p)}{p}} = o\lrp{|c_0(A_0^{XY})| p^{\max\{\delta_1,\delta_2\}-1}}$ under the assumption $\alpha(p) = o\lrp{p^{\max\{\delta_1,\delta_2\}}}$.
      
      It remains to think of $\cw_3$. Define $p\times p$ matrices $\Gamma_X = (\gamma_{X,j_1 j_2})_{p\times p}$ and $\Omega_X = (\omega_{X, j_1 j_2})_{p \times p}$, where
      \begin{equation*}
          \gamma_{X, j_1 j_2} = \cum(\tilde{x}_{1 j_1}^2 - \sigma_{X, j_1}^2, \tilde{x}_{1 j_2}),
          \qquad
          \omega_{X, j_1 j_2} = \cum(\tilde{x}_{1 j_1}^2 - \sigma_{X, j_1}^2, \tilde{x}_{1 j_2}^2 - \sigma_{X, j_2}^2).
      \end{equation*}
      Similarly, we define $\Gamma_Y$ and $\Omega_Y$. Let $\fone$ denote a $p$-dimensional vector whose elements are all ones. After some tedious calculations, we obtain that
      \BEqn
      & & \cw_3 \\
      &=& \rho^2 c_1^2(A_0^{X}) (A^{X})^{-2} 
          \big{(} \bbe{(|X_1-X_2|^2-A^{X})^2} \\
      & & \hspace{8em} - 2\bbe{(|X_1-X_2|^2-A^{X})(|X_1-X_3|^2-A^{X})} \big{)} \\
      & & + 2\rho(1-\rho) c_1^2(A_0^{XY}) (A^{XY})^{-2} 
            \big{(} \bbe{(|X_1-Y_1|^2-A^{XY})^2} \\
      & & \hspace{13em} - \bbe{(|X_1-Y_1|^2-A^{XY})(|X_1-Y_2|^2-A^{XY})} \\
      & & \hspace{13em} - \bbe{(|X_1-Y_1|^2-A^{XY})(|X_2-Y_1|^2-A^{XY})}\big{)} \\
      & & + (1-\rho)^2 c_1^2(A_0^{Y}) (A^{Y})^{-2} 
            \big{(} \bbe{(|Y_1-Y_2|^2-A^{Y})^2} \\
      & & \hspace{11em} - 2\bbe{(|Y_1-Y_2|^2-A^{Y})(|Y_1-Y_3|^2-A^{Y})} \big{)} \\ 
      & & + 2\rho^2(1-\rho) c_1(A_0^{X}) c_1(A_0^{XY}) (A^{X})^{-1} (A^{XY})^{-1} \\
      & & \hspace{2em} \times
          \big{(} \bbe{(|X_1-X_2|^2-A^{X})(|X_1-X_3|^2-A^{X})} \\
      & & \hspace{3em} - 2\bbe{(|X_1-X_2|^2-A^{X})(|X_1-Y_1|^2-A^{XY})} \\
      & & \hspace{3em} + \bbe{(|X_1-Y_1|^2-A^{XY})(|X_1-Y_2|^2-A^{XY})} \big{)} \\
      & & + 2\rho(1-\rho)^2 c_1(A_0^{Y}) c_1(A_0^{XY}) (A^{Y})^{-1} (A^{XY})^{-1} \\
      & & \hspace{2em} \times
          \big{(} \bbe{(|Y_1-Y_2|^2-A^{Y})(|Y_1-Y_3|^2-A^{Y})} \\
      & & \hspace{3em} - 2\bbe{(|X_1-Y_1|^2-A^{XY})(|Y_1-Y_2|^2-A^{Y})} \\
      & & \hspace{3em} + \bbe{(|X_1-Y_1|^2-A^{XY})(|X_2-Y_1 |^2-A^{XY})} \big{)} \\
      & & + 2\rho^2(1-\rho) c_1(A_0^{XY}) (A^{XY})^{-1}
          \lrp{c_1(A_0^{X})(A^{X})^{-1} - c_1(A_0^{XY})(A^{XY})^{-1}} \\
      & & \hspace{2em} \times
          \big{(} \bbe{(|X_1-X_2|^2-A^{X})(|X_1-X_3|^2-A^{X})} \\
      & & \hspace{3em} - \bbe{(|X_1-Y_1|^2-A^{XY})(|X_1-Y_2|^2-A^{XY})} \big{)} \\
      & & + 2\rho(1-\rho)^2 c_1(A_0^{XY}) (A^{XY})^{-1} 
          \lrp{c_1(A_0^{Y})(A^{Y})^{-1} - c_1(A_0^{XY})(A^{XY})^{-1}} \\
      & & \hspace{2em} \times
          \big{(} \bbe{(|Y_1-Y_2|^2-A^{Y})(|Y_1-Y_3|^2-A^{Y})} \\
      & & \hspace{2em} - \bbe{(|X_1-Y_1|^2-A^{XY})(|X_2-Y_1|^2-A^{XY})} \big{)} \\   
      & & + 2\rho^2(1-\rho) \lrp{c_1(A_0^{X})(A^{X})^{-1} - c_1(A_0^{XY})(A^{XY})^{-1}}^2 \\
      & & \hspace{2em} \times \bbe{(|X_1-X_2|^2-A^{X})(|X_1-X_3|^2-A^{X})} \\
      & & + 2\rho(1-\rho)^2 \lrp{c_1(A_0^{Y})(A^{Y})^{-1} - c_1(A_0^{XY})(A^{XY})^{-1}}^2 \\
      & & \hspace{2em} \times \bbe{(|Y_1-Y_2|^2-A^{Y})(|Y_1-Y_3|^2-A^{Y})} \\
      &=& 4c_1^2(A_0^{XY}) (A^{XY})^{-2}
          \lrp{\|\rho\Sigma_X+(1-\rho)\Sigma_Y\|_F^2
          + 2\rho(1-\rho) \Delta^{\top} (\rho\Sigma_X + (1-\rho)\Sigma_Y) \Delta} \\
      & & + 4\rho^2 \lrp{c_1^2(A_0^{X})(A^{X})^{-2} - c_1^2(A_0^{XY})(A^{XY})^{-2}} \|\Sigma_X\|_F^2 \\
      & & + 4(1-\rho)^2 \lrp{c_1^2(A_0^{Y})(A^{Y})^{-2} - c_1^2(A_0^{XY})(A^{XY})^{-2}} \|\Sigma_Y\|_F^2 \\
      & & - 8\rho^2(1-\rho) c_1(A_0^{XY})(A^{XY})^{-1} 
            \lrp{c_1(A_0^{X})(A^{X})^{-1} - c_1(A_0^{XY})(A^{XY})^{-1}}
            \fone^{\top} \Gamma_X \Delta \\
      & & + 8\rho(1-\rho)^2 c_1(A_0^{XY})(A^{XY})^{-1} 
            \lrp{c_1(A_0^{Y})(A^{Y})^{-1} - c_1(A_0^{XY})(A^{XY})^{-1}} 
            \fone^{\top} \Gamma_Y \Delta \\
      & & + 2\rho^2(1-\rho) 
            \lrp{c_1(A_0^{X})(A^{X})^{-1} - c_1(A_0^{XY})(A^{XY})^{-1}}^2 \fone^{\top} \Omega_X \fone \\
      & & + 2\rho(1-\rho)^2 
            \lrp{c_1(A_0^{Y})(A^{Y})^{-1} - c_1(A_0^{XY})(A^{XY})^{-1}}^2 \fone^{\top} \Omega_Y \fone.
      \EEqn
      Under Assumption \ref{Assumpt:uniform-kernel}\ref{Assumpt:uniform-kernel-1}, it naturally holds that for each $f^{(p)}$ and any $s\in D$ and $s_0 \in D_0$, we also have $f_1^{(p)}(s)= f_1^{(p)}(s_0) + f_2^{(p)}(\xi(s,s_0)) (s-s_0)$, and $f_2^{(p)}(s)= f_2^{(p)}(s_0) + f_3^{(p)}(\xi'(s,s_0)) (s-s_0)$ for some points $\xi(s,s_0),\xi'(s,s_0)$ between $s$ and $s_0$, thus the smoothness assumptions proposed in Lemma \ref{Lemma:PowerAnalysis-Prep-5} and Lemma \ref{Lemma:PowerAnalysis-Prep-6} are satisfied. By using the results derived in Lemma \ref{Lemma:PowerAnalysis-Prep-2}, Lemma \ref{Lemma:PowerAnalysis-Prep-5} and Lemma \ref{Lemma:PowerAnalysis-Prep-6}, we obtain under Assumption \ref{Assumpt:PowerAnalysis-2} that
      \BEqn
      & & \lrabs{\cw_3 - \lrp{4c_1^2(A_0^{XY}) (A^{XY})^{-2}
          \lrp{\|\rho\Sigma_X+(1-\rho)\Sigma_Y\|_F^2
          + 2\rho(1-\rho) \Delta^{\top} (\rho\Sigma_X + (1-\rho)\Sigma_Y) \Delta}}} \\
      &\le& C(\tilde{M},L_0,U_0,\rho) c_0^2(A_0^{XY}) \\
      & & \hspace{3em} \times 
          \lrp{\frac{|A^{X}-A^{XY}|+|A^{Y}-A^{XY}|}{p} + \frac{|A^{X}-A^{XY}|^2+|A^{Y}-A^{XY}|^2}{p^2}} 
          \lrp{\frac{\alpha(p)}{p}} \\
      &\le& C(\tilde{M},L_0,U_0,L_2,U_2,\rho) c_0^2(A_0^{XY}) p^{\max\{\delta_1-1,(\delta_2-1)/2\}} \lrp{\frac{\alpha(p)}{p}}.
      \EEqn
      Under Assumption\ref{Assumpt:uniform-kernel}\ref{Assumpt:uniform-kernel-4} and Assumption \ref{Assumpt:PowerAnalysis-2}, $|c_0(A_0^{XY})|$ and $|c_1(A_0^{XY})|$ are of the same order and the order of $\|\rho\Sigma_X+(1-\rho)\Sigma_Y\|_F^2$ is $\alpha(p)p$, thus we have that $4c_1^2(A_0^{X}) (A^{X})^{-2}  \|\rho\Sigma_X+(1-\rho)\Sigma_Y\|_F^2 + 8\rho(1-\rho) \Delta^{\top} (\rho\Sigma_X + (1-\rho)\Sigma_Y) \Delta$ is the leading term of $\cw_3$. 
      
      Therefore, the leading term of $\cv_k^2(Z)$ in this case is either $\cw_1$ or $\cw_3$, whichever that has the higher order. Specifically, when $\alpha(p) = o\lrp{p^{2\max\{\delta_1,\delta_2\}-1}}$, $\cw_1$ dominates $\cw_3$ in order and $\cw_3$ is the leading term of $\cv_k^2(Z)$ otherwise, which leads to the desired result when $p \geq p_0$ for some $p_0 = p_0(\tilde{M},\hat{M},U^{\ast},L_0,U_0,L_2,U_2,\rho)$.
      
      \item Following the analysis in the previous case, it suffices to compare the orders of $\cw_1$ and $\cw_3$. In this case, it holds that the order of $\cw_1$ is $c_0^2(A_0^{XY}) p^{2\delta_3-4}$ whereas the leading term of $\cw_3$ becomes $4c_1^2(A_0^{XY}) (A^{XY})^{-2} \|\rho\Sigma_X+(1-\rho)\Sigma_Y\|_F^2$. Note that it follows from the triangle inequality that
      \begin{equation*}
          \|\Sigma_Y\|_F - \rho\|\Sigma_X-\Sigma_Y\|_F 
      \le \|\rho\Sigma_X+(1-\rho)\Sigma_Y\|_F
      \le \rho\|\Sigma_X\|_F+(1-\rho)\|\Sigma_Y\|_F,
      \end{equation*}
      thus it holds under Assumption \ref{Assumpt:component-dept}\ref{Assumpt:component-dept-4} and Assumption \ref{Assumpt:PowerAnalysis-3} that, $\|\rho\Sigma_X+(1-\rho)\Sigma_Y\|_F^2$ has order $\frac{\alpha(p)}{p}$, which further implies that the order of $\cw_3$ is $c_1^2(A_0^{XY}) \lrp{\frac{\alpha(p)}{p}}$. Hence in this case, $\cw_3$ is always the leading term of $\cv_k^2(Z)$ when $p\geq p_0$ for some $p_0 = p_0(\tilde{M},\hat{M},U^{\ast},L_0,U_0,L_3,U_3,\rho)$, which completes the proof.

      \item Under Assumption \ref{Assumpt:PowerAnalysis-4}, we have $A^{X} = A^{XY} = A^{Y}$, which implies that $\cw_1=\cw_2=0$. It also hols in this case that $\Sigma_X=\Sigma_Y$, then based on the analysis in the previous cases, we obtain that
      \begin{equation*}
          \cw_3 = 4c_1^2(A_0^{XY})(A^{XY})^{-2} \|\Sigma_X\|_F^2,
      \end{equation*}
      which implies under Assumption \ref{Assumpt:component-dept}\ref{Assumpt:component-dept-4} that the order of $\cw_3$ is $c_1^2(A_0^{XY}) \lrp{\frac{\alpha(p)}{p}}$. Note that in this case, we have that
      \begin{equation*}
          |\cv_k^2(Z) - \cw_3| \le C(\tilde{M},U^{\ast},L_0,U_0,\rho) c_0^2(A_0^{XY}) \lrp{\frac{\alpha(p)}{p}}^2,
      \end{equation*}
      then under Assumption \ref{Assumpt:uniform-kernel}\ref{Assumpt:uniform-kernel-3} we obtain the desired results when $p \geq p_0$ for some $p_0 = p_0(\tilde{M},\hat{M},U^{\ast},L_0,U_0,\rho)$.
    \end{enumerate}  
\end{proof}

\begin{lemma}\label{Lemma:PowerAnalysis-3}
Suppose that $f=f^{(p)}$ satisfies Assumption \ref{Assumpt:uniform-kernel}\ref{Assumpt:uniform-kernel-1}-\ref{Assumpt:uniform-kernel-3}, and Assumption \ref{Assumpt:component-dept}\ref{Assumpt:component-dept-1}-\ref{Assumpt:component-dept-3} hold, then
\begin{enumerate}[label=(\roman*)]
    \item \label{Lemma:PowerAnalysis-3-1}
    it holds that $\bbe{k^2(Z_1,Z_2)} \le C(\tilde{M},U^{\ast},L_0,U_0,\rho) c_0^2(A_0^{XY})$.

    \item \label{Lemma:PowerAnalysis-3-2}
    it holds that
    \begin{equation*}
        \max\lrcp{\bbe{\tilde{k}^4(X_1,X_2)}, \bbe{\tilde{k}^4(X_1,Y_1)}, \bbe{\tilde{k}^4(Y_1,Y_2)}} 
    \le  C(\tilde{M},U^{\ast},L_0,U_0,\rho) c_0^4(A_0^{XY}) \lrp{\frac{\alpha(p)}{p}}^2. 
    \end{equation*}
    where $\tilde{k}$ denotes the centered version of $k$,

    \item \label{Lemma:PowerAnalysis-3-3}
    there exists some universal positive constant $C<\infty$, such that
    \BEqn
    & & \bbe{(h^k(X_1,X_2,Y_1,Y_2))^2} \\
    &\le& C \max\lrcp{C(\tilde{M},U^{\ast},L_0,U_0) c_0^2(A_0^{XY}) \lrp{\frac{\alpha(p)}{p}}, \bigp{\ce^k(X,Y)}^2}.
    \EEqn
\end{enumerate}
\end{lemma}
\begin{proof}
  \begin{enumerate}[label=(\roman*)]
    \item 
    By using Lemma \ref{Lemma:approx-2}, it holds that
    \BEqn
        k^2(X_1,Y_1)
    &=& \sum\limits_{i_1,i_2=0}^{6} c_{i_1}(A_0^{XY}) c_{i_2}(A_0^{XY}) (A^{XY})^{-i_1-i_2}
                                    \lrp{|X_1-Y_1|^2-A^{XY}}^{i_1+i_2} \\
    & & + 2R(X_1,Y_1)
          \sum\limits_{i=0}^{6} c_i(A_0^{XY}) (A^{XY})^{-i} \lrp{|X_1-Y_1|^2-A^{XY}}^i \\
    & & + R^2(X_1,Y_2),
    \EEqn
    and it follows from Lemma \ref{Lemma:Li-1} that
    \begin{equation*}
        \lrabs{\bbe{k^2(X_1,Y_1)} - c_0^2(A^{XY})} 
    \le C(\tilde{M},U^{\ast},L_0,U_0) c_0^2(X_0^{XY}) \lrp{\frac{\alpha(p)}{p}}.
    \end{equation*}
    We obtain the similar results for $\bbe{k^2(X_1,X_2)}$ and $\bbe{k^2(Y_1,Y_2)}$, and it thus follows from the definition of $Z$ and Lemma \ref{Lemma:PowerAnalysis-Prep-3} that
    \BEqn
    & & \lrabs{\bbe{k^2(Z_1,Z_2)} - \lrp{\rho^2 c_0^2(A_0^{X}) + 2\rho(1-\rho) c_0^2(A_0^{XY}) + (1-\rho)^2 c_0(A_0^{Y})}} \\
    &\le& C(\tilde{M},U^{\ast},L_0,U_0) c_0^2(X_0^{XY}) \lrp{\frac{\alpha(p)}{p}},
    \EEqn
    which implies that
    \BEqn
        \bbe{k^2(Z_1,Z_2)} 
    &\le& \lrabs{\rho^2 c_0^2(A_0^{X}) + 2\rho(1-\rho) c_0^2(A_0^{XY}) + (1-\rho)^2 c_0(A_0^{Y})} \\
    &   &  + C(\tilde{M},U^{\ast},L_0,U_0) c_0^2(X_0^{XY}) \lrp{\frac{\alpha(p)}{p}} \\
    &\le& C(\tilde{M},U^{\ast},L_0,U_0,\rho) c_0^2(A_0^{XY}).
    \EEqn
    \item The second step can be proved using the similar arguments used in Lemma \ref{Lemma:RateConvg-2}, and we thus spare the details.
    \item It follows from the definition of $h^k(X_1,X_2,Y_1,Y_2)$ and the $c_r$ inequality, we have
    \BEqn
    & & \bbe{(h^k(X_1,X_2,Y_1,Y_2))^2} \\
    &\le& 7\lrp{\bbe{\tilde{k}^2(X_1,X_2)} + \bbe{\tilde{k}^2(X_1,Y_1)} + \bbe{\tilde{k}^2(Y_1,Y_2)}} 
          + \lrp{\ce^k(X,Y)}^2.
    \EEqn
    Using the similar techniques as used to prove Lemma \ref{Lemma:RateConvg-2}, we can show that
    \begin{equation*}
        \max\{\bbe{\tilde{k}^2(X_1,X_2)}, \bbe{\tilde{k}^2(X_1,Y_1)}, \bbe{\tilde{k}^2(Y_1,Y_2}\} 
    \le C(\tilde{M},U^{\ast},L_0,U_0) c_0^2(A_0^{XY}) \lrp{\frac{\alpha(p)}{p}},
    \end{equation*}
    which completes the proof.
  \end{enumerate}
\end{proof}

\end{appendices}
\end{document}